\pdfoutput=1
\documentclass[11pt]{article}
\usepackage{color}
\usepackage{graphicx}
\usepackage{amsmath,amsfonts,amssymb,graphics,amsthm}
\usepackage{hyperref}
\usepackage{comment}
\usepackage{tabularx}
\usepackage[protrusion=true,expansion=true]{microtype}
\usepackage{enumerate}
\usepackage{bbm}
\usepackage{mathrsfs}
\usepackage{paralist}
\usepackage{enumitem}

\usepackage[margin=1in]{geometry}

\hypersetup{
    colorlinks=false,
    linktocpage,
    }

\numberwithin{equation}{section}

\newtheorem{theorem}{Theorem}[section]

\newtheorem{lemma}[theorem]{Lemma}

\newtheorem{prop}[theorem]{Proposition} 
\newtheorem{lem}[theorem]{Lemma} 
\newtheorem{cor}[theorem]{Corollary}

\theoremstyle{definition}
\newtheorem{definition}[theorem]{Definition}
\newtheorem{defn}[theorem]{Definition}

\theoremstyle{remark}
\newtheorem{remark}[theorem]{Remark}

\def\@rst #1 #2other{#1}
\newcommand\MR[1]{\relax\ifhmode\unskip\spacefactor3000 \space\fi
  \MRhref{\expandafter\@rst #1 other}{#1}}
\newcommand{\MRhref}[2]{\href{http://www.ams.org/mathscinet-getitem?mr=#1}{MR#2}}

\def\MR#1{\href{http://www.ams.org/mathscinet-getitem?mr=#1}{MR#1}}

\newcommand{\C}{\mathbbm{C}}
\newcommand{\D}{\mathbbm{D}}
\newcommand{\E}{\mathbbm{E}}
\newcommand{\N}{\mathbbm{N}}

\newcommand{\Z}{\mathbbm{Z}}
\newcommand{\R}{\mathbbm{R}}
\renewcommand{\P}{\mathbbm{P}}

\newcommand{\eps}{\varepsilon}
\newcommand{\ep}{\varepsilon}
\newcommand{\1}{\mathbf{1}}

\newcommand{\PE}{\P_\mathrm{E}}
\newcommand{\PQ}{\P_\mathrm{Q}}
\newcommand{\PW}{\P^k_\mathrm{W}}
\newcommand{\PI}{\P^k_\mathrm{I}}
\newcommand{\EE}{\E_\mathrm{E}}
\newcommand{\EQ}{\E_\mathrm{Q}}
\newcommand{\br}{\op{br}}
\newcommand{\sh}{\op{sh}}
\newcommand{\ti}{t}
\newcommand{\si}{s}

\newcommand{\scr}{\mathscr}

\DeclareMathOperator{\diam}{diam}

\DeclareMathOperator{\Cov}{Cov}
\DeclareMathOperator{\Var}{Var}

\def\cZ{\mathcal{Z}}
\def\cY{\mathcal{Y}}

\def\cW{\mathcal{W}}
\def\cV{\mathcal{V}}
\def\cU{\mathcal{U}}
\def\cT{\mathcal{T}}
\def\cS{\mathcal{S}}
\def\cR{\mathcal{R}}

\def\cN{\mathcal{N}}

\def\cL{\mathcal{L}}

\def\cJ{\mathcal{J}}
\def\cI{\mathcal{I}}
\def\cH{\mathcal{H}}
\def\cG{\mathcal{G}}
\def\cF{\mathcal{F}}
\def\cE{\mathcal{E}}
\def\cD{\mathcal{D}}
\def\cC{\mathcal{C}}

\def\cA{\mathcal{A}}

\newcommand{\aryb}{\begin{eqnarray*}}
\newcommand{\arye}{\end{eqnarray*}}
\def\alb#1\ale{\begin{align*}#1\end{align*}}
\newcommand{\eqb}{\begin{equation}}
\newcommand{\eqe}{\end{equation}}
\newcommand{\eqbn}{\begin{equation*}}
\newcommand{\eqen}{\end{equation*}}

\newcommand{\PhiM}{\Phi_{\op{M}}}
\newcommand{\PhiT}{\Phi_{\op{T}}}
\newcommand{\PhiTM}{\Phi_{\op{T}\to\op{M}} }

\newcommand{\BB}{\mathbbm}
\newcommand{\ol}{\overline}
\newcommand{\ul}{\underline}
\newcommand{\op}{\operatorname}

\newcommand{\frk}{\mathfrak}
\newcommand{\eqD}{\overset{d}{=}}

\newcommand{\rta}{\rightarrow}

\newcommand{\wt}{\widetilde}
\newcommand{\wh}{\widehat} 
\newcommand{\mcl}{\mathcal}

\newcommand{\tr}{\op{tr}}
\newcommand{\rmL}{\op{L}}

\newcommand{\rmD}{\op{D}}

\newcommand{\cJL}{{\cJ_{\op{L}}}}

\newcommand{\cJLp}{{\cJ^+_{\op{L}}}}
\newcommand{\cJRp}{{\cJ^+_{\op{R}}}}
\newcommand{\cJLm}{{\cJ^-_{\op{L}}}}
\newcommand{\cJRm}{{\cJ^-_{\op{R}}}}

\newcommand{\cJp}{{\cJ^+}}
\newcommand{\cJm}{{\cJ^-}}

\newcommand{\ET}{{\cE_{\op{T}}}}
\newcommand{\EG}{{\cE_{\op{G}}}}
\newcommand{\EM}{{\cE_{\op{M}}}}

\newcommand{\tL}{t_{\op{L}}}
\newcommand{\tR}{t_{\op{R}}}

\newcommand{\dZ}{{ d_{\cZ} }}

\newcommand{\triang}{T}
\newcommand{\tree}{\frk T}

\DeclareMathAlphabet{\mathpzc}{OT1}{pzc}{m}{it}

\begin{document}

\title{The Tutte embedding of tree-weighted planar maps converges to Liouville quantum gravity}
\author{Nina Holden\thanks{Courant Institute, New York University; nina.holden@nyu.edu.}}
\date{}

\maketitle

\begin{abstract}
	Suppose $M$ is a planar map with a macroscopic boundary decorated by a spanning tree $T$, such that $(M,T)$ has been sampled from the critical Boltzman measure. We show that $M$ converges to a Liouville quantum gravity disk with parameter $\gamma=\sqrt{2}$ under the Tutte embedding as the length of the boundary of $M$ goes to infinity. Furthermore, conditioned on $(M,T)$, the simple random walk on $M$ converges to a planar Brownian motion modulo reparametrization of time. %
\end{abstract}

\tableofcontents

\section{Introduction}

Planar maps are, roughly speaking, planar graphs embedded in the sphere, defined up to continuous deformations of the sphere. Planar maps are studied in many different branches of both mathematics and physics, including combinatorics, probability theory, geometry, string theory, and conformal field theory. In particular, in probability theory and theoretical physics, random planar maps serve as models for random surfaces. 

There are a number of natural probability measures one can put on the set of planar maps. One possibility is to consider the uniform measure, where all planar maps of a given size have the same probability. Another possibility is to consider a planar map with a statistical physics model such that the law of the planar map is reweighted by the partition function of the statistical physics model. In this paper, we will consider the latter setting where the statistical physics model is the spanning tree. We prove convergence of tree-weighted planar maps under the Tutte embedding to $\sqrt{2}$-LQG, together with a quenched invariance principle for the associated random walks on the embedded map.

A longstanding belief is that a number of random planar maps converge in the scaling limit to continuum random surfaces called Liouville quantum gravity (LQG) surfaces \cite{david-conformal-gauge,dk-qg,kpz-scaling}. This conjecture has been verified for several different planar map models and several notions or topologies of convergence. Let us describe three commonly considered notions of convergence. 

One gets a natural notion of convergence by viewing the random planar map as a metric space by equipping it with its graph distance and proving that this metric space, after appropriate rescaling, converges in law for the Gromov-Hausdorff topology. This notion of convergence has been established for a number of uniform planar map models, see e.g. \cite{legall-uniqueness,miermont-brownian-map}. The limiting surfaces are called Brownian surfaces, which are equivalent to LQG with parameter $\gamma=\sqrt{8/3}$ \cite{lqg-tbm3}.

An alternative is to use the mating-of-trees approach to LQG to define a notion of convergence. Several planar map models with a statistical physics model can be encoded bijectively by a 2D walk. Similarly, in the continuum a 2D Brownian motion encodes an LQG surface decorated by the random fractal curve known as the Schramm-Loewner evolution (SLE), and there is an explicit relationship between the parameter $\gamma\in(0,2)$ of the LQG and the correlation of the 2D Brownian motion. Convergence of the walk to 2D Brownian motion can be interpreted as a scaling limit result of the planar map to LQG. This notion of convergence has been established for a number of planar map models, see e.g. \cite{shef-burger} and the survey \cite{ghs-mating-survey}.

In this paper we will primarily be interested in a third notion of convergence which is called convergence under conformal embedding. One embeds the planar map into the complex plane $\BB C$ via a discrete conformal embedding and wants to prove that the measure induced by vertex counting measure on $\BB C$ converges under appropriate renormalization to the area measure of the LQG surface. There are a number of different discrete conformal embeddings which can be used, e.g.\ the Tutte embedding, circle packings, embedding via the uniformization theorem for Riemann surfaces, and embeddings which use observables of statistical physics models on the planar maps, such as the Cardy-Smirnov embedding. Convergence of planar maps under discrete conformal embedding has been stated as an open problem in several influential probabilistic papers, including \cite{shef-kpz,legall-sphere-survey,dkrv-lqg-sphere}.

Only a small number of convergence results for planar maps to LQG under discrete conformal embeddings have been proven in earlier literature.
 The case of uniform triangulations and the Cardy-Smirnov embedding was established in \cite{hs-cardy-embedding}. Furthermore, convergence was proven for two planar map models that can be viewed as coarse-grained LQG surfaces (namely, the mated-CRT map and the Poisson-Voronoi tessellation of an LQG surface) and embedding given by the Tutte embedding \cite{gms-tutte,gms-poisson-voronoi,afs-metric-ball} or the Smith embedding \cite{bgs-smith}.

Conformal embedding of planar maps has also been studied in e.g.\  \cite{curien-glimpse,gill-rohde-type,gjn-macroscopic-circles,gjn-dirichlet,chelkak-icm,chelkak-s-emb-periodic,clr-t-emb,klrr-circle-patterns} and subsequent works motivated by the relationship with LQG and Brownian surfaces.
Note that before the current work, the only combinatorially natural planar map model for which convergence to LQG had been established under conformal embedding was the case of uniform triangulations, and the methods of \cite{hs-cardy-embedding} are hard to adapt to non-uniform planar maps and to other embeddings.

In Section \ref{sec:intro-main} we will state our main result, namely the convergence of tree-weighted planar maps to LQG with $\gamma=\sqrt{2}$ under the Tutte embedding, along with an invariance principle for random walk on such embedded planar maps to 2D Brownian motion. In Section \ref{sec:intro-proof} we will give an outline of the proof and the remainder of the paper. In Section \ref{sec:intro-future} we will present future directions. In particular, we will explain that techniques developed in this paper can be used to prove counterparts of Theorem \ref{thm-tutte-conv0} for other discrete conformal embeddings and most likely also for other random planar map models.

\subsection{Main result}
\label{sec:intro-main}
A finite planar map $M$ is a proper embedding of a finite connected (multi)graph in the two-dimensional sphere $\BB S^2$, where we identify two planar maps if one can be obtained from the other via an orientation-preserving homeomorphism of $\BB S^2$. Let $\cV(M)$, $\cE(M)$ and $\cF(M)$ denote the set of vertices, edges and faces, respectively, of $M$. All the planar maps we consider will be rooted, meaning that there is a distinguished oriented edge $\frk e\in \cE(M)$. The tail vertex of $\frk e$ is called the \emph{root vertex} and the face to the right of $\frk e$ is called the \emph{root face}. 
A \emph{planar map with disk topology with boundary length $\ell\in\N$} is a finite planar map where the root face has simple boundary of length $\ell$. The vertices on the boundary of the root face are called \emph{boundary vertices} and are denoted by $\partial\cV(M)$ while the other vertices are called \emph{interior vertices}. 

For $n\in\N$ let $(M_n,\frk e_n)$ be a rooted planar map with boundary length $\lceil\sqrt{n}\rceil$. In particular, the numbers of vertices and edges are not fixed. 
 A collection of edges $\tree_n\subset \cE(M_n)$ is called a \emph{spanning tree with wired boundary conditions} if 
\begin{compactitem}
	\item for all vertices $v\in\cV(M_n)$ there is an edge in $\tree_n$ with $v$ as end-point,
	\item the edges of $\tree_n$ form exactly one cycle, which is equal to the boundary of the root face, and
	\item $\tree_n$ is connected.
\end{compactitem}
Notice that if we delete all the boundary edges and identify all the boundary vertices with a single vertex, then the remaining edges in $\frk T_n$ form a spanning tree. 
 We give the tuple $(M_n,\frk e_n, \tree_n)$ weight 
$$
w(M_n,\frk e_n, \tree_n) = 16^{-\# \cE(M_n)}. %
$$
Define a probability measure $\frk m$ on the set of tuples $(M_n,\frk e_n,\tree_n)$ such that for a normalizing constant $C_n>0$,
$$
\frk m(M_n,\frk e_n,\tree_n)=C_n^{-1}w(M_n,\frk e_n,\tree_n);
$$ 
as part of the proof of Theorem \ref{thm-tutte-conv0} below, we will see that $\frk m$ is well-defined.

Suppose $M_n$ has at least one interior vertex and let $\frk v\in \cV(M_n)\setminus\partial\cV(M_n)$ be some interior vertex. The {\bf Tutte embedding centered at $\frk v$} of $(M_n,\frk e_n)$ is a map $\phi:\cV(M_n)\to \ol{\BB D}$ defined as follows. Enumerate the boundary vertices of $M_n$ in counterclockwise order as $\{y_1,\dots,y_{\lceil\sqrt{n}\rceil}\}$, such that $y_1$ and $y_{\lceil\sqrt{n}\rceil}$ are the ending and starting vertices, respectively, of $\frk e_n$. For $j = 1,\dots,{\lceil\sqrt{n}\rceil}$ we let $\frk p(y_j)$ be the conditional probability given $M_n$ that a simple random walk started from $\frk v$ first hits the boundary at a vertex in the boundary arc $\{y_1,\dots,y_j\}$. The boundary vertices $y_1, y_2, \ldots, y_{\lceil\sqrt{n}\rceil}$ are then mapped in counterclockwise order to the unit circle $\partial \BB D$ such that $\phi(y_j):= e^{2\pi i \frk p(y_j)}$. Then the hitting probability of the random walk started from $\frk v$ approximates the uniform measure on the unit circle, provided the maximal hitting probability among the boundary vertices is small. We embed the interior vertices of $M_n$ via the discrete harmonic extension of the boundary values in the sense that for any interior vertex $v$, 
$$
\phi(v)= %
\frac{1}{\op{deg}_{\op{S}}(v)}
\sum_{w\sim v} %
\phi(w),
$$ 
where we write $w\sim v$ if $w$ and $v$ are connected by an edge and $w\neq v$, and $\op{deg}_{\op{S}}(v)$ is the number of vertices adjacent to $v$. With this convention we have unit conductance between any pair of adjacent vertices, even if there are multiple edges connecting them. %
Notice that the Tutte embedding is not necessarily injective. Also observe that the vertex $\frk v$ which we use to center the embedding is not necessarily mapped exactly to 0.

Sample a tuple $(M_n,\frk e_n,\tree_n)$ from the probability measure $\frk m$. If $M_n$ has at least one interior vertex let $\frk v\in \cV(M_n)$ be a uniformly sampled interior vertex and let $\phi$ be the Tutte embedding of $(M_n,\frk e_n)$ centered at $\frk v$. If no interior vertex exists, then define $\frk v$ and $\phi$ in an arbitrary way, for example let $\frk v$ be the starting vertex of the root edge and let $\phi$ be the map sending all vertices to 1.

 Let $\mu_n$ be the random measure on $\ol{\BB D}$ which assigns mass $n^{-1}$ to each embedded vertex of $\cV(M_n)$, i.e., for any Borel set $A\subset\ol{\D}$,
\eqb
\mu_n(A):= n^{-1}\cdot \#\{ w\in \cV(M_n)\,:\,\phi(w)\in A \}.
\label{eq:mu-ell}
\eqe
Similarly, let $\nu_n$ be the random measure on $\partial\BB D$ which assigns mass $n^{-1/2}$ to each embedded boundary vertex of $\cV(M_n)$, i.e., for any Borel set $A\subset \partial\D$,
\eqb
\nu_n(A):= n^{-1/2}\cdot \#\{ w\in \partial\cV(M_n)\,:\,\phi(w)\in A \}.
\label{eq:nu-ell}
\eqe

In probability theory, LQG is a theory of random surfaces. A $\gamma$-LQG surface for $\gamma\in(0,2)$ is, heuristically speaking, a random surface where the metric can be written in the form $e^{\gamma h}(dx^2+dy^2)$, where $h$ is a variant of the Gaussian free field (GFF) and $dx^2+dy^2$ is the standard Euclidean metric. One can define rigorously the area measure $\mu$, the boundary measure $\nu$, and the distance function (metric) of this surface. One can represent a $\gamma$-LQG surface by a pair $(D,h)$ called an \emph{embedding}, where $D\subset\C$ and $h$ is a variant of the GFF.

The unit boundary length $\gamma$-LQG disk is a canonical $\gamma$-LQG surface with disk topology. We refer to Section \ref{sec:lqg} for a description of the surface we are working with. In Theorem \ref{thm-tutte-conv0}, $\mu$ and $\nu$ refer to the area and boundary length measures on $\D$ and $\partial\D$, respectively, induced by $h$ in the special case $\gamma=\sqrt{2}$.
 The  embedding $(\D,h)$ is chosen such that a point sampled from $\nu$ is sent to 1 and a point sampled independently from $\mu/\mu(\D)$ is sent to 0.

For continuous curves $\beta_1:[0,T_{1}]\to\C$ and $\beta_2:[0,T_{2}]\to\C$, define 
	\begin{equation}\label{eq:curves-mod-reparam}
		d_{\rm CMP}(\beta_1,\beta_2):=\inf_\phi \sup_{t\in[0,T_{1}]}|\beta_1(t)-\beta_2(\phi(t))|,
	\end{equation}
	where the infimum is over all increasing homeomorphisms $\phi:[0,T_{1}]\to [0,T_{2}]$ and ``CMP'' stands for ``curves modulo parameterization". We allow curves  $\beta:\{0\}\to\C$ of duration 0 by identifying them with the constant curve $t\mapsto \beta(0)$ for $t\in[0,1]$. We identify curves at zero $d_{\rm CMP}$-distance and equip the resulting quotient space with the metric $d_{\rm CMP}$.

 For $v\in \cV(M_n)$ let $\cW: \{0,1,\dots,T \}\to \cV(M_n)$ denote the random walk on $M_n$ started from $v$ and stopped at the random time $T$ at which the walk hits a vertex on the boundary of $M_n$. We suppose the transition probabilities of the walk at any time are uniform on the set of adjacent vertices. Let $W: [0,T]\to\ol{\D}$ be the embedding of $\cW$, which is defined by $ W(t):=\phi(\cW(t))$ for $t\in\{0,1,\dots,T \}$ and by linear interpolation otherwise. Let $\sigma^{v,n}$ denote the law of $W$ conditioned on $(M_n,\frk e_n,\frk v)$. Let $\sigma^z$ denote the law of 2D Brownian motion started from $z$ run until it first exits $\D$. We view $\sigma^{v,n}$ and $\sigma^z$ as (random) probability measures on the space of curves in $\ol{\BB D}$ viewed modulo reparametrizations of time, see \eqref{eq:curves-mod-reparam}.

	Let $d_{\op{Prok}}$ denote the Prokhorov metric, with the underlying metric space depending on the context. We view $\mu_n$ and $\mu$ (resp.\ $\nu_n$ and $\nu$) as random variables with values in the space of finite measures on $\ol{\D}$ (resp.\ $\partial\D$) equipped with the Prokhorov metric. We use $\xrightarrow{d}$ %
	to denote convergence in law.
\begin{theorem}\label{thm-tutte-conv0}
	In the setting above, the probability measure $\frk m$ is well-defined,  and for an appropriate deterministic constant $c>0$ it holds that $(c\mu_n,\nu_n)$ converge in law as $n\rta\infty$ to the area and boundary length measures $(\mu,\nu)$ of the unit boundary length $\sqrt{2}$-LQG disk. 
 Furthermore, the conditional laws of the embedded random walks satisfy %
	$$
	\sup_{v\in\cV(M_n)}d_{\op{Prok}}(\sigma^{v,n},\sigma^{\phi(v)} ) 
	\xrightarrow[n\to\infty]{d} 
	0.
	$$
\end{theorem}

\subsection{Proof outline}
\label{sec:intro-proof}
By a variant of the Mullin bijection (see Proposition \ref{prop:mullin-wired-disk2}), the tree-decorated planar map $(M_n,\frk e_n,\frk T_n)$ can be encoded bijectively by a nearest-neighbor 2D lattice path $\cZ^n$ on $\N_0^2$ from $(\lceil\sqrt{n}\rceil-1,0)$ to $(0,0)$. 
 By mating of trees \cite{wedges,ag-disk}, a unit boundary length quantum disk of parameter $\gamma=\sqrt{2}$ decorated by a so-called clockwise space-filling SLE$_\kappa$ $\eta$ with $\kappa=8$ can be encoded by a 2D Brownian excursion $Z$ in $\R_+^2$ from $(1,0)$ to $(0,0)$. We say that $Z$ is the boundary length process of the SLE-decorated LQG surface. As $n\rta\infty$, the lattice path $\cZ^n$ converges in the scaling limit to $Z$. This gives us a first relationship between our tree-weighted planar map and $\sqrt{2}$-LQG.

In the Mullin bijection each vertex $v\in\cV(M_n)$ can be associated with some time $\frk t_v\in\R$ for the lattice path $\cZ^n$, see the text above Lemma \ref{prop-bij-steps} below. 
Suppose we have a coupling of $\cZ^n$ for $n\in\N$ and $Z$ such that $n^{-1/2}\cZ^n\rta Z$ as $n\rta\infty$, and  suppose $\frk t_v$ is expressed in the rescaled time coordinate for which this convergence holds.
 We get a natural embedding $\wt\phi:\cV(M_n)\to\ol{\D}$ of $M_n$ by setting $\wt\phi(v)=\eta(\frk t_v)$ for each $v\in\cV(M_n)$. 
 It is not hard to argue that the measure-convergence assertion of Theorem \ref{thm-tutte-conv0} (namely, the convergence of the two vertex counting measures to LQG) holds with $\phi$ replaced by $\wt\phi$. 
 In order to get the same result for $\phi$ it is sufficient to argue that $\|\wt\phi-\phi\|_\infty\rta 0$ as $n\rta\infty$. 
 By properties of the Tutte embedding one can argue that $\|\wt\phi-\phi\|_\infty\rta 0$ if random walk on the $\wt\phi$-embedded map converges to 2D Brownian motion, modulo reparametrization of time, provided some consistency and boundary conditions are also satisfied (Lemmas \ref{lem:rw-conv-implies-tutte-conv1} and \ref{lem:rw-conv-implies-tutte-conv2}). Therefore, in order to conclude the proof of Theorem \ref{thm-tutte-conv0} the main task is to establish this invariance principle.  
 Notice that $\wt\phi$, which defines the random environment of the walk, depends both on $Z$ and on $\cZ^n$.

In order to get the desired invariance principle, we will make use of a general invariance principle for  environments which, among other requirements, are assumed to be so-called translation invariant modulo scaling. This general invariance principle was first proven in \cite{gms-random-walk} and we slightly generalize it in Theorem \ref{thm-general-clt-uniform} below. We are given a so-called cell configuration $\cH^{\op{cell}}$, which is a collection of cells (i.e., compact connected subsets of $\C$) $H^{\op{cell}}$ which have union $\C$ and pairwise disjoint interior, along with a symmetric relation $\sim$ on $\cH^{\op{cell}}\times\cH^{\op{cell}}$ and pairwise conductances; see Definition \ref{def-cell-config}. The invariance principle concerns a random walk on $\cH^{\op{cell}}$ which is allowed to move between cells $H^{\op{cell}}_1,H^{\op{cell}}_2\in\cH^{\op{cell}}$ iff $H_1^{\op{cell}}\sim H_2^{\op{cell}}$. 
 The transition probabilities of the random walk are proportional to the conductances, which are always equal to 1 in our setting.

 The LQG surface known as the 0-quantum cone with a space-filling SLE satisfies some properties closely related to translation invariance modulo scaling. Our high-level strategy for establishing the desired invariance principle is the following: 
\begin{compactitem}
\item[(a)] couple the boundary length process $Z$ for the SLE-decorated 0-quantum cone with a (non-i.i.d.) nearest-neighbor walk $\cZ$ on $\Z^2$ (Section \ref{sec:kmt-0-cone}), 
\item[(b)] define a cell configuration $\cH^{\op{cell}}=\cH^{\op{cell}}(\cZ,Z)$ satisfying the conditions of Theorem \ref{thm-general-clt-uniform},\footnote{In fact, this theorem concerns two cell configurations $\cH^{\op{cell}}$ and $\cH^{\op{large}}$, but we focus on only $\cH^{\op{cell}}$ in this outline to simplify the presentation.} where in particular the cell configuration is isomorphic, viewed as a %
graph, 
 to the planar map encoded by the walk $\cZ$ via the Mullin bijection (Section \ref{sec:cell-conf}), and

\item[(c)] transfer the invariance principle to the planar map appearing in Theorem \ref{thm-tutte-conv0} (Section \ref{sec:tutte-conv}).

\end{compactitem}
The original invariance principle of \cite{gms-random-walk} was used in earlier works to establish convergence of planar maps to LQG under the Tutte embedding, more specifically for the mated-CRT map \cite{gms-tutte} and the Poisson-Voronoi tessellation \cite{gms-poisson-voronoi,afs-metric-ball}.

 The planar maps considered in these works are defined in an explicit way in terms of an SLE-decorated LQG surface and it is significantly easier to make use of \cite{gms-random-walk} to establish an invariance principle in these settings, and it is also significantly easier to transfer the invariance principle from the setting of the 0-quantum cone to the planar map of interest. In particular, these works do not need to make use of the Mullin bijection or couplings between random walk and Brownian motion, so the steps (a)-(c) above are either skipped or greatly simplified. Let us describe in the next few paragraphs (in chronological order) some of the new ideas that are introduced in this paper in order to overcome the difficulties that arise in the context of tree-weighted maps. These paragraphs also serve as brief overviews of Sections \ref{sec:kmt}, \ref{sec:kmt-0-cone}, \ref{sec:cell-conf}, and \ref{sec:tutte-conv}. We refer to the beginning of these sections (and associated subsections) for more details, and a discussion of further difficulties we encounter can also be found there.

\vskip5pt

\noindent {\bf Section \ref{sec:kmt}}: We construct a novel coupling between a higher-dimensional correlated Brownian motion $Z$ and a continuous-time random walk $\cZ$ with i.i.d.\ steps. The walk has a finite step size distribution and the jump times are given by an independent Poisson point process on $\R$. It is a KMT-type coupling such that $\sup_{s,t\in J}\|(Z_t-Z_s)-(\cZ_t-\cZ_s)\|_\infty$ is at most of order $\log |J|$ with high probability if $J\subset\R$ is a long interval of length $|J|$.

 Our coupling is significantly simpler than the coupling for higher-dimensional processes in \cite{zaitsev-kmt}. Furthermore, the coupling is translation invariant in law, unlike most other KMT-type couplings, where intervals of the form $[i2^k,(i+1)2^k]$ typically play a special role. %
  We prove that the walk restricted to a given interval is determined by the Brownian motion restricted to a quantitatively larger interval, plus some explicit additional randomness, with high probability. However, in order to determine $\cZ$ with probability 1 we need to observe the full process $Z$. 
 Several results in this section may have applications also in other contexts than planar maps and LQG.

\vskip5pt

\noindent 
{\bf Section \ref{sec:kmt-0-cone}}: The infinite-volume variant of the Mullin bijection gives a correspondence between a bi-infinite walk with i.i.d.\ increments in $\{(\pm 1,0),(0,\pm 1) \}$ and the so-called infinite tree-decorated planar map $(M,\frk e,\frk T)$. This walk can be coupled to an uncorrelated 2D Brownian motion via the coupling from Section \ref{sec:kmt}, and this 2D Brownian motion encodes a $\gamma$-quantum cone decorated with whole-plane space-filling SLE$_\kappa$ for $\gamma=\sqrt{2}$ and $\kappa=8$. 
 A first natural attempt to get the desired invariance principle for random walk on $M_n$ via Theorem \ref{thm-general-clt-uniform}, 
would be to define a cell configuration such that the cell associated with $v\in\cV(M)$ contains the point $\eta(\frk t_v)$. However, a fundamental difficulty with this approach is that a cell configuration satisfying this property is typically not translation invariant modulo scaling. The main issue is the $\gamma$-log singularity of $\gamma$-quantum cone at the origin. In order to get an environment which is translation invariant modulo scaling we need to instead consider a 0-quantum cone, which does not have such a singularity. However, the 0-quantum cone and the $\gamma$-quantum cone are singular due to different behavior near 0 and infinity. For an SLE-decorated 0-quantum cone one can define a boundary length process $Z$, but this process will be singular with respect to a 2D Brownian motion nearby the times 0 and $\pm\infty$. Therefore it is not even clear how to define an appropriate coupling between the SLE-decorated 0-quantum cone and a nearest-neighbor random walk on $\Z^2$; we need such a walk in order to define a suitable cell configuration $\cH^{\op{cell}}=\cH^{\op{cell}}(\cZ,Z)$. 

Our idea for constructing an appropriate coupling, and showing that some associated cell configuration satisfies the assumptions of Theorem \ref{thm-general-clt-uniform}, is to use that the laws of the 0-quantum cone and the $\gamma$-quantum cone can be explicitly compared in the following sense. If we are given a 0-quantum cone and sample a point $z_0$ from the $\gamma$-LQG area measure close to the origin, then the surface with marked points at $z_0$ (instead of 0) and $\infty$ will be absolutely continuous with respect to a $\gamma$-quantum cone restricted to any bounded set, with an explicit Radon-Nikodym derivative (Lemma \ref{lem:rooted}). 
 This allows us to bound probabilities for the SLE-decorated 0-quantum cone in terms of probabilities for the SLE-decorated $\gamma$-quantum cone. 
Even after the comparison result for the 0-quantum cone and the $\gamma$-quantum cone has been proven, there are major difficulties with proving the existence of a suitable coupling between the SLE-decorated 0-quantum cone and a 2D random walk with steps $(\pm 1,0),(0,\pm 1)$. 
 The issue is that the coupling between the 2D Brownian motion $Z$ and the bi-infinite i.i.d.\ walk $\cZ$ is non-local in the sense that one needs to observe the full process $Z$ in order to define $\cZ|_I$ for any bounded interval $I\subset\R$. The comparison results for the 0-quantum cone and the $\gamma$-quantum cone only apply for events which depend on bounded parts of the surface, and the moments of the Radon-Nikodym derivative term diverge fairly quickly (namely, polynomially with a possibly large power) as the diameter of the bounded set goes to infinity. This divergence rate is too large to directly apply the estimates from Section \ref{sec:kmt} discussed above.

 In order to establish a suitable coupling between the SLE-decorated 0-quantum cone and a 2D random walk, we prove that the explicit construction of a 2D random walk $\cZ$ from the boundary length process $Z$ which was developed in Section \ref{sec:kmt} still can be carried out in the setting of the 0-quantum cone. We need to use particular properties of this coupling (including Theorem \ref{prop-kmt}, Lemma \ref{prop:split-rand} and Lemma \ref{prop:add-one-jump}) and properties of the SLE-decorated 0-quantum cone, e.g.\ translation invariance modulo scaling. Key properties of the resulting coupling are stated in Theorem \ref{prop:zero-vs-gamma-walk} and Corollary \ref{cor:zero-vs-gamma}.

\vskip5pt

\noindent {\bf Section \ref{sec:cell-conf}}: Let $Z$ be the boundary length process of the SLE-decorated 0-quantum cone and let $\cZ$ be the associated random walk on $\Z^2$. In this section we construct a cell configuration $\cH^{\op{cell}}=\cH^{\op{cell}}(\cZ,Z)$  satisfying the assumptions of Theorem \ref{thm-general-clt-uniform}. In the Mullin bijection, each vertex $v\in\cV(M)$ can be associated with an interval $I_v\subset\R$ such that the intervals $I_v$ have union $\R$ and pairwise disjoint interior, and such that $\frk t_v\in I_v$. A first attempt to define a cell configuration would be to let the cell $H^{\op{cell}}_v\subset\C$ associated with $v$ be given by $\eta(I_v)$, and we say that $H^{\op{cell}}_v\sim H^{\op{cell}}_w$ iff $v,w$ are connected by an edge. However, this does not define a cell configuration, since in order to be a cell configuration it is required that if $H^{\op{cell}}_v\sim H^{\op{cell}}_w$ then $H^{\op{cell}}_v\cap H^{\op{cell}}_w\neq\emptyset$. To get around this issue, we 
associate each $v\in\cV(M)$ with a set $\cT^{\op{cell}}_v\subset\cV(M)\cup\cV(M_*)$, where $M_*$ is the dual map of $M$. We then define $H^{\op{cell}}_v$ to be the union of suitable subsets of $\eta(I_w)$ for $w\in\cT^{\op{cell}}_v$, where $I_w$ and $\frk t_w$ for $w\in\cV(M_*)$ are defined in a similar way as for primal vertices. In order to define sets $\cT^{\op{cell}}_v$ with the desired properties we need to study how adjacency of two vertices $v,w\in\cV(M)$ is related to the condition that $\eta(I_v)\cap\eta(I_w)\neq\emptyset$, which requires a careful study of combinatorial properties of the Mullin bijection.

Theorem \ref{thm-general-clt-uniform} requires, among other properties, that the cell configuration satisfies a certain moment condition. In order to establish this, we need various quantitative estimates for the non-i.i.d.\ 2D walk $\cZ$ coupled with the SLE-decorated 0-quantum cone, which we a priori know very little about. We get sufficient estimates for this walk by using Theorem \ref{prop:zero-vs-gamma-walk} and Corollary \ref{cor:zero-vs-gamma} (which compare our walk to a walk with i.i.d.\ increments), estimates for SLE and LQG (building on e.g.\ Section \ref{sec:0-cone-BL} and Appendix \ref{app2}), and combinatorial properties of the Mullin bijection. We get in particular bounds for the size of $\cT^{\op{cell}}_{v_0}$ and the diameter of $\{\frk t_w\,:\,w\in \cT^{\op{cell}}_{v_0} \}$, where $v_0$ is defined such that $0\in I_{v_0}$.
 Theorem \ref{thm-general-clt-uniform} also requires the cell configuration to satisfy a particular line connectivity property, which has been weakened as compared to the analogous property in \cite{gms-random-walk}, and which requires a study of delicate combinatorial and geometric properties of $\cT^{\op{cell}}_{v_0}$.

\vskip5pt

\noindent {\bf Section \ref{sec:tutte-conv}}: By Theorem \ref{thm-general-clt-uniform} and the results of Section \ref{sec:cell-conf}, we get an invariance principle for random  walk on our cell configuration $\cH^{\op{cell}}=\cH^{\op{cell}}(\cZ,Z)$ to 2D Brownian motion, where $Z$ is the boundary length process of the SLE-decorated 0-quantum cone and $\cZ$ is the associated (non-i.i.d.) 2D random walk. In order to conclude the proof, we need to transfer this invariance principle to the setting of a quantum disk, where the boundary length process is a Brownian cone excursion and the associated random walk encodes the planar map in Theorem \ref{thm-tutte-conv0} via the Mullin bijection. There are major difficulties with transferring the invariance principle from the 0-quantum cone case to the quantum disk case due to the non-local nature of the coupling we are using.
 More precisely, by using the absolute continuity result described above for the 0-quantum cone and the (recentered at another point) $\gamma$-quantum cone, it is not hard to compare locally the boundary length processes for the SLE-decorated 0-quantum cone  and the SLE-decorated quantum disk. However, it is hard to get a comparison result for the walks from such a local result, especially since, if we want to determine $\cZ|_{[-n,n]}$ with high probability we need to observe $Z|_{[-an,an]}$ for $a=a(n)\rta\infty$ as $n\rta\infty$.  

The hardest step of the proof is to transfer the invariance principle to a setting where the 2D random walk has i.i.d.\ increments on $\{(\pm 1,0),(0,\pm 1) \}$, namely, we want to prove an invariance principle for the random walk on $\cH^{\op{cell}}=\cH^{\op{cell}}(\cZ,Z)$, where $Z$ is the 2D Brownian motion describing the boundary length process of the SLE-decorated $\gamma$-quantum cone and $\cZ$ is the associated i.i.d.\ walk on $\Z^2$. This step is achieved in Section \ref{sec:tutte-conv0}. One difficulty is that we know very little about the boundary length process $Z$ and the walk $\cZ$ in the case of the 0-quantum cone. For example, unlike in the Brownian case, it is hard to show approximate independence of $Z$ on different scales, e.g.\ that the random variables $X_J$ in \eqref{eq:X-def} for $J\subset\R$ a short interval are approximately independent of the random variables $X_J$ for $J\subset\R$ a long interval. In order to obtain approximate independence results of this type, we make extensive use of Markov-type properties and various conditional laws for the GFF, along with the geometric interpretation of the random variables $X_J$ in terms of LQG lengths in Section \ref{sec:proof-prop10}. In Section \ref{sec:proof-prop8} we further make use of various exact laws and independence results for the coupling in Section \ref{sec:kmt}. We refer to Section \ref{sec:tutte-conv00} for a more extensive discussion of the proof in Section \ref{sec:tutte-conv0}.

\subsection{Future directions}
\label{sec:intro-future}
\begin{enumerate}
	\item The author and Pu Yu will use results from this paper and \cite{hy-cp-ru} to establish the counterpart of Theorem \ref{thm-tutte-conv0} for the circle packing embedding and the Riemann uniformization embedding. It would also be of interest to extend the theorem to other discrete conformal embeddings, including the Smith embedding, square tiling, s-embedding, t-embedding/circle pattern/Coulomb gauge, and to other planar map models. 
 For random planar maps which can be encoded with a mating-of-trees-type bijection by an i.i.d.\ walk \cite{bhs-site-perc,kmsw-bipolar,lsw-schnyder-wood}, most of the methods in this paper carry through. The main work in this case is to redo the combinatorial arguments in Section \ref{sec:cell-conf} and, in the case where $\gamma\in(\sqrt{2},2)$, address the fact that space-filling SLE segments are no longer simply connected. Significant novel ideas are needed for other random planar maps.
	\item The coupling between random walk and Brownian motion in Section \ref{sec:kmt} might be helpful for proving convergence to LQG of tree-weighted (and other non-uniform) planar maps in Gromov-Hausdorff topology. A tree-weighted map can be coupled together with a variant of the so-called mated-CRT map \cite{gms-tutte} such that there is a natural bijective correspondence between the two vertex sets.  
	Results in Section \ref{sec:cell-conf} imply that if we pick a uniform pair of adjacent vertices in the tree-weighted map, then the associated vertices in the mated-CRT map have distance $O(1)$ with high probability, and vice versa.	From Sections \ref{sec:kmt} and \ref{sec:cell-conf} we have a good understanding of what causes the exceptional events that give larger discrepancies in the distances. One could attempt to deduce from this that macroscopic distances for the two metrics differ by at most a constant multiple. If one has established a convergence result for the mated-CRT map to LQG for the Gromov-Hausdorff topology via a tightness and uniqueness argument, one might then be able to transfer it to the tree-weighted map if the uniqueness argument also holds in this setting.
	\item We believe one can deduce from results in this paper that the Peano curve of the spanning tree on $M_n$ in Theorem \ref{thm-tutte-conv0} converges to SLE$_8$ in a quenched sense conditioned on the planar map and its Tutte embedding.
 Annealed convergence to SLE$_8$ can be established with little work from the methods of this paper.
 Quenched convergence to SLE$_8$ can likely be established by using Wilson's algorithm and by proving quenched convergence of the loop-erased random walk to SLE$_2$, where the latter result can likely be derived from the analogous quenched result for the (non-loop-erased) random walk in Theorem \ref{thm-tutte-conv0}. In order to carry out this proof, one would need to extend our quenched invariance principle to the setting where conductances are proportional to edge multiplicities since this is the walk used by Wilson's algorithm; this should be possible to do via small modifications of the argument in the current paper.
	\item One can try to use estimates from this paper to study other statistical physics models on tree-weighted planar maps, such as diffusion-limited aggregation (DLA) \cite{gp-dla}, internal DLA \cite{bg-idla} and loop soups. One can study further properties of random walk on planar maps, such as convergence as a parameterized curve to the Liouville Brownian motion \cite{berestycki-gwynne-lbm}. One can also study other observables, such as the curvature of the planar map \cite{cp-gwynne-curvature}, the Laplacian determinant \cite{apps-central-charge}, other spectral properties (see e.g.\ \cite{bw-weyl,berestycki-spectral,bk-spectral} for recent continuum results), or correlation functions of statistical physics models, and one can consider other aspects of the discrete Laplacian and other differential operators on the planar map. As seen in the references, some of these questions have already been studied in the context of the mated-CRT map.
	\item All convergence results for planar maps to LQG in current literature use the mating-of-trees approach to LQG.\footnote{Convergence of uniform planar maps to Brownian surfaces does not use mating of trees, but the equivalence of Brownian surfaces and $\sqrt{8/3}$-LQG relies on mating of trees.} An interesting question is whether this can be avoided. It would also be of interest to establish convergence of planar maps to LQG via softer techniques, for example by comparing some given random planar map model to a random planar map model that is known to converge to LQG.
\end{enumerate}

\vskip10pt

\noindent {\bf Notation}. 
We let $\N_0=\{0,1,2,\dots \}$. 
For a set $I\subset\R$ we let $I_\Z:=I\cap \Z$ denote the set of integers in $I$. For $k\in\N$ we let $[k]=\{1,\dots,k \}$ denote the set of natural numbers between 1 and $k$. For $z\in\C$ and $r\geq 0$ we let $B(z,r)=\{w\in\C\,:\,|z-w|<r \}$ denote the open disk centered at $z$ with radius $r$. We let $\D=B(0,1)$ denote the unit disk.

\vskip10pt

\noindent {\bf Statement}. I have worked on this paper for several years and I am currently in the process of revising it. I was prompted to post the paper today since several closely related results were released by Open AI yesterday (October 6). In particular, the results of this paper are close to  ``A linear clock for random walk on tree-weighted planar maps'' by Open AI, which builds on ``Metric-measure limits of subcritical FK and spanning-tree planar maps'' and ``Canonical conformal limits of subcritical FK planar maps'', both also by Open AI.

The ``Future directions'' section above was written before Open AI posted these and related papers on October 6. The Open AI papers solve several of the open problems discussed in 1., 2.\ and 4.\ above.

\vskip10pt

\noindent {\bf Statement on AI usage}.  I used ChatGPT Astra for proofreading and revising the first draft of the paper. AI tools were not used to develop mathematical ideas or in the writing of the manuscript. All suggested changes were reviewed by the author, who takes full responsibility for the content of the paper. 

\vskip10pt

\noindent {\bf Acknowledgments}. 
I thank Pu Yu for many helpful comments on the draft. This research was supported by NSF grants DMS-2246820 and DMS-2553117, along with the Simons Collaboration on Probabilistic Paths to Quantum Field Theory.

\section{Preliminaries}

\subsection{Liouville quantum gravity and mating of trees}
\label{sec:mot}

\subsubsection{Liouville quantum gravity}
\label{sec:lqg}
Let $D\subset\C$ be a domain and let $h$ be a Gaussian free field (GFF) on $D$; see e.g.\ \cite{shef-gff,wp-gff} or the introductory sections of \cite{ig1,ig4,ss-contour,shef-zipper} for background on the GFF. We can define the Liouville area measure $\mu_h$ on $D$ as the regularized limit of $e^{\gamma h}$ times Lebesgue area measure on $D$, see e.g.\  \cite{shef-kpz,berestycki-gmt-elementary,kahane,rhodes-vargas-review}. Similarly, one can define a Liouville length measure on $\partial D$ and on SLE-type curves on $D$ which are independent from $h$.

Let $\gamma\in(0,2)$, $Q=2/\gamma+\gamma/2$, and $k,m\in\N_0$. A $\gamma$-LQG surface with $k$ marked points and $m$ curves is an equivalence class of tuples $(D,h,z_1,\dots,z_k,\eta_1,\dots,\eta_m)$, where $D\subset\BB C$ is a domain, $h$ is a distribution on $D$, $z_1,\dots,z_k\in D\cup\partial D$, and $\eta_j:I_j\to D\cup\partial D$ for $k=1,\dots,m$ are (continuous) curves for intervals $I_j\subseteq\R$. Two such tuples $(D,h,z_1,\dots,z_k,\eta_1,\dots,\eta_m)$ and $(\wt D,\wt h,\wt z_1,\dots,\wt z_k,\wt\eta_1,\dots,\wt\eta_m)$ are defined to be equivalent if there is a conformal map $\phi:\wt D\to D$ such that
\eqb
\wt h = h\circ \phi + Q\log|\phi'|,\qquad 
\phi(\wt z_j)=z_j\,\,\text{for}\,\,j=1,\dots,k,\qquad
\wt\eta_j=\phi^{-1}\circ\eta_j\,\,\text{for}\,\,j=1,\dots,m.
\label{eq:coord-change}
\eqe
The formula \eqref{eq:coord-change} implies that $\phi_*\mu_{\wt h}=\mu_h$, where $_*$ denotes pushforward of measures, and a similar transformation rule holds for the Liouville length measures. An \emph{embedding} of a $\gamma$-LQG surface refers to a representative $(D,h,z_1,\dots,z_k,\eta_1,\dots,\eta_m)$ of the equivalence class. Sometimes we will abuse notation and call $(D,h,z_1,\dots,z_k,\eta_1,\dots,\eta_m)$ a quantum surface, and we are then referring to the equivalence class that it defines.

The three main $\gamma$-LQG surfaces that we encounter in this paper are the $\alpha$-quantum cone for $\alpha\in\{0,\gamma \}$ and the quantum disk, which were first defined in \cite{wedges}. The $\alpha$-quantum cone, which can be defined for any $\alpha<Q$, is a $\gamma$-LQG surface $(D,h,0,\infty)$ with two marked points. Its $\gamma$-LQG measure has infinite total mass, but assigns finite mass to every bounded measurable
subset of $\BB C$. Roughly speaking, the $\alpha$-quantum cone is obtained by starting with a GFF plus $-\alpha\log |\cdot-z_0|$ for some fixed interior point $z_0$ of the domain and then ``zooming in'' near $z_0$ while rescaling the areas so that the quantum area of the unit disk remains of constant order \cite[Proposition 4.13(ii)]{wedges}. 

To define this surface more precisely, let $B$ be a standard linear Brownian motion and let $\wh B$ be an independent standard linear Brownian motion conditioned so that $\wh B_t+(Q-\alpha)t>0$. Let $A:\R\to\R$ be the process such that $A_t = B_t + \alpha t$ for $t\geq 0$ and 
$A_t =\wh B_{-t} +\alpha t$ for $t<0$. Let $h$ be the random distribution such that
if $h_r(0)$ denotes the circle average of $h$ on $\partial B_r(0)$ 
(see e.g.\ \cite[Section 3.1]{shef-kpz} for the definition and basic properties of the circle average), 
then $t\mapsto h_{e^{-t}}(0)$ has the same law as the process A and 
$h-h_{|\cdot|}(0)$ 
is independent from $h_{|\cdot|}(0)$ and has the same law as the analogous process for a whole-plane GFF. 
The tuple $(\C,h,0,\infty)$ is what we call the circle-average embedding of the $\alpha$-quantum cone. One reason this embedding is convenient, is that $h|_{\D}$ has the same law as the restriction to $\D$ of
a whole-plane GFF plus $-\alpha \log |\cdot|$, normalized so that its circle average over $\partial\D$ is 0.

We will use a special scale invariance property of the $\alpha$-quantum cone which we now describe.
Let $h$ be the circle-average embedding of an $\alpha$-quantum cone and for $b > 0$ let 
\eqb \label{eqn-mass-hit-time}
\frk R(b) := \sup\left\{ r > 0 : h_r(0) + Q \log r = \frac{1}{\gamma} \log b \right\}. 
\eqe 
In other words, $\frk R(b)$ is the largest radius $r> 0$ so that if we scale spatially by the factor $r$ and apply the coordinate change formula~\eqref{eq:coord-change}, then the average of the resulting field on $\partial \BB D$ is equal to $\gamma^{-1} \log b$. It is not hard to argue (see e.g.~\cite[Proposition 4.13(i)]{wedges}) that for each fixed $b>0$, 
\eqb \label{eqn-cone-scale}
h \eqD h(\frk R(b) \cdot) + Q \log \frk R(b) -  \frac{1}{\gamma} \log b .
\eqe
By~\eqref{eq:coord-change}, if we let $h^b$ be the field on the right side of~\eqref{eqn-cone-scale}, then a.s.\ $\mu_{h^b}(A) = b^{-1} \mu_h(\frk R(b) A)$ for each Borel set $A\subset \BB C$. In particular, typically $\mu_h(B_{\frk R(b)}) \asymp b$.

Given $r>0$ it will sometimes be convenient to work with what we call the $r$-embedding of the $\alpha$-quantum cone $(\C,h,0,\infty)$. This embedding is given by $(\C,\wt h,0,\infty)$ for $\wt h=h\circ\phi_r+Q\log|\phi_r'|$, where $\phi_r:\C\to\C$ is given by $\phi_r(z)=z\frk R(r^{\gamma Q})/r$ and $\frk R(r^{\gamma Q})$ is defined by \eqref{eqn-mass-hit-time}. Restricted to $B(0,r)$, $\wt h$ has the law of $-\alpha\log|r^{-1} \cdot|$ plus a whole-plane GFF with additive constant fixed so that it has average 0 on $\partial B(0,r)$.

Let $R>1$ be fixed and suppose $\alpha=0$. By the proof of \cite[Proposition 3.2]{gms-poisson-voronoi}, if $z_0$ is a point sampled from Lebesgue measure in $B(0,\frk R(b))$ then for appropriately chosen random $C>0$ we have the following convergence in total variation distance as $b\rta\infty$,   
\eqb
(h(C\cdot+z_0)+Q\log C)|_{B(0,R)}\rta h|_{B(0,R)}.
\label{eq50}
\eqe

A quantum disk is a finite-mass quantum surface $(\D, h)$ which is defined in \cite[Definition 4.21]{wedges}. We will omit the precise definition as it will not be needed. In this paper we will consider the variant of the quantum disk with unit boundary length (i.e., $\nu_h(\partial \D)=1$) and one marked boundary point. Conditional on $(\D,h)$, the marked point $z_0$ has the law of a point sampled from the boundary measure $\nu_h$. In Theorem \ref{thm-tutte-conv0} we work with an embedding such that $z_0=1$ and $w_0=0$, where $w_0$ has the law of a point sampled from the LQG area measure $\mu_h$ rescaled to be a probability measure.

\subsubsection{The Schramm-Loewner evolution and mating of trees}
\label{sec:sle}
The Schramm-Loewner evolution (SLE) is a random fractal curve introduced in \cite{schramm0} which describes the scaling limit of curves in critical statistical mechanics models. SLE curves have a parameter $\kappa>0$. In this paper we will consider the case $\kappa=8$, but the constructions we will describe work in the same way for other $\kappa\geq 8$ (in fact, some constructions even work for all $\kappa>4$) and we will work with general $\kappa\geq 8$ in the remainder of this section.

The most classical variant of SLE$_\kappa$ is the chordal SLE, which connects two given boundary points $a\neq b$ of a simply connected domain $D$. When $\kappa\geq 8$ this curve is space-filling, i.e., the curve hits every point of $D$. If we let $b\rta a$ from the counterclockwise (resp.\ clockwise) direction then in the limit we get a space-filling curve in $D$ which starts and ends at $a$. This variant of SLE$_\kappa$ is called the \emph{clockwise (resp.\ counterclockwise) space-filling SLE$_\kappa$ starting and ending at $a$}.

A third variant of SLE$_\kappa$ is the \emph{whole-plane space-filling SLE$_\kappa$}. It is a space-filling curve in $\C$ starting and ending at $\infty$. It can be constructed by considering a chordal SLE$_\kappa$ and zooming in nearby a deterministic point $z_0\in D$. An alternative construction can be given with imaginary geometry. The basic idea of the construction is to first construct the outer boundary of the curve stopped when it hits each fixed point $z\in\BB Q^2$ using a pair of flow lines of a GFF and then use these outer boundary curves to define a partial order on $\BB Q^2$ which can be extended to a continuous curve. 

Let $\ul\kappa  = 16/\kappa \in (0,2]$ and $\chi^{\op{IG}} := 2/\sqrt{\ul\kappa} -\sqrt{\ul\kappa}/2$.  Let $h^{\op{IG}}$ be a whole-plane GFF viewed modulo a global additive multiple of $2\pi\chi^{\op{IG}}$, as in~\cite{ig4}, where here ``IG'' stands for ``imaginary geometry".

For $z\in \BB C$ and $\theta \in (0,2\pi)$, we can define the \emph{flow line} of $h^{\op{IG}}$ started from $z$ with angle $\theta$ as in~\cite[Theorem~1.1]{ig4}. This flow line is a whole-plane SLE$_{\ul\kappa}(2-\ul\kappa)$ curve from $z$ to $\infty$ which is a.s.\ determined by $h^{\op{IG}}$, where we recall that a whole-plane SLE$_{\ul\kappa}(2-\ul\kappa)$ is a variant of SLE$_{\ul\kappa}$ which is defined in~\cite[Section~2.1]{ig4}. For $z\in\BB Q^2$, let $\eta_z^L$ and $\eta_z^R$ be the flow lines started from $z$ with angles $\pi/2$ and $-\pi/2$, respectively. These curves will be the left and right boundaries of $\eta$ stopped upon hitting $z$. 

For distinct $z,w \in \BB Q^2$, the flow lines $\eta_z^L$ and $\eta_w^L$ a.s.\ merge upon intersecting, and similarly with $R$ in place of $L$. The two flow lines $\eta_z^L$ and $\eta_z^R$ started at the same point a.s.\ do not touch each other except at their starting point~\cite[Theorem~1.7]{ig4}. 

We define a total order on $\BB Q^2$ by declaring that $z$ comes before $w$ if and only if $w$ is in a connected component of $\BB C\setminus (\eta_z^L\cup \eta_z^R)$ which lies to the right of $\eta_z^L$ (equivalently, to the left of $\eta_z^R$).  The whole-plane analog of~\cite[Theorem~4.12]{ig4} (see~\cite[Footnote 4]{wedges}) shows that there is a well-defined continuous curve $\eta : \BB R\rta \BB C$ which traces the points of $\BB Q^2$ in the above order. The curve $\eta$ is called the \emph{whole-plane space-filling SLE$_{\kappa}$ from $\infty$ to $\infty$} associated with $h^{\op{IG}}$. For any $a < b$ the set $\eta([a,b])$ has the topology of a closed disk.  

\begin{figure}
	\centering
	\includegraphics[scale=1]{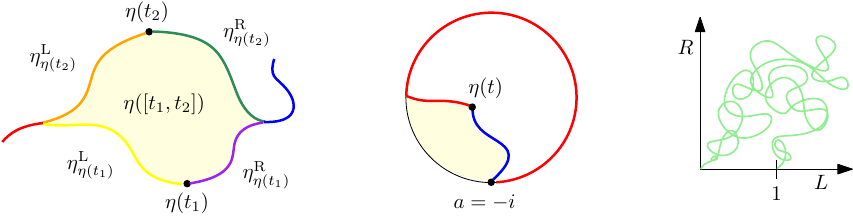}
	\caption{\label{fig:sle-mot} 
		\textbf{Left:} The left (resp.\ right) outer boundary of $\eta((-\infty,t_1])$, which we denote by $\eta_{\eta(t_1)}^{\mathrm{L}}$ (resp.\ $\eta_{\eta(t_1)}^{\mathrm{R}}$ ), is 
		the union of the red and yellow (resp.\ blue and purple) curves and the left (resp.\ right) outer boundary of $\eta((-\infty,t_2])$ is the union of the red and orange (resp.\ blue and green) curves. 
		The increment $L_{t_2}-L_{t_1}$ (resp.\ $R_{t_2}-R_{t_1}$) is the length of the orange (resp.\ green) curve minus the length of the yellow (resp.\ purple) curve. 
		\textbf{Middle:} A space-filling SLE$_\kappa$ loop in $\D$ from $a=-i$ to $a$ for $\kappa\geq 8$. The quantity $L_t$ (resp.\ $R_t$) is the LQG boundary length of the red (resp.\ blue) curve. 
		\textbf{Right:} The process $(L,R)$ encoding the SLE-decorated quantum disk in the middle figure is a correlated planar Brownian motion started from $(1,0)$ and conditioned to stay in the first quadrant  until reaching the origin. 
	}
	\label{fig:bdy-length-process}
\end{figure}

For $\gamma\in(0,\sqrt{2}]$ and $\kappa=16/\gamma^2\geq 8$ let $(\C,h,0,\infty,\eta)$ be the curve-decorated $\gamma$-LQG surface given by a $\gamma$-quantum cone with an independent whole-plane space-filling SLE$_\kappa$ $\eta$. It follows from \cite{ig4,wedges} that we may \emph{parameterize $\eta$ by $\mu_h$}, which means that $\eta(0)=0$ and that $\mu_h(\eta([s,t]))=t-s$ for any $-\infty<s<t<\infty$. Let $L_t$ (resp.\ $R_t$) denote the length of the left (resp.\ right) boundary arc of $\eta((-\infty,t])$ relative to the left (resp.\ right) boundary arc of $\eta((-\infty,0])$; see Figure \ref{fig:bdy-length-process} for an illustration. By \cite[Theorem 1.9]{wedges}, $(L,R)$ has the law of a two-sided Brownian motion such that the correlation between the coordinates is given by $-\cos(\pi\gamma^2/4)$ for $\kappa\leq 8$. It was proven in \cite{kappa8-cov} that this formula also holds for $\kappa>8$. Notice that in the special case $\gamma=\sqrt{2}$ the coordinates are independent. We will assume our Liouville boundary measure has been multiplied by an appropriate deterministic constant such that $\Var(L_1)=\Var(R_1)=0.5$.

Now let $(\D,h,-i,\eta)$ be the curve-decorated $\gamma$-LQG surface given by a unit boundary length $\gamma$-LQG disk with a clockwise space-filling SLE$_\kappa$ from $-i$ to $-i$, sampled independently of $h$ and then parameterized by $\mu_h$. Then one can define a similar boundary length process $(L,R)$, see Figure \ref{fig:bdy-length-process}. In this setting, $(L,R)$ has the law of a Brownian cone excursion in the first quadrant from $(1,0)$ to $(0,0)$ with the same variances and covariances as in the whole-plane case \cite{ag-disk}. The duration of the Brownian cone excursion is random and equal to $\mu_h(\D)$.

\subsection{The Mullin bijection}
\label{sec:bijection}
The Mullin bijection~\cite{mullin-maps} is a bijection between a spanning-tree decorated map and a nearest-neighbor walk in $\BB Z^2$. We first present the infinite-volume version of the bijection in Section \ref{sec:mullin-infinite} and then we give some finite volume versions of the bijection in Section \ref{sec:bijection-disk}. We believe Lemma \ref{prop:mullin-wired-disk} and Proposition \ref{prop:mullin-wired-disk2} have not been stated in earlier literature, while we refer to \cite{mullin-maps,bernardi-maps,shef-burger,chen-fk,ghs-mating-survey,ghs-map-dist} for further details on the other material. 

\subsubsection{Mullin bijection for bi-infinite walks}
\label{sec:mullin-infinite}
Let $(M, e_0, \tree)$ be the \emph{uniform infinite spanning-tree decorated planar map}, which is the Benjamini-Schramm~\cite{benjamini-schramm-topology} limit of uniform triples consisting of a planar map with an oriented root edge and a distinguished spanning tree. This infinite-volume limit is proven to exist in~\cite{shef-burger,chen-fk}.  

Let $M_*$ be the dual map of $M$ and let $\tree_*$ be the dual spanning tree of $\tree$, i.e., $\tree_*$ is the set of edges of $M_*$ which do not cross edges of $\tree$. Let $\triang$ be the triangulation with vertex set $\cV(M)\cup \cV(M_*)$ and edge set given by the union of $\tree\cup \tree_*$ and edges of the form $\{v,v^* \}$, where $v\in \cV(M)$ is incident to the face that is dual to $v^*\in \cV(M)$. Following the convention of %
\cite{ghs-mating-survey} we declare that the root edge of $\triang$ is the edge $\{v,v^* \}$ (oriented from $v$ to $v^*$) whose primal endpoint $v$ coincides with the terminal endpoint of $e_0$ and which is the first edge in counterclockwise order after $e_0$ with this property.

There is a unique path $\frk{tr} : \BB Z \rta \cF(\triang)$ in the dual of $\triang$ which hits each face of $\triang$ exactly once, does not cross $\tree$ or $\tree_*$ in the sense that $\frk{tr}(\ti)$ and $\frk{tr}(\ti-1)$ share an edge not belonging to $\tree\cup \tree_*$ for each $\ti\in \BB Z$, and it holds that $\frk{tr}(0)$ (resp.\ $\frk{tr}(1)$) is the face to the right (resp.\ left) of the root edge. We refer to e.g.\ \cite[Figure 6]{ghs-map-dist} for an illustration and to Section \ref{sec:bijection-disk}, \cite[Section 4.1]{shef-burger} and \cite{ghs-mating-survey} for illustrations of the finite volume case. %

Let $\cZ=(\cL,\cR)$ be a standard nearest-neighbor simple random walk on $\BB Z^2$. Let $\cJLp\subset\Z$ (resp.\ $\cJLm\subset\Z$) be the set of times at which $\cL$ goes up (resp.\ down), i.e., 
$$
\cJLp = \{ \ti\in\Z\,:\,\cL_{\ti}>\cL_{\ti-1} \}
\qquad\text{and}\qquad
\cJLm = \{ \ti\in\Z\,:\,\cL_{\ti}<\cL_{\ti-1} \}.
$$
Define $\cJRp$ and $\cJRm$ in the same manner. Observe that the four sets $\cJLp,\cJLm,\cJRp,\cJRm$ are disjoint with union $\BB Z$ and define $\cJp:=\cJLp\cup\cJRp$ and $\cJm:=\cJLm\cup\cJRm$.

\begin{figure}[ht!]
	\begin{center}
		\includegraphics[scale=.8]{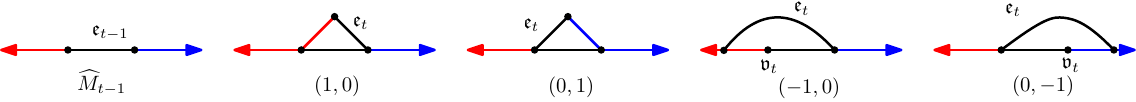} 
		\caption{\label{fig-mullin-sewing} Illustration of the sewing procedure. To the left you see the outer boundary of the map $\wh M_{\ti-1}$, which is a subgraph of the triangulation $\triang$. The external face lies above the line. The other four figures show how to construct $\wh M_{\ti}$ from $\wh M_{\ti-1}$ given the walk increment $\cZ_{\ti}-\cZ_{\ti-1}$. The face added in the four rightmost figures is $\frk{tr}(\ti)$.}
	\end{center}
\end{figure}

The walk $\cZ$ is naturally associated to an instance of $(M,e_0,\tree)$ such that $(M,e_0,\tree)$ is constructed from $\cZ$ via the so-called sewing procedure. The sewing procedure is an inductive relationship between certain submaps of the triangulation $\triang$. See Figure~\ref{fig-mullin-sewing}.  
For $\ti \in\BB Z$, let $\wh M_{\ti}$ be the submap of $\triang$ consisting of the vertices and edges on the boundaries of the triangles in $\frk{tr}((-\infty,\ti]_{\BB Z})$. 
For $\ti\in\BB Z$, let $\frk e_{\ti}$ be the edge shared by the triangles $\frk{tr}(\ti)$ and $\frk{tr}(\ti+1)$. 
Then $\wh M_{\ti}$ is an infinite triangulation with a simple boundary such that the external face has infinitely many edges lying to the left and right of $\frk e_{\ti}$. We can recover $(\wh M_{\ti} , \frk e_{\ti})$ from $(\wh M_{\ti-1 } , \frk e_{\ti-1 })$ by attaching the triangle $\frk{tr}(\ti)$, which has one edge equal to $\frk e_{\ti}$ and one edge equal to $\frk e_{\ti-1}$. If $\ti\in\cJLp$ (resp.\ $\ti\in\cJRp$) this triangle has two edges in the external face of $\wh M_{\ti-1}$ and $\frk e_{\ti}$ is the right (resp.\ left) one of these two edges if we face outward toward the external face. If $\ti\in\cJLm$ (resp.\ $\ti\in\cJRm$), then one of the edges of $\frk{tr}(\ti)$ coincides with the edge of $\wh M_{\ti-1}$ immediately to the left (resp.\ right) of $\frk e_{\ti-1}$, and the remaining edge lies in the external face of $\wh M_{\ti-1}$ and equals $\frk e_{\ti}$. In the latter case let $\frk v_{\ti}$ be the common end-point of this edge and $\frk e_{\ti-1}$. Note that $\ti$ is the first time for which $\frak v_{\ti}$ is a non-boundary vertex of the map $\wh M_{\ti}$.

The edges of the tree $\tree$ (resp.\ $\tree_*$) are precisely the edges of $\triang$ which are not crossed by $\frk{tr}$ and which lie to the left (resp.\ right) of the path $\frk{tr}$.
For each edge $e_*\in\tree_*$ its dual is the edge $e=\{v_1,v_2 \}\in \cE(M)$ such that $v_1,v_2\in \cV(M)$ lie on the boundary of the two faces $f_1,f_2\in \cF(T)$ that are incident to $e_*$. Using these observations we can recover $(M,e_0,\tree)$ from $\triang$ and $\frk{tr}$ by 
letting $\cV(M)$ be the set of vertices to the left of $\tr$ and by 
letting $\cE(M)$ be the union of $\tree$ and the dual edges to $\tree_*$. Furthermore, $e_0$ is uniquely characterized by the procedure described above which defines $\frk e_0$ from $e_0$.

From the description of the sewing procedure we see that, after omitting the steps at which $\cL$ stays constant, $\cL$ is a contour function for the infinite tree $\tree$ in the sense that it gives the relative distance from the various vertices to the root of the tree $\tree$ at infinity. Similarly, $\cR$ is a contour function for the infinite tree $\tree_*$.

The sewing procedure is uniquely characterizing the coupling between $\cZ$ and $(M,e_0,\tree)$. In fact, the procedure for constructing $(M,e_0,\tree)$ works for any bi-infinite walk $\cZ = (\cL, \cR)$ in $\BB Z^2$ with nearest neighbor steps such that 
\eqb
\liminf_{\ti\rta -\infty} \cL_{\ti} 
= \liminf_{\ti\rta  \infty} \cL_{\ti} 
= \liminf_{\ti\rta -\infty} \cR_{\ti} 
= \liminf_{\ti\rta \infty} \cR_{\ti} 
= -\infty .
\label{rmk:other-walks}
\eqe
We can define functions $\PhiM$ and $\PhiT$ on the set of nearest neighbor bi-infinite walks satisfying \eqref{rmk:other-walks} such that $\Phi_{\op{M}}(\cZ)=(M,e_0,\tree)$ and $\Phi_{\op{T}}(\cZ)=(T,\frk{tr})$. We also define a map $\Phi_{\op{T}\to\op{M}}$ on the range of $\Phi_{\op{T}}$ such that $\Phi_{\op{M}}=\Phi_{\op{T}\to\op{M}} \circ\Phi_{\op{T}}$.

There is a natural correspondence between the vertex and edge sets $\cV(\triang)$, $\cV(M)$, $\cE(\triang)\setminus(\frk T\cup\frk T_*)$ and the steps of the walk $\BB Z$. More precisely, the map $\ti\mapsto \frk e_{\ti}$ is a bijection between $\BB Z$ and $\cE(\triang)\setminus(\frk T\cup\frk T_*)$ and the map $\ti\mapsto \frk v_{\ti}$ is a bijection between $\cJm$ (resp.\ $\cJLm$) and $\cV(\triang)$ (resp.\ $\cV(M)$). For $v\in \cV(T)$ let $\frk t_v \in \cJm$ be such that the map $v\mapsto \frk t_v$ is the inverse of $\ti\mapsto \frk v_{\ti}$. Since we view $\cV(M)$ as a subset of $\cV(T)$, the time $\frk t_v$ is also well-defined for $v\in \cV(M)$.

\begin{figure}[ht!]
	\begin{center}
		\includegraphics[scale=.8]{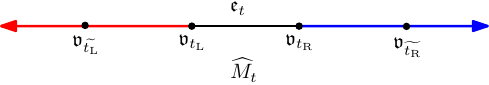} 
		\caption{\label{fig-mullin-sewing2} Illustration of Lemma \ref{prop-bij-steps}.}
	\end{center}
\end{figure}

The following lemma is immediate by the sewing procedure and the geometric meaning of $\cL$ and $\cR$ as contour functions. By the left (resp.\ right) end-point of $\frk e_t$ we mean the end-point lying on its left (resp.\ right) side in the sewing procedure, assuming the external face of $\wh M_t$ lies above its boundary as in Figures \ref{fig-mullin-sewing} and \ref{fig-mullin-sewing2}.
\begin{lem} 
	If $t\in\Z$ then the left (resp.\ right) end-point of $\frk e_t$ equals $\frak v_{\tL}$ (resp.\ $\frak v_{\tR}$), where 
	$$
	\tL=\inf\{s>t\,:\, \cL_s=\cL_t-1 \},\qquad
	\tR=\inf\{s>t\,:\, \cR_s=\cR_t-1 \}.
	$$ 
	Furthermore, if $\frk v_{\wt{\tL}}$ (resp.\ $\frk v_{\wt{\tR}}$) is the vertex on the frontier of $\wh M_t$ to the left (resp.\ right) of $\frk v_{{\tL}}$ (resp.\ $\frk v_{{\tR}}$), then
	$$
	\wt{\tL}=\inf\{s>t\,:\, \cL_s=\cL_t-2 \},\qquad
	\wt{\tR}=\inf\{s>t\,:\, \cR_s=\cR_t-2 \}.
	$$ 
	\label{prop-bij-steps}
\end{lem}

\subsubsection{Mullin bijection for finite walks and planar maps with disk topology}
\label{sec:bijection-disk}
The previous section describes a function $\PhiM=\PhiTM\circ\PhiT$ from bi-infinite nearest-neighbor walks on $\BB Z^2$ satisfying \eqref{rmk:other-walks} to planar maps with a spanning tree. In this section we will extend this map so that the input walk $\cZ=(\cZ_{\ti})_{\ti\in\{ 0,\dots,n \}}$ may be an arbitrary nearest-neighbor walk on $\BB Z^2$ of finite length. Special cases of $\PhiM$ give bijections between certain classes of finite-length walks and natural families of tree-decorated planar maps; see Propositions \ref{prop:mullin} and \ref{prop:mullin-wired-disk2} and Lemma \ref{prop:mullin-wired-disk} below. Our definition of $\PhiM(\cZ)$ for general finite-length $\cZ$ will be needed in Section \ref{sec:tutte-conv}, e.g.\ to define the map $M_n$ above Definition \ref{def:rw-conv}, while the special case in Proposition \ref{prop:mullin-wired-disk2} is relevant since it concerns the planar map appearing in Theorem \ref{thm-tutte-conv0}.

\begin{figure}%
	\centering
		\includegraphics[scale=.8]{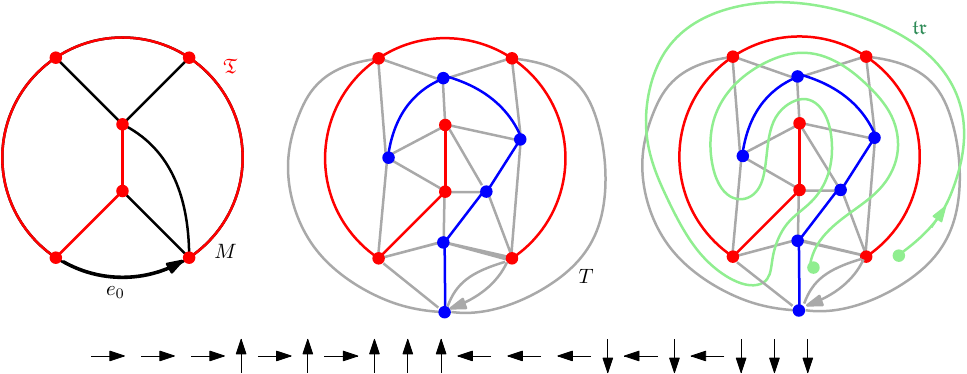} 
		\caption{
		Illustration of $(M,e_0,\frk T)$, $(T,\frk{tr})$ and $\cZ$ in the context of Lemma \ref{prop:mullin-wired-disk}, which is a special case of Proposition \ref{prop:mullin}. 
		{\bf Left}: A planar map $(M,e_0)$ with disk topology. The set $\frk T$ is a spanning tree of $M$ and $\frk T\cup\{ e_0\}$ defines a wired spanning tree on $M$. 
		{\bf Middle}: The planar map $T$. 
		{\bf Right}: The path $\frk{tr}$ is shown in green. The two green dots indicate the starting and ending triangles.
		{\bf Bottom}: The steps of the walk associated with $(M,e_0,\frk T)$ via the bijection in Lemma \ref{prop:mullin-wired-disk}. Observe that we recover the triangulation in the middle figure from the walk by doing the sewing procedure.}
\label{fig-mullin-disk} 
\end{figure}

First we consider the case where $\cZ=(\cZ_{\ti})_{\ti\in\{ 0,\dots,n \}}$ is such that for any $\si\in\{1,\dots,n \}$,
\eqb
\cL_{\si}\geq \cL_{0},\qquad
\cR_{\si}\geq \cR_{0}.
\label{eq:walk-pos}
\eqe
In this case we can apply the sewing procedure exactly as in Section \ref{sec:mullin-infinite}. See Figure %
\ref{fig-mullin-disk}.
More precisely, we initialize the bijection from a single edge $\frk e_0$; then, when reading the first step $\cZ_1-\cZ_0$ of the walk we get a map $\wh M_1$ by adding a triangle $\frk{tr}(1)$ as in the sewing procedure, and so on. The condition \eqref{eq:walk-pos} guarantees that every time we read a step $(-1,0)$ (resp.\ $(0,-1)$) there is a red (resp.\ blue) edge on our left (resp.\ right) side which can serve as an edge for the triangle we add. After we have read the $\ti$th step of the walk we have a map $\wh M_{\ti}$ where all the faces have three edges except a distinguished face with boundary length $(\cL_{\ti}-\cL_0)+(\cR_{\ti}-\cR_0)+2$. The two edges $\frk e_0$ and $\frk e_{\ti}$ lie on this face, in addition to $\cL_{\ti}-\cL_0$ (resp.\ $\cR_{\ti}-\cR_0$) edges which are part of the primal tree $\frk T$ (resp.\ the dual tree $\tree_*$). We refer to these $\cL_{\ti}-\cL_0$ (resp.\ $\cR_{\ti}-\cR_0$) edges as the left (resp.\ right) boundary of $\wh M_{\ti}$. Set $\PhiT(\cZ)=T=\wh M_n$.

If $\cZ$ satisfies $\cZ_0=\cZ_n=(0,0)$ in addition to \eqref{eq:walk-pos}, we get $(M,e_0,\tree)=\PhiM(\cZ)=\PhiTM(T,\frk{tr})$ from $T$ and $\frk{tr}$ similarly as before, namely by letting $\cV(M)$ be the set of vertices in $\tree$ and by letting $\cE(M)$ be the union of $\tree$ and the edges which are dual to $\tree_*$. Define $e_0$ as in the infinite-volume case by choosing the first edge of $M$ in clockwise order around the starting vertex of the root edge of $T$. We define $\PhiM$ and $\PhiTM$ for general finite-length walks $\cZ$ below Lemma \ref{prop:mullin-wired-disk}.

The following lemma follows from the Mullin bijection as presented in e.g.\ \cite{shef-burger}, although we remark that we use a slightly different convention for the assignment of a root to $T$. See Figure \ref{fig-mullin-disk}. 
\begin{prop}
	For any $n\in\N$ the map $\PhiM$ is a bijection between the following sets:
	\begin{compactitem}
		\item The set of nearest-neighbor walks $\cZ=(\cZ_{\ti})_{\ti\in\{0,\dots,2n \}}$ on $\BB Z^2$ satisfying \eqref{eq:walk-pos} such that $\cZ_0=\cZ_{2n}=(0,0)$. 
		\item The set of tuples $(M,e_0,\frk T)$, where $(M,e_0)$ is a rooted planar map with $n$ edges and $\frk T$ is a spanning tree on $M$.
	\end{compactitem}
\label{prop:mullin}
\end{prop}
In the special case considered in the proposition we can recover $\cZ$ from $(M,e_0,\frk T)$ by first constructing a map $T$, a dual tree $\frk T_*$ and a path $\frk{tr}$ in the dual map similarly as in Section \ref{sec:mullin-infinite}, and then letting $\cZ$ be the contour functions associated with the two trees, parametrized in the same way as $\frk{tr}$. More precisely, $T$ is the map whose vertex set is $\cV(M)\cup\cV(M_*)$ and edge set given by the union of $\tree\cup\tree_*$ and edges of the form $\{v,v^* \}$ where $v\in\cV(M)$ is incident to the face that is dual to $v^*\in\cV(M)$. We include one copy of each edge, except the edge $\{v,v^*\}$ whose primal endpoint $v$ coincides with the terminal endpoint of $e_0$ and which is the first edge in counterclockwise order after $e_0$ with this property; of this edge we include two copies and we declare that the second of the two edges in counterclockwise order after $e_0$ is the root of $T$, oriented from $v$ to $v^*$. Observe that $T$ is a triangulation except that the face to the right of its root edge only has degree 2. In other literature one typically does not include two copies of the root edge, which means that $T$ is a true triangulation. We decide to change the convention slightly so that $T$ is the map which arises in the sewing procedure for arbitrary finite-length nearest-neighbor walks on $\BB Z^2$. See Figure \ref{fig-mullin-disk} for an example.

The following lemma is a special case of the setting in Proposition \ref{prop:mullin} and illustrates that the restriction of $\PhiM$ to a special set of walks gives the tree-decorated planar maps considered in Theorem \ref{thm-tutte-conv0}. 
\begin{lem}
	For $\ell\in\N$ and $n\in\{2,3,\dots \}$ the map $\PhiM$ is a bijection between the following sets:
	\begin{compactitem}
		\item The set of nearest-neighbor walks $(\cZ_{\ti})_{\ti\in\{0,\dots,2n \}}$ on $\BB Z^2$ which satisfy the following: 
		\begin{compactitem}
			\item[(a)] $\cZ_{0}=\cZ_{2n}=(0,0)$,
			\item[(b)] $\cZ_{\ti}-\cZ_{\ti-1}=(1,0)$ for $\ti=1,\dots\ell$, 
			\item[(c)] $\cZ_{\ell+1}-\cZ_\ell=(0,1)$ and $\cZ_{2n}-\cZ_{2n-1}=(0,-1)$, 
			\item[(d)] $\cL_{\ti}\geq 0$ for $\ti=\ell+1,\dots,2n-1$, and
			\item[(e)] $\cR_\ti\geq 1$ for $\ti=\ell+1,\dots,2n-1$.
		\end{compactitem}
		\item The set of tuples $(M,e_0,\frk T)$, where $(M,e_0)$ is a planar map with disk topology of boundary length $\ell+1$ with $n$ edges and $\frk T$ is a spanning tree on $M$ such that $\partial \cE(M)\setminus\{ e_0\}\subset \frk T$.
	\end{compactitem}
\label{prop:mullin-wired-disk}
\end{lem}

\begin{figure}
	\centering
		\includegraphics[scale=1]{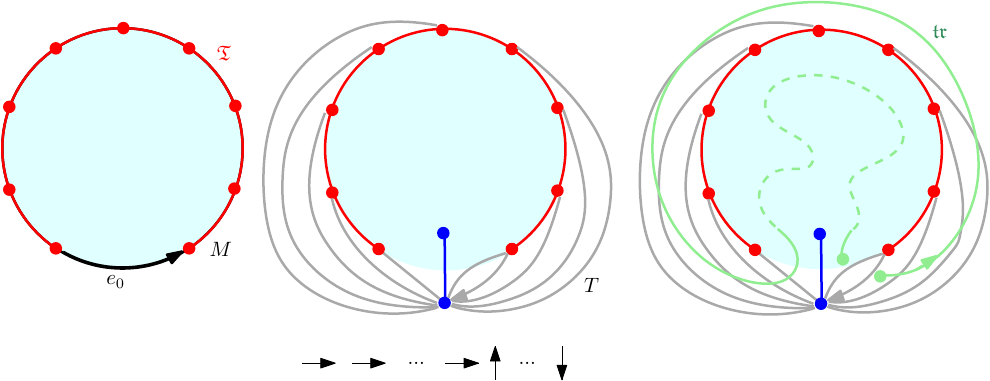} 
		\caption{Illustration of Lemma \ref{prop:mullin-wired-disk}. The domains with light blue background contain non-specified edges and vertices. The blue edge in the middle and right figures represents the edge of $\tree_*$ which is dual to $e_0$. In the walk at the bottom of the figure the first ``$\cdots$'' corresponds to $(1,0)$ steps such that the walk has $\ell$ initial $(1,0)$ steps, while the second ``$\cdots$'' corresponds to some general walk satisfying condition (d) and (e) of the lemma.}
	\label{fig-mullin-disk2}
\end{figure}

\begin{proof}
	Given a tuple $(M,e_0,\frk T)$ as in the second bullet point above we get a walk $\cZ$ as in the first bullet point by the procedure described below Proposition \ref{prop:mullin}; see Figures \ref{fig-mullin-disk} and \ref{fig-mullin-disk2}. Conditions (a) and (d) are satisfied since this lemma considers a special case of the setting in Proposition \ref{prop:mullin}. The steps described in (b) correspond to the faces of $T$ lying in the exterior face of $M$. The faces $\frk{tr}(\ell+1)$ and $\frk{tr}(2n)$ are incident to the blue edge in Figure \ref{fig-mullin-disk2}; this leads to conditions (c) and (e). Conversely, applying the sewing procedure shows that walks as in the first bullet point produce tree-decorated maps as in the second bullet point. 
\end{proof}

\begin{figure}
	\centering
	\includegraphics[scale=.8]{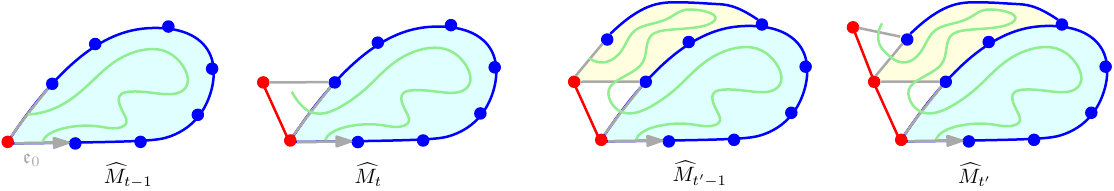} 
	\caption{Illustration of $\PhiT$ for a walk which satisfies the first condition of \eqref{eq:walk-pos2} at the times $\ti$ and $\ti'$. Edges and vertices in the interior of the blue and yellow regions are not shown. The two triangles on white background are $\frk{tr}(\ti)$ and $\frk{tr}(\ti')$. 
	}
\label{fig:bijection-general-walks} 
\end{figure}

Finally we will extend the functions $\PhiM$ and $\PhiT$ to finite nearest-neighbor walks $\cZ=(\cZ_{\ti})_{\ti\in\{0,\dots,n \}}$ on $\BB Z^2$ which do not necessarily satisfy the constraint \eqref{eq:walk-pos}.  We will do this construction with the sewing procedure as before, but we need to introduce a new kind of step for times $\ti\in\{1,\dots,n \}$ such that 
\eqb
\cL_{\ti} < \min_{0\leq \si< \ti} \cL_{\si}
\qquad\text{or}\qquad
\cR_{\ti} < \min_{0\leq \si< \ti} \cR_{\si}.
\label{eq:walk-pos2}
\eqe
See Figure \ref{fig:bijection-general-walks} for an illustration. Suppose $\ti\in\{1,\dots,n \}$ is such that the first condition of latter display holds, i.e., 
$\cL_\ti=\min_{0\leq \si<\ti}\cL_{\si}-1=\cL_{\ti-1}-1$. When doing the $\ti$th step of the sewing procedure we add a triangle $\frk{tr}(\ti)$ such that the edges of $\frk{tr}(\ti)$ are $\frk e_{\ti-1}$ and two edges in the external face of $\wh M_{\ti-1}$; the left of these edges will be part of $\tree$ and the right of these edges equals $\frk e_{\ti}$. In other words, we do the same as in the fourth panel, from the left, of Figure \ref{fig-mullin-sewing}, except that the red edge on the boundary of $\frk{tr}(\ti)$ is not on the boundary of $\wh M_{\ti-1}$; instead, this edge is in the external face of $\wh M_{\ti-1}$. We proceed similarly if the second condition of \eqref{eq:walk-pos2} is satisfied. When the sewing procedure has finished we have a map $T=\wh M_n$ and we set $\PhiT(\cZ)=T$. The left boundary of $T$ has length 
$\cL_0+\cL_n-2\inf_{0\leq t\leq n}\cL_t$ and the right boundary of $T$ has length 
$\cR_0+\cR_n-2\inf_{0\leq t\leq n}\cR_t$. 

\begin{figure}
	\centering
	\includegraphics{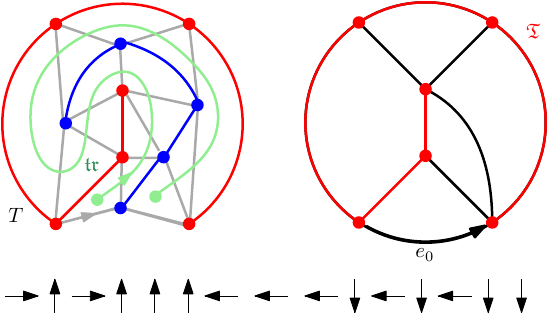}
	\caption{Illustration of Proposition \ref{prop:mullin-wired-disk2}. The figure shows the same tree-decorated planar map $(M,e_0,\frk T)$ as in Figure \ref{fig-mullin-disk}, but $(T,\frk {tr})$ and the walk $\cZ$ are different since we have applied the bijection from Proposition \ref{prop:mullin-wired-disk2} instead of the bijection from Lemma \ref{prop:mullin-wired-disk}.}
	\label{fig:T-M}
\end{figure}

We define $\PhiM$ and $\PhiTM$ similarly as before. Namely, we get $(M,e_0,\tree)=\PhiM(\cZ)=\PhiTM(T,\frk{tr})$ from $T$ and $\frk{tr}$ by letting $\cV(M)$ be the set of vertices in $\tree$ and letting $\cE(M)$ be the union of the following three sets:
    (i) $\tree$, 
	(ii) the edges which are dual to $\tree_*$, excluding the edges on the right boundary of $T$, and 
	(iii) if the left boundary has length at least one we add a new edge connecting the two end-points of $\frk e_0$ and $\frk e_{n}$ which are part of $\frk T$. The reason we decide to exclude the edges on the right boundary of $T$ in (ii) is that these edges have no naturally defined dual edge. The reason we add the new edge in (iii) is in order for Proposition \ref{prop:mullin-wired-disk2} below to hold. If we do add an edge in (iii) then this will be the root edge $e_0$ on $M$, oriented from the end-point of $\frk e_0$ to the end-point of $\frk e_{n}$ (see e.g.\ Figure \ref{fig:T-M}); otherwise we define $e_0$ as in the infinite-volume case by choosing the first edge of $M$ in clockwise order around the starting vertex of the root edge of $T$ (see e.g.\ Figure \ref{fig-mullin-disk}). The boundary of the face to the right of $e_0$ is called the boundary of the map.

The definition of the sewing procedure in the previous paragraph is motivated by the following lemma, whose proof is immediate from the construction. See also the argument in \cite[Section 3.1]{ghs-map-dist}, which is using closely related ideas. For some walk $\cZ=(\cZ_{\ti})_{\ti\in \{\ti_0,\dots,\ti_0+n \}}$ we define $\op{shift}(\cZ)$ to be the walk indexed by $\{0,\dots,n \}$ given by $\ti\mapsto \cZ_{\ti_0+\ti}-\cZ_{\ti_0}$. In other words, $\op{shift}(\cZ)$ is a shifted and translated version of $\cZ$ which starts from $(0,0)$ at time 0. For arbitrary finite nearest-neighbor walks $\cZ$ on $\Z^2$ define $\PhiT(\cZ):=\PhiT(\op{shift}(\cZ))$.

\begin{lem} %
	Suppose one of the following holds:
	\begin{compactitem}
		\item[(a)] $\cZ$ is a bi-infinite nearest-neighbor walk on $\BB Z^2$ satisfying \eqref{rmk:other-walks} and $I=\{t_0,t_0+1,\dots,t_0+n \}\subset\Z$ for some $t_0\in\Z$ and $n\in\N$, or
		\item[(b)] $\cZ$ is a nearest-neighbor walk on $\BB Z^2$ indexed by some finite set $I_0\subset \Z$ and $I=\{t_0,t_0+1,\dots,t_0+n \}\subset I_0$ for some $t_0\in\Z$ and $n\in\N$.
	\end{compactitem}
	Then $\PhiT(\cZ|_I)$ is the submap of $\PhiT(\cZ)$ with edge set and vertex set given precisely by the set of edges and vertices, respectively, lying on the boundary of some triangle $\frk{tr}(\ti)$ for $\ti\in I\setminus\{t_0 \}$.
	\label{prop:general-walk}
\end{lem}

Finally we present the bijection which will be used to prove Theorem \ref{thm-tutte-conv0} in Section \ref{sec:tutte-conv}. The set of tree-decorated planar maps is the same as considered in Lemma \ref{prop:mullin-wired-disk}, but the set of walks is different. See Figure \ref{fig:T-M}. 
In particular, we see that the map $\PhiM$ is not injective. Notice that in Theorem \ref{thm-tutte-conv0} all the boundary edges are assumed to be part of $\frk T$ and that we get to the set of tree-decorated maps described below by removing the root edge from $\frk T$.
\begin{prop}
	For any $\ell,n\in\N$ the map $\PhiM$ defined above is a bijection between the following sets:
	\begin{compactitem}
		\item The set of nearest-neighbor walks $(\cZ_{\ti})_{\ti\in\{\ell+1,\dots,2n-1 \}}$ on 
		$\BB N^2_0$ from $(\ell,0)$ to $(0,0)$.
		\item The set of tuples $(M,e_0,\frk T)$ as described in the second bullet point of Lemma \ref{prop:mullin-wired-disk}.
	\end{compactitem}
\label{prop:mullin-wired-disk2}
\end{prop}
\begin{proof}
	There is a natural bijection $\Psi$ 
	between 
	the walks in the first bullet point of Lemma \ref{prop:mullin-wired-disk} and 
	the walks in the first bullet point of this proposition, namely the map which restricts the  walks in Lemma \ref{prop:mullin-wired-disk} to $\{\ell+1,\dots,2n-1 \}$ and subtracts $(0,1)$. By case (b) of Lemma \ref{prop:general-walk} and the definition of $\PhiTM$ we see that $\PhiM(\cZ)=\PhiM\circ\Psi(\cZ)$ for walks $\cZ$ as in Lemma \ref{prop:mullin-wired-disk}, which concludes the proof. 
\end{proof}

\subsection{A variant of the Gwynne-Miller-Sheffield invariance principle} 
\label{app1}

In this section we state and prove a variant of~\cite[Theorem 3.10]{gms-random-walk}, which gives a condition under which a random walk on the adjacency graph of a random collection of cells in the plane converges to 2D Brownian motion. The following definition makes precise what we mean by an ``adjacency graph of cells". 

\begin{defn} \label{def-cell-config}
A \emph{cell configuration} on $\BB C$ consists of the following objects.
\begin{enumerate}
	\item A locally finite collection $\mcl H$ of compact connected subsets of $\BB C$ (called ``cells'') with non-empty interiors whose union is all of $\BB C$ and such that the intersection of any two elements of $\mcl H$ has zero Lebesgue measure.
	\item A symmetric relation $\sim$ on $\mcl H\times \mcl H$ such that  
	\eqb
	\text{if $H\sim H'$, then $H\cap H'\not=\emptyset$ and $H\not=H'$}.
	\label{eq:cell-conf-adjacency}
	\eqe 
	\item A function $\frk c$ (the conductances) from the set of pairs $(H,H')$ with $H\sim H'$ to $(0,\infty)$ such that $\frk c(H,H') = \frk c(H',H)$. 
\end{enumerate}
\end{defn}
We view $\mcl H$ as a weighted graph whose vertices are the cells of $\mcl H$ and whose edge set is $\{ \{H,H'\} \in \mcl H\times \mcl H : H\sim H' \}$, such that the edge $\{H,H'\}$ is assigned weight $\frk c(H,H')$.  For any $A\subset\BB C$ define
\eqb \label{eqn-cell-restrict}
\mcl H(A) := \left\{H\in\mcl H : H\cap A \not=\emptyset \right\}
\eqe 
and then define a metric on the space of cell configurations by
\eqbn
\begin{split}
&\BB d^{\op{CC}}(\mcl H,\mcl H') := \int_0^\infty e^{-r} \wedge \inf_{f_r} \big\{ \sup_{z\in \BB C} |z - f_r(z)|   \notag \\
&\qquad \qquad   + \max \{ |\frk c(H_1,H_2) - \frk c'(f_r(H_1) , f_r(H_2) ) |\,:\, H_1\sim H_2, H_j\in \cH( B(0,r)) \} \big\}  \,dr,
\end{split}
\eqen
where each of the infima is over all homeomorphisms $f_r : \BB C\rta \BB C$ such that $f_r$ takes each cell in $\mcl H$ intersecting $B(0,r)$ to a cell in $\mcl H'$ intersecting $B(0,r)$ and preserves the adjacency relation, and $f_r^{-1}$ does the same with $\mcl H$ and $\mcl H'$ reversed. We interpret the maximum of an empty set of nonnegative discrepancies as $0$, and an infimum over no admissible homeomorphisms as $\infty$. When we work with random cell configurations we view them as random variables with values in the space of cell configurations equipped with the Borel $\sigma$-algebra generated by this metric. 
We also define the primal and dual stationary measures on $\mcl H$ by
\eqb \label{eqn-stationary-measure}
\pi(H) := \sum_{\substack{H'\in\mcl H : H'\sim H}} \frk c(H,H') \quad \op{and} \quad
\pi^*(H) := \sum_{\substack{H'\in\mcl H : H'\sim H}} \frac{1}{ \frk c(H,H') }. 
\eqe

We will work with random cell configurations satisfying the following.  
\begin{enumerate}
\item \textbf{Translation invariance modulo scaling.} There is a (possibly random and $\mcl H$-dependent) increasing sequence of open sets $U_j \subset \BB C$, each of which is either a square or a disk, whose union is all of $\BB C$ such that the following is true. Conditional on $\mcl H$ and $U_j$, let $z_j$ for $j\in\BB N$ be sampled uniformly from Lebesgue measure on $U_j$. Then the shifted cell configurations $\mcl H - z_j$ converge in law to $\mcl H$ modulo scaling as $j\rta\infty$, i.e., there are random numbers $C_j >0$ (possibly depending on $\mcl H$ and $z_j$) such that $ C_j(\mcl H-z_j) \rta \mcl H$ in law with respect to the metric $\BB d^{\op{CC}}$. Several equivalent definitions are given in \cite[Definition 1.2]{gms-random-walk}. \label{item-resampling}
\item \textbf{Ergodicity modulo scaling.} Let $F = F(\mcl H)$ be a real-valued measurable function which is invariant under translation and scaling, i.e., it holds a.s.\ that
\eqb
F(C(\mcl H-z)) = F(\mcl H),\qquad\forall z\in\BB C,C>0. 
\label{eq:transl3}
\eqe 
Then $F$ is a.s.\ equal to a deterministic constant.
\label{item-ergodic} 
\item \textbf{Finite expectation.} With $H_0$ the cell containing 0 and $\pi$ and $\pi^*$ as in~\eqref{eqn-stationary-measure},  \label{item-moment}
\eqb\label{eqn-moment}
\BB E\left[ \frac{\op{diam}(H_0)^2}{\op{area}(H_0)} \pi(H_0)     \right] <\infty  \quad\op{and} \quad 
\BB E\left[ \frac{\op{diam}(H_0)^2}{\op{area}(H_0)} \pi^*(H_0)     \right] <\infty 
\eqe 
where $\op{diam}$ and $\op{area}$ denote Euclidean diameter and Lebesgue measure, respectively. 
\end{enumerate}

Several equivalent formulations of hypothesis~\ref{item-resampling} are given in~\cite[Definition 1.2]{gms-random-walk}, but we state only the formulation which we will use in the present paper. In order to state the next theorem we need to define a variant of \eqref{eq:curves-mod-reparam} for curves of infinite duration. Let $\beta_1,\beta_2:[0,\infty)\to\C$ be curves and for $r>0$ and $j=1,2$ let $T_{j,r}$ be the exit time of $\beta_j$ from the ball $B(0,r)$ (or 0 if the curve starts outside $B(0,r)$). If $\beta_1,\beta_2$ are such that $T_{1,r},T_{2,r}<\infty$ for all $r$, set
\eqb 
d^{\rm loc}_{\rm CMP}(\beta_1,\beta_2):=\int_1^\infty e^{-r}\Big( 1 \wedge
d_{\rm CMP}(
\beta_1|_{[0,T_{1,r}]},
\beta_2|_{[0,T_{2,r}]}
)
\Big) \,dr,
\label{eq:curves-mod-reparam-inf}
\eqe
where we define $d_{\rm CMP}(
\beta_1|_{[0,T_{1,r}]},
\beta_2|_{[0,T_{2,r}]}
)=\infty$ if 
$\beta_1|_{[0,T_{1,r}]}$ and/or $\beta_2|_{[0,T_{2,r}]}$ have infinite duration. We identify two curves with zero $d_{\rm CMP}^{\rm loc}$ distance.
\begin{theorem}\label{thm-general-clt-uniform}
Let $\cH^{\op{cell}}$ and $\cH^{\op{large}}$ be random cell configurations satisfying the above three hypotheses in addition to the following a.s.:
\begin{enumerate}
\item[(i)] There is a bijection $\phi:\cH^{\op{cell}}\to\cH^{\op{large}}$ such that $\phi(H)\sim \phi(H')$ if and only if $H\sim H'$. Furthermore, if $H\sim H'$ then $\frk c_{\op{cell}}(H,H')=\frk c_{\op{large}}(\phi(H),\phi(H')) $.
\item[(ii)] The convex hull of $H$ is contained in the convex hull of $\phi(H)$ for each $H\in\cH^{\op{cell}}$.
\item[(iii)] For sufficiently large $r>1$, for each horizontal or vertical line segment $\ell\subset\C$ of length at least $r$ and at distance at most $r$ from the origin the following holds. For any two cells $H,H'\in\cH^{\op{cell}}$ intersecting $\ell$ we can find a path $H\sim H_1\sim\dots\sim H_k\sim H'$ such that $\phi(H),\phi(H_{1}),\phi(H_2),\dots,\phi(H_k),\phi(H')$ all intersect $\ell$.
\end{enumerate}
For $z\in\BB C$, let $Y^z:\N_0\to\cH^{\op{cell}}$ denote the random walk on $\cH^{\op{cell}}$ with conductances $\frk c_{\op{cell}}$ started from a cell containing $z$, chosen arbitrarily when there is more than one such cell. For $j\in\BB N_0$, let $\wh Y_j^z$ be an arbitrarily chosen point of the cell $Y_j^z$ and extend $\wh Y^z$ from $\BB N_0$ to $[0,\infty)$ by piecewise linear interpolation. 
There is a deterministic covariance matrix $\Sigma$ with $\det\Sigma\not=0$ such that the following is true. For each fixed compact set $A \subset \BB C$, it is a.s.\ the case that as $\ep \rta 0$, the maximum over all $z\in A$ of the Prokhorov distance between the conditional law of $\ep \wh Y^{z/\ep}$ given $\cH^{\op{cell}}$ and the law of Brownian motion started from $z$ with covariance matrix $\Sigma$, with respect to the topology on curves modulo time parameterization (as defined in \eqref{eq:curves-mod-reparam-inf}), tends to 0.
\end{theorem}

The difference as compared to \cite[Theorem 3.10]{gms-random-walk} is that they do not consider $\cH^{\op{large}}$ and instead they require $\cH^{\op{cell}}$ to satisfy the following property (called ``connectivity along lines''): Almost surely, for each horizontal or vertical line segment $\ell \subset \BB C$, the subgraph of $\mcl H^{\op{cell}}$ induced by the set of cells which intersect $\ell$ is connected. Note that \cite[Theorem 3.10]{gms-random-walk} is a special case of the above theorem since we recover \cite[Theorem 3.10]{gms-random-walk} by setting $\cH^{\op{large}}=\cH^{\op{cell}}$.

The proof of Theorem \ref{thm-general-clt-uniform} is almost identical to the proof of \cite[Theorem 3.10]{gms-random-walk}. The only place where line connectivity is used is in the proof of \cite[Lemma 2.19]{gms-random-walk}. In our setting we get from Assumptions (i)-(iii) that if $\ell$ is as in Assumption (iii), %
$H_1,H_2\in\cH^{\op{cell}}(\ell)$, and $g:\cH^{\op{large}}\to\R$, then
$$
|g(\phi(H_1))-g(\phi(H_2))| \leq \sum_{ H,H'\in\cH^{\op{large}}(\ell)\,:\,H\sim H' } |g(H)-g(H')|.
$$
The sum on the right side can be bounded exactly as in \cite{gms-random-walk} since $\cH^{\op{large}}$ satisfies the assumptions of that paper (except the line connectivity, which is not needed). Note that it is sufficient to consider line segments $\ell$ as described in Assumption (iii) due to the particular type of line segments considered in the proof of \cite[Lemma 2.19]{gms-random-walk} and since this lemma concerns an a.s.\ limit as the line lengths go to infinity.

\section{Strong coupling of Brownian motion and random walk}
 \label{sec:kmt}
In this section we will construct a strong coupling between higher-dimensional Brownian motion and a continuous-time random walk. The main result is Theorem \ref{prop-kmt} below. The construction is an infinite-volume higher-dimensional variant of a construction of Pollard \cite{pollard}. For any $n\in\N$, Pollard gives a coupling of the following two processes: 
\begin{compactitem}
	\item a Brownian bridge $(\check Y_t)_{t\in[0,1]}$, i.e.\ a standard Brownian motion restricted to $[0,1]$ conditioned to satisfy $\check Y_1=0$, and
	\item a process $(\check\cY_t)_{t\in [0,1]}$ of the form $\check\cY_t=n^{-1/2}(\#([0,t]\cap A)-nt)$, where $A$ is a set of $n$ uniform and independent points in $[0,1]$.
\end{compactitem}
Let us briefly review some key ideas of Pollard's coupling; see Section \ref{sec:kmt1} for more details. Given a realization of $\wt Y$ he constructs $\wt{\cY}$ as follows. He determines the number of points $\#(A\cap J_{i,k})$ in dyadic intervals $J_{i,k}=[i2^{k},(i+1)2^{k})$ iteratively for decreasing $k\in\Z_\leq=\{\dots,-2,-1,0  \}$. Before starting the iteration he is given $\#(A\cap J_{0,0})=\#A=n$. In the first step of the iteration he determines $\#(A\cap J_{0,-1})$ and $\#(A\cap J_{1,-1})$, then he determines 
$\#(A\cap J_{0,-2})$, $\#(A\cap J_{1,-2})$, $\#(A\cap J_{2,-2})$ and $\#(A\cap J_{3,-2})$, and so on. By sending $k\rta-\infty$ he will eventually know $A$, and $A$ determines $\check{\cY}$. Notice that $J_{i,k}=J_{2i,k-1}\cup J_{2i+1,k-1}$ for all $k\in\Z_\leq$ and $i\in\{0,2^{-k}-1\}$.  
Given $\#(A\cap J_{i,k})$ he determines $\#(A\cap J_{2i,k-1})$ and $\#(A\cap J_{2i+1,k-1})$ (which much sum to $\#(A\cap J_{i,k})$ and have the law of binomial random variables of parameters $\#(A\cap J_{i,k})$ and $0.5$ conditioned on $\#(A\cap J_{i,k})$) by considering the increments  
$\wt Y_{(2i+1)2^{k-1}}-\wt Y_{i2^k}$ 
and 
$\wt Y_{(i+1)2^{k}}-\wt Y_{(2i+1)2^{k-1}}$ of $\wt Y$ during the intervals $J_{2i,k-1}$ and $J_{2i+1,k-1}$, respectively. The error in the coupling is bounded by using Tusn\'ady's lemma (see Lemma \ref{prop:tus} below), which gives a coupling of a binomial random variable and a normal random variable of the same mean and variance.

In Section \ref{sec:kmt1} we present a variant of Pollard's coupling in the case where $Y$ is a two-sided Brownian motion and $(\cY_t)_{t\in\R}=(P_t-t)_{t\in\R}$ for $(P_t)_{t\in\R}$ a two-sided Poisson process (see below for the precise definition). Our construction can be viewed as the local limit of Pollard's construction when we zoom in near any fixed time at the same time as $n\rta\infty$. This construction is not straightforward since we cannot initialize the iterative construction of Pollard due to 
the fact that the number of jumps of $P$ in any given dyadic interval is random. A key property of the coupling is that $\cY$ can be determined almost locally by $Y$ in the sense that  $\cY|_{[-\Delta,\Delta]}$ is determined by $Y|_{[-\Delta^{1+\eps},\Delta^{1+\eps}]}$ with high probability for  any fixed $\eps>0$ and $\Delta\gg 1$; see Proposition \ref{prop-kmt0} below for a stronger and more precise statement. Another important property stated in the cited proposition is translation invariance in law.

Earlier literature contains many variants of the KMT coupling \cite{kmt}, i.e., couplings between a Brownian motion and a random walk. In many of these couplings the random walk is indexed by the non-negative integers. One key difference in our setting is that our discrete process jumps at the jump times of an independent Poisson point process, i.e., the times between consecutive jump times are given by i.i.d.\ exponential random variables instead of being equal to 1. This makes it possible to extend the coupling to the higher-dimensional case where the steps of the walk are sampled independently from some finite step distribution. We extend the coupling to this case in Section \ref{sec:kmt3}.

\subsection{One-dimensional process on compact interval}
\label{sec:kmt0}

Let $k_0\in\Z$, $n\in\N$ and $p>0$. In this section we present a coupling between the following processes:
\begin{compactitem}
	\item A process $\wt{Y}:[0,2^{k_0}]\to\R$ given by 
	a variance $p$ Brownian bridge on $[0,2^{k_0}]$. 
	\item A process $\wt{\cY}:[0,2^{k_0}]\to\R$ of the form 
	$$
	\wt{\cY}_t= \#(A\cap[0,t])-pt,
	$$
	where $A$ has the law of $n$ uniform points sampled from $[0,2^{k_0}]$. 
\end{compactitem}
We will typically have $n\approx p2^{k_0}$ and in this case our coupling is such that $\|\wt Y-\wt{\cY}+t2^{-k_0}(n-p2^{k_0})\|_{\infty}$ is small with high probability, where $\|\cdot\|_\infty$ denotes the supremum norm.
\begin{remark}
		For comparison, Pollard considered the case $k_0=0$ and $p=1$. He proved that in this context $\|\wt Y-n^{-1/2}(\wt{\cY} + (1-n)t )\|_\infty$ is typically not larger than a constant times $n^{-1/2}\log n$. His bound is the main input to the proof of Lemma \ref{prop:pol} below. One key difference between Pollard's setting and our setting is that he keeps $k_0$ constant and sends $n\rta\infty$, while we typically send $n,k_0\rta\infty$ simultaneously such that $n\approx p2^{k_0}$. Another key difference is that he chooses the drift of the discrete process so it ends at 0, while we let the discrete process end at $n-p2^{k_0}$ since this will be more convenient in Section \ref{sec:kmt1} below. However, as you will see in Lemma \ref{prop:pol}, our coupling is equivalent to the one of Pollard and the discrete processes in the two settings can be related via a scaling and the addition of a linear drift term; see \eqref{eq41}.
	\label{rmk:pol} 
\end{remark}

When describing the coupling between $\wt Y$ and $\wt{\cY}$ we will need the following classical result due to Tusn\'ady. See e.g.\ \cite[Chapter 10.5, Lemma 19]{pollard} for a statement. 
\begin{lem}[Tusn\'ady's lemma]
	There exists a (deterministic) function $\tau:\N_0\times\R\to\N_0$ such that the following hold:
	\begin{compactitem}
		\item For any $n\in\N_0$ the function $x\mapsto \tau(n,x)$ is increasing.
		\item If $X\sim\cN(0,1)$ is a standard normal random variable then $\xi:=\tau(n,X)$ has $\op{Binomial}(n,1/2)$ distribution.
		\item The random variable $\xi$ satisfies the inequalities
		$$
		|\xi - \wh X| \leq 1+\frac{X^2}{8},\qquad |\xi-n/2| \leq 1+0.5\sqrt{n}|X|, 
		$$
		where $\wh X=n/2+X\sqrt{n/4}$, which has a $\cN(n/2,n/4)$ distribution.
	\end{compactitem}
	\label{prop:tus}
\end{lem}
Let $\wt Y$ be as above. For $k=\{\dots,k_0-1,k_0 \}$ and $i=\{0,1,\dots,2^{k_0-k}-1 \}$ define $X_{i,k}$ by 
\eqb
X_{i,k}:=2^{-k/2}(2\wt Y_{(2i+1)2^{k-1}}-\wt Y_{i2^{k}}-\wt Y_{(i+1)2^{k}}).
\label{eq:xik}
\eqe
The following lemma is a consequence of \cite{pollard}.  
\begin{lem}
	Suppose $p^{-1}2^{-k_0}n\geq 1/2$. 
	Then there is a coupling of $\wt Y$ and $\wt{\cY}$ such that if $\Delta:=2^{k_0/2}(1-p^{-1}2^{-k_0}n)$,
	\eqb
	\P\big[\|\wt{\cY}-\wt Y-t2^{-k_0}(n-p2^{k_0})\|_\infty>
	4p\Delta^2 
	+ C_1|\Delta| 2^{-k_0/2}\sqrt{nx} 
	+ C_1 \sqrt{p2^{k_0}/n} (x+\log n)\big] \leq C_0\exp(-x),
	\label{eq:pol}
	\eqe 
	where $C_0,C_1>1$ are constants which depend on $p$ but are independent of $k_0$ and $n$. 
	In this coupling, $\wt{\cY}$ is determined by $\wt Y$ and for all $k\in\{ \dots,k_0-2,k_0-1,k_0\},i\in\{0,1,\dots,2^{k_0-k}-1 \}$ we have
	\eqb
	\#(A\cap J_{2i,k-1})=\tau(
	\#(A\cap J_{i,k}),
	p^{-1/2}X_{i,k}), 
	\label{eq:tau-ite}
	\eqe
	where $\tau$ is as in Lemma \ref{prop:tus}.
	\label{prop:pol}
\end{lem}
\begin{proof}
	By \cite[Chapter 10, Theorem 26]{pollard} (see also Remark \ref{rmk:pol}), there is a coupling between $\wt Y$ and the process $\check{\cY}:[0,2^{k_0}]\to\R$ given by 
	\eqb
	\check{\cY}_t
	:= (p^{-1}2^{-k_0}n)^{-1/2}(\#(A\cap[0,t]) -  2^{-k_0}nt )
	= (p^{-1}2^{-k_0}n)^{-1/2}(\wt{\cY}_t + (p-2^{-k_0}n)t)
	\label{eq41}
	\eqe 
	such that
	\eqb
	\P\Big[\|\check {\cY} - \wt Y \|_\infty >C_1(p^{-1}2^{-k_0}n)^{-1/2}(x+\log n)\Big] \leq C_0\exp(-x).
	\label{eq:pol0}
	\eqe
	To be precise, Pollard considered only the case $k_0=0$ and $p=1$ but the case of general $k_0$ and $p$ is immediate via scaling since the processes $(p^{-1/2}2^{-k_0/2}\wt Y_{2^{k_0}t})_{t\in[0,1]}$ and $(p^{-1/2}2^{-k_0/2}\check{\cY}_{2^{k_0}t})_{t\in[0,1]}$ fit into Pollard's setting of $k_0=0$ and $p=1$. In Pollard's coupling $\check{\cY}$ is determined by $\wt Y$ and the set of jump times $A$ of $\check{\cY}$ satisfies \eqref{eq:tau-ite}.
	
	Since $\check{\cY}$ and $\wt{\cY}$ are related by \eqref{eq41}, it is immediate that $\wt{\cY}$ is also determined by $\wt Y$. Furthermore, since the set of jump times $A$ is the same for $\check{\cY}$ and $\wt{\cY}$ we get \eqref{eq:tau-ite}. In order to complete the proof we need to establish \eqref{eq:pol}. By \eqref{eq:pol0} and the triangle inequality the probability on the left side of \eqref{eq:pol} is bounded above by 
	\eqb
	\P\Big[\big\|\wt{\cY}-\check{\cY}-t2^{-k_0}(n-p2^{k_0})\big\|_\infty > 4p\Delta^2 +C_1|\Delta|2^{-k_0/2}\sqrt{nx} \Big]+C_0\exp(-x).
	\label{eq65}
	\eqe
	By the definition of $\check {\cY}$, 
	\eqbn
	\begin{split}
	\|\wt{\cY}-\check{\cY}-t2^{-k_0}(n-p2^{k_0})\|_\infty 
	&= \|(1-(p^{-1}2^{-k_0}n)^{-1/2})\wt{\cY} + t(p-n2^{-k_0})(1-(p^{-1}2^{-k_0}n)^{-1/2})\|_\infty\\
	&\leq 2|\Delta|2^{-k_0/2}\|\wt {\cY}\|_\infty
	+ 2p\Delta^2,	
	\end{split}
	\label{eq:pol1}
	\eqen
	where we use in the second step that $|(p^{-1}2^{-k_0}n)^{-1/2}-1|=|(1-2^{-k_0/2}\Delta)^{-1/2}-1|\leq 2^{1-k_0/2}|\Delta|$. 
	By considering the law of $\wt {\cY}$ we get $\P[\|\wt {\cY} + (p-2^{-k_0}n)t\|_\infty> 0.5C_1 \sqrt{nx}]<\exp(-x)$ for sufficiently large $x$, which implies that with probability at least $1-\exp(-x)$ we have $\|\wt {\cY}\|_\infty<0.5C_1 \sqrt{nx} + |\Delta|p2^{k_0/2}$. In combination with the last display this gives that with probability at least $1-\exp(-x)$, 
	\eqbn
	\begin{split}
		\|\wt{\cY}-\check{\cY}-t2^{-k_0}(n-p2^{k_0})\|_\infty 
		&\leq 4p\Delta^2 +C_1|\Delta|2^{-k_0/2}\sqrt{nx}.
	\end{split}
	\eqen
	This concludes the proof of \eqref{eq:pol} in combination with \eqref{eq65}.
\end{proof}
Notice that \eqref{eq:tau-ite} uniquely characterizes the coupling of $\wt Y$ and $\wt{\cY}$ since it tells us how to define $\wt{\cY}$ as a function of $\wt Y$. More explicitly, we first use $\wt Y$ to determine $A$ and then we use that $\wt{\cY}$ is determined by $A$. In order to get $A$ from $\wt Y$ we determine $\#(A\cap J_{i,k})$ iteratively for decreasing $k$ and $i\in\{0,1,\dots,2^{k_0-k}-1 \}$. Namely, we are given that $\#(A\cap J_{0,k_0})=\# A=n$, we then get $\#(A\cap J_{0,k_0-1})$ and $\#(A\cap J_{1,k_0-1})$ by using \eqref{eq:tau-ite}, then we get 
$\#(A\cap J_{0,k_0-2})$,
$\#(A\cap J_{1,k_0-2})$,
$\#(A\cap J_{2,k_0-2})$, 
$\#(A\cap J_{3,k_0-2})$ by applying \eqref{eq:tau-ite} again, and so on. 

For an interval $J\subset\R$ let $C(J,\R)$ (resp.\ $D(J,\R)$) denote the set of continuous (resp.\ c\`adl\`ag) real-valued functions on $J$. Define the function spaces
$$
\wt{\cC}:=\bigcup_{k_0\in\Z} C([0,2^{k_0}],\R),\qquad
\wt{\cD}:=\bigcup_{k_0\in\Z} D([0,2^{k_0}],\R).
$$
Let $y\in\wt{\cC}$, $y:[0,2^{k_0}]\to\R$, be an arbitrary continuous function. For $k\in\{\dots,k_0-1,k_0\}$ and $i\in \{0,\dots,2^{k_0-k}-1\}$ define
\eqb
x_{i,k}:=2^{-k/2}(2y_{(2i+1)2^{k-1}}-y_{i2^{k}}-y_{(i+1)2^{k}}).
\label{eq:xik-0}
\eqe 
Define $\wt{\cC}_{\op{good}}\subset \wt{\cC}$ such that $y\in \wt{\cC}_{\op{good}}$ if and only if for every $n\in\N$ there is a set $A=A(y,n)\subset[0,2^{k_0})$ such that for any $k\in\{\dots,k_0-1,k_0 \}$ and $i\in\{0,1,\dots,2^{k_0-k}-1 \}$ we have\footnote{Given $y$ one can attempt to construct an appropriate set $A$ by determining $\#(A\cap J_{i,k})$ iteratively for decreasing $k$ as described in the previous paragraph. An appropriate set $A$ may fail to exist if, for example, there is some $t_0\in[0,2^{k_0}]$ such that, according to the iterative construction, we should have $\#(A\cap J_{i,k})\geq 2$ for all $i,k$ such that $t_0\in J_{i,k}$. An appropriate set $A$ may also fail to exist if there is a dyadic $t_0\in[0,2^{k_0}]$ such that, according to the iterative construction, we should have $\#(A\cap J_{i,k})\geq 1$ for all $i,k$ such that $t_0$ is a right end-point of $J_{i,k}$; there is no set $A$ satisfying this since the intersection of all sets $J_{i,k}$ with $t_0$ as its right end-point is empty. However, Brownian curves of the kind we consider in this paper are elements of $\wt{\cC}_{\op{good}}$, see \eqref{eq:tau-ite}.}  
\eqb
\#(A\cap J_{2i,k-1})= \tau(
\#(A\cap J_{i,k}),
p^{-1/2} x_{i,k}),\qquad \#A=n.
\label{eq:A}
\eqe
Define the function
$$
\wt F_p:\wt{\cC}\times\N_0\to\wt{\cD}
$$
such that if $y\in\wt{\cC}_{\op{good}}$ then $\wt F_p(y,n)$ is the real-valued process on $[0,2^{k_0}]$ given by
$$
t\mapsto \#(A\cap[0,t]) - pt,
$$
while if $y\in \wt{\cC}\setminus\wt{\cC}_{\op{good}}$ then $\wt F_p(y,n)$ is defined in the exact same way but with $A=\emptyset$. 

Observe that in the context Lemma \ref{prop:pol}, $\wt F_p(\wt Y,n)=\wt{\cY}$ a.s. Also observe that if $z=\wt F_p(Y,n)$ then 
\eqb
z|_{[0,2^{k_0-1}]} = \wt F_p(y|_{[0,2^{k_0-1}]},s),\qquad\text{where}\qquad 
s:= \#\{ t\in [0,2^{k_0-1}]\,:\, z_t\neq z_{t^-} \}.
\label{eq43}
\eqe

Let $\wh J\subset J$ be arbitrary intervals such that $\wh J$ is bounded, and suppose $y:J\to\R$ is continuous. Let $t_1<t_2$ denote the end-points of the interval $\wh J$.  
Define $\br(y,\wh J):\wh J\to\R$ by 
\eqb
\br(y,\wh J):=(y+h)|_{\wh J},
\label{eq:br}
\eqe
 where $h:\wh J\to\R$ is the unique linear function $h(t)=a_1+ta_2$ chosen such that %
$\br(y,\wh J)_{t_1}=\br(y,\wh J)_{t_2}=0$. %
We suppose $y$ is such that the limits $\lim_{t\downarrow t_1} y(t)$ and $\lim_{t\uparrow t_2} y(t)$ are well-defined, which is a necessary requirement in order for $\br(y,\wh J)$ to be well-defined. 

We see from \eqref{eq:xik-0} that if we add a linear function to $y$ then $x_{i,k}$ does not change. Therefore, for any function $h:[0,2^{k_0}]\to\R$ of the form $h(t)=a_1+ta_2$ for $a_1,a_2\in\R$,
$$
\wt F_p(y+h,n)=\wt F_p(y,n).
$$
In particular,
\eqb
\wt F_p(y,n)=\wt F_p(\br(y,[0,2^{k_0}]),n).
\label{eq:br-meas}
\eqe

The following is an important property of the function $\tau$ which we will need later.
\begin{lem} Let $\tau$ be the function appearing in Lemma \ref{prop:tus}. 
	Then for any $x\in\R$,
	\eqb
	\tau(n+1,x)\geq\tau(n,x)\qquad\text{and}\qquad (n+1)-\tau(n+1,x)\geq n-\tau(n,x).
	\label{eq:add-one-jump}
	\eqe
	For $k_0\in\Z$, $y: [0,2^{k_0}]\to \R$ satisfying $y\in\wt{\cC}_{\op{good}}$, and $n\in\N$, let $A(y,n)\subset[0,2^{k_0})$ satisfy \eqref{eq:A}. Then $A(y,n)\subset A(y,n+1)$.
	\label{prop:add-one-jump}
\end{lem}
\begin{proof}
	For fixed $n\in\N_0$ the function $x\mapsto \tau(n,x)$ is stepwise constant since it is increasing and it takes only integer values.  It is uniquely determined (except for the values at the jump points) by the fact that $\tau(n,x)=k$ for $k\in\N_0$ and $x\in(x_{n,k},x_{n,k+1})$, where $x_{n,k}$ is defined by
	\eqb
	\P[\xi_n<k]=\P[X<x_{n,k}],
	\quad\xi_n\sim \op{Binomial}(n,1/2),
	\quad X\sim \cN(0,1).
	\label{eq:xi-def}
	\eqe
	By properties of binomial random variables, 
	$\P[\xi_{n+1}<k]<\P[\xi_n<k]$. Using this in combination with \eqref{eq:xi-def} we get $x_{n+1,k}<x_{n,k}$, 
	which implies the first part of \eqref{eq:add-one-jump}.
	
	To prove the second assertion of the lemma, observe that the function $\wh\tau:\N_0\times\R\to\N_0$ defined by $\wh\tau(n,x)=n-\tau(n,-x)$ is increasing in the second variable and that $\wh\tau(n,X)$ has Bin$(n,1/2)$ distribution. Therefore $\wh\tau$ and $\tau$ are identical, except that the two functions may take different values at the points of discontinuity. Using the first assertion we get that away from the discontinuities,
	$$
	(n+1)-\tau(n+1,x) = \wh\tau(n+1,-x) \geq \wh\tau(n,-x) = n-\tau(n,x), 
	$$
	which gives the second part of \eqref{eq:add-one-jump}.
	
	In order to get the last assertion of the lemma we will first prove by induction on $k$ (with $k$ decreasing) that we have 
	\eqb
	\# (A(y,n)\cap J_{i,k})\leq \# (A(y,n+1)\cap J_{i,k})\leq \# (A(y,n)\cap J_{i,k})+1, \quad\forall i\in \{0,\dots,2^{k_0-k}-1\}.
	\label{eq:ind}
	\eqe 
	The induction hypothesis holds for $k=k_0$. We suppose the induction hypothesis holds for some $k$ and want to show that it also holds with $k$ replaced by $k-1$. It is sufficient to argue that for each $i\in \{0,\dots,2^{k_0-k}-1\}$,
	$$
	\#(A(y,n)\cap J_{2i,k-1})\leq \#(A(y,n+1)\cap J_{2i,k-1}),\quad
	\#(A(y,n)\cap J_{2i+1,k-1})\leq \#(A(y,n+1)\cap J_{2i+1,k-1}).
	$$
	The first inequality of the last display is equivalent to  
	$
	\tau\big(\#\big(A(y,n)\cap J_{i,k}\big),p^{-1/2} x_{i,k}\big)
	\leq \tau\big(\#\big(A(y,n+1)\cap J_{i,k}\big), p^{-1/2}x_{i,k}\big)
	$, which holds by the first part of \eqref{eq:add-one-jump} and the induction hypothesis. The second inequality of the last display holds for a similar reason using the second part of \eqref{eq:add-one-jump}. We conclude by induction that \eqref{eq:ind} holds for all $k\in\{\dots,k_0-1,k_0 \}$. 
	
	For $k$ sufficiently small we have $\# (A(y,n+1)\cap J_{i,k})\in\{0,1 \}$ for all $i\in \{0,\dots,2^{k_0-k}-1\}$. For such $k$ there will be exactly $n$ $i$'s such that $\# (A(y,n+1)\cap J_{i,k})= \# (A(y,n)\cap J_{i,k}) =1$ and exactly one $i$ such that $\# (A(y,n+1)\cap J_{i,k})=\# (A(y,n)\cap J_{i,k})+1=1$. For the former $n$ $i$'s it follows from the iterative construction that we will have $\# (A(y,n+1)\cap J_{i',k'})= \# (A(y,n)\cap J_{i',k'})$ for all $i',k'$ such that $k'\leq k$ and $J_{i',k'}\subseteq J_{i,k}$, which implies $A(y,n+1)\cap J_{i,k}= A(y,n)\cap J_{i,k}$ and further the last assertion of the lemma.
\end{proof}

\subsection{One-dimensional process on $\R$}
\label{sec:kmt1}
Let $p>0$ and let $Y^\pm:\R\to\R$ be two independent variance $p$ Brownian motions. Define $(Y_t)_{t\in\R}$ by
\eqb
Y_t = Y^-_{(-t)^-}\1_{t<0}+Y^+_{t}\1_{t>0},\qquad t\in\R.
\label{eq:Y}
\eqe
Then $Y$ has the law of a two-sided variance $p$ Brownian motion. 
Let $(P_t)_{t\geq 0}$ be a rate $p$ Poisson process and define $(\cY^+_t)_{t\geq 0}$ by $\cY^+_t=P_t-pt$. Let $\cY^-$ be an independent copy of $\cY^+$. In this section we will define a coupling of $Y$ and the process 
\eqb
\big(-\cY^-_{(-t)^-}\1_{t<0}+\cY^+_{t}\1_{t>0}\big)_{t\in\R}.
\label{eq:cY-law}
\eqe 
The main result is Proposition \ref{prop-kmt0} below. Lemma \ref{prop:br-indep} will also play an important role in later sections.

It is possible to argue that the coupling between $Y$ and $\cY$ that we give in this section is the local limit of the coupling in Section \ref{sec:kmt0} as we zoom in nearby a uniform time in $[0,2^{k_0}]$ at the same time as $n=p2^{k_0}\rta\infty$. However, we will not be proving this result and our actual construction of the coupling in this section takes a somewhat different route, although it does build heavily on results from Section \ref{sec:kmt0}. Let us explain some key ideas of our construction. We are given $Y$ and want to determine a process of the form \eqref{eq:Y} such that $Y$ and $\cY$ are uniformly close with high probability when restricted to compact intervals. For $k_0\gg 1$ suppose $J\subset\R$ is an interval of length $2^{k_0}$ containing 0 in the interior. We start by constructing a coupling of $Y|_J$ and a process $\cY^{k_0}:J\to\R$ of the form $t\mapsto \#(A\cap (-\infty,t])-pt - \#(A\cap (-\infty,0])$ by the method in Section \ref{sec:kmt0}. In fact, once we are given the number of jumps that $\cY^{k_0}$ makes in $J$ our function $\cY^{k_0}$ will be defined precisely as the function $\wt{\cY}$ in Section \ref{sec:kmt0}, except that we do a translation of the time interval (so the domain of $\cY^{k_0}$ is $J$ instead of $[0,2^{k_0}]$) and we add a constant to $\cY^{k_0}$ (so that $\cY^{k_0}_0=0$); see \eqref{eq:cYm}. We want that $\|Y|_J-\cY^{k_0}\|_\infty$ is small. Then it is natural to require that $\cY^{k_0}$ makes approximately $p2^{k_0}+\Delta_J Y$ jumps during $J$, where $\Delta_JY:=Y_b-Y_a$ for $b$ (resp.\ $a$) equal to the right (resp.\ left) end-point of $J$, since this guarantees that the two functions have (approximately) the same increment during $J$; namely, on the event that the number of jumps is exactly equal to 
$\lceil p2^{k_0}+\Delta_J Y\rceil$,
\eqb
|\Delta_J \cY^{k_0} - \Delta_J Y|
= | (\lceil p2^{k_0}+\Delta_J Y\rceil - p2^{k_0}) - \Delta_J Y|  < 1.
\label{eq44}
\eqe
We define a function $\cY^{k_0}$ as described here for each $k_0\in\N$ (and some appropriate choice of the intervals $J$) and we argue that $\cY^{k_0}|_{I}$ converges a.s.\ as $k_0\rta\infty$ for any fixed bounded interval $I$. We define a function $\cY$ of the form \eqref{eq:cY-law} such that $\cY|_I$ is equal to this limit for any $I$; see \eqref{eq:def-walk0}. 

Before stating our results we need to introduce some conventions related to dyadic intervals. 
For some interval $[a,b)\subset\R$ we say that the two intervals $[a,(a+b)/2)$ and $[(a+b)/2,b)$ are the {\emph{dyadic children}} of $[a,b)$. We say that $[2a-b,b)$ and $[a,2b-a)$ are the left and right, respectively, {\emph{dyadic parents}} of $[a,b)$. 
We say that $J$ is an {\emph{ancestor}} of $I$ if there exists $k\in\N$ and intervals $I_0,\dots,I_k$ such that $I=I_0$, $I_j$ is a dyadic parent of $I_{j-1}$ for $j=1,\dots,k$, and $I_k=J$.

For $\theta\in[0,1)$ and $V_k\in\{-1,1\}$ for $k\in\N$ let $\frk D=\frk D(\theta,V_1,V_2,\dots)$ be the smallest collection of intervals satisfying the following. 
\begin{compactitem}
	\item $[\theta,\theta+1)\in\frk D$.
	\item If $I\in\frk D$ then both the dyadic children of $I$ are also in $\frk D$.
	\item If $I\in\frk D$, $0\in I$, $|I|=2^{k-1}$ and $V_k=-1$ (resp.\ $V_k=1$) for some $k\in\N$ then the left (resp.\ right) dyadic parent of $I$ is also contained in $\frk D$.
	\item If $[a,a+2^k)\in\frk D$ for some $a\in\R$ and $k\in\Z$ then $[a-2^k,a),[a+2^k,a+2^{k+1})\in\frk D$.
\end{compactitem}
Let $\Omega_{\frk D}=\{ \frk D(\theta,V_1,V_2,\dots)\,:\,\theta\in[0,1),V_1,V_2,\dots\in\{-1,1 \} \}$.
An element of $\Omega_{\frk D}$ is called a \emph{dyadic  decomposition} of $\R$. 

Given any dyadic decomposition $\frk D$ of $\R$ and $k\in\Z$ let $I(k)$ denote the unique element of $\frk D$ which has length $2^k$ and which contains the number 0. Given any dyadic decomposition $\frk D$ of $\R$ and $t_0\in\R$ let $\frk D-t_0$ denote the dyadic decomposition
\eqb
\frk D-t_0 :=\{ [t-t_0,s-t_0)\,:\,[t,s)\in\frk D \}.
\label{eq:D-shift}
\eqe

If $\theta\in[0,1)$ has been sampled uniformly at random and $(V_k)_{k\in\N_0}$ are i.i.d.\ and uniform on $\{-1,1 \}$ then we call the associated dyadic decomposition $\frk D(\theta,V_1,V_2,\dots)$ a \emph{uniform dyadic decomposition of $\R$.} We assume in the remainder of this subsection that $\frk D$  is a uniform dyadic decomposition independent of $Y$.

\begin{figure}
	\includegraphics[]{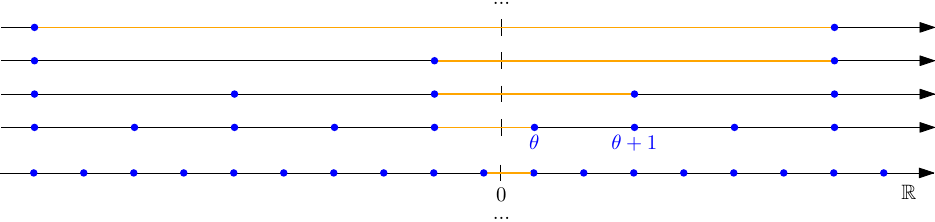}
	\caption{An illustration of the dyadic decomposition $\frk D(\theta,V_1,V_2,...)$ for $V_1=1,V_2=1,V_3=-1,\dots$. The intervals between two consecutive blue dots are elements of the dyadic decomposition. The intervals in orange are, from bottom to top, $I(-1),I(0),I(1),I(2),I(3)$. The bottom (resp.\ top) three dots ($...$) indicate that there are additional shorter (resp.\ longer) intervals in the dyadic decomposition that are not shown.}
\end{figure}

Recall the function $\wt F_p$ from Section \ref{sec:kmt0}. We will now extend this function by enlarging its domain. More precisely, we will define a function
$$
F_p: \cC\times\R\to\cD
$$
such that $F_p|_{\wt{\cC}\times\N_0}=\wt F_p$, where
$$
\cC:=\bigcup_{k_0\in\Z,a\in [-2^{k_0},0]} C([a,a+2^{k_0}],\R),\qquad
\cD:=\bigcup_{k_0\in\Z,a\in [-2^{k_0},0]} D([a,a+2^{k_0}],\R).
$$
Notice that the domain of $F$ is larger than the domain of $\wt F$ in two ways: The second input argument is an arbitrary real number $r\in\R$ (rather than $r\in\N_0$) and we allow $a\in[-2^{k_0},0]$ instead of requiring $a=0$. Roughly speaking, we define $F_p$ by doing a shift by $a$, approximating $r$ by an element of $\N_0$, and then applying $\wt F_p$. 
More precisely, given any $y\in\cC$ with domain $[a,a+2^{k_0}]$ define $\wt y\in\wt{\cC}$ by $\wt y(t):=y(t+a)$. Then let $F_p(y,r):[a,a+2^{k_0}]\to\R$ be the function
\eqb
t\mapsto \wt F_p(\wt y,0\vee \lceil r \rceil)_{t-a} - \wt F_p(\wt y,0\vee \lceil r \rceil)_{-a}.
\label{eq68}
\eqe 
Then $F_p(y,r)$ is a function of the form $t\mapsto \#(A\cap(-\infty,t]) - pt + b$ for a set $A\subset[a,a+2^{k_0}]$, where the constant $b\in\R$ is chosen such that $F_p(y,r)_0=0$. If $\wt y\in\wt{\cC}_{\op{good}}$, the number of jumps is $\#A=0\vee\lceil r\rceil$; otherwise $\#A=0$ by the definition of $\wt F_p$.

The random variables $X_J$ introduced in the below lemma will play a central role throughout the remainder of the paper.
\begin{lem}
	Conditioned on any realization of $\frk D$
	the following is a collection of independent variance $p$ normal random variables 
	\eqb
	X_{J}:= (b-a)^{-1/2}(2Y_{0.5(a+b)}-Y_a-Y_b),\qquad
	J=[a,b)\in\frk D.
	\label{eq:X-def}
	\eqe
	For any $I\in\frk D$ let $J^I_0,J^I_1,\dots\in\frk D$ be such that $J^I_0=I$ and $J^I_{j}$ is the dyadic parent of $J^I_{j-1}$ for $j\in\N$, and set $V^I_j=1$ (resp.\ $V^I_j=-1$) if $J^I_j$ is the right (resp.\ left) dyadic parent of $J^I_{j-1}$. Then
	\eqb
	\Delta_{I}Y = \sum_{j\in\N} |I|^{1/2} 2^{-j/2}V^I_j X_{J^I_j}. 
	\label{eq:Z-brZ-X}
	\eqe
	\label{prop:br-indep}
\end{lem}
\begin{proof}
	A computation of variances and covariances gives that each $X_J$ is a variance $p$ normal random variable and that the $X_J$'s are pairwise uncorrelated. This implies the first assertion of the lemma since the $X_J$'s are jointly Gaussian. The last assertion of the lemma follows since, by the definition in \eqref{eq:X-def} (more precisely, by writing $X_{J^I_j}$ as a sum of $\pm\Delta_J Y$ for disjoint intervals $J$ of length $|I|$ with union $J_j^I$) and $|J^I_j|=|I|2^j$, the right side of \eqref{eq:Z-brZ-X} is equal to
	\eqbn
	\begin{split}
		\sum_{j\in\N} \Big(
		\sum_{J\in\frk D,\:\,|J|=|I|,J\subset J^I_{j-1}}  
		2^{-j} \Delta_J Y  -
		\sum_{J\in\frk D,\:\,|J|=|I|,J\subset J^I_j\setminus J^I_{j-1}} 
		2^{-j} \Delta_J Y \Big)
		= \sum_{j\in\N} 2^{-j}\Delta_I Y = \Delta_I Y,
	\end{split}
	\eqen
	Restricting the sum to $j\leq N$ and swapping the order of summation we get
	$$
	\sum_{j=1}^{N}2^{1-j}
	\left(\Delta_{J_{j-1}^I}Y-\frac12\Delta_{J_j^I}Y\right)
	=\sum_{j=0}^{N-1}2^{-j}\Delta_{J_j^I}Y
	-\sum_{j=1}^{N}2^{-j}\Delta_{J_j^I}Y\\
	=\Delta_IY-2^{-N}\Delta_{J_N^I}Y.
	$$
	The second term on the right side converges a.s.\ to zero as $N\rta\infty$, which justifies \eqref{eq:Z-brZ-X}.
\end{proof}

Let $Y$ be as in \eqref{eq:Y} and let $\frk D$ be an independent uniform dyadic decomposition of $\R$. For $m\in\N_0$ and $S\in\R$ define $\cY^{m,S}:I(m)\to\R$, $\cY^{m}:I(m)\to\R$, and $S(m)\in\R$ by
\eqb
\cY^{m,S} := F_p(Y|_{I(m)}, S),\quad
\cY^{m} := \cY^{m,S(m)},\quad
S(m):= 2^m p + \sum_{j\in\{1,\dots,2m \}} 2^{(m-j)/2} 
V_{m+j} X_{I(m+j)}.
\label{eq:cYm}
\eqe
Whenever $F_p$ is applied to a continuous path restricted to a half-open dyadic interval, first extend the restriction continuously to the closed interval, apply $F_p$, and restrict the result back to the half-open interval. 
Notice that $\cY^{m}$ is measurable with respect to $Y|_{I(3m)}$. 
Let $(\cY_t)_{t\in\R}$  
be the process satisfying the following, provided such a process exists and the limit is well-defined for each $k\in\N$
\eqb
\cY|_{I(k)}=\lim_{m\rta\infty} \cY^m|_{I(k)}.
\label{eq:def-walk0} 
\eqe
See the second paragraph of this section for intuition and motivation for our definition of $\cY$ and $\cY^m$. 

The following lemma shows that $S(m)$ is close to $2^{m}p+\Delta_{I(m)}Y$ as desired. 
\begin{lem}
	We have that $S(m)\in\sigma(\frk D,\sigma(Y|_{I(3m)}))$. Furthermore $\P[|S(m)-2^mp-\Delta_{I(m)}Y|>x]\leq C\exp(-c2^{m}x^2)$.
	\label{prop:Sm}
\end{lem} 
\begin{proof}
	The first assertion is immediate. To get the second assertion observe that by \eqref{eq:Z-brZ-X} the considered probability is equal to
	$$
	\P\Big[\Big|\sum_{j=2m+1 }^\infty 2^{(m-j)/2} V_{m+j} X_{I(m+j)}\Big|>x \Big],
	$$
	where the bound follows since the considered sum is a centered normal random variable with variance of order $2^{-m}$.
\end{proof}
The reader might ask why we define $S(m)$ by \eqref{eq:cYm} instead of setting $S(m)=\lceil 2^{m}p+\Delta_{I(m)}Y \rceil\vee 0$. We do this due to properties of the higher-dimensional coupling we study in Section \ref{sec:kmt3} later. More precisely, we do it in order to guarantee that $S(m)$ is determined by the higher-dimensional Brownian motion $Z$ restricted to $I(3m)$ (along with $\frk D$ and an independent random variable $\frk U$) also in the higher-dimensional setting.

The following is the main result of this section.
\begin{prop} 
	The function $\cY:\R\to\R$ in \eqref{eq:def-walk0} is well-defined a.s.\ and is equal in law to the process \eqref{eq:cY-law}. Furthermore, the following hold. 
	\begin{compactitem}
		\item[(i)] $(Y,\cY,\frk D)$
		is translation invariant in law in the sense that for any fixed $t_0\in\R$,
		$$
		( (Y_{t+t_0}-Y_{t_0})_{t\in\R},  (\cY_{t+t_0}-\cY_{t_0})_{t\in\R},
		\frk D-t_0 ) \eqD (Y,\cY,\frk D).
		$$
		\item[(ii)] There is a $c\in(0,1)$ such that for any $k\in\N_0$ and $m\in\N\cap[2\log_2(k\vee 1),\infty)$, 
		\eqbn
		\begin{split}
			\P\big[ \cY|_{I(k)}\neq \cY^{k+m}|_{I(k)} \big]&<\exp(-cm).
		\end{split}
		\eqen
		\item[(iii)] For any $C_1>1$ we can find $C'>1$ such that for any $k,m\in\N$,
		\eqb
		\P\Big[ \sup_{t\in I(m)} |Y_t-\cY^{m+k}_t|>C'm \Big] \leq 2^{-C_1m}
		\label{eq:err3-0}
		\eqe
		and
		\eqb
		\P\Big[ \sup_{t\in I(m)} |Y_t-\cY_t|>C'm \Big] \leq 2^{-C_1m}.	
		\label{eq:err3-1}
		\eqe
	\end{compactitem}
	\label{prop-kmt0}
\end{prop}

The remainder of this section is devoted to the proof of Proposition \ref{prop-kmt0}.  Before starting the proof let us give a brief outline. In order to argue that the limit in \eqref{eq:def-walk0} is well-defined we will argue that for each fixed $m\in\N$ and for large $n\in\N$ it holds with high probability that $\cY^n|_{I(m)}=\cY^{n+1}|_{I(m)}$ (see Lemma \ref{prop-sqrtk-E} below). By Lemma \ref{prop43} it is sufficient to argue the following two results: 
(i) $\wt S(n)\approx S(n)$, where $\wt S(n)$ is as defined in Lemma \ref{prop43} and satisfies $\cY^{n+1}|_{I(m)}=\cY^{n,\wt S(n)}|_{I(m)}$, and 
(ii) for all $S\approx S(n)$ it holds that $\cY^{n,S}|_{I(m)}=\cY^n|_{I(m)}$. The latter two results are proven in Lemmas \ref{prop:err3} and \ref{prop:same-E}, respectively. The bounds appearing in (ii) and (iii) of Proposition \ref{prop-kmt0} are proven by combining \eqref{eq:pol} and estimates from the proof of existence of $\cY$.

The next lemma is an immediate consequence of \eqref{eq43}. Before stating it we will introduce a function $S$ which will be used several times in the remainder of the paper. For intervals $I,J\subseteq \R$ with $J\subseteq I$ and a c\`adl\`ag function $y:I\to\R^m$, let $S(y,J)\in\N_0\cup\{\infty \}$ denote the number of jumps made by $y$ during $J$, i.e., 
\eqb
S(y,J) := \#\{ t\in J\,:\,\Delta_t y\neq 0 \},
\qquad\text{where}\qquad 
\Delta_t y:=y(t)-\lim_{s\uparrow t} y(s).
\label{eq:defS}
\eqe
\begin{lem}
	Let $n\in\N$ and set $\wt S(n):= S(\cY^{n+1},I(n))$. Then 
	$\cY^{n+1}|_{I(n)}=\cY^{n,\wt S(n)}$ a.s.
	\label{prop43}
\end{lem}

The following lemma establishes (ii) from the paragraph below Proposition \ref{prop-kmt0}. It is a consequence of Lemma \ref{prop:add-one-jump} and the translation invariance in law of the set of jump times of $\cY^{m+n,S}$.
\begin{lem}
		For all $C'>1$ there is a $c\in(0,1)$ such that for all $m,n\in\N$ such that $n\geq 2\log_2(m)$,
		\eqb
		\P[ \exists S_1,S_2 \in \cI\,:\, \cY^{m+n,S_1}|_{I(m)} \neq \cY^{m+n,S_2}|_{I(m)}  ]<2^{-cn},
		\label{eq8-E}
		\eqe
		where for $S(j)$ as in \eqref{eq:cYm},
		$$
		\cI:=[
		S(m+n)-C'(m+n), 
		S(m+n)+C'(m+n) ]\cap\N_0,\qquad 
		S(j):=p2^{j}+\Delta Y_{I(j)}.
		$$
		\label{prop:same-E}
	\end{lem}
	\begin{proof}
		It is sufficient to prove the lemma for large $n$. 
		We will argue that for all $S\in\N_0$ (possibly random, but chosen independently of $\br(Y,I(m+n))$), for sufficiently large $n$ and $n\geq 2\log_2(m)$,
		\eqb
		\P[\cY^{m+n,S}|_{I(m)} \neq \cY^{m+n,S+1}|_{I(m)}  ] \leq 2^{-n}.
		\label{eq:one-S}
		\eqe
		This is sufficient to conclude the proof since then a union bound, \eqref{eq:one-S}, and the last assertion of Lemma \ref{prop:br-indep} imply that the left side of \eqref{eq8-E} is bounded above by 
		$$
		(2C'(m+n)+1) 2^{-n} \leq 2^{-cn}
		$$
		for an appropriate $c\in(0,1)$. Let $A_S\subset I(m+n)$ (resp.\ $A_{S+1}\subset I(m+n)$) denote the set of times at which  $\cY^{m+n,S}$ (resp.\ $\cY^{m+n,S+1}$) makes a jump. Then $\# A_S=S$ and $\# A_{S+1}=S+1$. Furthermore, by Lemma \ref{prop:add-one-jump} there is a $t_0\in I(m+n)$ such that $A_{S+1}=A_S\cup\{t_0 \}$. Therefore
		$$
		\{\cY^{m+n,S}|_{I(m)} \neq \cY^{m+n,S+1}|_{I(m)}\}
		= \{t_0\in I(m) \}.
		$$ 
		Conditioned on $S$ and $I(m+n)$ it holds that $A_S$ (resp.\ $A_{S+1}$) has the law of a collection of $S$ (resp.\ $S+1$) independent uniform random variables in $I(m+n)$, which implies that $\P[t_0\in I(m)]=2^{-n}$ and gives \eqref{eq:one-S}.
	\end{proof}
	
	The following lemma establishes (i) from the paragraph below Proposition \ref{prop-kmt0}. The key input is the error bound in Tusn\'ady's lemma. 
	\begin{lem}
		For any $C_0>0$ it holds for $C'>1$ sufficiently large and for all $j\in\N$ that
		$$
		\P[| \wt S(j)%
		-S(j)|>C'j] < 2^{-C_0 j},
		$$
		where $\wt S(j)$ is as defined in Lemma \ref{prop43}.
		\label{prop:err3}
	\end{lem}
	\begin{proof}
		We assume throughout the proof that $j$ is sufficiently large. 
		For concreteness we assume throughout the proof that $I(j+1)$ is the right dyadic parent of $I(j)$, but the other case (with ``left'' instead of ``right'') can be treated in the same way. Consider the independent normal random variables
		$$
		X_1 := 2^{-j/2}\Delta_{I(j)} Y, 
		\qquad 
		X_2 := 2^{-j/2}(\Delta_{I(j+1)} Y-\Delta_{I(j)} Y).
		$$
		Throughout the proof we work on the following high probability event, which we may do since its complement has probability smaller than $2^{-C_0j}$ by Lemma \ref{prop:Sm}
		$$
		|\lceil S(j+1) \rceil - \lceil p2^{j+1}+ 2^{j/2} (X_1+X_2)\rceil|\leq 1,\qquad
		|S(j)-(p2^j+2^{j/2}X_1)|<1.
		$$
		We will apply Tusn\'ady's lemma (Lemma \ref{prop:tus}) with
		$$
		n := \lceil S(j+1)\rceil\vee 0,\qquad
		X := (2p)^{-1/2}(X_1-X_2).
		$$
		Define $\wh X=n/2+X\sqrt{n/4}$ and $\xi=\xi(n,X)$ and observe that $\wt S(j)=\xi$ by \eqref{eq:tau-ite}. Suppose $S(j+1)>0$, which happens with probability at least $1- 2^{-C_0j-1}$. Then for some constant $C_p>1$ depending on $p$,
		\eqbn
		\begin{split}
			| \wt S(j) 
			-S(j)|
			&\leq |\xi-(p2^j+2^{j/2}X_1)|+1
			\leq |\xi-\wh X | + |\wh X-(p2^j+2^{j/2}X_1)|+1\\
			&\leq 3+X^2/8 + |X_1^2-X_2^2|/(4p)+p^{-1/2}2^{-(j+1)/2}|X|
			\leq C_p(1+X_1^2+X_2^2),
		\end{split}
		\eqen
		where we applied Tusn\'ady's lemma, $|n-p2^{j+1}-2^{j/2}(X_1+X_2)|\leq2$, and
		\[
		|\sqrt n-\sqrt{p2^{j+1}}|\leq p^{-1/2} 2^{-(j+1)/2} |n-p2^{j+1}|.
		\] 
		The lemma now follows from the Gaussian tail bound.
	\end{proof}
The following lemma implies in particular that the process $\cY$ defined by \eqref{eq:def-walk0} is well-defined. It is a consequence of Lemmas \ref{prop:same-E} and \ref{prop:err3} above.
\begin{lemma}
	For $j_1\leq j_2$ let $B_{j_1,j_2}$ be the event that
	\eqb
	\cY^{j_2}|_{I(j_1)} \neq
	\cY^{j_2+1}|_{I(j_1)}.
	\label{eq:B-E}
	\eqe
	Then for some $C_0>1$ it holds for all $m,n\in\N_0$ such that $n\geq 2\log_2(m\vee 1)$,
	\eqb
	\P\left[ \bigcup_{j\geq m+n} B_{m,j} \right] \leq 2^{-n/C_0}.
	\label{eq:same-E}
	\eqe
	In particular, $\cY$ is well-defined a.s.\ and for all $m,n\in\N$ such that $n\geq 2\log_2(m)$,
	$$
	\P[\cY|_{I(m)}\neq
	\cY^{m+n}|_{I(m)}
	] < 2^{-n/C_0}.
	$$
	\label{prop-sqrtk-E}
\end{lemma}
\begin{proof}
	It is sufficient to prove the lemma for large $n$ and we assume throughout the proof that $n$ is large. 
	First observe that for any constant $C'>1$ we have $B_{m,j}\subset E_j\cup \wt B_{m,j}$, where
	\eqbn
	\begin{split}
	E_{j} &:= \{ |\wt S(j)-S(j)|>C'j \},\\
	\wt B_{m,j}&:=
	\{ \exists S_1,S_2 \in \cI\,:\, \cY^{j,S_1}|_{I(m)} \neq \cY^{j,S_2}|_{I(m)}  \},\qquad
	\cI:=[
	S(j)-C'j, 
	S(j)+C'j ]_{\Z}.
\end{split}
	\eqen
	By choosing $C'$ sufficiently large, Lemma \ref{prop:err3} gives that
	$\P[E_{j}] \leq 2^{-j}.$
	By Lemma \ref{prop:same-E} and $j-m\geq n\geq 2\log_2 m$, there exists $c\in(0,1)$ such that
	$\P[\wt B_{m,j}] \leq 2^{-c(j-m)}$. 
	The latter two estimates and a union bound %
	give $\P[B_{m,j}]\leq 2^{-j}+2^{-c(j-m)}$, which implies \eqref{eq:same-E}. 
	On the complement of the union in \eqref{eq:same-E} the limit on the right side of \eqref{eq:def-walk0} exists if we replace $k$ and $m$ by $m$ and $\ell$, respectively. It follows that $\cY$ is well-defined and that 
	$\cY|_{I(m)}=\cY^{m+n}|_{I(m)}$ on the complement of the union in \eqref{eq:same-E}, 
	which gives the last assertion of the lemma.
\end{proof}
Next we will establish the error bounds in \eqref{eq:err3-0} and \eqref{eq:err3-1} in Lemma \ref{prop:eucl-error-E} building on  Lemma \ref{prop:pol2}.
The following elementary lemma will be applied in the proof of Lemma \ref{prop:pol2} and we leave the proof for the reader. 
\begin{lem}
	Suppose $b\geq a\geq 0$, $h:[0,b]\to\R$ is a linear function, and $f_1:[-a,b-a]\to \R$ and $f_2:[0,b]\to\R$ are functions such that $f_2(t)=f_1(t-a)+h(t)$ for all $t\in[0,b]$, $f_2(0)=f_2(b)=0$, $|f_1(-a)-f_1(b-a)|\leq 2$, and $f_1(0)=0$. Then $\|f_1\|_\infty \leq 2\|f_2\|_\infty+2$.
	\label{prop:func-diff-linear}
\end{lem}

The next lemma follows from \eqref{eq:pol}.
\begin{lem}
	For any $C_0>1$ we can find $C'>1$ such that for all $m\in\N$,
	$$
	\P\Big[\sup_{t\in I(m)}|Y_t-\cY^{m}_t|>C'm\Big]
	\leq \exp(-C_0m).
	$$
	\label{prop:pol2}
\end{lem}
\begin{proof}
	Let $K>1$ be such that $K^2/(2p)>2(C_0+1)$. Defining the event $E=\{ |S(m)-p2^m|<K2^{m/2}\sqrt{m};|\lceil S(m)\rceil -2^mp-\Delta_{I(m)}Y|<2\}$,
	\eqb
	\begin{split}
		\P\Big[\sup_{t\in I(m)}|Y_t-\cY^{m}_t|>C'm\Big]
		\leq \E\Big[ \P\Big[ \sup_{t\in I(m)}|Y_t-\cY^{m}_t|>C'm\,|\, \sigma( S(m),\Delta_{I(m)}Y ) \Big]\1_E\Big] + \P[E^c].
	\end{split}
	\label{eq:pol2}
	\eqe
	We have $\P[E^c]\leq \exp(-2(C_0+1)m)$ for sufficiently large $m$ by Lemma \ref{prop:Sm} and since it follows from \eqref{eq:cYm} that $S(m)-p2^{m}$ is a centered Gaussian of variance slightly smaller than $p2^{m}$. To bound the first term on the right side of \eqref{eq:pol2} we will apply Lemma \ref{prop:pol} to argue that for sufficiently large $C'>1$, on the event $E$, 
	\eqb
	\P\Big[ \sup_{t\in I(m)}|Y_t-\cY^{m}_t|>C'm\,\,\Big|\,\, \sigma( S(m),\Delta_{I(m)}Y ) \Big] < \exp(-2C_0m),
	\label{eq22}
	\eqe
	which is sufficient to conclude the proof of the lemma. 
	
	Let us explain how we apply Lemma \ref{prop:pol} to get \eqref{eq22}. There are three obstacles to applying the lemma, namely 
	(a) $C'm$ in \eqref{eq22} is different from the bound in the lemma, 
	(b) the function $Y-\cY^{m}$ does not have exactly the same form as the function $\wt{\cY}-\wt Y-t2^{-k_0}(n-p2^{k_0})$ appearing in the lemma, and 
	(c) we consider a conditional expectation in \eqref{eq22}. To address (a), recall that $\cY^m$ makes $\lceil S(m)\rceil$ jumps during $I(m)$ and  observe that on the event $E$, if we set $\Delta=2^{m/2}(1-p^{-1}2^{-m}\lceil S(m)\rceil)$ then $|\Delta|\leq p^{-1}(1+K)\sqrt{m}$. Therefore, if we pick $C'$ sufficiently large, on the event $E$ and with $x_m=4C_0m$,
	\eqb
	4p\Delta^2 +C_1|\Delta|2^{-m/2}\sqrt{\lceil S(m)\rceil x_m} +C_1 \sqrt{p 2^m/\lceil S(m)\rceil } (x_m+\log \lceil S(m)\rceil ) < (C'm-2)/2,
	\label{eq24}
	\eqe
	where we observe that the left side of the last display is of a similar form as the bound in Lemma \ref{prop:pol} for $x=x_m,n=\lceil S(m) \rceil,k_0=m$. 
	
	To address (b) and (c) we define auxiliary functions $\wt Y,\wt{\cY}:[0,2^m]\to\R$ by 
	$$
	\wt Y_t := \br(Y,I(m))_{t-a} - \br(Y,I(m))_{-a},\qquad
	\wt{\cY}_t := \cY^m_{t-a} - \cY^m_{-a}, 
	$$
	where $-a<0$ is the left end-point of the interval $I(m)$. The functions $\wt Y,\wt{\cY}$ are of the same form as the functions in Lemma \ref{prop:pol} and are coupled together as in that lemma by the definition of $\cY^m$ in \eqref{eq:cYm}. 
	By Lemma \ref{prop:func-diff-linear}, %
	\eqb
	\|(Y-\cY^m)|_{I(m)}\|_\infty \leq 2\|\wt Y-\wt{\cY} + t2^{-m}(\lceil S(m) \rceil -p2^m) \|_\infty+2,
	\label{eq23}
	\eqe
	where we observe that the assumptions of the lemma are exactly satisfied by the function $f_1=(Y-\cY^m)|_{I(m)}$ and $f_2=\wt Y-\wt{\cY} + t2^{-m}(\lceil S(m)\rceil-p2^{m})$; in particular, we have 
	$$
	f_2(2^m) = 0 - \wt{\cY}_{2^m} + (\lceil S(m)\rceil - p2^m) = 0
	$$
	since the process $\cY^{m,S}$ defined in \eqref{eq:cYm} makes $\lceil S\rceil \vee 0$ jumps during $I(m)$ by definition,
	and, furthermore,
	$$
	|f_1(-a)-f_1(2^m-a)|
	=|\Delta_{I(m)}Y-\Delta_{I(m)}\cY^m| 
	= |\Delta_{I(m)}Y - (\lceil S(m)\rceil -p2^m )|
	\leq 2
	$$ 
	by the same considerations as in \eqref{eq44} and by the definition of $E$.  
	By Lemma \ref{prop:br-indep}, $\wt{Y}$ is independent of $S(m)$ and $\Delta_{I(m)} Y$ so we can treat $S(m)$ and $\Delta_{I(m)} Y$ as constants when we bound the right side of \eqref{eq23} given $S(m)$; this addresses issue (c) mentioned above. 
	Combining the above observations, we get \eqref{eq22} since on the event $E$,
	\eqb
	\begin{split}
		&\P\Big[ \|(Y-\cY^{m})|_{I(m)}\|_{\infty}>C'm\,\,\Big|\,\, \sigma( S(m),\Delta_{I(m)}Y ) \Big]\\
		&\qquad\leq\P\Big[ \| \wt Y-\wt{\cY} + t2^{-m}(\lceil S(m)\rceil-p2^{m}) \|_\infty > (C'm-2)/2\,|\, \sigma( S(m),\Delta_{I(m)}Y ) \Big] 
		\leq \wt C_0\exp(-x_m), %
	\end{split}
	\eqe
	where we apply Lemma \ref{prop:pol}, \eqref{eq24}, and independence of $\wt Y$ and $S(m)$ in the last step and $\wt C_0>1$ is a constant. The right side of the last display is bounded above by  $\exp(-2C_0m)$ for sufficiently large $m$.
\end{proof}

The next lemma follows from Lemma \ref{prop-sqrtk-E} and Lemma \ref{prop:pol2}.
\begin{lem}
	For any $C_1>1$ we can find $C'>1$ such that for all $m,k\in\N$,
	$$\P[\sup_{t\in I(m)}|Y_t-\cY_t|>C' m]\leq 2^{-C_1 m},\qquad
	\P[\sup_{t\in I(m)}|Y_t-\cY^{m+k}_t|>C' m]\leq 2^{-C_1 m}
	$$
	\label{prop:eucl-error-E}
\end{lem}
\begin{proof}
	For any $k\in\N$, 
	\eqbn
	\begin{split}
		&\P[\sup_{t\in I(m)}|Y_t-\cY_t|>C'm]
		\leq \P[\sup_{t\in I(m)}|Y_t-\cY^{m+k}_t|>C'm]
		+ \P[\cY^{m+k}|_{I(m)}\neq\cY|_{I(m)}].
	\end{split}
	\eqen
	The second term on the right side is smaller than $2^{-k/C_0}$ by Lemma \ref{prop-sqrtk-E} provided $k\geq 2\log_2 m$; we pick $k=\lceil (C_1+1) C_0 m\rceil$. To bound the first term on the right side we apply Lemma \ref{prop:pol2} to get that for $C'$ sufficiently large,
	$$
	\P[\sup_{t\in I(m)}|Y_t-\cY^{m+k}_t|>C'm]
	\leq \P\Big[\sup_{t\in I(m+k)}|Y_t-\cY^{m+k}_t|> \frac{C'}{1+C_0C_1} (m+k)\Big] \leq \exp(-C_1(m+k)),
	$$
	which gives the first assertion of the lemma. To establish the second assertion of the lemma, first observe that if $k\leq \wt Cm$ for any given constant $\wt C>1$ then the estimate is an immediate consequence of Lemma \ref{prop:pol2}. For $k\geq \wt Cm$ we have
	$$
	\P\Big[\sup_{t\in I(m)}|Y_t-\cY^{m+k}_t|>C' m\Big]
	\leq \P\Big[\sup_{t\in I(m)}|Y_t-\cY_t|>C' m\Big]
	+ \P\Big[\cY|_{I(m)}\neq\cY^{m+k}|_{I(m)}\Big]
	$$
	By Lemma \ref{prop-sqrtk-E} the second term on the right side is bounded by $2^{-k/C_0}\leq 2^{-C_1 m-1}$ if $\wt C$ is sufficiently large. By the first assertion of the current lemma, the first term on the right side is bounded above by $2^{-C_1m-1}$ for $C'$ sufficiently large.
\end{proof}

\begin{proof}[Proof of Proposition \ref{prop-kmt0}]
	Most assertions are proven in Lemmas \ref{prop-sqrtk-E} and \ref{prop:eucl-error-E}. It remains to argue translation invariance in law as in Assertion (i) and that $\cY$ has the desired law \eqref{eq:cY-law}. It is immediate from \ref{eq:def-walk0} and the law of $\cY^m$ that $\cY$ has the desired law. To get translation invariance in law, first observe that $Y$ and $\frk D$ are translation invariant in the sense of Assertion (i), and they are jointly translation invariant by independence. We now get Assertion (i) by using that $\cY$ is a measurable function of $Y$ and $\frk D$. 
\end{proof}

We end the section by stating some properties of the coupling which will be used to prove Lemma \ref{prop:def-cZn} later. 
In the next lemma we let $Y,\cY:\R\to\R$ be as in Proposition \ref{prop-kmt0} and we recall the definition of $X_J$ for $J\in\frk D$ from \eqref{eq:X-def}, the definition of $\tau$ from Lemma \ref{prop:tus}, and the definition of $S(\cdot,\cdot)$ from \eqref{eq:defS}.
\begin{lem}
	The coupling in Proposition \ref{prop-kmt0} satisfies the following.
	\begin{compactitem}
		\item[(i)] For any $I,J\in\frk D$ such that $I$ is the right (resp.\ left) dyadic parent of $J$ it holds that
		$$
		S(\cY,J) = \tau(  S(\cY,I), p^{-1/2} X_I )\qquad
		(\text{resp.}\,\,\, S(\cY,J) = S(\cY,I)-\tau(  S(\cY,I), p^{-1/2} X_I )).
		$$
		\item[(ii)] For any $m\in\Z$ it holds that $S( \cY, I(m) )$ is a measurable function of $(V_j,X_{I(j)})_{j\geq m+1}$, where for $j=\{0,-1,-2,\dots \}$ we set $V_j=1$ (resp.\ $V_j=-1$) if $I(j)$ is the right (resp.\ left) dyadic parent of $I(j-1)$.
		\item[(iii)] For any $m\in\Z$, $\cY|_{I(m)}$ is measurable with respect to $\br(Y,I(m))$ and $S(\cY,I(m))$.
	\end{compactitem}
\label{prop:coupling-properties}
\end{lem}
\begin{proof}
	The first part of Assertion (i), namely the case when $I$ is the right dyadic parent of $J$, follows from \eqref{eq43}, the definition of $F_p(\cdot,\cdot)$ in \eqref{eq68}, the definition of $\cY^m$ in \eqref{eq:cYm}, and the definition of $\cY$ in \eqref{eq:def-walk0}. The second part of Assertion (i) follows from the first part along with the fact that $S(\cY,J)+S(\cY,I\setminus J)=S(\cY,I)$.
	
	By \eqref{eq:def-walk0}, in order to prove Assertion (ii) it is sufficient to argue that for any $n\in\{\max\{1, m+1\},\dots \}$ it holds that $S(\cY^n, I(m))$ is a measurable function of $(V_j,X_{I(j)})_{j\geq m+1}$. By iterating the result in (i) we see that $S(\cY^n, I(m))$ is a measurable function of $S(\cY^n,I(n)),(V_j,X_{I(j)})_{m+1\leq j\leq n}$. By \eqref{eq:cYm} it holds that $S(\cY^n,I(n))=\lceil S(n)\rceil\vee 0\in\sigma\big( (V_j,X_{I(j)})_{n+1\leq j\leq 3n} \big)$, which completes the proof.
	
	The process $\cY|_{I(m)}$ is measurable with respect to its set of jumps $A\cap I(m)$. By Assertion (i), the set $A\cap I(m)$ is determined by $S( \cY,I(m) )$ along with $(X_J)_{J\in\frk D,J\subseteq I(m)}$. The latter set is determined by $\br(Y,I(m))$ by \eqref{eq:X-def}, which implies Assertion (iii).
\end{proof}

\renewcommand{\d}{d}

\subsection{Higher-dimensional process on $\R$}
\label{sec:kmt3}
Let $\d,m\in\N$. In this section we extend the coupling of the previous section to the case of $m$-dimensional processes $(\cZ_t)_{t\in\R}$ and $(Z_t)_{t\in\R}$. The process $Z$ is a two-sided $m$-dimensional correlated Brownian motion, while the process $\cZ$ is a continuous time random walk, as in the following definition.
\begin{definition}
	Let $v_1,\dots,v_d\in\R^m$. A stochastic process $(\cZ_t)_{t\in\R}$ with values in $\R^m$ is a \emph{continuous time random walk with steps $\{ v_1,\dots,v_d\}$} if the trajectories $t\mapsto \cZ_t$ are almost surely piecewise constant and right-continuous with a finite number of jumps in any bounded interval, such that all the jumps of the process take their value in $\{v_1,\dots,v_d \}$. We say that such a process has \emph{rate} $r>0$ if the set of times at which the process jumps has the law of a rate $r$ Poisson point process.
	\label{def:cont-time-rw}
\end{definition}

Let $v_1,\dots,v_\d\in\R^m\setminus\{0 \}$ and $p_1,\dots,p_\d\in(0,1]$ be such that $p_1+\dots+p_\d=1$ and $p_1v_1+\dots+p_\d v_\d=0$. Suppose the $v_\ell\neq v_{\wt\ell}$ if $\ell\neq\wt\ell$. Let $(\frk Z_j)_{j\in\N_0}$ be the $m$-dimensional centered random walk such that $\frk Z_0=0$ and the increments of  $(\frk Z_j)_{j\in\N_0}$ are independent with 
$$
\P[\frk Z_{j+1}-\frk Z_j=v_\ell] = p_\ell,\qquad \ell=1,\dots,\d,\quad j\in\N_0.
$$
Let $(P_t)_{t\geq 0}$ be an independent unit rate Poisson process and let $(\cZ^+_t)_{t\geq 0}$ be the $\R^m$-valued stochastic process defined by $\cZ^+_t=\frk Z_{P_t}$. Let $\cZ^-$ be an independent copy of $\cZ^+$ and define $(\cZ_t)_{t\in\R}$ by
$$
\cZ_t = -\cZ_{(-t)^-}^-\1_{t<0}+\cZ_{t}^+\1_{t>0}.
$$

Let us now construct the coupling between $Z$ and $\cZ$ by applying the coupling from Section \ref{sec:kmt1} $d$ times, once for each of the allowed jumps. For $t\in\R$ define $\Delta_t\cZ:=\cZ_t-\cZ_{t^-}$, so in particular $\Delta_t\cZ\neq 0$ if and only if $\cZ$ jumps at time $t$. 
For $\ell=1,\dots,\d$ and $t\in\R$ let 
$$
P^\ell_t = 
\1_{t>0}\cdot\#\{s\in[0,t]\,:\,\Delta \cZ_s=v_\ell  \}
- \1_{t<0}\cdot\#\{s\in(t,0]\,:\,\Delta \cZ_s=v_\ell  \}.
$$
In other words, $(P^\ell_t)_{t\in\R}$ is the unique piecewise constant right-continuous real-valued process such that $P^\ell_0=0$ and $P^\ell$ makes an upward jump of 1 exactly at the times at which $\cZ$ makes a jump of $v_\ell$. 
Then $P^1,\dots,P^\d$ have the law of two-sided independent Poisson processes and the rate of $P^\ell$ is equal to $p_\ell$, see e.g.\, \cite[Example 3.6.7]{durrett}. 
For each $\ell\in[\d]$ define $\cY^\ell_t=P^\ell_t-p_\ell t$ and observe that
\eqb
\cZ=\sum_{\ell=1}^\d P^\ell v_\ell
=\sum_{\ell=1}^\d (P^\ell-p_\ell t)v_\ell 
= \sum_{\ell=1}^\d \cY^\ell v_\ell.
\label{eq97}
\eqe
Couple each process $\cY^\ell$ with a one-dimensional Brownian motion $Y^\ell$ of variance $p_\ell$ as in Proposition \ref{prop-kmt0}. 
Suppose the same uniform dyadic decomposition $\frk D$ is used for all $\ell\in[\d]$ and that $Y^1,\dots,Y^d,\frk D$ are jointly independent. 
Define $(Z_t)_{t\in\R}$ by
\eqb
Z_t:=\sum_{\ell=1}^\d Y^\ell_t v_\ell.
\label{eq69}
\eqe
Then $(Z_t)_{t\in\R}$ is a two-sided $m$-dimensional Brownian motion such that 
$$
\Cov(Z^i_1,Z^j_1)=\sum_{\ell=1}^\d v^i_\ell v^j_\ell p_\ell=\Cov(\frk Z^i_1,\frk Z^j_1),\qquad i,j\in[m],
$$ 
where we write $Z=(Z^1,\dots,Z^m)$ and  $v_\ell=(v_\ell^1,\dots,v_\ell^m)$. This concludes our construction of the coupling of $Z$ and $\cZ$.

Before stating our main result (Theorem \ref{prop-kmt}) we need to introduce some notation. Similarly as in \eqref{eq:X-def}, for each $J=[a,b)\in\frk D$ and $\ell\in[d]$ define
\eqb
X_{J}:= (b-a)^{-1/2}(2Z_{0.5(a+b)}-Z_a-Z_b),\qquad
X^\ell_{J}:= (b-a)^{-1/2}(2Y^\ell_{0.5(a+b)}-Y^\ell_a-Y^\ell_b).
\label{eq:X-def2}
\eqe
For each $\ell\in[\d]$  and $n\in\N_0$ define the process $\cY^{\ell,n}:I(n)\to\R$ by
\eqb
\cY^{\ell,n} := F_{p_\ell}(Y^\ell|_{I(n)}, S_\ell(n)),\qquad
S_\ell(n):= 2^n p_\ell + \sum_{j\in\{1,\dots,2n \}} 2^{(n-j)/2} V_{n+j} X^\ell_{I(n+j)}.
\label{eq:cYm2}
\eqe
Recall that for all $k\in\N$, a.s.,
\eqb
\cY^\ell|_{I(k)}=\lim_{n\rta\infty} \cY^{\ell,n}|_{I(k)}.
\label{eq:def-walk}
\eqe
The following theorem is an immediate consequence of Proposition \ref{prop-kmt0}.
\begin{theorem}
	In the setting above the following hold. 
	\begin{compactitem}
		\item[(i)] $(Z,\cZ,\frk D)$
		is translation invariant in law in the sense that for any fixed $t_0\in\R$,
		$$
		( (Z_{t+t_0}-Z_{t_0})_{t\in\R},  (\cZ_{t+t_0}-\cZ_{t_0})_{t\in\R},
		\frk D-t_0 ) \eqD (Z,\cZ,\frk D).
		$$
		\item[(ii)] For any $k\in\N_0$, $m\in\N\cap[2\log_2(k\vee 1),\infty)$, and $\ell\in\{1,\dots,d \}$, for some constant $c>0$,
		\eqbn
		\begin{split}
			\P\big[ \cY^\ell|_{I(k)}\neq \cY^{\ell,k+m}|_{I(k)} \big]&<\exp(-cm).
		\end{split}
		\eqen
		\item[(iii)] For any $C_1>1$ we can find $C'>1$ such that for any $k,n\in\N$ and $\ell\in\{1,\dots,d \}$,
		\eqb
		\P\Big[ \sup_{t\in I(n)} |Y^\ell_t-\cY^{\ell,n+k}_t|>C'n \Big] \leq 2^{-C_1n}
		\label{eq:err3}
		\eqe
		and
		\eqb
		\P\Big[ \sup_{t\in I(n)} |Z_t-\cZ_t|>C'n \Big] \leq 2^{-C_1n}.
		\label{eq:err4}	
		\eqe
	\end{compactitem}
\label{prop-kmt}
\end{theorem}
By the definition in \eqref{eq:def-walk0}, for each $\ell\in[d]$ the function $\cY^\ell$ can be written as a measurable function of $Y^\ell$ and $\frk D$. In contrast, $\cZ$ is generally \emph{not} determined by $Z$ and $\frk D$ since $(Y^1,\dots,Y^d)$ is generally not determined by $Z$ and $\frk D$. In the next lemma we will see that $\cZ$ is determined by $Z$, $\frk D$ and some independent randomness denoted by $\frk U$. It will be convenient for us in later sections to represent $\cZ$ in this way.
\begin{lem}
	It is possible to define random variables $\frk U=(U_J)_{J\in\frk D},Z,(Y^1,\dots,Y^{d})$ on the same probability space, along with measurable functions $\cU$ and $\wt {\cU}$, such that $(Z,(Y^1,\dots,Y^{d})\big)$ has the same law as above and the following hold.
	\begin{compactitem}
		\item[(i)] Conditioned on $\frk D$, the random variable $\frk U=(U_{J})_{J\in\frk D} $ is a collection of independent standard $\dZ$-dimensional normal random vectors indexed by $J\in\frk D$. Furthermore, $\frk U$ and $Z$ are independent. 
		\item[(ii)] $(Y^1,\dots,Y^{d}) = {\cU}(Z,\frk U)$ a.s.
		\item[(iii)] $(X_J^1,\dots,X_J^{d})=\wt {\cU}(X_J,U_J)$ a.s.\ for all $J\in\frk D$.
	\end{compactitem}
	\label{prop:split-rand}
\end{lem}
\begin{proof}
	Let $Z,\frk D$ be independent and have the desired marginal laws, and conditional on $Z$ and $\frk D$ let $\frk U$ be a collection of independent standard $d$-dimensional normal random variables indexed by $\frk D$. In order to conclude, it is sufficient to argue that 
	we can find functions $\cU$ and $\wt{\cU}$ such that (ii) and (iii) are satisfied, 
	$(Y^1,\dots,Y^{d})$ has the desired law, and $(Y^1,\dots,Y^{d})$ is coupled with $Z$ and $\frk D$ in the desired way.
	
	Let $A$ denote the matrix with columns $v_\ell$, let $D$ be the diagonal matrix with entries $p_1,\dots,p_{d}$ on the diagonal and set $\Sigma=ADA^T$. Then define $\wt{\cU}(x,u)$ such that
	$$
	\wt{\cU}(x,u)^T = DA^T\Sigma^+x + D^{1/2}(I-D^{1/2}A^T\Sigma^+AD^{1/2})u,
	$$
	where $\Sigma^+$ is the pseudoinverse of $\Sigma$. Then $\wt{\cU}(X_J,U_J)$ gives the conditional law of $(X_J^1,\dots,X_J^d)$; this can be checked by verifying that if we define $(X_J^1,\dots,X_J^d)$ as in Assertion (iii), then $(X_J^1,\dots,X_J^d)$ has the desired marginal law and the desired joint law with $X_J$, and additionally that $X_J= X^1_J v_1+\dots+X^\d_J v_\d$.
	
	Define $(X_J^1,\dots,X_J^{d})=\wt {\cU}(X_J,U_J)$ for all $J\in\frk D$ and observe that $\{ X_J^\ell\,:\,J\in\frk D,\ell\in[d] \}$ is a collection of independent normal random variables. Define $Y^1,\dots,Y^{d}$ from the $X_J^\ell$'s by using \eqref{eq:Z-brZ-X} repeatedly along with continuity of $Y^1,\dots,Y^{d}$. Equivalently, we can define
	\eqb
	\cU((z_t)_{t\in\R}, (u_J)_{J\in\frk D})^T = DA^T\Sigma^+z_t+D^{1/2}(I-D^{1/2}A^T\Sigma^+AD^{1/2})w_t,
	\label{eq:def-cU}
	\eqe
	where $(w_t)_{t\in\R}$ is the Brownian motion such that
	\eqref{eq:X-def} holds with $X_J$ and $Y$ replaced by $u_J$ and $w$. 
	This shows that we can find an appropriate function $\cU$ satisfying Assertion (ii) such that $(Y^1,\dots,Y^{d})$ has the desired law.
	Finally, $(Y^1,\dots,Y^{d})$ and $(Z,\frk D,\frk U)$ are coupled together as desired since \eqref{eq69} is satisfied.

\end{proof}

\begin{remark}
	The random variable $\frk U$ is a collection of random variables indexed by $\frk D$. This means in particular that $\frk D$ is measurable with respect to $\frk U$. This measurability is necessary in order for Assertion (ii) of the above lemma to hold and will also be important later when we say that certain random variables can be written as a measurable function of $\frk U$ and other random variables. There are two sources of randomness used to define $\frk U$: The randomness of $\frk D$ and, conditioned on $\frk D$, we associate each element of $\frk D$ with an independent $d$-dimensional standard normal random vector.
	\label{rmk:U-determines-D}
\end{remark}

\begin{remark}
	Lemma \ref{prop:split-rand} implies that $\cZ$ can be written as a natural measurable function of $(Z,\frk U)$. Namely, first define $(Y^1,\dots,Y^d):={\cU}(Z,\frk U)$, then define $(\cY^1,\dots,\cY^d)$ by \eqref{eq:def-walk0} with $(Y,\cY)$ replaced by $(Y^\ell,\cY^\ell)$ for $\ell=1,\dots,d$, and finally define $\cZ:\R\to\R^m$ by $\cZ=\sum_{\ell=1}^{d}\cY^\ell v_\ell$.
	\label{rmk:split-rand}
\end{remark}

For a dyadic system $\frk D$ and $t_0\in\R$ recall the definition of $\frk D - t_0$ from \eqref{eq:D-shift}. We will now define $\frk U-t_0$ in a similar fashion for $\frk U=(U_J)_{J\in\frk D}$. Let $\frk U - t_0$ denote the set $(\wt U_J )_{J\in\frk D-t_0}$ where $\wt U_J:=U_{J+t_0}$. In particular, $\frk U-t_0$ is indexed by $\frk D-t_0$. From this definition it is immediate that with notation as in Lemma \ref{prop:split-rand}, if $(\wt Y^1,\dots,\wt Y^{d}):=\cU(Z_{t+t_0}-Z_{t_0},\frk U-t_0)$, for all $t\in\R$,
\eqb
(\wt Y^1_t,\dots,\wt Y^{d}_t) = (Y^1_{t+t_0},\dots, Y^{d}_{t+t_0})-(Y^1_{t_0},\dots,Y^{d}_{t_0}).
\label{eq:U-translate}
\eqe 

The following two lemmas will be helpful in later sections.
\begin{lem}
	For any $\ell\in[d]$ and $k\in\N$ the random variable $S_\ell(k)$ and the process $\cY^{\ell,k}$ are measurable with respect to $Z|_{I(3k)}$ and $\frk U$.
	\label{eq:cYell-measurable}
\end{lem}
\begin{proof}
	First observe that $(X_J)_{J\in\frk D,J\subset I(3k)}$ is measurable with respect to $Z|_{I(3k)}$ by \eqref{eq:X-def2}, which implies that $(X^\ell_J)_{J\in\frk D,J\subset I(3k)}$ is measurable with respect to $Z|_{I(3k)},\frk U$ by Assertion (iii) of Lemma \ref{prop:split-rand}. The lemma now follows from the definition \eqref{eq:cYm2} of $S_\ell(k)$ and $\cY^{\ell,k}$.
\end{proof}

The next lemma is an immediate consequence of \eqref{eq:Z-brZ-X}.
\begin{lem}
	For any $j\in\Z$,
	\eqb
	\Delta_{I(j)}Z = \sum_{i=1}^\infty  2^{(j-i)/2}V_{i}^{I(j)} X_{I(j+i)}.
	\label{eq:Z-brZ-X2}
	\eqe
\end{lem}

\section{Strong coupling with random walk for the 0-quantum cone}
\label{sec:kmt-0-cone}
In Section 3.3 we gave a strong coupling between a possibly higher-dimensional Brownian motion $Z$ and a random walk $\cZ$. We now specialize to the two-dimensional case, so that $Z=(Z_t)_{t\in\R}$ is a correlated two-sided two-dimensional Brownian motion and $\cZ=(\cZ_t)_{t\in\R}$ is a continuous-time random walk with the same variances and covariances as $Z$. 
More precisely, we suppose that $\Var(Z_1^1)=\Var(Z_1^2)=1/2$ and $\operatorname{Cov}(Z_1^1,Z_1^2)=-\tfrac12\cos(\pi\gamma^2/4)$ for some fixed $\gamma\in(0,2)$. 
 Recall that the set of times at which $\cZ$ makes a jump is given by a unit rate Poisson point process on $\R$ and the steps of $\cZ$ are i.i.d.\ conditioned on this Poisson point process such that $\Delta_t\cZ$ takes $\dZ\in\N$ distinct values.\footnote{We write $\dZ$ instead of $d$ in the remainder of the paper.}

By mating-of-trees theory, $Z$ encodes a $\gamma$-quantum cone with a space-filling SLE (see Section \ref{sec:mot}). For several different choices of $\cZ$, this walk encodes a natural infinite planar map model, see Section \ref{sec:bijection} for the special case of a tree-weighted planar map and the walk whose steps are uniformly distributed on $\{(1,0),(-1,0),(0,1),(0,-1) \}$. Therefore, the coupling mentioned in the previous paragraph also gives a coupling between an SLE-decorated $\gamma$-quantum cone and various infinite planar map models.

In this subsection we will extend this coupling to the case of the SLE-decorated 0-quantum cone (instead of the SLE-decorated $\gamma$-quantum cone) and establish a general result which allows us to compare the probability of general events in the two settings (the 0-quantum cone and the $\gamma$-quantum cone). The main results of this subsection are Theorem \ref{prop:zero-vs-gamma-walk},  Corollary \ref{cor:zero-vs-gamma}, and Proposition \ref{prop:zero-vs-gamma-ball}. A priori it is not obvious how one would proceed to extend the coupling to the 0-quantum cone case since the $\gamma$-quantum cone and the 0-quantum cone are singular with respect to each other. The idea of our approach is to use that in the $\gamma$-quantum cone case, $\cZ$ can be expressed as an explicit function of $Z$ (and an independent random variable $\frk U$; see Remark \ref{rmk:split-rand}) and that we get a natural coupling in the 0-quantum cone case by applying the exact same function, which turns out to be a.s.\ well-defined also in the 0-quantum cone case. Furthermore, we are able to compare the two settings by using that if we have a 0-quantum cone $(\C,0,\infty)$ and sample a point $z_0$ from the LQG area measure then $(\C,z_0,\infty)$ is absolutely continuous with respect to a $\gamma$-quantum cone when restricted to appropriate bounded subsets of $\C$. Using the latter observation we get Proposition \ref{prop:zero-vs-gamma-ball} relatively easily, but proving existence of $\cZ$ still requires significantly more work; see the beginning of Section \ref{sec:0-cone-walk-exists} for a more detailed proof outline for this result. 

The following definition will be used frequently in the remainder of the paper.
\begin{definition}\label{def:PQ-PE}
	Let $\PE$ be the law of $(h,\eta,\frk U)$ where $(\C,h,0,\infty)$ is an (unspecified) embedding of the 0-quantum cone, $\eta$ is a whole-plane space-filling SLE$_\kappa$, 
	$\frk U$ is as in Lemma \ref{prop:split-rand}, and the three random variables $h,\eta,\frk U$ are all independent. Let $\PQ$ be defined in the exact same way, except that $(\C,h,0,\infty)$ is an (unspecified) embedding of the $\gamma$-quantum cone instead of the 0-quantum cone. Under either of these probability measures we define the boundary length process $Z$ of $(h,\eta)$  as in Section \ref{sec:mot}.
	\label{PE-PQ}
\end{definition}
Note that in the above definition we consider an arbitrary embedding of the 0-quantum cone and the $\gamma$-quantum cone, and most of the results we state will hold for any choice of embedding of the surface. However, sometimes we will need to fix the embedding; this will be specified when needed. 

Let us explain more explicitly how we define the random walk $\cZ$ in the 0-quantum cone case. Let $(h,\eta,\frk U)$ be sampled from $\PE$ and let $Z$ be the boundary length process of $(h,\eta)$. Define $$(Y^1,\dots,Y^{\dZ})={\cU}(Z,\frk U),$$ where ${\cU}$ is as in Lemma \ref{prop:split-rand}. Define $X_I$ and $X^\ell_I$ for $I\in\frk D$ and $\ell\in[\dZ]$ by \eqref{eq:X-def2}. For $m,S\in\N_0$ and $\ell\in[\dZ]$ define $\cY^{\ell,m,S}:I(m)\to\R$, $\cY^{\ell,m}:I(m)\to\R$, and $S_\ell(m)\in\R$ as in \eqref{eq:cYm2}. Recall Lemma \ref{eq:cYell-measurable}, which says that $\cY^{\ell,m}$ is  measurable with respect to $Z|_{I(3m)}$ and $\frk U$. Let $\cY^\ell:\R\to\R$  be the process satisfying \eqref{eq:def-walk}, %
provided such a process exists and the limit is well-defined for $k\in\N$. We say that 
\eqb
\cZ:=\sum_{\ell\in[\dZ]} \cY^\ell v_\ell
\label{eq:cZ}
\eqe 
is well-defined if $\cY^\ell$ is well-defined for each $\ell\in[\dZ]$. 
\begin{theorem} 
	Under the probability measure $\PE$ the walk $\cZ$ given by \eqref{eq:cZ} is well-defined. 
	Let $G:D(\R,\R^2)\times\Omega_{\frk U}\to\{0,1 \}$ be a measurable function and define $A$ by 
	$$
		A = \{ \exists s\in I(0) \,:\, G(\cZ(\cdot+s)-\cZ(s),\frk U-s)=1 \}.
		$$
	We require that $A\in \sigma(\cZ|_{I(k)},I(k))$. 
	Then there exists a constant $C_0>1$ such that for all $k\in\N$,
	\eqb
	\PE[A] \leq \PQ[A]^{1/4} 2^{C_0k} + 2^{-k/C_0}.
	\label{eq:zero-vs-gamma-walk}
	\eqe
	Furthermore, for any $c>0$ there is a $C_1>1$ such that for any $m\in\N$ and $k\geq cm$ we have
	$$
	\PE\big[\exists n\in [m+k,\infty)_{\Z}  ,
	\ell\in[\dZ]
	\text{\,\,s.t.\,\,}
	\cY^\ell|_{I(m)} \neq
	\cY^{\ell,n}|_{I(m)}\big]
	\leq 2^{-k/C_1},
	$$
	where we recall that $\cY^{\ell,n}$ is defined by \eqref{eq:cYm2}.
	Moreover, for any $C_2>1$ we can find $\wt C>1$ such that for any $m\in\N$, 
	\eqb
	\PE\bigg[\sup_{t\in[-2^m,2^m]}|Y_t^\ell-\cY_t^\ell|>\wt C m\bigg] \leq 2^{-C_2 m}.
	\label{eq56}
	\eqe
	Finally, there is a measurable function $g$ such that $\cZ=g(h,h^{\op{IG}},\frk U)$ and this function is translation invariant modulo scaling in the sense that a.s.,  
		\eqb
		(\cZ_{t+\tau_z}-\cZ_{\tau_z})_{t\in\R} = g( h(C\cdot+z)+Q\log C,h^{\op{IG}}(C\cdot +z),\frk U-\tau_z )\qquad\forall C>0, z\in\C,
		\label{eq:transl0}
		\eqe
		where we refer to \eqref{eq:U-translate} for the definition of $\frk U-\tau_z$.
	\label{prop:zero-vs-gamma-walk}
\end{theorem}
Let us describe the event $A$ in the above theorem briefly in words. The function $G$ is typically chosen such that $G(\cZ,\frk U)=1$ if $(\cZ,\frk U)$ satisfies some undesirable and rare event. The event $A$ occurs if this bad event occurs for $(\cZ(\cdot+s)-\cZ(s),\frk U-s)$ for some $s\in I(0)$. The estimate in \eqref{eq:zero-vs-gamma-walk}, which will play a crucial role in the paper, is natural due to the fact that, as discussed earlier in this section, it is possible to find some point $z_0$ close to the origin such that a 0-quantum cone is locally absolutely continuous with respect to a $\gamma$-quantum cone if recentered at the point $z_0$. The factor $2^{C_0k}$ on the right side of this estimate represents the Radon-Nikodym derivative, while $2^{-k/C_0}$ is an error term. A more detailed proof outline can be found above Lemma \ref{lem:trans-inv}. 

The estimate \eqref{eq:zero-vs-gamma-walk} is powerful, but one difficulty with applying it, is that we require $A\in \sigma(\cZ|_{I(k)})$. For many events we will consider, this measurability result does not hold for any deterministic $k\in\N$. On the other hand, it often holds that our event $A$ is measurable with respect to $\cZ|_{I(k)}$ for large $k$ with high probability. In such situations we can instead apply the following result, which is a corollary of Theorem \ref{prop:zero-vs-gamma-walk}.
\begin{cor}
	Let $G,G_1,G_2,\dots:D(\R;\R^2)\times\Omega_{\frk U}\to\{0,1 \}$ and define
	\eqb
	\begin{split}
		A&=\{\exists s\in I(0)\,:\,G(\cZ(\cdot+s)-\cZ(s),\frk U-s)=1 \},\\
		A_k&=\{\exists s\in I(0)\,:\,G_k(\cZ(\cdot+s)-\cZ(s),\frk U-s)=1 \},\quad k\in\N.
	\end{split}
	\label{eq:Ak}
	\eqe	
	Suppose $G,G_1,G_2,\dots$ are such that $A_k\in\sigma(\cZ|_{I(k)},I(k))$ for any $k\in\N$, $A\in\sigma(\cZ,\frk U)$, and $A_k\subset A_{k+1}\subset A$ for any $k\in\N$.
	Also suppose $\PE[\1_{A_k}\neq \1_A]\leq \wt C 2^{-\wt ck}$ for some $\wt c\in(0,1)$, $\wt C>1$ and all $k\in\N$. Then there is an $\eps>0$ independent of $(\wt c,\wt C)$ such that
	$$
	\PE[A] \leq \wt C(\wt c\eps)^{-1}\PQ[A]^{\wt c \eps}.
	$$  
	\label{prop:eq0}
	\label{cor:zero-vs-gamma}
\end{cor}
The following proposition is a variant of Theorem \ref{prop:zero-vs-gamma-walk} where the considered event depends on $h$ and $h^{\op{IG}}$ restricted to some ball instead of $\cZ|_{I(k)}$. In a similar spirit as above, the event $A(K)$ means that some bad event occurs for $h(\cdot-z),h^{\op{IG}}(\cdot-z)$ for some $z\in K$. 
\begin{prop} 
	Let $k\in\N$ and suppose now that the quantum cones under the probability measures $\PE$ and $\PQ$ have the $2^{2k}$-embedding. 
	Let $G: D'(\C)^2\to\{0,1 \}$ be a measurable function, let $h^{\op{IG}}$ denote the GFF generating $\eta$ (see Section \ref{sec:mot}),
	and for a compact set $K\subseteq \ol{B(0,2^{-k})}$ define $A(K)$ by 
	$$
		 A(K) = \{ \exists z\in K \,:\, G(h(\cdot-z),h^{\op{IG}}(\cdot-z))=1 \}.
	$$
	We require that if $K$ is a closed and possibly random ball  contained in $\ol{B(0,2^{-k})}$ then $A(K)\in \sigma(h|_{B(0,2^{2k}-2^{-k+1})},h^{\op{IG}}|_{B(0,2^{2k}-2^{-k+1})})$.
	Then there exists a constant $C_0>1$ such that for all $k\in\N$,
	\eqb
	\PE[A(\{0 \} )] \leq \PQ\Big[A\Big( \ol{B(0,2^{-k})} \Big)\Big]^{1/2} 2^{C_0k}.
	\label{eq:zero-vs-gamma}
	\eqe
	\label{prop:zero-vs-gamma-ball}
\end{prop}
The latter proposition is significantly easier to prove than Theorem \ref{prop:zero-vs-gamma-walk} and Corollary \ref{cor:zero-vs-gamma} above. It can be proven via a brief argument once we have established Lemma \ref{lem:rooted} right below. Lemma \ref{lem:rooted} makes precise the result discussed above which says that one can find a point $z_0$ close to the origin such that a 0-quantum cone is locally absolutely continuous with respect to a $\gamma$-quantum cone if recentered at $z_0$. This allows us to say that if some event $A(0)$ is likely to happen for the 0-quantum cone then it is also fairly likely to happen for the $\gamma$-quantum cone \emph{after recentering the $\gamma$-quantum at some point of our choice in a small ball around the origin}. The factor $2^{C_0k}$ on the right side of \eqref{eq:zero-vs-gamma} comes from the Radon-Nikodym derivative in the next lemma.

\begin{lem}
	For $R,k>1$ let $h$ be the $2^R$-embedding of a 0-quantum cone. Let $U=B(0,2^{-k})$ and let $z_0$ be sampled from $\mu_h|_{U}$ renormalized to be a probability measure. Conditioned on $z_0$, the following two laws on fields are the same:
	\begin{compactitem}
		\item the law of $h(\cdot+z_0)|_{B(0,2^{R} -2^{-k} )}$, weighted by $\mu_h(U)\E[\mu_h(U)]^{-1}$ and
		\item the law of the $2^R$-embedding of a $\gamma$-quantum cone restricted to $B(0,2^{R}-2^{-k})$.
	\end{compactitem}
	Furthermore, for any $p>0$ there is a constant $C_p>1$ such that
	\eqb
	\E[\mu_h(U)^{-p}] \leq \exp(C_p (R+k))
	\label{eq:annulus-moment}
	\eqe
	Finally, for any $p\in(0,4/\gamma^2)$ there are $C>1$ and $c\in(0,1)$ such that  $\E[\mu_h(U)^p]\leq C\exp(CR)$.
	\label{lem:rooted}
\end{lem}
The proof of the lemma is deferred to Appendix \ref{app2}. Using the lemma we can now complete the proof of Proposition \ref{prop:zero-vs-gamma-ball}.
\begin{proof}[Proof of Proposition \ref{prop:zero-vs-gamma-ball}]
	Enlarge the probability space associated with $\PE$ so $\PE$ is a probability measure on tuples $(h,\eta,\frk U,z_0)$, where $z_0$ is sampled from the ball $U$ as in Lemma \ref{lem:rooted} with $R=2k$.  Define $\wh{\P}_Q$ by $\xi d\wh{\P}_Q=d\PE$ for $\xi=C_1\mu(U)^{-1}$ and  
$C_1=\E_{\op{E}}[\mu(U)]>0$ 
a normalizing constant chosen such that $\wh{\P}_Q$ is a probability measure. Then  
\eqb
\begin{split}
	\PE[A(0)] 
	=\wh\E_{\op{Q}}[A(0);\xi] 
	\leq
	\wh{\P}_{\op{Q}}[A(0)]^{1/2} \wh\E_{\op{Q}}[\xi^2]^{1/2}.
\end{split}
\label{eq11}
\eqe
Therefore it is sufficient to argue that the following hold for some $C_0>1$ 
\eqb
\begin{split}
	(i)\,\,\wh{\P}_{\op{Q}}[A(0)]\leq \PQ[A(\ol{B(0,2^{-k})})],\qquad
	(ii)\,\,\wh{\E}_{\op{Q}}[\xi^2] \leq 2^{C_0k}.
\end{split}
\eqe
We get (ii) by using the definitions of $\wh \E_{\op{Q}}$ and $\xi$ and then Lemma \ref{lem:rooted}
\eqb
\wh\E_{\op{Q}}[\xi^2]
= \EE[\mu_h(U)^{-1}]\EE[\mu_h(U)] %
\leq 2^{C_0 k}.
\label{eq14}
\eqe
Let $(h,h^{\op{IG}},\frk U,z_0)$ be sampled from $\wh{\P}_{\op{Q}}$ and let $(h^Q,h^{\op{IG,Q}},\frk U^Q)$ be sampled from $\PQ$. By Lemma \ref{lem:rooted}, we have $h(\cdot+z_0)|_{B(0,2^{2k}-2^{-k})}\eqD h^Q|_{B(0,2^{2k}-2^{-k})}$. Since $(h,z_0)$ is independent of $h^{\op{IG}}$ and $h^{\op{IG}}$ is translation invariant in law we have $h^{\op{IG}}(\cdot+z_0)\eqD h^{\op{IG,Q}}$ and independence of $h(\cdot+z_0)$ and $h^{\op{IG}}(\cdot+z_0)$. This implies
\eqb
	(h(\cdot+z_0), h^{\op{IG}}(\cdot+z_0))|_{B(0,2^{2k}-2^{-k})}
	\eqD 
	(h^Q,h^{\op{IG,Q}})|_{B(0,2^{2k}-2^{-k})}.
	\label{eq:law-4tuple-h}
\eqe
We now obtain (i) by
\eqbn
\begin{split}
	\wh\P_{\op{Q}}[ A(0)]
	&= \wh\P_{\op{Q}}[G(h,h^{\op{IG}})=1]
	\leq \wh\P_{\op{Q}}[\exists z\in \ol{B(z_0,2^{-k})}\,\,\text{s.t.}\,\,G(h(\cdot+z),h^{\op{IG}}(\cdot+z))=1  ]\\
	&\leq\PQ[\exists z\in \ol{B(0,2^{-k})}\,\,\text{s.t.}\,\,G(h(\cdot+z),h^{\op{IG}}(\cdot+z))=1]=\PQ[A(\overline{B(0,2^{-k})})],
\end{split}
\eqen
where the first inequality follows since $0\in B(z_0,2^{-k})$ and the second inequality follows from \eqref{eq:law-4tuple-h}.	
\end{proof}

\begin{lem}
	Under $\PE$, the process $Z$ satisfies \eqref{eq:Z-brZ-X2} and the processes $Y^1,\dots,Y^{\dZ}$ satisfy \eqref{eq:Z-brZ-X}.
\end{lem}
\begin{proof}
	Throughout the proof we let $c\in(0,1)$ denote a constant which might decrease from line to line. 
	As in the proof of Lemma \ref{prop:br-indep}, we need to argue that $2^{-N}\Delta_{J_N^I}Y\rta 0$ a.s.\ as $N\rta\infty$. By first doing a union bound and then using Proposition \ref{prop:zero-vs-gamma-ball} (assuming we work with the $2^{(1-2\delta)N/(\gamma Q)}$-embedding), the definition of a uniform dyadic system and Lemma \ref{prop:contained1}, for arbitrary $\eps>0$,
	\eqb
	\begin{split}
		&\PE[ \Delta_{J_N^I}Y>\eps 2^N ]
		\leq 
		\PE[\Delta_{J_N^I}Y>\eps 2^N;J_N^I\subset I((1-\delta)N)\subset B(0,2^{(1-2\delta)N/(\gamma Q)})] \\
		&\qquad\qquad\qquad+\PE[J_N^I\not\subset I((1-\delta)N)]
		+\PE[I((1-\delta)N)\not\subset B(0,2^{(1-2\delta)N/(\gamma Q)})]\\
		&\qquad\leq 2^{C_0(1-2\delta)N/(\gamma Q)/2}\PQ[\Delta_{J_N^I}Y>\eps 2^N;J_N^I\subset I((1-\delta)N)\subset B(0,2^{(1-2\delta)N/(\gamma Q)})]^{1/2} + 2^{-cN}.
	\end{split}
	\eqe
	By using that $Y$ has the law of a Brownian motion under $\PQ$, we get that the probability on the right side is bounded above by $\exp(-c\eps^22^N)$, which means that the right side of the last display is bounded above by $2^{-cN}$.
\end{proof}

The next two lemmas will be used in the proof of Lemma \ref{prop:zero-vs-gamma0} below.
\begin{lem}
	Let $(h,h^{\op{IG}},\frk U)$ be sampled from $\PE$ or $\PQ$ and let $z_0\in\C$ be sampled in an arbitrary way such that $(h,\eta,z_0)$ is independent of $\frk U$. 
	Let $\tau_{z_0}:=\inf\{t\in\R\,:\,\eta(t)=z_0 \}$ be the first time at which $\eta$ hits $z_0$. 
	Then $\frk U-\tau_{z_0}\eqD\frk U$ and the random variables $(h,\eta,z_0)$ and $\frk U-\tau_{z_0}$ are independent.
	\label{lem:trans-inv}
\end{lem}
\begin{proof}
	Recall that we can write $\frk D=\frk D(\theta,V_1,V_2,\dots)$ for $\theta\in[0,1)$ uniform and $V_1,V_2,\dots$ i.i.d.\ and uniform on $\{-1,1 \}$. Note that $\theta+\tau-\lfloor \theta+\tau \rfloor$ is uniform on $[0,1]$ and independent of $\tau$ for any real-valued random variable $\tau$ independent of $\theta$. 
	This implies that $\frk D-\tau_{z_0}\eqD \frk D$ and that the random variables
	$(h,\eta,z_0)$ and $\frk D-\tau_{z_0}$ are independent. We get the lemma by using that 
	the standard normal random vectors defining $\frk U$ are independent of all other randomness.
\end{proof}

\begin{lem}
	Let $(h,h^{\op{IG}},\frk U)$ be sampled from $\PE$ or $\PQ$. Let $m\in\N$ and $t_0\in\R$. For each $j\in\N$ let $\wt I(j)$ be the unique length $2^j$ interval in $\frk D-t_0$ containing 0 and define for $\ell\in[\dZ]$ and $a<b$,
	$$
	(\wt Y^1,\dots,\wt Y^{\dZ}) = {\cU}( (Z_{t+t_0}-Z_{t_0})_{t\in\R}, \frk U-t_0 ),\quad 
	\wt{\cY}^{\ell,m}:=F_{p_\ell}(\wt Y^\ell|_{\wt I(m)}, \wt S_\ell(m) ),
	$$
	$$
	\wt X^\ell_{[a,b)}:= (b-a)^{-1/2}(2\wt Y^\ell_{0.5(a+b)}-\wt Y^\ell_a-\wt Y^\ell_b),\qquad
	\wt S_\ell(n):= 2^n p_\ell + \sum_{j\in\{1,\dots,2n \}} 2^{(n-j)/2} \wt V_{n+j} \wt X^\ell_{\wt I(n+j)},
	$$
	where $\wt V_1,\wt V_2,\dots$ are defined such that $\frk D-t_0=\frk D(\wt\theta,\wt V_1,\wt V_2,\dots)$ for some $\wt\theta\in[0,1)$.
	Then it holds a.s.\ simultaneously for all $t_0\in I(m)$ that
	\eqb
	\wt{\cY}^{\ell,m}=(\cY^{\ell,m}_{t+t_0}-\cY^{\ell,m}_{t_0})_{t\in I(m)-t_0}.
	\label{eq:cZ-shift}
	\eqe
\end{lem}
\begin{proof}
	We saw in the proof of Lemma \ref{prop:split-rand} that the function $\cU$ can be expressed explicitly by using the function $\wt{\cU}$ and \eqref{eq:Z-brZ-X}. It follows from this and $\wt X^\ell_{[a,b)}=X^\ell_{[a+t_0,b+t_0)}$ that $\wt Y^\ell = (Y^\ell_{t+t_0}-Y^\ell_{t_0})_{t\in\R}$ a.s.\ simultaneously for all $t_0$. We get further from the definition of $F_{p_\ell}$ and $\wt S_\ell(n)$ that \eqref{eq:cZ-shift} holds for $t_0\in I(m)$.
\end{proof}

We are now ready to prove a first variant of Theorem \ref{prop:zero-vs-gamma-walk} by using similar ideas as in the proof of Proposition \ref{prop:zero-vs-gamma-ball}. The main difference between Proposition \ref{prop:zero-vs-gamma-ball} and Lemma \eqref{prop:zero-vs-gamma0} is that our event depends on a walk $\cZ^{2k}$ instead of the fields $(h,h^{\op{IG}})$. The main estimate \eqref{eq:zero-vs-gamma0} is of a similar form as \eqref{eq:zero-vs-gamma}, but we also need an error term $\PE[B_0]$ due to potential abnormal behavior of the SLE $\eta$ or the LQG surface. The key difference between Theorem \ref{prop:zero-vs-gamma-walk} and the below lemma is that below we consider the walk $\cZ^{2k}$, which is measurable with respect to $(Z|_{I(6k)}, \frk U)$, while in Theorem \ref{prop:zero-vs-gamma-walk} we consider the walk $\cZ$, which we do not yet know is well-defined and which we need the full boundary length process $Z$ to define. Let  $\Omega_{\frk U}$ denote the measurable space in which $\frk U$ takes its value.
\begin{lem}
	Let $(h,h^{\op{IG}},\frk U)$ be sampled from $\PQ$ or $\PE$. Let $k\in\N$ and let $\wt{\cZ}^{2k}=\sum_{\ell\in[\dZ]} v_\ell\cY^{\ell,2k}$, see \eqref{eq:cYm2} for the definition of the right side. Let $\wt G:D(\R,\R^2)^2\times\Omega_{\frk U}\to\{0,1 \}$  be a measurable function and define the event $\wt A$ by 
	$$
		\wt A = \{ \exists s\in I(0) \,:\,  \wt G(\wt{\cZ}^{2k}(\cdot+s)-\wt{\cZ}^{2k}(s),Z(\cdot+s)-Z(s),\frk U-s)=1 \}.
	$$
	We require that $\wt A\in \sigma(\wt{\cZ}^{2k},Z|_{I(6k)},\frk U)$. Suppose that $h$ has the $2^{6k}$-embedding.  
	Then for some constant $C_0>1$ depending only on $\gamma$,
	\eqb
	\begin{split}
		\PE[\wt A] &\leq  \PQ[\wt A]^{1/2}2^{C_0k} + \PE[B_0],
	\end{split}
\label{eq:zero-vs-gamma0}
	\eqe
	where 
	$$
	B_0:= \{ \eta[I(6k)]\not\subset B(0,2^{6k}-2^{1-k}) \}\cup \{ B(0,2^{-k})\not\subset \eta[I(0)]\}
	$$
	satisfies $\PE[B_0]\leq 2^{-k/C_0}$.
	\label{prop:zero-vs-gamma0}
\end{lem}
\begin{proof}
	Defining $\xi$ and $\wh\E_{\op{Q}}$ as in the proof of Proposition \ref{prop:zero-vs-gamma-ball} we have 
		\eqb
		\begin{split}
			\PE[\wt A] 
			&\leq
			\PE[\wt A;B_0^c] + \PE[B_0]
			=\wh\E_{\op{Q}}[\wt A;B_0^c;\xi] + \PE[B_0] \\
			&\leq
			\wh{\P}_{\op{Q}}[\wt A;B_0^c]^{1/2} \wh\E_{\op{Q}}[\xi^2]^{1/2} +\PE[B_0].
		\end{split}
		\eqe
		Therefore it is sufficient to argue that the following hold for some $C_0>1$ in order to get the lemma
		\eqb
		\begin{split}
			(i)\,\,\wh{\P}_{\op{Q}}[\wt A;B_0^c]\leq \PQ[\wt A],\qquad
			(ii)\,\,\wh{\E}_{\op{Q}}[\xi^2] \leq 2^{C_0k},\qquad
			(iii)\,\,\PE[B_0]\leq 2^{-k/C_0},
			\label{eq12}
		\end{split}
		\eqe
	We established (ii) in the proof of Proposition \ref{prop:zero-vs-gamma-ball} and
	Lemmas \ref{prop:contained1} and \ref{prop:contained2} imply (iii) so it remains to prove (i). 	
	
	Let $(h,h^{\op{IG}},\frk U)$ be sampled from $\wh{\P}_{\op{Q}}$ and let $(h^Q,h^{\op{IG,Q}},\frk U^Q)$ be sampled from $\PQ$. Let $\eta$ be the whole-plane space-filling SLE$_\kappa$ associated with $h^{\op{IG}}$ and set $\tau_{z_0}=\inf\{t\in\R\,:\,\eta(t)=z_0 \}$. Lemma \ref{lem:trans-inv} implies that 
		\eqb
		\begin{split}
			&(h(\cdot+z_0)|_{B(0,2^{6k}-2^{-k})}, h^{\op{IG}}(\cdot+z_0)|_{B(0,2^{6k}-2^{-k})}, \frk U-\tau_{z_0})\\
			&\qquad\qquad\qquad\qquad\qquad
			\eqD 
			(h^Q|_{B(0,2^{6k}-2^{-k})},h^{\op{IG,Q}}|_{B(0,2^{6k}-2^{-k})},\frk U^Q).
		\end{split}
		\label{eq:law-4tuple}
		\eqe
	Letting $H$ be a measurable function such that $H(\wt{\cZ}^{2k},Z,\frk U)=\1_{\wt A}$ we have
	\eqbn
	\begin{split}
		\wh\P_{\op{Q}}[\wt A;B_0^c]
		&= \wh\E_{\op{Q}}[H(\wt{\cZ}^{2k},Z,\frk U) \1_{B_0^c}]\\
		&=\wh\E_{\op{Q}}[H(
		\wt{\cZ}^{2k}(\cdot+\tau_{z_0})-\wt{\cZ}^{2k}(\tau_{z_0}),
		Z(\cdot+\tau_{z_0})-Z(\tau_{z_0}), 
		\frk U-\tau_{z_0}) \1_{B_0^c}] \\
		&\leq \E_{\op{Q}}[H(\wt{\cZ}^{2k},Z,\frk U)]
		= \PQ[\wt A],
	\end{split}
	\eqen
	where the second equality follows since on the event $B_0^c\subset\{ B(0,2^{-k})\subset \eta[I(0)]\}$ we have $\tau_{z_0}\in I(0)$, which implies by the definition of $\wt A$ that
	$$
	H(\wt{\cZ}^{2k},Z,\frk U) =H(
	\wt{\cZ}^{2k}(\cdot+\tau_{z_0})-\wt{\cZ}^{2k}(\tau_{z_0}),
	Z(\cdot+\tau_{z_0})-Z(\tau_{z_0}), 
	\frk U-\tau_{z_0}),
	$$
	and the last inequality follows from \eqref{eq:law-4tuple}, \eqref{eq:cZ-shift}, $\cZ^{2k}\in\sigma(Z|_{I(6k)},\frk U)$ and the fact that $B_0^c\subset\{ \eta(I(6k))\subset B(0,2^{6k}-2^{1-k}) \}$.  
\end{proof}

\begin{proof}[Proof of Corollary \ref{prop:eq0} assuming Theorem \ref{prop:zero-vs-gamma-walk}]
	If $\PQ[A]=1$ the result is immediate. If $\PQ[A]=0$, apply Theorem \ref{prop:zero-vs-gamma-walk} to each $A_k$ at scales $k+n$ and send $n\to\infty$. This gives $\PE[A_k]=0$ and hence $\PE[A]=0$. We may therefore assume $0<\PQ[A]<1$. 
	Let $E_{0}=\emptyset$ and $E_k=\{\1_{A_{k}}=\1_{A}\}$ and define $x:=-0.25\log_2\PQ[A]>0$. Then
	\eqbn
	\begin{split}
		\PE[ A ]
		&=\sum_{k=1}^\infty \min\{ \PE[E_k\setminus E_{k-1}]; \PE[A;E_k\setminus E_{k-1}] \}
		\leq\sum_{k=1}^\infty \min\{ \PE[E_k\setminus E_{k-1}]; \PE[A_k] \}.
		\end{split}
	\eqen
	By Theorem \ref{prop:zero-vs-gamma-walk} and since $A_k\in\sigma(\cZ|_{I(k)},I(k))\subset\sigma(\cZ|_{I(k+n)},I(k+n))$, we get that the right side of the last display is bounded above by the following, where $\eps>0$ is some constant independent of $\wt c$ and $\wt C$ which may change in each step
	\eqbn
	\begin{split}
		&\leq\sum_{k=1}^\infty \min\{ \PE[E_k\setminus E_{k-1}]; \inf_{n\in\N_0}\PQ[A_k]^{1/4}2^{C_0(k+n)}+2^{-(k+n)/C_0} \}\\
		&\leq\sum_{k=1}^\infty \min\{ \PE[E_k\setminus E_{k-1}]; \inf_{n\in\N_0}\PQ[A]^{1/4}2^{C_0(k+n)}+2^{-(k+n)/C_0} \}\\
		&\leq \sum_{k=1}^\infty \min\{ \wt C 2^{-\wt c(k-1)}; \inf_{n\in\N_0} 2^{-x+C_0(k+n)}+2^{-(k+n)/C_0} \}\\
		&\leq  \frac {\wt C}{\eps\wt c} 2^{-\wt c(C_0+C_0^{-1})^{-1}x}+ \sum_{k=1}^{\lfloor(C_0+C_0^{-1})^{-1}x\rfloor} \inf_{n\in\N_0} 2^{-x+C_0(k+n)}+2^{-(k+n)/C_0} \\
		&\leq \frac {\wt C}{\eps\wt c} 2^{-\wt c(C_0+C_0^{-1})^{-1}x}
		+3\cdot2^{-C_0^{-1}(C_0+C_0^{-1})^{-1}x}\cdot (C_0+C_0^{-1})^{-1}x
		\lesssim \frac {\wt C}{\eps\wt c} 2^{-\eps \wt c x},
	\end{split}
	\eqen
	where in the last line we picked $n\in\N_0$ such that $|(C_0+C_0^{-1})(k+n)-x|<(C_0+C_0^{-1})$ and we see that an appropriate $n\in\N_0$ exists since $k\leq \lfloor(C_0+C_0^{-1})^{-1}x\rfloor$. 
\end{proof}

\subsection{Proof of Theorem \ref{prop:zero-vs-gamma-walk}: Existence of the 0-quantum cone walk}
\label{sec:0-cone-walk-exists}

In this section we will complete the proof of Theorem \ref{prop:zero-vs-gamma-walk}. The starting point of the proof is Lemma \ref{prop:zero-vs-gamma0} proved above. The main difference between the two results is that the event in Lemma \ref{prop:zero-vs-gamma0} depends on the walk $\cZ^{2k}\in\sigma(Z|_{I(6k)},\frk U)$ while the walk in Theorem \ref{prop:zero-vs-gamma-walk} depends on $\cZ$. The latter walk depends on the full process $Z$ (not just $Z|_{I(6k)}$), which makes the theorem significantly harder since the absolute continuity result for the 0-quantum cone and the $\gamma$-quantum cone (Lemma \ref{lem:rooted}) is only for bounded regions. Furthermore, we do not yet know that the process $\cZ$ (defined in \eqref{eq:cZ}) is well-defined under $\PE$ and a large part of this section will be devoted to proving this.

The high-level idea for proving that $\cZ$ is well-defined is similar as in the case when $Z$ is a Brownian motion, which was addressed in Theorem \ref{prop-kmt}. 
Namely, we will argue that for each fixed $m\in\N$ and $\ell\in[\dZ]$ and for large $n\in\N$ it holds with high probability that $\cY^{\ell,n}|_{I(m)}=\cY^{\ell,n+1}|_{I(m)}$ (see Lemma \ref{prop-sqrtk} below). To get this it is sufficient to prove that 
(i) $\wt S_\ell(n)\approx S_\ell(n)$, where $\wt S_\ell(n)$ satisfies $\cY^{\ell,n+1}|_{I(m)}=\cY^{\ell,n,\wt S_\ell(n)}|_{I(m)}$, and 
(ii) for all $S\approx S_\ell(n)$ it holds that $\cY^{\ell,n,S}|_{I(m)}=\cY^{\ell,n}|_{I(m)}$. The result (i) will follow from a short argument combining Lemmas \ref{prop:zero-vs-gamma0} and \ref{prop:err3}, see \eqref{eq:E} below. The result (ii), on the other hand, will require significant work and is harder than in the Brownian case addressed in Lemma \ref{prop:same-E}. The proof of the cited lemma relies in an essential way on the translation invariance in law of the set of jump times on $\cY^{m+n,S}$, which no longer holds under the measure $\PE$. Instead, in Lemma \ref{prop:same} below we will use translation invariance modulo scaling as in Lemma \ref{prop:translation}, along with a number of estimates for LQG and SLE.

After proving that $\cZ$ is well defined we will establish \eqref{eq56} and \eqref{eq:transl0} in Lemmas \ref{prop:eucl-error} and \ref{prop:transl0}, respectively. For the former result we use Lemma \ref{prop:zero-vs-gamma0} along with an estimate from the Brownian case and from this subsection. In order to complete the proof of the theorem it remains to establish \eqref{eq:zero-vs-gamma-walk}, where the key difficulty is to treat events with $\cZ$ instead of $\cZ^{2k}$. This is, along with Lemma \ref{prop:same}, the hardest proof in this subsection and the key idea is to prove absolute continuity of $\cZ|_{I(k)}$ and $\cZ^{2k}|_{I(k)}$ along with an upper bound on the Radon-Nikodym derivative. 
\begin{lem}
	Let $(h,h^{\op{IG}},\frk U)$ be sampled from $\PE$, 
	let $m\in\N$, let $w_0$ be sampled from Lebesgue measure on $\eta(I(m))$, and let $\tau_{w_0}=\inf\{ t\in\R\,:\,\eta(t)=w_0\}$ be the time at which $\eta$ first hits $w_0$. Then for some random constant $C>0$,
	\eqbn
	\begin{split}
		&(h,h^{\op{IG}},\frk U,Z,\wt{\cZ}^{m})
		\eqD
		(h(C\cdot+w_0)+Q\log C,h^{\op{IG}}(C\cdot+w_0),\frk U-\tau_{w_0},\\
		&\qquad\qquad\qquad\qquad\qquad\qquad\qquad\qquad\qquad(Z_{t+\tau_{w_0}}-Z_{\tau_{w_0}})_{t\in\R},
		(\wt{\cZ}^{m}_{t+\tau_{w_0}}-\wt{\cZ}^{m}_{\tau_{w_0}})_{t\in (I(m)-\tau_{w_0})  } ).
	\end{split}
	\eqen
	\label{prop:translation}
\end{lem}
\begin{proof}
	For concreteness, we assume that $h$ has the circle average embedding. 
	Let $R>1$ and $\eps\in(0,1)$ be arbitrary. We will prove the following, where $d_{\op{tv}}$ denotes total variation distance and $C$ is independent of $\eps$ and $R$
	\eqb
	\begin{split}
		&d_{\op{tv}}\Big( (h|_{B(0,R)},h^{\op{IG}}|_{B(0,R)},\frk U),\\
		&\qquad\qquad((h(C\cdot+w_0)+Q\log C)|_{B(0,R)},h^{\op{IG}}(C\cdot+w_0)|_{B(0,R)},\frk U-\tau_{w_0}) \Big)\leq\eps.
	\end{split}
	\label{eq70}
	\eqe
	Let us first explain why it is sufficient to prove this. Since $R$ and $\eps$ were arbitrary, \eqref{eq70} implies the same result with $R=\infty$ and $\eps=0$. By \eqref{eq:coord-change} and the coordinate change formula  in imaginary geometry (see \cite[equation (1.2) and the paragraph below Theorem 1.2]{ig4}), we see that the boundary length process associated with $(h(C\cdot+w_0)+Q\log C,h^{\op{IG}}(C\cdot+w_0))$ is given by $(Z_{t+\tau_{w_0}}-Z_{\tau_{w_0}})_{t\in\R}$. Finally, by  \eqref{eq:cZ-shift} it holds that the associated walk is given by $(\wt{\cZ}^{m}_{t+\tau_{w_0}} - \wt{\cZ}^{m}_{\tau_{w_0}})_{t\in (I(m)-\tau_{w_0})}$, which concludes the proof that it is sufficient to establish \eqref{eq70}. 		
	We will assume without loss of generality that $h$ has the $\wt R$-embedding for some large $\wt R$ depending on $R$ and $\eps$ which will be chosen below. Let $v_0\in\C$ be sampled uniformly from Lebesgue measure in $B(0,\wt R/3)$. Then 
	\eqbn
	(h|_{B(0,\wt R/3)},h^{\op{IG}}|_{B(0,\wt R/3)},\frk U)
	\eqD
	(h(\cdot+v_0)|_{B(0,\wt R/3)},h^{\op{IG}}(\cdot+v_0)|_{B(0,\wt R/3)},\frk U-\tau_{v_0}).
	\eqen
	Let $\wt I\in\frk D$ denote the length $2^m$ interval of $\frk D$ containing $\tau_{v_0}$ and let $x_0$ be sampled from Lebesgue area measure on $\eta(\wt I)$. By choosing $\wt R$ sufficiently large we can guarantee that for some suitable random constant $C_1>0$,
	\eqbn
	\begin{split}
		&d_{\op{tv}}\Big( (h(C\cdot+w_0)+Q\log C)|_{B(0,R)},h^{\op{IG}}(C\cdot+w_0)|_{B(0,R)},\frk U-\tau_{w_0},\\
		&\qquad\qquad( (h(C_1\cdot+x_0)+Q\log C_1)|_{B(0,R)},h^{\op{IG}}(C_1\cdot+x_0)|_{B(0,R)},\frk U-\tau_{x_0}) \Big)<\eps/2,
	\end{split}
	\eqen 
	where we remark that the reason for the positive (but small) total variation distance is that $B(x_0,C_1 R)$ and/or $B(w_0,C R)$ might go outside of $B(0,\wt R)$. By choosing $\wt R$ sufficiently large, the total variation distance between $x_0$ and $v_0$ conditioned on the environment is $<\eps/10$ with probability at least $1-\eps/10$, which implies that for random $C_2>0$,

	\eqbn
	\begin{split}
		&d_{\op{tv}}\Big( (h(C_2\cdot+v_0)+Q\log C_2)|_{B(0,R)},h^{\op{IG}}(C_2\cdot+v_0)|_{B(0,R)},\frk U-\tau_{v_0},\\
		&\qquad\qquad( (h(C_1\cdot+x_0)+Q\log C_1)|_{B(0,R)},h^{\op{IG}}(C_1\cdot+x_0)|_{B(0,R)},\frk U-\tau_{x_0}) \Big)<\eps/2,
	\end{split}
	\eqen
	Combining the last three displays gives \eqref{eq70}.
\end{proof}

\begin{lem}
	Let $c_0\in(0,1)$ and $C'>1$, and choose $m,n\in\N$ such that $n\geq c_0 m$. There exists a constant $c>0$ depending on $c_0,C'$ such that for each $\ell\in[\dZ]$,
	\eqb
	\PE[ \exists S_1,S_2 \in \cI\,:\, \cY^{\ell,m+n,S_1}|_{I(m)} \neq \cY^{\ell,m+n,S_2}|_{I(m)}  ]<2^{-cn},
	\label{eq8}
	\eqe
	where
	$$
	\cI:=[ S_\ell(m+n)-C'(m+n), 
	S_\ell(m+n)+C'(m+n)]\cap\N_0.
	$$
	\label{prop:same}
\end{lem}
\begin{proof}
	First observe that the probability we want to bound is bounded by  
	\eqbn
	\begin{split}
		&\PE[\exists S\in\cI \text{\,s.t.\,} \cY^{\ell,m+n,S}|_{I(m)} \neq \cY^{\ell,m+n,S+1}|_{I(m)}]
	\end{split}
	\eqen
	Let $w_0$ be sampled from Lebesgue measure on $\eta(I(m+n))$, let $\tau_{w_0}=\inf\{t\in\R\,:\,\eta(t)=w_0 \}$ be the hitting time of $w_0$, and let $I(m,w_0)\in\frk D$ denote the length $2^{m}$ interval in $\frk D$ such that $\tau_{w_0}\in I(m,w_0)$.  
	By Lemma \ref{prop:translation},
	\eqbn
	\begin{split}
		&\PE[\exists S\in\cI \text{\,s.t.\,} \cY^{\ell,m+n,S}|_{I(m)} \neq \cY^{\ell,m+n,S+1}|_{I(m)}]\\
		&\qquad=
		\PE[\exists S\in\cI \text{\,s.t.\,} 
		(\cY^{\ell,m+n,S}_{t+\tau_{w_0}} - \cY^{\ell,m+n,S}_{\tau_{w_0}} )_{t\in I(m,w_0)-\tau_{w_0}} \neq (\cY^{\ell,m+n,S+1}_{t+\tau_{w_0}}-\cY^{\ell,m+n,S+1}_{\tau_{w_0}})_{t\in I(m,w_0)-\tau_{w_0}}].
	\end{split}
	\eqen
	Letting $A_S:=\{t\in I(m+n)\,:\,\Delta_t \cY^{\ell,m+n,S}\neq 0\}$ denote the set of times at which the function $\cY^{\ell,m+n,S}$ jumps, it follows from Lemma \ref{prop:add-one-jump} and the definition of $\wt F_p$ that $A_S\subseteq A_{S+1}$ and $A_{S+1}\setminus A_S$ has at most one point. If $A_{S+1}=A_S$ we have equality of the two functions considered in the last indented equation so let us assume that $A_{S+1}=A_S\cup\{t_S \}$ for some $t_S\in I(m+n)$. Up to null sets, 
	\eqbn
	\begin{split}
		&\{ \exists S\in\cI \text{\,s.t.\,} 
		(\cY^{\ell,m+n,S}_{t+\tau_{w_0}} - \cY^{\ell,m+n,S}_{\tau_{w_0}} )_{t\in I(m,w_0)-\tau_{w_0}} \neq (\cY^{\ell,m+n,S+1}_{t+\tau_{w_0}}-\cY^{\ell,m+n,S+1}_{\tau_{w_0}})_{t\in I(m,w_0)-\tau_{w_0}} \}\\
		&\qquad\qquad\qquad\qquad = \{\exists S\in\cI \text{\,s.t.\,}t_S\in I(m,w_0)\}
		= \{\exists S\in\cI \text{\,s.t.\,}w_0\in \eta(J(m,t_S)) \},
	\end{split}
	\eqen
	where $J(m,t_S)$ is the unique length $2^m$ interval in $\frk D$ containing $t_S$. 
	Letting $\cF$ denote the $\sigma$-algebra generated by all randomness except $w_0$, we get from the definition of $w_0$ and the last display that
	\eqbn
	\begin{split}
		&\PE[\exists S\in\cI \text{\,s.t.\,} (\cY^{\ell,m+n,S}_{t+\tau_{w_0}}
		- \cY^{\ell,m+n,S}_{\tau_{w_0}})_
		{t\in I(m,w_0)-\tau_{w_0}} 
		\neq (\cY^{\ell,m+n,S+1}_{t+\tau_{w_0}}
		-\cY^{\ell,m+n,S+1}_{\tau_{w_0}})_{t\in I(m,w_0)-\tau_{w_0}}\,|\,\cF]\\
		&\qquad\qquad\qquad\qquad\qquad\qquad\qquad\qquad\qquad\qquad\qquad\qquad\qquad\qquad
		\leq \sum_{S\in\cI} \frac{ \op{Leb}(\eta(J(m,t_S))) }{ \op{Leb}(\eta(I(m+n))) }.
	\end{split}
	\eqen
	Suppose now that $h$ has the $2^{(1+2\eps)\frac{m+n}{\gamma Q}}$-embedding. Define the event $B_3$ by
	\eqbn
	\begin{split}
		B_3= & \{ B(0,2^{(1-\eps)\frac{m+n}{\gamma Q}}) \not\subset \eta(I(m+n)) \}\cup 
		\{ \eta(I(m+n))\not\subset B(0,2^{(1+\eps)\frac{m+n}{\gamma Q}})\}.
	\end{split}
	\eqen
	By Lemmas \ref{prop:contained1} and \ref{prop:contained2}, $\PE[B_3]\leq 2^{-c(m+n)}$. Taking expectation in the last displayed estimate and using that $\#\cI<3C'(m+n)$,

	\eqbn
	\begin{split}
		&\PE[\exists S\in\cI \text{\,s.t.\,} (\cY^{\ell,m+n,S}_{t+\tau_{w_0}}
		- \cY^{\ell,m+n,S}_{\tau_{w_0}})_
		{t\in I(m,w_0)-\tau_{w_0}} 
		\neq (\cY^{\ell,m+n,S+1}_{t+\tau_{w_0}}
		-\cY^{\ell,m+n,S+1}_{\tau_{w_0}})_{t\in I(m,w_0)-\tau_{w_0}}]\\
		&\quad\leq 
		3C'(m+n)
		\EE\Big[\frac{\max \{ \op{Leb}(\eta(J))\,:\,J\in\frk D,J\subset I(m+n),|J|=2^{m} \} }{\op{Leb}(\eta(I(m+n)))}\Big] \\
		& \quad\leq 
		3C'(m+n) \E_{\op{E}}\Big[\frac{\max \{ \op{Leb}(\eta(J))\,:\,J\in\frk D,J\subset I(m+n),|J|=2^{m} \} }{\op{Leb}(\eta(I(m+n)))}\1_{B_3^c}\Big]+3C'(m+n)\PE[B_3]\\
		&\quad\leq
		3C'(m+n)
		\EE\Big[\frac{\max \{ \op{Leb}(\eta(J))\,:\,\eta(J)\subset B(0,2^{(1+\eps)\frac{m+n}{\gamma Q}}),J=[a,a+2^m]\text{ for }a\in\R \} }{2^{2(1-\eps)\frac{m+n}{\gamma Q} }}\Big]
		\\
		&\qquad\qquad\qquad\qquad\qquad\qquad\qquad\qquad\qquad\qquad\qquad\qquad\qquad\qquad\qquad\qquad\qquad\qquad+2^{-c(m+n)} \\
		&\quad\leq 
		3C'(m+n)
		2^{8\eps(m+n)/(\gamma Q)}(2^{-2q(n+2\eps m+2\eps n)}+ 2^{2-\alpha(q,\gamma)(n+2\eps m+2\eps n+1)+g(n+2\eps m+2\eps n)})+2^{-c(m+n)},
	\end{split}
	\eqen
	where we use the definition of $B_3$ to get the second inequality and we use the last assertion of Lemma \ref{prop:max-cell-diam} to get the last inequality. The right side of the last display is bounded above by $2^{-cn}$ since we can first pick an arbitrary $q\in(0,2(2+\gamma)^{-2})$ (observing that any such $q$ gives $\alpha(q,\gamma)>0$) and then pick $\eps$ sufficiently small.
\end{proof}

\begin{lemma}
	For $j_1\leq j_2$ let $B_{j_1,j_2}$ be the event that for some
	$\ell\in[\dZ]$, 
	\eqb
	\cY^{\ell,j_2}|_{I(j_1)} \neq
	\cY^{\ell,j_2+1}|_{I(j_1)}.
	\label{eq:B}
	\eqe
	Then for any $c_0>0$ there is a  $C_0>1$ such that for all $m,n\in\N_0$ for which $n\geq c_0m$,
	\eqb
	\PE\left[ \bigcup_{j\geq m+n} B_{m,j} \right] \leq 2^{-n/C_0}.
	\label{eq:same}
	\eqe
	In particular, $\cZ$ is well-defined and for all $m,n\in\N$ such that $n\geq c_0m$,
	\eqb
	\PE[\exists\ell\in[\dZ]\text{\,\,such\,\,that\,\,}\cY^\ell|_{I(m)}\neq
	\cY^{\ell,m+n}|_{I(m)}
	] < 2^{-n/C_0}.
	\label{eq:same2}
	\eqe
	\label{prop-sqrtk}
\end{lemma}
\begin{proof}
	By a union bound,
	\eqb
	\PE\left[B_{m,j} \right] \leq 
	\PE[ E_{j}  ]+\PE[\wt B_{m,j}],
	\label{eq:Bsplit}
	\eqe
	where for some constant $C'>1$ to be chosen later and $S(\cdot,\cdot)$ as in \eqref{eq:defS},
	$$
	E_{j} = \bigcup_{\ell\in[\dZ]}\{ |S(\cY^{\ell,j+1},I(j)) 
	-S_\ell(j)|>C'j \}
	$$
	and	
	$$
	\wt B_{m,j}=
	\bigcup_{\ell\in[\dZ]}\{ \exists S_1,S_2 \in \cI\,:\, \cY^{\ell,j,S_1}|_{I(m)} \neq \cY^{\ell,j,S_2}|_{I(m)}  \},\qquad
	\cI:=[0,\infty)\cap [
	S_\ell(j)-C'j, 
	S_\ell(j)+C'j ]_{\Z}.
	$$
	By Lemma \ref{prop:same}, $\PE[\wt B_{m,j}] \leq 2^{-c(j-m)}$. 
	We will now apply Lemma \ref{prop:zero-vs-gamma0} with $\wt G=\1_{E_{j}}$ and $k=j$ (in fact, any $k$ such that $ck\leq j\leq 2k-1$ for some $c>0$ would work). Observe that the event $\wt A$ in the cited lemma is then given by $\wt A=E_j$. Furthermore, $\wt A\in \sigma(Z|_{I(3j+3)},\frk U)\subset \sigma(\wt{\cZ}^{2j},Z|_{I(6j)},\frk U)$ by Lemma \ref{eq:cYell-measurable}.  
	We get
	\eqb
	\PE[E_{j}] 
	\leq 
	\PQ[E_{j}]^{1/2}2^{C_0j}+2^{-j/C_0}.
	\label{eq:E}
	\eqe
	By choosing $C'$ sufficiently large, Lemma \ref{prop:err3} gives that
	$$
	\PQ[E_{j}] \leq 2^{-3C_0j}.
	$$
	Combining this with the other above estimates and inserting into \eqref{eq:Bsplit} we get $\PE[B_{m,j}]\leq 2^{-cj}+2^{-c(j-m)}$, which implies \eqref{eq:same}. 
	
	On the complement of the event in \eqref{eq:same} the limit on the right side of \eqref{eq:def-walk} exists for $I=I(m)$ and
	$$
	\cY^\ell|_{I(m)}= \cY^{\ell,m+n} |_{I(m)},
	$$
	which gives the last assertion of the lemma.
\end{proof}

\begin{lem}
	For any $C_3>1$ we can find $C'>1$ such that
	$\PE[\sup_{t\in I(m)}|Y^\ell_t-\cY^\ell_t|>C' m]\leq 2^{-C_3 m}$ for $\ell=1,\dots,\dZ$
	\label{prop:eucl-error0}
\end{lem}
\begin{proof}
	For any $k\in\N$, 
	\eqbn
	\begin{split}
		&\PE\Big[\sup_{t\in I(m)}|Y^\ell_t-\cY^\ell_t|>C'm\Big]
		\leq \PE\Big[\sup_{t\in I(m)}|Y^\ell_t-\cY^{\ell,m+k}_t|>C'm\Big]
		+ \PE[\cY^{\ell,m+k}|_{I(m)}\neq\cY^\ell|_{I(m)}].
	\end{split}
	\eqen
	The second term on the right side is smaller than $2^{-k/C_0}$ by Lemma \ref{prop-sqrtk} as long as $k\geq m$. We bound the first term on the right side by using Lemma \ref{prop:zero-vs-gamma0} with 
	$$
	\wt G(\cZ^{\ell,m+k},Z,\frk U)= \1_{\wt A},\qquad 
	\wt A=\{\sup_{t,s\in I(m)}|\Delta_{s,t}|>C'm\},\qquad
	\Delta_{s,t} := (Y^\ell_s-\cY^{\ell,m+k}_s)-(Y^\ell_t-\cY^{\ell,m+k}_t),
	$$
	$k$ replaced by $m+k$, and $\wt A=\{\wt G(\cZ,Z,\frk U)=1 \}$, which gives
	$$
	\PE\Big[\sup_{s,t\in I(m)}|\Delta_{s,t}|>C'm\Big]
	\leq 
	\PQ\Big[\sup_{s,t\in I(m)}|\Delta_{s,t}|>C'm\Big]^{1/2}2^{C_0(m+k)} + 2^{-(m+k)/C_0}.
	$$
	By \eqref{eq:err3}, given arbitrary $C_1>1$ we can find $C'>1$ such that
	$$
	\PQ\Big[\sup_{s,t\in I(m)}|\Delta_{s,t}|>C'm\Big] \leq 2^{-C_1m}. 
	$$
	Setting $k=2C_0C_3m$ and choosing $C_1>1$ sufficiently large such that $-C_1/2+C_0+2C^2_0C_3<-2C_3$ we get
	$$
	\PE\Big[\sup_{t\in I(m)}|Y^\ell_t-\cY^\ell_t|>C'm\Big]
	\leq 2^{-C_1m/2+C_0(m+k)}+2^{-(m+k)/C_0}+2^{-k/C_0} < 3\cdot 2^{-2C_3m}.
	$$
\end{proof}

\begin{lem}
	For any $C_3>1$ we can find $C'>1$ such that
	 $\PE[\sup_{t\in [-2^m,2^m]}|Y^\ell_t-\cY^\ell_t|>C' m]\leq 2^{-C_3 m}$ for $\ell=1\dots,\dZ$.
	\label{prop:eucl-error}
\end{lem}
\begin{proof}
	For any $a>1$,
	$$\PE[\sup_{t\in [-2^m,2^m]}|Y^\ell_t-\cY^\ell_t|>C' m]
	\leq 
	\PE[ \sup_{t\in I(am) } |Y^\ell_t-\cY^\ell_t|>C' m ]
	+ \PE[ [-2^m,2^m]\not\subset I(am) ].
	$$
	By choosing $a$ sufficiently large, the second probability on the right side is bounded above by $2^{-2C_3m}$. It then follows from Lemma \ref{prop:eucl-error0} that by choosing $C'$ sufficiently large (depending on $a$), the first probability on the right side is also bounded above by $2^{-2C_3m}$.
\end{proof}

\begin{lem}
	There is a measurable function $g$ such that $\cZ=g(h,h^{\op{IG}},\frk U)$ and this function is translation invariant modulo scaling in the sense that \eqref{eq:transl0} a.s.\ holds.  
	\label{prop:transl0}
\end{lem}
\begin{proof}
	The first assertion is an immediate consequence of the definition of $\cZ$ since $\cZ$ is determined by $\cY^\ell$ by \eqref{eq:cZ}, $\cY^\ell$ is determined by $Y,\frk U$ by Lemma \ref{eq:cYell-measurable} and \eqref{eq:def-walk}, and $Y$ is determined by $h,h^{\op{IG}}$. We get the second assertion by using \eqref{eq:cZ-shift} along with the fact that the boundary length process associated with the fields in the right side of \eqref{eq:transl0} is given by $(Z_{t+\tau_z}-Z_{\tau_z})_{t\in\R}$, a.s.\ simultaneously for all $C>0$ and $z\in\C$.
\end{proof}

In order to complete the proof of Theorem \ref{prop:zero-vs-gamma-walk} it remains to establish \eqref{eq:zero-vs-gamma-walk}. The first strategy that comes to mind is to combine \eqref{eq:zero-vs-gamma0} with the fact that $\cZ^{2k}|_{I(k)}=\cZ|_{I(k)}$ with high probability. We will make use of these results in the proof; however, we also need significant additional ideas since the failure probability in \eqref{eq:same2} is too large to rely only on this strategy. The key additional idea is to prove that, although $\cZ^{2k}|_{I(k)}$ and $\cZ|_{I(k)}$ might not be identical, their laws are absolutely continuous with a good Radon-Nikodym derivative. We prove the latter result by using the explicit laws and certain independence properties of the processes $\cZ^{2k}|_{I(k)}$ and $\cZ|_{I(k)}$, along with the next lemma.
\begin{lem}
	Let $V\sim\op{Pois}(n)$, let $\wt V\sim\cN(n,n+\eps_n)$ for $|\eps_n|<10/n$, and let $\wt\xi$ denote the Radon-Nikodym derivative of the law of $\lceil\wt V\rceil\vee 0$ with respect to the law of $V$. Then for any fixed $C_0>1$, 
	$$
	\lim_{n\rta\infty}\sup_{|x|<C_0\log(n)} \wt\xi(\lceil n+\sqrt{n}x\rceil ) = 1
	= \lim_{n\rta\infty}\inf_{|x|<C_0\log(n)} \wt\xi(\lceil n+\sqrt{n}x\rceil ).
	$$
	In particular, there is a constant $C_1>0$ such that $\E[\wt\xi(V)^2\1_{|V-n|<C_0\sqrt{n}\log(n)}]<C_1$ for all $n\in\N$.
	\label{prop:RN-poisson}
\end{lem}
\begin{proof}
	For any $x\in[-C_0\log(n),C_0\log(n)]$ such that $n+\sqrt{n}x$ is an integer,
	\eqbn
	\begin{split}
			\P[V= n+\sqrt{n}x] &= \frac{n^{ n+\sqrt{n}x} \exp(-n)}{ (n+\sqrt{n}x)!}
			\sim 
			\frac{n^{ n+\sqrt{n}x} \exp(\sqrt{n}x)}{ \sqrt{2\pi  (n+\sqrt{n}x)}  (n+\sqrt{n}x)^{ n+\sqrt{n}x }}
			\sim 
			\frac{\exp(-x^2/2 )}{ \sqrt{2\pi n} },
		\end{split}
		\eqen
		where we apply Stirling's formula in the second step and $\sim$ means that the ratio of the left and right sides go to 1 as $n\rta\infty$ and $|x|<C_0\log(n)$. Combining this with the following estimate gives the lemma
		\eqbn
		\P\big[\lceil\wt V\rceil= n+\sqrt{n}x\big] \sim
		\frac{\exp(-x^2/2 )}{\sqrt{2\pi n}}.
		\eqen
	\end{proof}
	
\begin{proof}[Proof of \eqref{eq:zero-vs-gamma-walk} in Theorem \ref{prop:zero-vs-gamma-walk}] 
		We apply Lemma \ref{prop:zero-vs-gamma0} with 
		$
		\wt G=G\circ\op{proj}_{1,3},
		$ 
		where proj$_{1,3}$ denotes projection onto the first and third coordinates. 
		Set 
		$B_1 = \{\exists\ell\in[\dZ]\text{\,s.t.\,} \cY^{\ell,2k}|_{I(k)}\neq\cY^\ell|_{I(k)} \}$. 
		Notice that $A\cap B^c_1=\wt A\cap B^c_1$. We get

		$$
		\PE[A] \leq \PE[\wt A] + \PE[B_1]
		\leq \PQ[\wt A]^{1/2}2^{C_0k} + \PE[B_0]+ \PE[B_1].
		$$	
		We bound $\PE[B_0]+\PE[B_1]$ by $2^{-k/C_0}$ by applying Lemmas \ref{prop:zero-vs-gamma0} and \ref{prop-sqrtk}. In order to bound $\PQ[\wt A]$, we will first prove the following:
		\begin{compactitem}
			\item[(a)] There is a measurable function $H$ with values in $\{ 0,1\}$ such that 
			\eqb
			H(\frk D,\wt{\cS},\frk X)=\1_{\wt A}
			\quad\text{and}\quad
			H(\frk D,\cS,\frk X)=\1_{A},
			\label{eq:H}
			\eqe
			where
			\eqb
			\begin{split}
				&S_\ell=S( \cY^{\ell},I(2k) ),\quad
				\wt S_\ell=S( \cY^{\ell,2k},I(2k) ),\quad
				\cS=(S_\ell)_{\ell\in[\dZ]},\quad
				\wt{\cS}=(\wt S_\ell)_{\ell\in[\dZ]},\\
				&\frk X=(X_{I}^\ell)_{\ell\in[\dZ],I\in\frk D,I\subseteq I(2k)} 
			\end{split}
			\label{eq:HA}
			\eqe
			and we recall the definition of $S(\cdot,\cdot)$ from \eqref{eq:defS} and the definition of $X_I^\ell$ from \eqref{eq:X-def2}.
			\item[(b)] Conditioned on $\frk D$, the random variables $\cS$ (resp.\ $\wt{\cS}$) and $\frk X$ are independent.
			\item[(c)] $\wt{\cS}$ is absolutely continuous with respect to $\cS$ conditioned on $\frk D$ and the Radon-Nikodym derivative $\xi$ satisfies $\E[\xi(\cS)^2\1_{\cS\in K}\,|\,\frk D]\leq C_1$ for some constant $C_1>1$, where $K=\{ (x_1,\dots,x_{\dZ})\in\R^{\dZ}\,:\, |x_\ell-2^{2k}p_\ell|\leq k2^k\sqrt{p_\ell}\}$. %
		\end{compactitem}
		 
		Proof of (a): By assumption we have $A\in\sigma(\cZ|_{I(k)},I(k))$, which implies by the definition of $\wt G$ that $\wt A\in\sigma(\cZ^{2k}|_{I(k)},I(k))\subset\sigma(\cY^{1,2k},\dots,\cY^{\dZ,2k},I(k) )$. Furthermore, for each $\ell\in[\dZ]$ we have $\cY^{\ell,2k}\in\sigma(\wt S_\ell,\br(Y^\ell,I(2k)))$ by the definition in \eqref{eq:cYm2} along with \eqref{eq:br-meas}. By using that $\br(Y^\ell,I(2k))\in \sigma(I(2k),(X^\ell_{I})_{I\in\frk D,I\subseteq I(2k)} )\subset\sigma(I(2k),\frk X)$ due to \eqref{eq:Z-brZ-X} we conclude that a function $H$ satisfying the first identity in \eqref{eq:H} exists. To get the second identity in \eqref{eq:H} we use that $\cY^{\ell}|_{I(2k)}=F(Y^\ell|_{I(2k)},S_\ell)$ by  \eqref{eq:cYm2} and \eqref{eq:def-walk}, 
		and we note that the right side of this formula is identical to the definition of $\cY^{\ell,2k}$ except that $\wt S_\ell$ has been replaced by $S_\ell$. \eqref{eq:cYell-measurable} and Lemma \ref{prop:coupling-properties}

		Proof of (b): By \eqref{eq:Z-brZ-X} we see that $\wt{\cS}$  
		is measurable with respect to $\cG:=\sigma\big(\frk D, (X_{I}^\ell)_{\ell\in[\dZ], I\in\frk D, I(2k)\subsetneq I}\big)$. It follows from Lemma \ref{prop:br-indep} that, 
		conditioned on $\frk D$, the random variables $\wt{\cS}$  and $\frk X$ are independent. We argue in the same way for $\cS$, except that we use 
		Assertion (ii) of Lemma \ref{prop:coupling-properties} to argue that $\cS$ is measurable with respect to $\cG$.

		Proof of (c): By the explicit law of the processes $P^1,\dots,P^{\dZ}$ from the beginning of Section \ref{sec:kmt3}, the random variables $S_1,\dots,S_{\dZ}$ are independent and $S_\ell\sim\op{Pois}(2^{2k}p_\ell)$ for $\ell\in[\dZ]$. 
		By \eqref{eq:cYm2} and since $\wt S_\ell=\lceil S_\ell(2k)\rceil\vee 0$ by \eqref{eq43} and \eqref{eq68}, the random variables $\wt{S}_1,\dots\wt S_{\dZ}$ are also independent and $S_{\ell}(2k)\sim\cN(2^{2k}p_\ell,2^{2k}p_\ell-2^{-2k}p_\ell)$. It follows that $\wt{\cS}$ is absolutely continuous with respect to $\cS$. This also holds if we condition on $\frk D$ since this does not change the laws of $\cS$ and $\wt{\cS}$. Let $\xi$ denote the Radon-Nikodym derivative of the law of $\wt{\cS}$ with respect to the law of $\cS$ given $\frk D$. Since the random variables $S_\ell$ (resp.\ $\wt S_\ell$) are independent given $\frk D$, $\xi(\cS)=\prod_{\ell\in[\dZ]}\xi_\ell(S_\ell)$, where $\xi_\ell$ is the Radon-Nikodym derivative of the law of $\wt S_\ell$ with respect to the law of $S_\ell$. The bound in (c) now follows from the last assertion of
		Lemma \ref{prop:RN-poisson}.
		
		We have for sufficiently large $k$,
		$$
		\PQ[\wt A]
		\leq 
		\PQ[\wt A\cap \{\wt\cS\in K\}]+\P[\wt\cS\not\in K] 
		\leq \PQ[\wt A\cap \{\wt\cS\in K\}]+\exp(-0.1k^2),
		$$
		where the second inequality follows from the tail behavior of Gaussian random variables. We have further 
		\eqb
		\begin{split}
			&\PQ[\wt A\cap \{\wt\cS\in K\}]
			=\EQ[\EQ[H(\frk D,\wt{\cS},\frk X)\1_{\wt\cS\in K}\,|\,\frk D]]\\
			&\qquad=\EQ[\EQ[H(\frk D,\cS,\frk X)\xi(\cS)\1_{\cS\in K}\,|\,\frk D]]\\
			&\qquad\leq \EQ[\EQ[H(\frk D,\cS,\frk X)\,|\,\frk D]^{1/2}\EQ[\xi(\cS)^2\1_{\cS\in K}\,|\,\frk D]^{1/2}]\\
			&\qquad\leq C_1^{1/2}
			\EQ[\EQ[H(\frk D,\cS,\frk X)\,|\,\frk D]^{1/2}]
			\leq C_1^{1/2}\EQ[H(\frk D,\cS,\frk X)]^{1/2}=C_1^{1/2}\PQ[A]^{1/2}.
		\end{split}
		\eqe
		where we use (b) to get the second equality and we use (c) to get the second inequality. 
	\end{proof}

\subsection{The boundary length process of the 0-quantum cone with SLE}
\label{sec:0-cone-BL}
This section collects various results on the boundary length process associated with the SLE-decorated 0-quantum cone. Most lemmas are fairly easy to derive from other results, with the exception of Lemma \ref{prop:eucl-inf-lower}, which requires a more involved argument.
\begin{lem}
	$\cZ$ under the measure $\PE$ has a.s.\ finitely many jumps in any compact interval.
	\label{prop-finite-jumps}
\end{lem}
\begin{proof}
	This is immediate from the observation that if the limit in \eqref{eq:def-walk} exists then $\cY^\ell$ makes finitely many jumps in $I(k)$ for any $k\in\N$.
\end{proof}

The following scaling property of $Z$ is immediate from the invariance of the 0-quantum cone under adding a constant.
\begin{lem}
If $Z$ is the boundary length process under the measure $\PE$ and $a>0$ is a deterministic constant then
$Z \eqD (a^{-1/2}Z_{at})_{t\in\R}$.	\label{prop:E-scaling}
\end{lem}
\begin{proof}
	If we add the constant $\gamma^{-1}\log a^{-1}$ to the field $h$ then the LQG area is multiplied by $a^{-1}$ and the LQG boundary lengths are multiplied by $a^{-1/2}$, so the new boundary length process is $(a^{-1/2}Z_{at})_{t\in\R}$. This boundary length process is equal in law to $Z$ since the law of a $0$-quantum cone is invariant under adding a deterministic constant to the field, see \cite[Proposition 4.13]{wedges}. 
\end{proof}

The following upper bound for $Z$ follows from the analogous bound in the $\gamma$-quantum cone case combined with Lemma \ref{prop:zero-vs-gamma0}.
\begin{lem}
	There is a constant $c>0$ such that for any sufficiently large $x>1$,
	$$
	\PE\bigg[\sup_{t\in I(0) } |L_t| > x \bigg] \leq 2^{-cx^2}.
	$$
	\label{prop:eucl-inf-upper}
\end{lem}
\begin{proof}
	We apply Lemma \ref{prop:zero-vs-gamma0} with $k$ equal to a sufficiently small constant multiple of $x^2$ and get that for some $\wt c,c>0$,
	$$
	\PE\Big[\sup_{s,t\in I(0) } |L_t-L_s| > x \Big] 
	\leq 
	\PQ\Big[\sup_{s,t\in I(0) } |L_t-L_s| > x \Big]^{1/2} 2^{C_0k} + 2^{-k/C_0}
	\leq 2^{-\wt c x^2 + C_0k} + 2^{-k/C_0} \leq 2^{-c x^2}.
	$$
	This concludes the proof since $L_0=0$, so the probability in the lemma is bounded above by the probability in the above display.
\end{proof}

\begin{figure}
	\centering
	\includegraphics{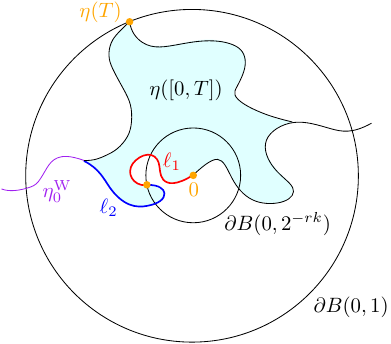}
	\caption{Illustration of the proof of Lemma \ref{prop:eucl-inf-lower}.}
	\label{fig-euclideanprocess}
\end{figure}

The following lower bound for the infimum of $L$ cannot be established by the same method as the previous lemma. 
The issue is that the probability of the considered event decays too slowly in the Brownian case to offset the factor $2^{C_0k}$ from the Radon--Nikodym estimate.
Our goal is to bound from above the probability of the event that a particular imaginary geometry flow line segment started from 0 on a 0-quantum cone has a very small length. A key idea of the proof is to use that, inside the unit disk away from 0, the fields of the 0-quantum cone and the $\gamma$-quantum cone can be coupled so that their difference is bounded, and for the $\gamma$-quantum cone we can use estimates for Brownian motion. A number of estimates for LQG and SLE are needed as part of the proof, including those in Lemmas \ref{prop-area-until-hit} and \ref{prop:cone-moment}.
\begin{lem}
	For $x\in\R$ define
	$
	\tau_x = \inf\{ t\geq 0\,:\,L_t=x \}.
	$
	Then there is a constant $c>0$ such that for any $k>0$,
	$\PE[\tau_{-1}>2^k]\leq 2^{-ck}$.
	\label{prop:eucl-inf-lower}
\end{lem}
\begin{proof} 
	We assume throughout the proof that the field $h$ has the circle average embedding and observe that the lemma does not depend on the choice of embedding.  
	Furthermore, throughout the proof let $c>0$ be some small constant which may vary from line to line. Define $a=0.4$, $r=3a/(2+\gamma^2)$, and let $\delta\in(0,1)$ be sufficiently small such that
	\eqb
	\begin{split}
		2a>\delta+\gamma^2r,\qquad [(1-\delta)(2-\gamma^2/2)+\gamma^2]r>\delta+2a.	
	\end{split}
	\label{eq:a-r-delta}
	\eqe
	Define
	$$
	T=\inf\{ t\geq 0\,:\,\eta([0,t])\not\subset B(0,1) \}.
	$$
	Recall from Section \ref{sec:sle} that $\eta_0^{\op{L}}$  denotes the angle $\pi/2$ flow line started from 0 of the imaginary geometry field $h^{\op{IG}}$ coupled together with $\eta$. The curve $\eta_0^{\op{L}}$ is also equal to the left boundary arc of $\eta((-\infty,0])$. Parameterize $\eta_0^{\op{L}}$ in an arbitrary way, let $\tau_r$ be the last exit time from $B(0,2^{-rk})$, i.e.,
	$$
	\tau_r = \sup\{t\geq 0\,:\,\eta_0^{\op{L}}(t)\in B(0,2^{-rk}) \},
	$$
	and set
	$$
	\ell_1=\eta_0^{\op{L}}([0,\tau_r]) \cap \eta([0,T]),
	\qquad 
	\ell_2 = \Big(\eta_0^{\op{L}}([0,\infty)) \cap \eta([0,T])  \Big) \setminus \ell_1
	.
	$$
	See Figure \ref{fig-euclideanprocess} for an illustration. We may have $\ell_2=\emptyset$ but we will see below that this is a low probability event. Observe that $\ell_1\cup\ell_2=\eta_0^{\op{L}}([0,\infty)) \cap \eta([0,T])$ is a boundary segment of $\eta([0,T])$ of LQG length $\nu(\ell_1\cup\ell_2)=-\inf_{t\in [0,T]}L_t$, where $\nu$ denotes the $\gamma$-LQG length measure on the SLE$_{\gamma^2}$-type curve $\ell_1\cup\ell_2$. By first applying Lemma \ref{prop:E-scaling} and then doing a union bound we have
	\eqbn
	\PE[\tau_{-1}>2^{k+2a k}]
	\leq \PE[ \tau_{-2^{-a k}}> 2^{k} ]
	\leq \PE[\mu(B(0,1)) > 2^{k} )]
	+ \PE[\nu(\ell_1\cup\ell_2 )< 2^{-a k} ].
	\eqen
	The first term on the right side is bounded above by a constant multiple of $2^{-k}$ by Markov's inequality. In the remainder of the proof we will bound the second term on the right side. 
	
	Let the measure $\wt \nu$ be defined exactly as  $\nu$ except that we replace the field $h$ by $h+\gamma\log|\cdot|^{-1}$. Note that the law of $\wt \nu$ under $\PE$ is equal to the law of  $\nu$ under $\PQ$. Also observe that since $\ell_2$ lies outside the ball $B(0,2^{-rk})$ we have $\wt\nu(\ell_2)\leq 2^{rk\gamma^2/2}\nu(\ell_2)$. Using this we get
	\eqb
	\begin{split}
		\PE[\nu(\ell_1\cup\ell_2 )< 2^{-a k} ] 
		&\leq \PE[\nu(\ell_2)< 2^{-a k} ]
		\leq \PQ[\nu(\ell_2)<2^{-a k+rk\gamma^2/2}]\\
		&\leq \PQ[\nu(\ell_1\cup\ell_2)<2^{-a k+1+rk\gamma^2/2}]
		+ \PQ[\nu(\ell_1)>2^{-a k+rk\gamma^2/2}]
	\end{split}
\label{eq21}
	\eqe
	To bound the first term on the right side of \eqref{eq21}, 
	\eqbn
	\begin{split}
	 &\PQ[\nu(\ell_1\cup\ell_2)<2^{-a k+1+rk\gamma^2/2}]
	 =\PQ[-\inf_{t\in[0,T]} L_t<2^{-a k+1+rk\gamma^2/2}]\\
	 &\qquad\leq 
	 \PQ[ \mu(\eta[0,T])<2^{-2a k+2+rk\gamma^2+\delta k}]
	 +\PQ\Big[ -\inf_{t\in[0,2^{-2a k+2+rk\gamma^2+\delta k}]} L_t < 2^{-a k+1+rk\gamma^2/2}\Big]\leq 2^{-ck},
	\end{split}
	\eqen	
	where we get the last estimate by using Lemma \ref{prop-area-until-hit}, the first inequality of \eqref{eq:a-r-delta}, and the fact that $L$ has the law of a Brownian motion.
	To bound the second term on the right side of \eqref{eq21}, first define
	$$
	T_+=\inf\{t\geq 0\,:\,\eta([0,t])\not\subset B(0,2^{-(1-\delta)rk}) \},\qquad
	T_-=\sup\{t\leq 0\,:\,\eta([t,0])\not\subset B(0,2^{-(1-\delta)rk}) \}.
	$$
	By a union bound,
	\eqbn
	\begin{split}
		&\PQ[\nu(\ell_1)>2^{-a k + kr\gamma^2/2}] \leq\, 
		\PQ[B(0,2^{-rk})\not\subset\eta([T_-,T_+])]\\
		&\qquad\qquad\qquad+
		\PQ\big[\mu(B(0,2^{-(1-\delta)rk}))>2^{-(1-\delta)(Q\gamma-\gamma^2) rk+\delta k} \big] \\
		&\qquad\qquad\qquad+ 
		\PQ\Big[-\min_{t\in[0,2^{-(1-\delta)(Q\gamma-\gamma^2) rk+\delta k}]} L_{t} > 2^{-a k + kr\gamma^2/2}\Big]\leq 2^{-ck},
	\end{split}
	\eqen
	where we get the last bound by using \cite[Proposition 6.2]{hs-euclidean} to bound the first term, 
	we use Lemma \ref{prop:cone-moment} to bound the second term, 
	and we use the second inequality of \eqref{eq:a-r-delta} and that $L$ has the law of a Brownian motion to bound the third term.
\end{proof}

\begin{remark}
	The conclusion of Lemmas \ref{prop:E-scaling}, \ref{prop:eucl-inf-upper}, and \ref{prop:eucl-inf-lower} also hold if we replace $\PE$ by $\PQ$. This is immediate by properties of Brownian motion and since $L$ has the law of a two-sided Brownian motion under $\PQ$.
	\label{rmk:Q-inf}
\end{remark}

The last two lemmas of this section were used in the proof of Lemma \ref{prop:eucl-inf-lower} above.

\begin{lem}
	For some $c\in(0,1)$ we have
	$\PQ[ \mu(\eta[0,T])< x ] \leq x^c$, where we suppose the field $h$ used to define $\mu$ has the circle average embedding and $T=\inf\{t\geq 0\,:\,\eta([0,t])\not\subset B(0,1) \}$ is as defined in the proof of Lemma \ref{prop:eucl-inf-lower}.
	\label{prop-area-until-hit}
\end{lem}
\begin{proof}
	Let $s\in(0,1/(\gamma Q))$. 
	By \cite[Lemma 3.6]{ghm-kpz} there is a constant $c\in(0,1)$ such that the set $\eta([0,T])$ contains a ball $B$ of radius $x^s$ except on an event of probability $\exp(-cx^{-s})$. By \eqref{eq:ball-upper} with $\eps=x^s$, $\delta=x^{1-\gamma Q s}$ and $\sigma^2=\log x^{-s}$ we get the lemma. 
\end{proof}

\section{The cell configurations}
\label{sec:cell-conf}
In this section we are given a random infinite tree-decorated planar map $(M,e_0,\tree)$ and a process $(\cL_t,\cR_t)_{t\in\R}$ encoding it via the correspondence presented in Section \ref{sec:bijection}. We will define two cell configurations (Definition \ref{def-cell-config}) $(\cH^{\op{cell}},\sim_{\op{cell}},\frk c_{\op{cell}})$ and $(\cH^{\op{large}},\sim_{\op{large}},\frk c_{\op{large}})$ satisfying two properties:
(i) The cell configurations (viewed as multi-graphs) are isomorphic to the planar map $M$, 
(ii) the cell configurations satisfy the hypotheses of Theorem \ref{thm-general-clt-uniform}. Theorem \ref{thm-general-clt-uniform} then gives an invariance principle for simple random walk on $M$ embedded according to the cell configuration (see Proposition \ref{prop:bm-conv-E} below), which is a key step to establishing the main result of this paper, namely Theorem \ref{thm-tutte-conv0}.

Let us first explain what probability measure we put on $(M,e_0,\tree)$. Recall the two probability measures $\PQ$ and $\PE$ from Definition \ref{def:PQ-PE}, each of which induce a law on continuous-time two-dimensional random walks $\cZ$. We now specialize to the setting where %
\eqb
\dZ=4,\quad p_1=p_2=p_3=p_4=0.25,\quad v_1=(1,0),\, v_2=(-1,0),\, v_3=(0,1), \,v_4=(0,-1).
\label{eq:p-v}
\eqe
We suppose throughout this section that $\cZ$ is sampled either from $\PE$ or $\PQ$ in this setting.
Under $\PQ$ the steps of $\cZ$ are then uniformly distributed on $\{(1,0),(-1,0),(0,1),(0,-1)\}$ and $(M,e_0,\tree)$ has the law of the uniform infinite tree-decorated planar map. We will use the measure $\PE$ when applying Theorem \ref{thm-general-clt-uniform} as described in the previous paragraph, but in our proofs we will need to also use $\PQ$, and estimates are typically easier for $\PQ$ since the laws of various random variables are more explicit.

We will now explain the high-level idea for the construction of the cell configurations. Let $T$ denote the infinite triangulation associated with $(M,e_0,\frk T)$. Recall that $\cV(M)$ and $\cV(\triang)$ are in bijection with the time sets $\cJLm$ and $\cJm$, respectively, via the map $t\mapsto\frk v_t$; see Section \ref{sec:mullin-infinite} for notation. A natural candidate for a cell configuration satisfying (i) above would be the one consisting of cells $(\eta(I_t))_{t\in\cJLm}$, where $I_t$ is the closure of the component of $\R\setminus \cJLm$ having $t$ as its right end-point (see \eqref{eq:It} below). However, this does not define a cell configuration since property \eqref{eq:cell-conf-adjacency} is typically violated. Our idea for constructing the cell configurations is to associate each $t\in\cJLm$ with a finite set $\cT^{\op{cell}}_t\subset \cJm$ (resp.\ $\cT^{\op{large}}_t\subset \cJm$) such that the cell $H^{\op{cell}}_t$ (resp.\ $H^{\op{large}}_t$) associated with $t$ is the union of appropriate subsets of ${\eta}(I_s)$ for $s\in\cT^{\op{cell}}_t$ (resp.\ $s\in\cT^{\op{large}}_t$). 

Much care needs to be taken when defining $\cT^{\op{cell}}_t$ and $\cT^{\op{large}}_t$ such that the assumptions of Theorem \ref{thm-general-clt-uniform} are satisfied. In order to make sure that $H^{\op{cell}}_t$ is connected for each $t$ we require that the union of $\ol{\eta}(I_s)$ for $s\in \cT^{\op{cell}}_t$ is connected for each $t\in\cJLm$. In order to make sure that \eqref{eq:cell-conf-adjacency} holds for $\cH^{\op{cell}}$ we require that if $\{s,t\}\in \cE(M)$ then $s\in \cT^{\op{cell}}_t$ and $t\in \cT^{\op{cell}}_s$. In order to make sure that the moment condition \eqref{eqn-moment} by $\cH^{\op{cell}}$ is satisfied we need to bound $\#\cT^{\op{cell}}_t$ %
along with the shape of each individual set $\eta(I_s)$ for $s\in \cT^{\op{cell}}_t$. We also verify these properties for $\cH^{\op{large}}$ (instead of $\cH^{\op{cell}}$) by the same method. Finally, in order to make sure Assertion 3.\ of Theorem \ref{thm-general-clt-uniform} is satisfied we need the sets $\cT^{\op{large}}_t$ to be sufficiently much larger than the sets $\cT^{\op{cell}}_t$ in a very precise sense. Ergodicity modulo scaling (Property \ref{item-ergodic}) will be automatically satisfied by our construction due to the fact that the 0-quantum cone satisfies a variant of this property. Since the encoding of $(M,e_0,\frk T)$ in terms of $\cZ$ is via the triangulation $\triang$ we will often make use of $\triang$ when verifying that our sets $\cT^{\op{cell}}_t$ and $\cT^{\op{cell}}_t$ satisfy all the required properties.

Observe that when verifying the properties in the previous paragraph we need to keep careful track of two different notions of connectivity: We need to know whether $\eta(I_t)$ and $\eta(I_s)$ touch (as subsets of $\C$) for $s,t\in\cJm$ and we need know whether $\frak v_t$ and $\frak v_s$ are adjacent in $M$ for $s,t\in\cJLm$. Additionally, we will consider whether $\frak v_t$ and $\frak v_s$ are adjacent in $T$ for $s,t\in\cJm$. These three notions of connectivity give rise to the three graphs $(\cJm,\EG)$, $(\cJm,\ET)$ and $(\cJLm,\EM)$ considered below. The latter two graphs are determined by $\cZ$ while (by the mating-of-trees construction; see Section \ref{sec:sle}) the former graph can be defined in a rather similar manner as $(\cJm,\ET)$  using $Z$ instead of $\cZ$. We need to compare the different notions of connectivity and due to this it is important to bound the error $\sup_{t\in I}\|Z_t-\cZ_t\|$ for appropriate intervals $I\subset\R$.

When checking that our cell configurations satisfy the moment bound \eqref{eqn-moment}, we will use Theorem \ref{prop:zero-vs-gamma-walk}, Corollary \ref{cor:zero-vs-gamma}, and Proposition \ref{prop:zero-vs-gamma-ball}, since estimates are typically easier to obtain under $\PQ$ than under $\PE$. However, several of the events appearing in Theorem \ref{prop:zero-vs-gamma-walk} and Corollary \ref{cor:zero-vs-gamma} are measurable with respect to $\cZ|_{I(k)}$ for a deterministic $k\in\N$, while the sets 
$\cT^{\op{cell}}_t$ and $\cT^{\op{large}}_t$ depend on the full process $\cZ$ in addition to $Z$. In order to get into a setting where we may apply Theorem \ref{prop:zero-vs-gamma-walk} and Corollary \ref{cor:zero-vs-gamma}, we will define variants $\cT^{h,k}_t$ of these sets for $h>1$ and $k\in\N$ which depend only on $\cZ|_{I(k)}$; see Section \ref{sec:cells2}. When applying the cited theorem and corollary, it will also be important to control the maximum and minimum element of the various time sets and this is achieved in Section \ref{sec:max-min-time-sets}. We will also use estimates from Sections \ref{sec:kmt} and \ref{sec:0-cone-BL} along with Appendix \ref{app2}.

In Section \ref{sec:cells} we will define the cell configurations and verify that they satisfy all requirements of Theorem \ref{thm-general-clt-uniform} except possibly the moment condition. In Section \ref{sec:moment} we verify that the moment condition is also satisfied.

\subsection{Time sets for the cell configurations}
\label{sec:cells}
Suppose $(\cZ_t)_{t\in\R}$ is the continuous-time two-dimensional random walk sampled according to $\PE$ or $\PQ$ in the setting of \eqref{eq:p-v}. Recall that the steps of $\cZ$ are uniformly and independently distributed on $(1,0),(-1,0),(0,1),(0,-1)$ when sampled from $\PQ$.
We define time sets $\cJLm,\cJRm,\cJm,\cJLp,\cJRp,\cJp\subset\R$ for $\cZ$ as in Section \ref{sec:bijection}, e.g.\ $\cJLm=\{t\in\R\,:\,\Delta\cL_t<0 \}$.\footnote{Note that this is a slight abuse of notation since the walk in Section \ref{sec:bijection} is  indexed by $\Z$ while the process $\cZ$ in this section is indexed by $\R$, but the definitions are identical.} We also define $\frk v_t$ for $t\in\cJm$ and $\frk e_t$ for $t\in\cJ:=\cJm\cup\cJp$ as in that section.

In the remainder of this section we do the following. In Section \ref{sec:cells1} we associate each $t\in\cJm$ and $d\in\N$ with two subsets of $\cJm$, namely $\cT^{\op{C},d}_t$ and $\cT^{\op{D},d}_t$. In Section \ref{sec:cells2} we define variants of $\cT^{\op{D},d}_t$ which are locally defined by $\cZ$, i.e., they can be defined as a function of $\cZ|_{I(m)}$ for fixed $m\in\N$. We need to work with these locally defined variants of the time sets when applying Theorem \ref{prop:zero-vs-gamma-walk} and Corollary \ref{cor:zero-vs-gamma}. Bounds for the minimal and maximal element of the various time sets are found in Section \ref{sec:max-min-time-sets}.  
In Section \ref{sec:cells3} we define time sets $\cT^{\op{cell}}_t$ and $\cT^{\op{large}}_t$ using $\cT^{\op{C},1}_t$ and $\cT^{\op{C},3}_t$, respectively, 
and we use these time sets to define the two cell configurations which will be used when applying Theorem \ref{thm-general-clt-uniform}. 

\subsubsection{Time sets for cell configurations}
\label{sec:cells1}
We will associate each $t\in\cJm$ and $d\in\N$ with two subsets of $\cJm$, denoted by $\cT^{\op{C},d}_t$ and $\cT^{\op{D},d}_t$. Recall that there is a bijection between $\cJm$ and the vertices of the triangulation $\triang$, and the sets $\cT^{\op{C},d}_t$ and $\cT^{\op{D},d}_t$ should be viewed as vertex sets associated with the vertex $\frak v_t$. The variable $d\in\N$ indicates how large the sets are. The only result from this section which will be used outside of this section and Section \ref{sec:cells2} is  
Proposition \ref{prop:time-sets} right below and the reader can skip the rest of the section on the first reading if desired.

For $t\in\cJm$ define
\eqb
I_{t}:=[t^-,t], \quad\text{where}\quad
t^- := 
\begin{cases}
	\sup\{s\in\cJLm\,:\,s<t \} &\text{if\,} t\in\cJLm,\\
	\sup\{s\in\cJRm\,:\,s<t \} &\text{if\,} t\in\cJRm.
\end{cases}
\label{eq:It}
\eqe 
Notice that $\{ I_{t}\,:\, t\in\cJLm \}$ is a collection of intervals with union $\R$ and pairwise disjoint interior, and that the same holds with $\cJRm$ instead of $\cJLm$.

Let us define two graphs with vertex set $\cJm$, namely $(\cJm,\ET)$ and $(\cJm,\EG)$, where $\ET,\EG\subset\{ \{s,u\}\,:\,s,u\in\cJm \}$ are the collection of edges. The letter ``T'' (resp.\ ``G'') is chosen since it is the first letter of the word ``triangulation'' (resp.\ ``geometric''). 
For $s,u\in\cJm$ we have $\{s,u \}\in\ET$ if and only if the vertices $\frak v_s$ and $\frak v_u$ are adjacent in the triangulation $\triang$ encoded by $\cZ$.   
For $s,u\in\cJm$ we have $\{s,u \}\in \EG$ if and only if $s\neq u$ and at least one of the three following properties hold:
\begin{compactitem}
	\item[(i)] $s,u\in\cJLm$ and there are $\wt s\in I_s$ and $\wt u\in I_u$ such that
	\eqb
	L_{\wt s}=L_{\wt u}=\inf_{t\in[\wt s\wedge \wt u,\wt s\vee\wt u]} L_t,
	\label{eq:line-under}
	\eqe
	\item[(ii)] the property in (i) holds with $R$ and R instead of $L$ and L, respectively, or
	\item[(iii)] $I_s\cap I_u\neq\emptyset$.
\end{compactitem} 
We say that $J\subset\cJm$ is T-connected (resp.\ G-connected) if it is connected viewed as a subset of the vertices of the graph $(\cJm,\ET)$ (resp.\ $(\cJm,\EG)$). 

Our motivation for the definition of $\EG$ is that, as explained in \cite[Section 8.2]{wedges}, if $u,s\in\cJLm$ (resp.\ $u,s\in\cJLm$) then property (i) (resp.\ property (ii)) is satisfied if and only if $\eta(I_u)$ and $\eta(I_s)$ have a non-empty intersection. Furthermore, property (iii) implies that $\eta(I_u)$ and $\eta(I_s)$ have a non-empty intersection. For $J\subset\cJm$ define
\eqb
\ol\eta(J):=\bigcup_{s\in J} \eta(I_s)
\label{eq:ol-eta}
\eqe
and observe that if $J$ is G-connected then $\ol\eta(J)$ is a connected subset of $\BB C$.

\begin{figure}
	\centering
	\includegraphics[scale=1]{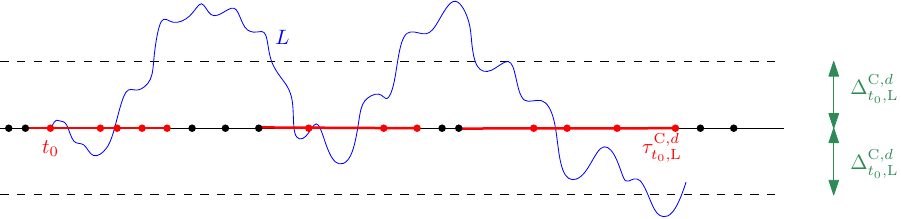}
	\caption{The dots represent the times in $\cJLm$, i.e., the times at which $\cL$ makes a downward jump. The red dots represent times in $\cT_{t_0,\op{L}}^{\op{C},d}$ and the red intervals correspond to intervals $I_t$ for $t\in \cT_{t_0,\op{L}}^{\op{C},d}$. The black dots represent times in $\cJLm\setminus \cT_{t_0,\op{L}}^{\op{C},d}$.}
	\label{fig:Tcell}
\end{figure}

The main goal of this section is to prove the following proposition. Observe that the proposition is an immediate consequence of Lemmas \ref{prop-time-set-cell} and \ref{prop-time-set-graph} below. The letters ``C'' and ``D'' are chosen since they are the first letter of the words ``continuum'' and ``discrete''. As you can see below, $\cT^{\op{C},d}_s$ and $\cT^{\op{D},d}_s$ are defined in a similar way, except that the roles of $\cZ$ and $Z$ have been swapped, and the continuum process $Z$ (resp., the discrete process $\cZ$) plays a more prominent role when defining  $\cT^{\op{C},d}_s$ (resp.\  $\cT^{\op{D},d}_s$). %
\begin{prop} For each $s\in\cJm$ and $d\in\N$ there are sets 
	$\cT^{\op{C},d}_s,\cT^{\op{D},d}_s\subset\cJm$
	such that the following holds for any $s,u\in\cJm$:
	\begin{enumerate}
		\item {\bf Inclusion}: $s\in \cT^{\op{C},d}_s \subset \cT^{\op{D},d+1}_s$ and $s\in\cT^{\op{D},d}_s\subset \cT^{\op{C},d+1}_s$.
		\item {\bf Symmetry}: $s\in \cT^{\op{C},d}_u$ if and only if $u\in\cT^{\op{C},d}_s$. Similarly, 
		$s\in \cT^{\op{D},d}_u$ if and only if $u\in\cT^{\op{D},d}_s$.
		\item {\bf $T$-adjacency}: If $\{s,u \}\in\ET$ then $u\in\cT^{\op{C},d}_s$.
		\item {\bf $G$-adjacency}:If $\{s,u \}\in\EG$ and $d\geq 2$ then $u\in\cT^{\op{D},d}_s$.
		\item {\bf Connectedness}: $\cT^{\op{C},d}_s$ is G-connected and, if $d\geq 2$, $\cT^{\op{D},d}_s$ is T-connected. 
	\end{enumerate}
\label{prop:time-sets}
\end{prop}
We first define $\cT^{\op{C},d}_s$; see Figure \ref{fig:Tcell}. For an interval $[t_0,t_1]$ define $N_{\op{L}}[t_0,t_1]$ to be the maximal number of consecutive upsteps for $\cL|_{[t_0,t_1]}$ without any downsteps, i.e., 
$$
N_{\op{L}}[t_0,t_1] := \max \{ \#(J\cap\cJLp) \,:\,J\subset[t_0,t_1] \op{ interval},J\cap\cJLm=\emptyset  \},
$$
and then set
$$
\Delta_{\op{L}}[t_0,t_1] 
= \sup_{t\in[t_0,t_1]} |(L_t-\cL_t)-(L_{t_0}-\cL_{t_0})|+N_{\op{L}}[t_0,t_1] +2.
$$
We will use multiple times below that $\Delta_{\op{L}}[t_0,t_1] 
> \sup_{t\in[t_0,t_1]} |(L_t-\cL_t)-(L_{t_0}-\cL_{t_0})|$ since this inequality is crucial when comparing $L$ and $\cL$. The second term in the definition of $\Delta_{\op{L}}[t_0,t_1]$ (namely, $N_{\op{L}}[t_0,t_1]$) will only play a role at a single place, namely when proving the second part of Proposition \ref{prop:time-sets}, Assertion 1., (i.e., $\cT^{\op{D},d}_s\subset \cT^{\op{C},d+1}_s$).

Now fix $d\in\N$ and for $t_0\in\cJ$ define %
\eqb
\tau^{\op{C},d}_{t_0,\op{L}} = \inf\{ t\in\cJLm\,:\, t>t_0 \text{ and } 
\inf_{s\in [t_0,t]} (L_s-L_{t_0}) < -  d\cdot \Delta_{\op{L}}[t_0,t]
\},
\label{eq36}
\eqe
where the infimum of an empty set is infinity, and 
$$
\Delta^{\op{C},d}_{t_0,\op{L}}= d\cdot\Delta_{\op{L}}[t_0,\tau^{\op{C},d}_{t_0,\op{L}} ] .
$$
It follows from Section \ref{sec:max-min-time-sets} below that $\tau^{\op{C},d}_{t_0,\op{L}}$ and $\Delta^{\op{C},d}_{t_0,\op{L}}$ are both finite a.s., but note that we do not need a proof of this in order to define the next objects. Let 
\eqb
 \cT_{t_0,\op{L}}^{\op{C},d} = \{ t\in\cJLm\cap [t_0,\tau^{\op{C},d}_{t_0,\op{L}}] \,:\, \{ L_s\,:\,s\in I_t\} \cap [L_{t_0}-\Delta^{\op{C},d}_{t_0,\op{L}},L_{t_0}+\Delta^{\op{C},d}_{t_0,\op{L}}]\neq\emptyset \}.
 \label{eq:TCd0L}
\eqe

Define $N_{\op{R}}[t_0,t_1],\Delta_{\op{R}}[t_0,t_1],\Delta^{\op{C},d}_{t_0,\op{R}},
	\tau^{\op{C},d}_{t_0,\op{R}},\cT^{\op{C},d}_{t_0,\op{R}}$ in the exact same way by using $\cR$ and $R$ instead of $\cL$ and $L$. 
The following lemma is a  consequence of Lemma \ref{prop-bij-steps} and the above definitions, in particular the choice of $\Delta^{\op{C},d}_{t_0,\op{L}}$ and $\Delta^{\op{C},d}_{t_0,\op{R}}$.
\begin{lem}
	Suppose $t_0\in\cJ$ and let $v_{\op{L}}$ and $v_{\op{R}}$ denote the left and right, respectively, end-point of the edge $\frk e_{t_0}\in\cE(T)$. Let $\wt v_{\op{L}}$ (resp.\  $\wt v_{\op{R}}$) denote the vertex to the left (resp.\ right) of $v_{\op{L}}$ (resp.\  $v_{\op{R}}$ on the frontier of the map $M_{t_0}$ in the sewing procedure.
	
	Then  
	$\frk t_{v_{\op{L}}}, \frk t_{\wt v_{\op{L}}} \in \cT_{t_0,\op{L}}^{\op{C},d}$ 
	and 
	$\frk t_{v_{\op{R}}}, \frk t_{\wt v_{\op{R}}} \in \cT_{t_0,\op{R}}^{\op{C},d}$ for any $d\in\N$ a.s. Furthermore, if $t_0\in\cJm$ then $t_0\in \cT_{t_0,\op{L}}^{\op{C},d}\cup \cT_{t_0,\op{R}}^{\op{C},d}$ a.s.
	\label{prop-edge-end-points}
\end{lem}
\begin{proof}
	We first prove that $\frk t_{v_{\op{L}}}\in\cT_{t_0,\op{L}}^{\op{C},d}$ for $t_0\in\cJp$ and then argue how to modify the proof to get $\frk t_{\wt v_{\op{L}}}\in\cT_{t_0,\op{L}}^{\op{C},d}$. This is sufficient to conclude the proof since  $\frk t_{v_{\op{R}}},\frk t_{\wt v_{\op{R}}} \in\cT_{t_0,\op{R}}^{\op{C},d}$ is proven by an identical argument, while the last assertion of the lemma (namely, $t_0\in \cT_{t_0,\op{L}}^{\op{C},d}\cup \cT_{t_0,\op{R}}^{\op{C},d}$) is an immediate consequence of the definitions. The first assertion of Lemma \ref{prop-bij-steps} gives
	$\frk t_{v_{\op{L}}} = \inf\{ s>t_0\,:\, \cL_s=\cL_{t_0}-1 \}$.
	Therefore, for any $t\in[t_0,\frk t_{v_{\op{L}}}]$ we have $\cL_t-\cL_{t_0}\geq -1$, so 
	\eqb
	L_t - L_{t_0} \geq \cL_t-\cL_{t_0} - \sup_{u\in[t_0,t]} |L_u - L_{t_0} - \cL_u+\cL_{t_0}| 
	\geq- d\cdot\Delta_{\op{L}}[t_0,t],
	\label{eq31}
	\eqe
	which implies that $\frk t_{v_{\op{L}}}< \tau^{\op{C},d}_{t_0,\op{L}}$. We get further that $\Delta^{\op{C},d}_{t_0,\op{L}}$ is greater than or equal to the absolute value of the right side of \eqref{eq31} and that $\frk t_{v_{\op{L}}}\in \cT_{t_0,\op{L}}^{\op{C},d}$ since 
	$$
	|L_{t_{v_{\op{L}}}}-L_{t_0}|
	\leq |\cL_{t_{v_{\op{L}}}}-\cL_{t_0}|
	+ \sup_{t\in[t_0,\frk t_{v_{\op{L}}}]}|L_t - L_{t_0} - \cL_t+\cL_{t_0}|\leq \Delta^{\op{C},d}_{t_0,\op{L}},
	$$
	where we used that $|\cL_{\frk t_{v_{\op{L}}}}-\cL_{t_0}|=1$ to get the second inequality.

	To prove $\frk t_{\wt v_{\op{L}}}\in\cT_{t_0,\op{L}}^{\op{C},d}$ we proceed in the exact same way, except that now Lemma \ref{prop-bij-steps} gives
	$\frk t_{\wt v_{\op{L}}} = \inf\{ s>t_0\,:\, \cL_s=\cL_{t_0}-2 \}$, which implies that for any $t\in[t_0,\frk t_{\wt v_{\op{L}}}]$ we have $\cL_t-\cL_{t_0}\geq -2$. The indented equations of the previous paragraph with $\wt v_{\op{L}}$ instead of $v_{\op{L}}$ follow by the same argument as before.
\end{proof}

For $t_0\in\cJ$ and $s_0\in\cJm$ set
	\eqb
	\cT^{\op{C},d,0}_{t_0} = \cT^{\op{C},d}_{t_0,\op{L}}\cup \cT^{\op{C},d}_{t_0,\op{R}},\qquad
	\cT^{\op{C},d}_{s_0} = \bigcup_{t_0\in\cJ\,:\,s_0\in \cT_{t_0}^{\op{C},d,0} } \cT^{\op{C},d,0}_{t_0}.
	\label{eq-def-Tcell}
	\eqe
\begin{lem}
	For each $s_0,u_0\in\cJm$ the following hold a.s.
	\begin{enumerate}
		\item $s_0\in \cT^{\op{C},d}_{s_0}$.
		\item $s_0\in\cT^{\op{C},d}_{u_0}$ if and only if $u_0\in\cT^{\op{C},d}_{s_0}$.
		\item If $\{s,u \}\in\ET$ then $u\in\cT^{\op{C},d}_{s}$.
		\item $\cT^{\op{C},d}_{s_0}$ is G-connected.
	\end{enumerate}
	\label{prop-time-set-cell}
\end{lem}
\begin{proof}
	When proving Assertion 1.\ we may assume by symmetry in L and R that $s_0\in\cJLm$ (instead of $s_0\in\cJRm$). Then we have by \eqref{eq:TCd0L} that $s_0\in\cT_{s_0,\op{L}}^{\op{C},d}\subset \cT_{s_0}^{\op{C},d}$, which gives Assertion 1. Assertion 2.\ follows from \eqref{eq-def-Tcell}. More precisely, by \eqref{eq-def-Tcell} we have $u_0\in\cT^{\op{C},d}_{s_0}$ if and only if there is a $t_0\in\cJ$ such that $s_0,u_0\in\cT^{\op{C},d,0}_{t_0}$, and the latter condition is symmetric in $s_0$ and $u_0$. 
	
	To prove Assertion 3., let $e\in\cE(T)$ be the edge  connecting $\frk v_s$ and $\frk v_u$. 
	If $e\in\cE(T)\setminus(\frk T\cup\frk T_*)$ let $\frk t_e$ denote the time such that $e=\frk e_{\frk t_e}$. By Lemma \ref{prop-edge-end-points}, $s,u\in\cT^{\op{C},d,0}_{\frk t_e}$, which implies the assertion. Now suppose $e\in \frk T\cup\frk T_*$, assume without loss of generality that $e\in\frk T$, and let $t_0$ denote the time at which the edge $e$ is added in the sewing procedure. Notice that in step  $t_0$ of the sewing procedure we added a triangle with the two new edges $e$ and $e_{t_0}$, with $e$ lying to the left of $e_{t_0}$. In the notation of Lemma \ref{prop-edge-end-points} with this choice of $t_0$, we have $\{s,u \}=\{\frk t_{v_{\op{L}}},\frk t_{\wt v_{\op{L}}} \}$, so the lemma implies that again that $s,u\in\cT^{\op{C},d,0}_{\frk t_0}$.

	To prove Assertion 4., first note that each set $\cT^{\op{C},d}_{t_0,\op{L}}$ is G-connected by \eqref{eq:line-under} and, similarly, $\cT^{\op{C},d}_{t_0,\op{R}}$ is also G-connected. 
	The set $\cT^{\op{C},d,0}_{t_0}$ is G-connected since we can find $s_0\in \cT^{\op{C},d}_{t_0,\op{L}}$ (resp.\  $u_0\in\cT^{\op{C},d}_{t_0,\op{R}}$) such that $t_0\in I_{s_0}$ (resp.\ $t_0\in I_{u_0}$), which implies that $\{u_0,s_0 \}\in\EG$ since (iii) in the definition of $\EG$ is satisfied. 
	Finally $\cT^{\op{C},d}_{s_0}$ is G-connected since all the sets $\cT^{\op{C},d,0}_{t_0}$ in the definition of $\cT^{\op{C},d}_{s_0}$ (equation \eqref{eq-def-Tcell}) contain $s_0$.
\end{proof}

\begin{figure}
	\centering
	\includegraphics[scale=1]{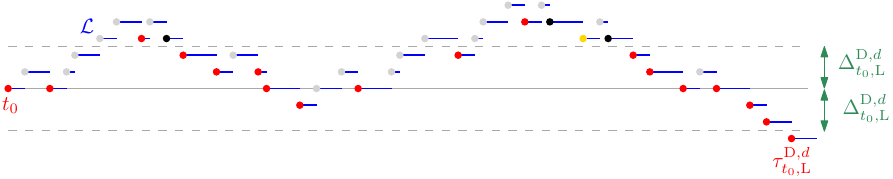}
	\caption{Illustration of $\cT_{t_0,\op{L}}^{\op{D},d}$ and $\wh{\cT}_{t_0,\op{L}}^{\op{D},d}$. The blue graph is $\mathcal L$, 
	the red dots are points $(t,\cL_t)$ for $t\in\cT_{t_0,\op{L}}^{\op{D},d}$,
	the orange dot is the point $(t,\cL_t)$ for $t\in\wh{\cT}_{t_0,\op{L}}^{\op{D},d}\setminus \cT_{t_0,\op{L}}^{\op{D},d}$, 
	the black dots are points $(t,\cL_t)$ for $t\in\cJLm\setminus \wh{\cT}_{t_0,\op{L}}^{\op{D},d}$, and the gray dots are points $(t,\cL_t)$ for $t\in\cJLp$. 
	}
	\label{fig:Tgraph}
\end{figure}

We now define $\cT^{\op{D},d}_{t_0}$; see Figure \ref{fig:Tgraph}. For $t_0\in\cJ$ define 
\eqbn
\begin{split}
	\tau^{\op{D},d}_{t_0,\op{L}} &:= \inf\{ t\in\cJLm\,:\, t>t_0\text{ and }\inf_{s\in[t_0,t]} (\cL_s-\cL_{t_0})<- d\cdot \Delta_{\op{L}}[t_0,t]  \}\\
	&= \inf\{ t\in\cJLm\,:\, t>t_0\text{ and }\cL_t-\cL_{t_0}<
	- d\cdot \Delta_{\op{L}}[t_0,t]  \},
\end{split}
\eqen
where the second equality holds since $\cL$ is piecewise constant with downward jumps only at the times in $\cJLm$. Also define
$$
\Delta^{\op{D},d}_{t_0,\op{L}}= d \cdot \Delta_{\op{L}}[t_0, \tau^{\op{D},d}_{t_0,\op{L}}].
$$
Let 
$$
\cT_{t_0,\op{L}}^{\op{D},d} = \{ t\in\cJLm\cap [t_0,\tau^{\op{D},d}_{t_0,\op{L}}] \,:\, \{ \cL_s\,:\,s\in I_t\} \cap [\cL_{t_0}-\Delta^{\op{D},d}_{t_0,\op{L}},\cL_{t_0}+\Delta^{\op{D},d}_{t_0,\op{L}}]\neq\emptyset \}.
$$

\begin{figure}
	\centering
	\includegraphics[scale=1]{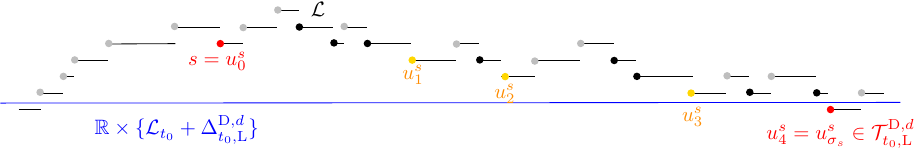}
	\caption{Illustration of the times $u^s_0,u^s_1,\dots$ which are used when defining $\wh{\cT}_{t_0,\op{L}}^{\op{D},d}$. The colors red, orange, black and gray have the same meaning as in Figure \ref{fig:Tgraph}. %
	}
	\label{fig:hat-time-set}
\end{figure}

Next we will define a set $\wh{\cT}_{t_0,\op{L}}^{\op{D},d}$ which is slightly larger than $\cT_{t_0,\op{L}}^{\op{D},d}$. See Figure \ref{fig:hat-time-set}. We will need to use the set $\wh{\cT}_{t_0,\op{L}}^{\op{D},d}$ instead of $\cT_{t_0,\op{L}}^{\op{D},d}$ in \eqref{eq48} below since only $\wh{\cT}_{t_0,\op{L}}^{\op{D},d}$ is T-connected, which is a property we will need when verifying Assertion 4.\ of Lemma \ref{prop-time-set-graph} below. For any $s\in \cT_{t_0,\op{L}}^{\op{D},d}$ define $u_0^s,u_1^s,\dots$ iteratively by $u^s_0=s$ and
$$
u^s_j = \inf\{u > u^s_{j-1}\,:\,\cL_u = \cL_{u^s_{j-1}}-1  \}.
$$
Then set $\sigma_s = \inf\{j\in \N\,:\,u^s_j\in \cT_{t_0,\op{L}}^{\op{D},d} \text{ or } u^s_j \geq  \tau^{\op{D},d}_{t_0,\op{L}} \}$ and
$$
\wh{\cT}_{t_0,\op{L}}^{\op{D},d} = \bigcup_{s\in \cT_{t_0,\op{L}}^{\op{D},d}} \{ u^s_0,\dots,u^s_{\sigma_s-1} \}.
$$  
Define $\Delta^{\op{D},d}_{t_0,\op{R}},\tau^{\op{D},d}_{t_0,\op{R}},
\cT_{t_0,\op{R}}^{\op{D},d}$ in the exact same way by using $\cR$ and $R$ instead of $\cL$ and $L$, and set
\eqb
\cT_{t_0}^{\op{D},d,0}=
\wh{\cT}_{t_0,\op{L}}^{\op{D},d}\cup \wh{\cT}_{t_0,\op{R}}^{\op{D},d}.
\label{eq48}
\eqe
Finally, for $s_0\in\cJm$ set
\eqb
\cT^{\op{D},d}_{s_0} = \bigcup_{t_0\in\cJ\,:\,s_0\in \cT_{t_0}^{\op{D},d,0} } \cT^{\op{D},d,0}_{t_0}.
\label{eq:Tgraph}
\eqe
\begin{lem}
	For each $t_0\in\cJm$ the set $\cT^{\op{D},d}_{t_0}$ satisfies the following a.s.
	\begin{enumerate}
		\item $t_0\in \cT^{\op{D},d}_{t_0}$.
		\item $s\in\cT^{\op{D},d}_{t_0}$ if and only if $t_0\in\cT^{\op{D},d}_s$.
		\item If $\{s,u \}\in\EG$ and $d\geq 2$ then $u\in\cT^{\op{D},d}_s$.
		\item $\cT^{\op{D},d}_{t_0}$ is T-connected for $d\geq 2$.
		\item $\cT^{\op{C},d}_{t_0}\subseteq\cT^{\op{D},d+1}_{t_0}$ and $\cT^{\op{D},d}_{t_0}\subseteq\cT^{\op{C},d+1}_{t_0}$.
	\end{enumerate}
	\label{prop-time-set-graph}
\end{lem}
\begin{figure}
	\centering
	\includegraphics{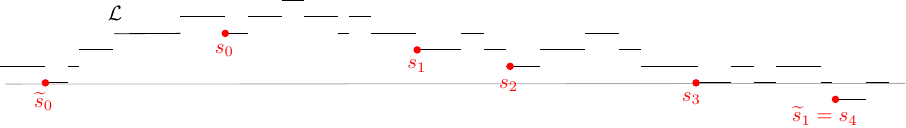}
	\caption{Illustration of the proof of Assertion 4.\ of Lemma \ref{prop-time-set-graph}. The red dots are the points $(t,\cL_t)$ for $t=\wt s_0,\wt s_1,s_0,s_1,\dots,s_{k_0}$ with $k_0=4$. The black graph is the graph of $\cL$. The gray line is $\R\times\{\wt s_0 \}$.}
	\label{fig:s0s1}
\end{figure}
\begin{proof}
	Assertions 1.\ and 2.\ are proven exactly as the corresponding assertions of Lemma \ref{prop-time-set-cell}.
	
	To prove Assertion 3., we first consider the case where (iii) in the definition of $\EG$ is satisfied. Suppose without loss of generality that $s<u$ and that $u\in\cJLm$. Then $u\in \cT^{\op{D},d}_{s,\op{L}}$, which implies $u\in \cT^{\op{D},d}_{s}$ as desired. The case where (ii) in the definition of $\EG$ is satisfied is proven exactly as case (i), so in order to conclude it is sufficient to treat the case where (i) in the definition of $\EG$ is satisfied, i.e., $s,u\in\cJLm$ and there are $\wt s\in I_s$ and $\wt u\in I_u$ such that \eqref{eq:line-under} is satisfied. Assume without loss of generality that $s<u$, which implies $\wt s<\wt u$ since $I_s$ and $I_u$ have disjoint interior. Let $\ul s=\sup(\cJ\cap(-\infty,\wt s] )$. Then $s\in \cT^{\op{D},d}_{\ul s,\op{L}}$ since $\ul s\in I_s$. We will argue now that $u\in \cT^{\op{D},d}_{\ul s,\op{L}}$, which is sufficient to conclude the proof of Assertion 3.\ since it implies that $u\in \cT^{\op{D},d}_{s}$.
	
	By the definition of $\ul s$ and since $\cL$ is piecewise constant with jumps only in $\cJL$, we have $\cL_{\wt s}=\cL_{\underline{s}}$. Observe that $\tau_{\underline{s},\rmL}^{\rmD,d}\geq u$ since for any $t\in (\underline{s},\wt u]$, 
		\begin{equation*}
			\begin{split}
				\cL_t-\cL_{\underline{s}}&=\cL_t-\cL_{\wt s}\geq L_t-L_{\wt s} - \sup_{r\in [\wt s,t]}|(\cL_{r}-\cL_{\wt s})-(L_{r}-L_{\wt s})|
				\geq -2\sup_{r\in [\underline s,t]}|(\cL_{r}-\cL_{\underline s})-(L_{r}-L_{\underline s})|,
			\end{split}
		\end{equation*}
		where in order to get the last inequality we use that $L_t\geq L_{\wt s}$ for $t\in [\wt s,\wt u]$ since \eqref{eq:line-under} is satisfied. By the same argument and using $L_{\wt u}=L_{\wt s}$,
		\begin{equation}\label{eq:pf-def-T-Dd-1}
			\begin{split}
				|\cL_{\wt u}-\cL_{\underline{s}}|\leq 2\sup_{r\in [\underline s,\wt u]}|(\cL_{r}-\cL_{\underline s})-(L_{r}-L_{\underline s})| < \Delta^{\op{D},d}_{\ul s,\op{L}},
			\end{split}
		\end{equation}
	which implies $u\in \cT_{\underline{s},\rmL}^{\rmD,d}$ as desired.

	We prove Assertion 5.\ before Assertion 4.\ since a partial result in the former proof is needed in the proof of the latter. To prove the first part of Assertion 5., first observe that $\tau^{\op{D},d+1}_{t_0,\op{L}}\geq \tau^{\op{C},d}_{t_0,\op{L}}$ since for any $s\in\cJLm$ satisfying $s<\tau^{\op{C},d}_{t_0,\op{L}}$ 
	\eqb
	\begin{split}
		\inf_{t\in [t_0,s]} (\cL_t-\cL_{t_0}) 
		&\geq 
		\inf_{t\in [t_0,s]} (L_t-L_{t_0})
		- \sup_{t\in[t_0,s]} |(L_t-\cL_t)-(L_{t_0}-\cL_{t_0})|\\
		&\geq 
		- d \Delta_{\op{L}}[t_0,s]
		- \sup_{t\in[t_0,s]} |(L_t-\cL_t)-(L_{t_0}-\cL_{t_0})| 
		\geq - (d+1)\Delta_{\op{L}}[t_0,\tau^{\op{C},d}_{t_0,\op{L}}],
	\end{split}
	\label{eq32}
	\eqe
	where we use the definition of $\tau^{\op{C},d}_{t_0,\op{L}}$ and the fact that $s<\tau^{\op{C},d}_{t_0,\op{L}}$ to get the second inequality. 
	This implies further that 
	$$
	\Delta^{\op{D},d+1}_{t_0,\op{L}}\geq  \Delta^{\op{C},d}_{t_0,\op{L}}
	+ \sup_{t\in[t_0,\tau^{\op{C},d}_{t_0,\op{L}}]}|(L_t-\cL_t)-(L_{t_0}-\cL_{t_0})|.
	$$ 
	If $t\in \cT^{\op{C},d}_{t_0,\op{L}}\setminus\{t_0 \}$ then we can find $r\in I_t$ such that $|L_r-L_{t_0}|\leq \Delta^{\op{C},d}_{t_0,\op{L}}$. Using that $r\in [t_0,\tau^{\op{C},d}_{t_0,\op{L}}]$, we get
	\eqb
	|\cL_r - \cL_{t_0}| \leq |L_r - L_{t_0}| + \sup_{t\in[t_0,\tau^{\op{C},d}_{t_0,\op{L}}]}|(L_t-\cL_t)-(L_{t_0}-\cL_{t_0})| \leq \Delta^{\op{D},d+1}_{t_0,\op{L}},
	\label{eq49}
	\eqe
	which implies $t\in \cT^{\op{D},d+1}_{t_0,\op{L}}$. Since $t$ was arbitrary, we get $\cT^{\op{C},d}_{t_0,\op{L}}\subset \cT^{\op{D},d+1}_{t_0,\op{L}}\subset \wh{\cT}^{\op{D},d+1}_{t_0,\op{L}}$. The latter result gives the first part of Assertion 5. 
	
	The second part of Assertion 5.\ is proven via a similar argument, except that some extra care is needed since we use $\wh{\cT}^{\op{D},d}_{t_0,\op{L}}$ and $\wh{\cT}^{\op{D},d}_{t_0,\op{L}}$ (instead of $\cT^{\op{D},d}_{t_0,\op{L}}$ and $\cT^{\op{D},d}_{t_0,\op{L}}$) when defining ${\cT}^{\op{D},d,0}_{t_0}$ in \eqref{eq48}. We first argue as above that 
	$\tau^{\op{C},d+1}_{t_0,\op{L}}\geq \tau^{\op{D},d}_{t_0,\op{L}}$ and
	$$
	\Delta^{\op{C},d+1}_{t_0,\op{L}}\geq  \Delta^{\op{D},d}_{t_0,\op{L}}
	+ \sup_{t\in[t_0,\tau^{\op{D},d}_{t_0,\op{L}}]}|(L_t-\cL_t)-(L_{t_0}-\cL_{t_0})|
	+ N_{\op{L}}[t_0, \tau^{\op{D},d}_{t_0,\op{L}}].
	$$
	If $t\in \cT^{\op{D},d}_{t_0,\op{L}}\setminus\{t_0\}$ then we can find $r\in I_t$ such that $|\cL_r-\cL_{t_0}|\leq \Delta^{\op{D},d}_{t_0,\op{L}}$. Therefore, by the definition of $N_{\op{L}}[t_0, \tau^{\op{D},d}_{t_0,\op{L}}]$, if $t\in \wh{\cT}^{\op{D},d}_{t_0,\op{L}} \setminus\{t_0\}$ then we can find $r\in I_t$ such that $|\cL_r-\cL_{t_0}|\leq \Delta^{\op{D},d}_{t_0,\op{L}}+N_{\op{L}}[t_0, \tau^{\op{D},d}_{t_0,\op{L}}]$, which implies by the last display and a similar argument as \eqref{eq49} that 
	$|L_r-L_{t_0}|\leq \Delta^{\op{C},d+1}_{t_0,\op{L}}$, so
	 $t\in \cT^{\op{C},d+1}_{t_0,\op{L}}$. We conclude as before.

	To prove Assertion 4.\ we will first argue that  $\wh{\cT}^{\op{D},d}_{t_0,\op{L}}$ is T-connected. 
	Consider arbitrary consecutive times $\wt s_0<s_0$ in $\wh{\cT}^{\op{D},d}_{t_0,\op{L}}$, i.e., $\wt s_0,s_0\in \wh{\cT}^{\op{D},d}_{t_0,\op{L}}$ are such that $\wt s_0<s_0$ and $(\wt s_0,s_0)\cap \wh{\cT}^{\op{D},d}_{t_0,\op{L}}=\emptyset$. We will argue that $\wt s_0$ and $s_0$ are connected by a T-path such that all vertices on the path are in $\wh{\cT}^{\op{D},d}_{t_0,\op{L}}$, which is sufficient to conclude that $\wh{\cT}^{\op{D},d}_{t_0,\op{L}}$ is T-connected since $\wt s_0$ and $s_0$ were arbitrary. Define 
	$$
	\wt s_1 = \inf\{ u>\wt s_0\,:\,\cL_u=\cL_{\wt s_0}-1 \}
	$$
	and for $k\in\N$ define $s_k$ iteratively by
	$$
	s_k = \inf\{ u>s_{k-1}\,:\,\cL_u=\cL_{s_{k-1}}-1 \}.
	$$
	See Figure \ref{fig:s0s1} for an illustration. 
	Observe that, by the definition of $\tau^{\op{D},d}_{t_0,\op{L}}$ and $\wh{\cT}^{\op{D},d}_{t_0,\op{L}}$, for any $s\in\wh{\cT}^{\op{D},d}_{t_0,\op{L}}$ such that $s\neq\tau^{\op{D},d}_{t_0,\op{L}}$ it holds that 
	\eqb
	\inf\{u>s\,:\,\cL_u=\cL_s-1 \}\in\wh{\cT}^{\op{D},d}_{t_0,\op{L}}.
	\label{eq:also-in-T0}
	\eqe 
	This implies that $\wt s_1\in \wh{\cT}^{\op{D},d}_{t_0,\op{L}}$. Since $\wt s_0$ and $s_0$ are consecutive elements of $\wh{\cT}^{\op{D},d}_{t_0,\op{L}}$ we get that $s_0\leq \wt s_1$. By the definition of $\wt s_1$ we get further that $\cL_{s_0}\geq \cL_{\wt s_1}$. Combining the two previous sentences and using the iterative definition of the $s_k$'s we get that $\wt s_1=s_{k_0}$ for some $k_0\in\N_0$. 
	By Lemma \ref{prop-bij-steps} we have that $\{\wt s_0,\wt s_1 \}\in\ET$ and that $\{s_{k-1},s_k \}\in\ET$ for any $k\in\N$, which implies that $s_0,s_1,\dots,s_{k_0},\wt s_0$ is a path in the graph $(\cJm,\ET)$. By using \eqref{eq:also-in-T0} and the iterative definition of the $s_k$'s we get that all vertices on this path are in $\wh{\cT}^{\op{D},d}_{t_0,\op{L}}$. 
	We conclude that $\wh{\cT}^{\op{D},d}_{t_0,\op{L}}$ is T-connected. By symmetry we get that $\wh{\cT}^{\op{D},d}_{t_0,\op{R}}$ is also T-connected. 
	
	Observe that Lemma \ref{prop-edge-end-points} still holds with $\wh{\cT}^{\op{D},d}_{t_0,\op{L}}$ and $\wh{\cT}^{\op{D},d}_{t_0,\op{R}}$ (with $d\geq 2$) instead of $\cT^{\op{C},d}_{t_0,\op{L}}$ and $\cT^{\op{C},d}_{t_0,\op{R}}$ (with $d\geq 1$) since, as argued in the proof of Assertion 5., $\cT^{\op{C},d-1,0}_{t_0,\op{L}}\subset \wh{\cT}^{\op{D},d}_{t_0,\op{L}}$ and $\cT^{\op{C},d-1,0}_{t_0,\op{R}}\subset \wh{\cT}^{\op{D},d}_{t_0,\op{R}}$ for $d\geq 2$. By using this observation along with the fact that $\wh{\cT}^{\op{D},d}_{t_0,\op{L}}$ and $\wh{\cT}^{\op{D},d}_{t_0,\op{R}}$ are T-connected, we get that $\cT^{\op{D},d,0}_{t_0}$ is also T-connected. Finally, $\cT^{\op{D},d}_{s_0}$ is T-connected since $s_0$ is a common element of all the sets $\cT^{\op{D},d,0}_{t_0}$ in \eqref{eq:Tgraph}.
	
\end{proof}

\subsubsection{Time sets determined locally by $\cZ$}\label{sec:cells2}

In order to define our cell configurations and argue that they satisfy the conditions of Theorem \ref{thm-general-clt-uniform} we need various estimates for the time sets $\cT^{\op{C},d}_t$ and $\cT^{\op{D},d}_t$ defined in Section \ref{sec:cells1}. Estimates of this type are typically easier to achieve under the probability measure $\PQ$ than the probability measure $\PE$ since, under $\PQ$, $Z$ has the law of a Brownian motion and the law of $\cZ$ is explicit. In order to get estimates under $\PE$ we will then use Theorem \ref{prop:zero-vs-gamma-walk}, Corollary \ref{cor:zero-vs-gamma}, and Proposition \ref{prop:zero-vs-gamma-ball}.

Several of the events appearing in Theorem \ref{prop:zero-vs-gamma-walk} and Corollary \ref{cor:zero-vs-gamma} are measurable with respect to $\cZ|_{I(k)}$ for a deterministic $k\in\N$. On the other hand, the time sets  $\cT^{\op{C},d}_t$ and $\cT^{\op{D},d}_t$ defined in Section \ref{sec:cells1} also depend on $Z$ and we need to observe the full process $\cZ$ (not only $\cZ$ restricted to $I(k)$ for some $k$) in order to define them. In order to get into a setting where we may apply Theorem \ref{prop:zero-vs-gamma-walk} and Corollary \ref{cor:zero-vs-gamma} we will define variants $\cT^{h,k}_t$ of $\cT^{\op{D},d}_t$ for $h>1$ and $k\in\N\cup\{\infty \}$ which depend only on $\cZ|_{I(k)}$, where we set $I(\infty)=\R$. These sets $\cT^{h,k}_t$ are defined just as $\cT^{\op{D},d}_t$, except for the following three changes: (i) In order for $\cT^{h,k}_t$ to not depend on $Z$ we replace $\Delta^{\op{D},d}_{t,\op{L}}$ by a deterministic constant $h>1$, (ii) we require that $\cT^{h,k}$ is a subset of $I(k)$, which will imply that we do not need to observe $\cZ$ outside of the interval $I(k)$ in order to define it, and 
(iii) we do not define the equivalents of 
$\wh{\cT}^{\op{D},d}_{t_0,\op{L}}$ and $\wh{\cT}^{\op{D},d}_{t_0,\op{R}}$, but use the equivalents of $\cT^{\op{D},d}_{t_0,\op{L}}$ and $\cT^{\op{D},d}_{t_0,\op{R}}$ instead.

For $t_0\in\cJ$ and $h\geq 1$ define
$$
\tau^{h}_{t_0,\op{L}} 
:= \inf\{ t\in\cJLm\cap[t_0,\infty)\,:\, \inf_{s\in[t_0,t]} (\cL_s-\cL_{t_0})<-h  \}
= \inf\{ t\in\cJLm\cap[t_0,\infty)\,:\,  \cL_t-\cL_{t_0}<-h  \},
$$
where the second equality holds since $\cL$ is piece-wise constant with downward jumps only at the times in $\cJLm$. For $t_0\in I(m)$ set
\eqb
\cT_{t_0,\op{L}}^{h,m} 
=
\{ t\in\cJLm\cap [t_0,\tau^{h}_{t_0,\op{L}}]\cap I(m) \,:\, \{ \cL_s\,:\,s\in I_t\} \cap [\cL_{t_0}-h,\cL_{t_0}+h]\neq\emptyset \},
\label{eq29}
\eqe
and set $\cT_{t_0,\op{L}}^{h,m} =\emptyset$ if $t_0\not\in I(m)$.
Define $\tau^{h}_{t_0,\op{R}}$ and $\cT^{h,m}_{t_0,\op{R}}$  in the exact same way by using $\cR$ instead of $\cL$, and set
\eqb
\cT_{t_0}^{h,m,0}=\cT_{t_0,\op{L}}^{h,m}\cup \cT_{t_0,\op{R}}^{h,m}.
\label{eq27}
\eqe
Finally, for $s_0\in\cJm$ set
\eqb
\cT^{h,m}_{s_0} = \bigcup_{t_0\in\cJ\,:\,s_0\in \cT_{t_0}^{h,m,0} } \cT^{h,m,0}_{t_0}.
\label{eq28}
\eqe

\begin{lem}\label{prop:localT}
	\begin{compactitem}
		\item[(i)] $s_0\in\cT^{h,m}_{s_0} \subset \cT^{h,m+1}_{s_0} \subset \cT^{h,\infty}_{s_0}$ for any $s_0\in\cJm$ such that $s_0\in I(m)$. Furthermore, for $s_0,t_0\in\cJm$ we have $s_0\in \cT^{h,m}_{t_0}$ if and only if $t_0\in \cT^{h,m}_{s_0}$.
		\item[(ii)] $\cT^{h,m}_{s_0}$ is a measurable function of $\cZ|_{I(m)}$.
		\item[(iii)] If $I_u\subset I(m)$ for all $u\in\cT^{h,\infty}_{s_0}$ then $\cT^{h,m}_{s_0}=\cT^{h,\infty}_{s_0}$.
		\item[(iv)] 
		\eqb
		\begin{split}
			&\{ d\sup_{s,t\in I(m)} |(L_t-\cL_t)-(L_{s}-\cL_{s})|
			+(d+1)N_{\op{L}}[I(m)]+2d <h \}
			\\&\qquad
			\cap 
			\{ d\sup_{s,t\in I(m)} |(R_t-\cR_t)-(R_{s}-\cR_{s})|
			+(d+1)N_{\op{R}}[I(m)]+2d <h \}\\
			&\qquad\qquad\subset \{ \cT^{\op{D},d}_{t}\subset \cT^{h,\infty}_{t}\,\,\forall t\in\cJm \text{\,s.t.\,}\cT^{\op{D},d}_{t}\subset I(m) \},
		\end{split}
		\label{eq18}
		\eqe
		where for an interval $J$ with end-points $t_1<t_2$ we define $N_{\op{L}}[J]:=N_{\op{L}}[t_1,t_2]$.
	\end{compactitem}
\end{lem}
\begin{proof}
	Assertions (i), (ii), and (iii) are immediate by the construction. To prove Assertion (iv), pick $s_0\in J$ such that $\cT^{\op{D},d}_{s_0}\subset I(m)$. Recall that $\cT^{\op{D},d}_{s_0}$ is the union of sets $\wh\cT_{t_0,\op{L}}^{\op{D},d}$ and $\wh\cT_{t_0,\op{R}}^{\op{D},d}$, which much also be contained in $I(m)$, so in particular $\tau^{\op{D},d}_{t_0,\op{L}}\in \wh\cT_{t_0,\op{L}}^{\op{D},d}\subset I(m)$ and $\tau^{\op{D},d}_{t_0,\op{R}}\in \wh\cT_{t_0,\op{R}}^{\op{D},d}\subset I(m)$. For such a $t_0$, on the first event on the left side of \eqref{eq18},
	\eqb
	\Delta_{t_0,\op{L}}^{\op{D},d} + N_{\op{L}}[t_0,\tau^{\op{D},d}_{t_0,\op{L}}] 
	\leq  
	\sup_{t\in[t_0,\tau^{\op{D},d}_{t_0,\op{L}}]} d|(L_t-\cL_t)-(L_{t_0}-\cL_{t_0})|+(d+1)N_{\op{L}}[I(m)]+2d
	<h, 
	\label{eq17}
	\eqe
	which implies that $\wh{\cT}^{\op{D},d}_{t_0,\op{L}}\subset \cT^{h,\infty}_{t_0,\op{L}}$. Similarly, $\wh{\cT}^{\op{D},d}_{t_0,\op{R}}\subset \cT^{h,\infty}_{t_0,\op{R}}$. It follows from \eqref{eq:Tgraph} and \eqref{eq28} that $\cT^{\op{D},d}_{s_0}\subset \cT^{h,\infty}_{s_0}$ as desired.
\end{proof}

\subsubsection{Bound for maximum and minimum element of time sets}
\label{sec:max-min-time-sets}
When proving estimates for the time sets $\cT^{\op{C},d}_t,\cT^{\op{D},d}_t,\cT^{h,m}_t$ (defined in Sections \ref{sec:cells1} and \ref{sec:cells2}) and the cell configurations $\cH^{\op{cell}}$ and $\cH^{\op{large}}$ (defined below in Section \ref{sec:cells3}) it will be essential to have estimates for the minimal and maximal element of the various time sets. In this section we obtain such estimates for $\cT^{\op{C},d}_t$ and $\cT^{h,\infty}_t$ under the measures $\PQ$ and $\PE$.

\begin{prop}\label{prop:Text-len-tail}
	\begin{compactitem}
		\item[(i)] Let $d\in\N$. 
		Let $J\subset\R$ denote a random interval containing 0 such that $\PE[J\not\subset I(k)]\leq 2^{-c_0 k}$ for some $c_0>0$ and for all $k\in\N_+$. Then there is a constant $c_1=c_1(d,c_0)$ such that for all $k\in\N_+$,
		\eqbn
		\begin{split}
			\PE\Big[\bigcup_{t\in J} \cT^{\op{C},d}_t \not\subset I(k) \Big] \leq 2^{-c_1 k}.
		\end{split}
		\eqen 
		\item[(ii)] Assertion (i) still holds if we replace $\PE$ by $\PQ$.
		\item[(iii)] 
		Suppose there is a constant $c_0>0$ and intervals $J_h\subset\R$ for $h>1$ such that for any $k\in\N$ and $h>1$ we have $\PE[ J_h\not\subset I(k) ]\leq h^{1/c_0}2^{-c_0k}$. Then for some $c_1>0$,
		\eqbn
		\begin{split}
			\PE\Big[\bigcup_{t\in J_h} \cT^{h,\infty}_t \not\subset I(k) \Big] \leq h^{c_1^{-1}}2^{-c_1 k}.
		\end{split}
		\eqen		
		\item[(iv)] Assertion (iii) still holds if we replace $\PE$ by $\PQ$.
	\end{compactitem}
\end{prop}

The proof for $\PQ$ can be done relatively easily since the laws of both $Z$ and $\cZ$ are explicit under this probability measure. For $\PE$ the proof is much harder since we a priori know little about both $Z$ and $\cZ$. Assertions (i) and (iii) are proven in a rather similar manner, so let us give a brief outline only for the proof of Assertion (i). Two key partial steps are to bound the value of $\tau^{\op{C},d}_{t,\op{L}}$ 
	(Lemma \ref{prop:tau-k-max}) 
and to bound the probability of the event that $\cT^{\op{C},d}_{t,\op{L}}$ or $\cT^{\op{C},d}_{t,\op{R}}$ for some $t\ll -1$ has a non-empty intersection with a given interval containing 0 
(Lemma \ref{prop11}); once these bounds are proven it is quite easy to conclude the proof using the definition of $\cT^{\op{C},d}_t$.  
When we prove Lemmas \ref{prop:tau-k-max} and \ref{prop11} we use 
properties of $Z$ under $\PE$ (Lemma \ref{prop:eucl-inf-upper} and Lemma \ref{prop:eucl-inf-lower}), Theorem 4.2 (which allows to bound $\PE$-probabilities by $\PQ$-probabilities), and the fact that $Z$ and $\cZ$ are uniformly close on bounded intervals with high probability (Lemma \ref{prop:eucl-error}). 
\begin{lem}
	For any $a_1>0$ and $a_2\in(0,a_1)$ there is a $c\in(0,1)$ such that for any $m\in\N$,
	$$
	\PE\Big[\inf_{t\in [2^m,2^{(1+a_1)m}]} L_t>-2^{0.5(1+a_2)m}\Big] < 2^{-cm}.
	$$
	The same holds for $\PQ$ instead of $\PE$.
	\label{prop:eucl-inf-lower2}
\end{lem}
\begin{proof}
	We have
	\eqbn
	\begin{split}
		&\PE\Big[\inf_{t\in [2^m,2^{(1+a_1)m}]} L_t>-2^{0.5(1+a_2)m}\Big] 
		= 
		\PE\Big[\inf_{t\in [0,2^{(1+a_1)m}]} L_t>-2^{0.5(1+a_2)m}\Big]\\
		&\qquad+ \PE\Big[\inf_{t\in [2^m,2^{(1+a_1)m}]} L_t>-2^{0.5(1+a_2)m}; \inf_{t\in [0,2^m]} L_t\leq -2^{0.5(1+a_2)m}\Big]
		\leq 2^{-cm},
	\end{split}
	\eqen
	where we apply Lemmas \ref{prop:eucl-inf-lower} and \ref{prop:eucl-inf-upper}, along with the scaling property of Lemma \ref{prop:E-scaling} in the last step. The proof for $\PQ$ can be done in the exact same way due to Remark \ref{rmk:Q-inf}. 
\end{proof}

\begin{lem}
	For any $a\in(0,1)$ there is a $c\in(0,1)$ such that for any $m\in\N$, $$\PE[\cJLm\cap[-2^{m},-(1-a)2^{m}]=\emptyset]<2^{-c m}.$$
	The same result holds for $\PQ$ instead of $\PE$.
	\label{prop:JLm-empty}
\end{lem}
\begin{proof}
	The result for $\PQ$ is immediate since $\cJLm$ has the law of a rate $1/4$ Poisson point process under $\PQ$; in fact, the considered probability is bounded above by $2^{-c2^{m}}$ in that case. To prove the estimate for $\PE$ we will apply Theorem \ref{prop:zero-vs-gamma-walk} with $k=2m$ and 
	$$
	G(\cZ,\frk D) = \1(\{ \cJLm\cap [-2^m,-(1-a)2^m] =\emptyset; [-2^m,-(1-a)2^m]\subset I(k) \}),
	$$
	which gives, in the notation of the theorem,
	$$
	A = 
	\bigcup_{u\in I(0)} \{ \cJLm\cap [-2^m+u,-(1-a)2^m+u] =\emptyset;  [-2^m+u,-(1-a)2^m+u]\subset I(k) \}.
	$$
	We choose this particular definition of $G$ since it fits into the context of the cited theorem at the same time as it can be compared to the event we want to bound the probability of. 
	Observe that $A\in\sigma(\cZ|_{I(k)})$. By Theorem \ref{prop:zero-vs-gamma-walk},
	$$
	\PE[A] \leq \PQ[A]^{1/4}2^{C_0 k} + 2^{-k/C_0}.
	$$
	By using that $\cJLm$ has the law of a rate $1/4$ Poisson point process under $\PQ$ we get that $\PQ[A]\leq\PQ[ \cJLm\cap [-2^{m}+1,-(1-a)2^{m}-1] =\emptyset ]\leq 2^{-c2^{cm}}$. Inserting this into the last indented equation we get $\PE[A]\leq 2^{-cm}$. This concludes the proof of the lemma since 
	$$
	\PE[\cJLm\cap[-2^{m},-(1-a)2^{m}]=\emptyset \leq \PE[ A ] + \PE[  [-2^m,-(1-a)2^m]\not\subset I(k) ] \leq 2^{-cm},
	$$
	where in the last step we use that $\frk D$ has the law of a uniform dyadic decomposition under $\PE$.
\end{proof}

\begin{lemma}
	There exist $C>1$ such that for all $m\in\N$ we have $\PE[N_{\op{L}}[-2^m,2^m] \geq Cm ]<2^{-m/C}$. The same result holds with $\PQ$ instead of $\PE$.
	\label{prop22}
\end{lemma}
\begin{proof}
	We proceed similarly as in the previous lemma, but now with
	\eqbn
	\begin{split}
		G(\cZ,\frk D) &= \1(\{ N_{\op{L}}[I(1.1m)]\geq Cm; [-2^m,2^m]\subset I(1.1m) \}),\\
		A &= \bigcup_{u\in I(0)}
		\{ N_{\op{L}}[I(1.1m)]\geq Cm;
		[-2^m+u,2^m+u]\subset I(1.1m) \},
	\end{split}
	\eqen 
	where we define $I(r):=I(\lceil r\rceil)$ for $r\in\R$. 
	The constant $c$ may vary from line to line throughout the proof. By Theorem \ref{prop:zero-vs-gamma-walk} applied with $k=1.1m$,
	\eqb
	\PE[A] \leq \PQ[A]^{1/4}2^{1.1C_0 m} + 2^{-1.1m/C_0}.
	\label{eq73}
	\eqe
	To bound the right side of the last display, observe that 
	\eqbn
	\begin{split}
		\PQ[A] 
		&\leq \PQ[N_{\op{L}}[I(1.1m)]\geq Cm]\\
		&\leq 
		\PQ[ \# (\cJL\cap I(1.1m) ) > 2^{1.1m+1} ]
		+ \PQ[\# (\cJL\cap I(1.1m) ) \leq 2^{1.1m+1}; N_{\op{L}}[I(1.1m)]\geq Cm ].
	\end{split}
	\eqen 
	The first probability is bounded above by $2^{-c2^{1.1m}}$ since $\# (\cJL\cap I(1.1m))$ is a Poisson random variable with parameter $2^{1.1m}$. The second probability is bounded above by $2^{-(C-1.1)m+1}$ since for a simple random walk on $\BB Z$ with $2^{1.1m+1}$ steps, the probability that it makes at least $Cm$ consecutive up-steps is bounded above by $2^{-Cm}\cdot 2^{1.1m+1}= 2^{-(C-1.1)m+1}$. Choosing $C$ sufficiently large and inserting into \eqref{eq73} we get $\PE[A]\leq 2^{-cm}$.
	Combining the results above concludes the proof of the lemma since
	$$
	\PE[N_{\op{L}}[-2^m,2^m] \geq Cm ]
	\leq \PE[A] + \PE[ [-2^m-1,2^m+1]\not\subset I(1.1m) ] \leq 2^{-cm}.
	$$
\end{proof}

\begin{lem}
	\begin{compactitem}
		\item[(i)] For any $a>0$ and $d>1$ there is a $c\in(0,1)$ such that for any $m\in\N$,
		\eqb
		\begin{split}
			&\PE\Big[  \max \{ \tau^{\op{C},d}_{t,\op{L}}\,:\, t\in[-2^m,2^m]\cap\cJ  \} >2^{(1+a)m} \Big]<2^{-cm}.
		\end{split}
		\label{eq39}
		\eqe
		The estimates also hold if we replace $\PE$ by $\PQ$.
		\item[(ii)] For any $a>0$ there is a $c\in(0,1)$ such that for any $m\in\N$ and $h<2^{0.5m}$,
		\eqb
		\begin{split}
			&\PE\Big[  \max \{ \tau^{h}_{t,\op{L}}\,:\, t\in[-2^m,2^m]\cap\cJ  \} >2^{(1+a)m} \Big]<2^{-cm}.
		\end{split}
		\label{eq39h}
		\eqe
		The estimates also hold if we replace $\PE$ by $\PQ$.
	\end{compactitem}
\label{prop:tau-k-max}
\end{lem}
\begin{proof} 
	We first prove Assertion (i). We start by arguing that for any $C>1$,
	\eqb
	\begin{split}
		&\Big\{\sup_{s_1,s_2\in[-2^m,
		2^{(1+a)m}]} |(L_{s_1}-\cL_{s_1}) - (L_{s_2} - \cL_{s_2})|<Cm \Big\}
		\cap \{ N_{\op{L}}[ -2^m,2^{(1+a)m} ]<Cm \}
		\\& 
		\cap\Big\{\inf_{t\in [-2^m,2^m]} L_t>-2^{0.5(1+0.5a)m}\Big\}
		\cap\Big\{\inf_{t\in [2^m,2^{(1+0.9a)m}]} L_t<-2^{0.5(1+0.5a)m}-2Cdm-2d\Big\}
		\\& 
		\cap\{ \cJLm\cap [2^{(1+0.9a)m},2^{(1+a)m}]\neq \emptyset \}
		\subset \{ \max \{ \tau^{\op{C},d}_{t,\op{L}}\,:\, t\in[-2^m,2^m]\cap\cJ  \} \leq 2^{(1+a)m} \}.
	\end{split}
\label{eq40}
	\eqe 
	Suppose that the five events on the left side occur. Let $t\in[-2^m,2^m]\cap\cJ$ and $t_0\in \cJLm\cap [2^{(1+0.9a)m},2^{(1+a)m}]$, where we see that a suitable $t_0$ exists due to the occurrence on the fifth event on the left side of the last display. By first using occurrence of the third and fourth events on the left side of \eqref{eq40}, and then using occurrence of the first and second events on the left side of \eqref{eq40},
	\eqb
	\begin{split}
		\inf_{s\in [t,t_0]} (L_s-L_{t}) 
		&< \big(-2^{0.5(1+0.5a)m}-2Cdm-2d\big) - (-2^{0.5(1+0.5a)m})\\ 
		&<-  
		d\sup_{s\in[t,t_0]} |(L_s-\cL_s)-(L_{t}-\cL_{t})|
		-dN_{\op{L}}[-2^m,2^{(1+a)m}]
		-2d\leq -d\Delta_{\op{L}}[t,t_0].
	\end{split}
	\eqe 
	This implies that $\tau^{\op{C},d}_{t,\op{L}}\leq t_0\leq 2^{(1+a)m}$. Since $t\in[-2^m,2^m]\cap\cJ$ was arbitrary, this implies that the event on the right side of \eqref{eq40} occurs.
	By choosing $C>1$ sufficiently large, the events on the left side of \eqref{eq40} each occur with $\PE$-probability at least $1-2^{-cm}$ by 
	Lemma \ref{prop:eucl-error},  
	Lemma \ref{prop22}, 
	Lemma \ref{prop:eucl-inf-upper},
	Lemma \ref{prop:eucl-inf-lower2}, and
	Lemma \ref{prop:JLm-empty},
	respectively. This concludes the proof of Assertion (i) for $\PE$. The case of $\PQ$ is handled in the exact same way, except that we use \eqref{eq:err4} (instead of Lemma \ref{prop:eucl-error}) and Remark \ref{rmk:Q-inf}. 
	In order to establish Assertion (ii), by similar considerations as in the proof of Assertion (i),
	\eqb
	\begin{split}
		&\Big\{\inf_{t\in [-2^m,2^m]} L_t>-2^{0.5(1+0.5a)m}\Big\}
		\cap\Big\{\inf_{t\in [2^m,2^{(1+0.9a)m}]} L_t<-2^{0.5(1+0.5a)m}-h\Big\}\\
		&\qquad 
		\cap\{ \cJLm\cap [2^{(1+0.9a)m},2^{(1+a)m}]\neq \emptyset \}
		\subset \{ \max \{ \tau^{h}_{t,\op{L}}\,:\, t\in[-2^m,2^m]\cap\cJ  \} \leq 2^{(1+a)m} \}.
	\end{split}
	\eqe 
	The events on the left side have $\PE$-probability at least $1-2^{-cm}$ by Lemmas \ref{prop:eucl-inf-upper},   \ref{prop:eucl-inf-lower2}, and \ref{prop:JLm-empty}, respectively, where we use that $h\leq 2^{0.5m}$ when bounding the second probability. This completes the proof of Assertion (ii) for $\PE$. Again the proof for $\PQ$ is identical. 

\end{proof}

\begin{lem}
	\begin{compactitem}
		\item[(i)] There is a constant $c>0$ such that for all $m\in\N$, %
		$$
		\PE[\exists t\in [-2^{m+1},-2^m]\cap\cJ\,:\,\cT^{\op{C},d}_{t,\op{L}}\cap[-2^{0.8m}, 2^{0.8m}]\neq\emptyset] \leq 2^{-cm}.
		$$
		The same holds for $\PQ$ instead of $\PE$.
		\item[(ii)] There is a constant $c>0$ such that for all $m\in\N$ and $h>1$ such that $h<2^{0.45m}$,
		$$
		\PE[\exists t\in [-2^{m+1},-2^m]\cap\cJ\,:\,\cT^{h,\infty}_{t,\op{L}}\cap[-2^{0.8m}, 2^{0.8m}]\neq\emptyset] \leq 2^{-cm}.
		$$
		The same holds for $\PQ$ instead of $\PE$.
	\end{compactitem}
	\label{prop11}
\end{lem}
\begin{proof}
	We first prove Assertion (i) for $\PE$. 
	We consider separately the case (a) where $t\in [-2^{m+1},-2^{m}]\cap\cJ$ satisfies $L_t\geq -2^{0.45m}$ and the case (b) where $t\in [-2^{m+1},-2^{m}]\cap\cJ$ satisfies $L_t<-2^{0.45m}$, starting with case (a). First we argue that 
	\eqb
	\begin{split}
		&\Big\{\sup_{s_1,s_2\in[-2^{m+1},2^{m+1}]} |(L_{s_1}-\cL_{s_1}) - (L_{s_2} - \cL_{s_2})|<Cm \Big\}
		\cap \{N_{\op{L}}[-2^{m+1},2^{m+1}]<Cm \}
		\\&
		\cap\{\inf_{t\in[-2^m,-2^{0.9m}]} L_t<-2^{0.45m}-2dCm-2d\} 
		\cap \{ \cJLm \cap (-2^{0.9m},-2^{0.8m})\neq\emptyset \}
		\\&
		\qquad
		\subset\{\not\exists t\in[-2^{m+1},-2^m]\cap\cJ\,:\, L_t \geq  -2^{0.45m} \text{\,\,and\,\,} \cT^{\op{C},d}_{t,\op{L}}\cap[-2^{0.8m}, 2^{0.8m}]\neq\emptyset\}. 
	\end{split}
	\label{eq16}
	\eqe
	Suppose that the events on the left side of \eqref{eq16} occur. Suppose $t\in[-2^{m+1},-2^{m}]\cap\cJ$ satisfies $L_t\geq -2^{0.45m}$. Choose $t'\in[-2^m,-2^{0.9m}]$ such that $L_{t'}<-2^{0.45m}-2dCm-2d$ and let $t''\in \cJLm \cap (-2^{0.9m},-2^{0.8m})$; note that we can find appropriate $t'$ and $t''$ due to occurrence of the third and fourth events on the left side of \eqref{eq16}. We see that $\tau^{\op{C},d}_{t,\op{L}}\leq t''<-2^{0.8m}$ since 
	\eqbn
	\begin{split}
		\inf_{s\in [t,t'']} (L_s-L_{t})
		&\leq L_{t'}-L_t 
		< (-2^{0.45m}  -2dCm-2d) - (-2^{0.45m})\\
		&<  
		-d\sup_{u\in [-2^{m+1},2^{m+1}]} |(L_u-\cL_u)-(L_{t}-\cL_{t})|
		- dN_{\op{L}}[-2^{m+1},2^{m+1}]
		-2d\\
		&=-d\Delta_{\op{L}}[-2^{m+1},2^{m+1}] 
		\leq 
		-d\Delta_{\op{L}}[t,t''].
	\end{split}
	\eqen
	Since $\tau^{\op{C},d}_{t,\op{L}}< 2^{0.8m}$ we get from \eqref{eq:TCd0L} that $\cT^{\op{C},d}_{t,\op{L}}\cap[-2^{0.8m},2^{0.8m}]=\emptyset$, which implies \eqref{eq16}.

	For $C>1$ sufficiently large, the four events on the left side of \eqref{eq16} each have probability above $1-2^{-cm}$ by 
	Lemma \ref{prop:eucl-error}, 
	Lemma \ref{prop22}, 
	Lemma \ref{prop:eucl-inf-lower2} and Lemma \ref{prop:JLm-empty}, respectively. By using this and \eqref{eq16} we get
	\eqb
	\PE[\exists t\in[-2^{m+1},-2^m]\,:\, L_t \geq  -2^{0.45m} \text{\,\,and\,\,} \cT^{\op{C},d}_{t,\op{L}}\cap[-2^{0.8m}, 2^{0.8m}]\neq\emptyset] \leq 2^{-cm}.
	\label{eq38}
	\eqe
	Next we will obtain the counterpart of \eqref{eq38} for case (b), namely the case where $L_t<-2^{0.45m}$. We first argue that 
	\eqb
	\begin{split}
		&\Big\{\sup_{s_1,s_2\in[-2^{m+1},2^{2m}]} |L_{s_1}-\cL_{s_1} - (L_{s_2} - \cL_{s_2})|<Cm\Big\}
		\cap 
		\{ N_{\op{L}}[-2^{m+1},2^{2m}]<Cm \}
		\\&
		\cap
		\Big\{\inf_{t\in[-2^{0.85m},2^{0.8m}]} L_t>-2^{0.45m}+2dCm+2d\Big\} 
		\cap\Big\{ \max \{ \tau^{\op{C},d}_{t,\op{L}}\,:\,t\in [-2^{m+1},-2^m] \} \leq 2^{2m} \Big\}\\
		\\&
		\cap \{ \cJLm \cap (-2^{0.85m},-2^{0.8m})\neq\emptyset \}
		\subset \big\{\not\exists t\in[-2^{m+1},-2^m]\,:\, L_t <  -2^{0.45m}, %
		\cT^{\op{C},d}_{t,\op{L}}\cap[-2^{0.8m}, 2^{0.8m}]\neq\emptyset\big\}. 
	\end{split}
	\label{eq37}
	\eqe
	Suppose the events on the left side of \eqref{eq37} occur and let $t\in[-2^{m+1},-2^m]$ be such that $L_t <  -2^{0.45m}$. In order to establish \eqref{eq37} we need to argue that $\cT^{\op{C},d}_{t,\op{L}}\cap[-2^{0.8m}, 2^{0.8m}]=\emptyset$. By occurrence of the first and fourth events on the left side of \eqref{eq37} we have 
	$\Delta^{\op{C},d}_{t,\op{L}} \leq 2Cdm+2d$. For any $s\in \cT^{\op{C},d}_{t,\op{L}} \cap[-2^{0.8m},2^{0.8m}]$ we have $I_s\subset (-2^{0.85m},2^{0.8m}]$ by occurrence of the fifth event on the left side of \eqref{eq37}. Combining this with $L_t<-2^{0.45m}$, $\Delta^{\op{C},d}_{t,\op{L}} \leq 2Cdm+2d$, the definition \eqref{eq:TCd0L} of $\cT^{\op{C},d}_{t,\op{L}}$ and occurrence of the third event on the left side of \eqref{eq37} we get that $s\not\in \cT^{\op{C},d}_{t,\op{L}}$. Since $s\in \cT^{\op{C},d}_{t,\op{L}} \cap[-2^{0.8m},2^{0.8m}]$ was arbitrary we can conclude that $\cT^{\op{C},d}_{t,\op{L}}\cap[-2^{0.8m}, 2^{0.8m}]=\emptyset$ and \eqref{eq37} follows.
	
	The five events on the left side of \eqref{eq37} all have $\PE$-probability at least $1-2^{-cm}$ by Lemma \ref{prop:eucl-error}, Lemma \ref{prop22}, Lemma \ref{prop:eucl-inf-upper}, Lemma \ref{prop:tau-k-max}, and Lemma \ref{prop:JLm-empty}, respectively. It follows from \eqref{eq37} that
	$$
	\PE\Big[\exists t\in[-2^{m+1},-2^m]\,:\, L_t <  -2^{0.45m} \text{\,\,and\,\,} \cT^{\op{C},d}_{t,\op{L}}\cap[-2^{0.8m}, 2^{0.8m}]\neq\emptyset\Big] \leq 2^{-cm}.
	$$
	Combining this with \eqref{eq38} we obtain Assertion (i) for $\PE$.
		
	The proof of Assertion (ii) for $\PE$ is done similarly, except that it is easier since we do not need to consider the first two events on the left side of \eqref{eq16} and \eqref{eq37}. 

	The proof of Assertions (i) and (ii) for $\PQ$ is identical, except that we apply \eqref{eq:err4} and Remark \ref{rmk:Q-inf} instead of some of the results cited above.
\end{proof}

\begin{proof}[Proof of Proposition \ref{prop:Text-len-tail}]
	We first prove Assertion (i). Throughout the proof $c_1>0$ denotes a constant satisfying the conditions of the lemma which may decrease throughout the proof. Setting $t_1=-2^{0.8k}$, $t_2=-2^{0.5k}$, $t_3=2^{0.5k}$, $t_4=2^{0.9k}$, we have 
	\eqb\begin{split}
		\PE\Big[\bigcup_{t\in J} \cT^{\op{C},d}_t \not\subset I(k)\Big] \leq&\, \PE[ [t_1,t_4]\not\subset I(k)] + \PE\Big[\bigcup_{t\in J} \cT^{\op{C},d}_t \not\subset [t_1,t_4]\Big] \\
		\leq 
		&\,\PE[ [t_1,t_4]\not\subset I(k)] 
		+  \PE[J\not\subset [t_2,t_3]]
		+ \PE[\exists t<t_1 \,:\, \cT^{\op{C},d}_{t,\op{L}}\cap[t_2,t_3]\neq\emptyset]\\
		&+\PE[\exists t<t_1 \,:\, \cT^{\op{C},d}_{t,\op{R}}\cap[t_2,t_3]\neq\emptyset]+ \PE[\exists t\in [t_1,t_3]\,:\,\cT^{\op{C},d}_{t,\op{L}}\not\subset[t_1,t_4]]\\
		&+ \PE[\exists t\in [t_1,t_3]\,:\,\cT^{\op{C},d}_{t,\op{R}}\not\subset[t_1,t_4]],
	\end{split}\label{eq15}\eqe
	where in the last step we use \eqref{eq-def-Tcell}. 
	The first term on the right side is bounded above by $2^{-c_1 k}$ by the definition of a uniform dyadic decomposition. 
	The second term on the right side is bounded above by $2^{-c_1 k}$ by the assumption of the lemma. 
	We bound the third term on the right side of \eqref{eq15} by observing that
	\eqbn
	\begin{split}
		&\PE[\exists t<t_1 \,:\, \cT^{\op{C},d}_{t,\op{L}}\cap[t_2,t_3]\neq\emptyset]\\
		&\qquad\leq \sum_{m=0.8k}^{\infty} \PE[\exists t\in [-2^{m+1},-2^m]\cap\cJ\,:\,\cT^{\op{C},d}_{t,\op{L}}\cap[t_2,t_3]\neq\emptyset] 
		\leq \sum_{m=0.8k}^{\infty} 2^{-c_1 m} \leq 2^{-c_1 k},
	\end{split}
	\eqen
	where we applied Lemma \ref{prop11} in the second step. 
	In order to bound the fifth term on the right side of \eqref{eq15}, observe that 
	\eqb
	\begin{split}
		&\PE[\exists t\in [t_1,t_3]\,:\,\cT^{\op{C},d}_{t,\op{L}}\not\subset[t_1,t_4]]
		\leq \PE[\exists t\in [t_1,t_3]\,:\,\tau^{\op{C},d}_{t,\op{L}} > t_4] \leq 2^{-c_1 k},
	\end{split}
	\label{eq20}
	\eqe
	where we use Lemma \ref{prop:tau-k-max} in the last step. The forth and sixth terms are bounded in a similar way as the third and fifth terms, respectively, by symmetry in L and R. 
	This completes the proof of Assertion (i). Assertion (ii) is proven exactly as Assertion (i) since all the applied lemmas also hold for $\PQ$. 
	
	To prove Assertion (iii), first observe that we may assume $h<2^{0.45 k}$ since for larger $h$ the bound we want to derive is trivially satisfied as the right-hand side $h^{c_1^{-1}}2^{-c_1 k}$ is larger than 1 for sufficiently small $c_1$. For $h<2^{0.45 k}$ we can proceed exactly as in the proof of Assertions (i) and (ii) with $\cT^{\op{C},d}_{t},\cT^{\op{C},d}_{t,\op{L}},\cT^{\op{C},d}_{t,\op{R}}$ replaced by $\cT^{h,\infty}_{t},\cT^{h,\infty}_{t,\op{L}},\cT^{h,\infty}_{t,\op{R}}$, and we get 
	\eqbn
	\begin{split}
		\PE\Big[\bigcup_{t\in J_h} \cT^{h,\infty}_t \not\subset I(k) \Big] \leq 2^{-c_1 k}.
	\end{split}
	\eqen		
	The right side is smaller than $h^{c_1^{-1}}2^{-c_1 k}$ since $h>1$, which concludes the proof of Assertion (iii). The proof of Assertion (iv) is identical.

\end{proof}

An immediate corollary of Proposition \ref{prop:Text-len-tail} is the following.
\begin{cor}
	The sets $\cT_{t_0}^{\op{C},d},\cT_{t_0}^{\op{D},d},\cT_{t_0}^{h,\infty}$ are finite a.s.\ under the probability measures $\PQ$ and $\PE$.
	\label{prop:time-sets-finite}
\end{cor}
\begin{proof}
	The result for $\cT_{t_0}^{\op{C},d}$ and $\cT_{t_0}^{h,d}$ follows by combining Lemma \ref{prop-finite-jumps} (or the corresponding result for $\PQ$, which holds since the jumps of $\cZ$ have the law of a Poisson point process under $\PQ$) and Proposition \ref{prop:Text-len-tail}. The result for $\cT_{t_0}^{\op{D},d}$ is obtained by using that $\cT_{t_0}^{\op{D},d}\subset \cT_{t_0}^{\op{C},d+1}$ according to Proposition \ref{prop:time-sets}.
\end{proof}

\subsection{The cell configurations}\label{sec:cells3}
\begin{figure}
	\centering
	\includegraphics{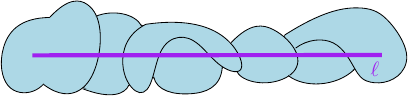}
	\caption{Illustration of $\cT^{\eta}_\ell$. The compact sets on blue background are the sets $\ol\eta(t)$ for $t\in \cT^{\eta}_\ell$.}
	\label{fig:Tell}
\end{figure}

Recall that $\cZ$ encodes a planar map $M$ with a spanning tree via the Mullin bijection and that there is a natural bijection $t\mapsto \frk v_t$ between the vertex set $\cV(M)$ of $M$ and the time set $\cJLm\subset\R$ for $\cZ$. Let $\EM \subset \{ \{s,u\}\,:\,s,u\in\cJLm \}$ be such that $\{s,u \}\in\EM$ if and only if the vertices $\frak v_s$ and $\frak v_u$ are adjacent in $M$ and $s\neq u$. Then $(\cJLm,\EM)$ is a graph isomorphic to $M$.\footnote{
	Here we view $M$ as a graph (rather than a multigraph) by deleting self-loops and collapsing all the multiple edges of $M$, i.e., for any distinct $v,w\in V(M)$ all the edges between $v$ and $w$ (if any) are collapsed to a single edge. 
} 
Notice that the vertex set of this graph is a subset of the vertex set of the graphs $(\cJm,\ET)$ and $(\cJm,\EG)$ considered in the previous subsection. We call a path in the graph $(\cJLm,\EM)$ (resp., $(\cJm,\ET)$, $(\cJm,\EG)$) an M-path (resp., T-path, G-path), i.e., an M-path is an ordered sequence $t_1,\dots,t_k\in\cJLm$ for some $k\in\N$ such that $\{t_j,t_{j+1} \}\in\EM$ for $j=1,\dots,k-1$. 
For $s\in\cJLm$ define
\eqb
\cT^{\op{cell}}_s := \bigcup_{u\in\cJm\,:\,\{u,s \}\in\ET\text{ or }u=s} \cT^{\op{C},1}_u,\qquad
\cT^{\op{large}}_s := \bigcup_{u\in\cJm\,:\,\{u,s \}\in\ET\text{ or }u=s} \cT^{\op{C},3}_u.
\label{eq:wtT} 
\eqe
The main motivation of this definition is that Proposition \ref{prop:conditions} right below holds, which is the main result of this section. We refer to \eqref{eq:ol-eta} for the definition of $\ol\eta$ and set
\eqb
\cT^{\op{cell,inv}}_s := \{u\in\cJLm\,:\, s\in \cT^{\op{cell}}_u \},\qquad
\cT^{\op{large,inv}}_s := \{u\in\cJLm\,:\, s\in \cT^{\op{large}}_u \}.
\label{eq74}
\eqe
\begin{prop}
	There is a cell configuration $(\cH^{\op{cell}},\sim_{\op{cell}},\frk c_{\op{cell}})$ with $\cH^{\op{cell}}=(H^{\op{cell}}_t)_{t\in\cJLm}$ such that for each $t\in\cJLm$, 
	\eqb
	H^{\op{cell}}_t \subset \bigcup_{s\in \cT^{\op{cell}}_t} \ol\eta(s),\qquad\text{and}\qquad
	\op{area}(H^{\op{cell}}_t) = \sum_{s\in \cT^{\op{cell}}_t } \frac{\op{area}(\ol\eta(s))}{2\cdot \# \cT^{\op{cell,inv}}_s }.
	\label{eq:cell-conf}
	\eqe
	Furthermore, $H^{\op{cell}}_t \sim_{\op{cell}} H^{\op{cell}}_s$ if and only if $\frk v_t$ and $\frk v_s$ are adjacent in $M$ and $s\neq t$, and in this case $\frk c_{\op{cell}}(H^{\op{cell}}_t,H^{\op{cell}}_s)=1$.  %
	There is also a cell configuration $(\cH^{\op{large}},\sim_{\op{large}},\frk c_{\op{large}})$ satisfying the exact same conditions with 
	``$\mathrm{cell}$'' replaced by ``$\mathrm{large}$''. 
	Finally, the cell configurations 
	$(\cH^{\op{cell}},\sim_{\op{cell}},\frk c_{\op{cell}})$ and 
	$(\cH^{\op{large}},\sim_{\op{large}},\frk c_{\op{large}})$
	satisfy all conditions of Theorem \ref{thm-general-clt-uniform}.
	\label{prop:conditions}
\end{prop}
The proposition is an immediate consequence of Lemmas \ref{prop25}, \ref{lem:transl-inv}, \ref{prop:erg}, and \ref{prop:i-iii} and Proposition \ref{prop:moment} below. For a line segment $\ell$ in $\BB C$ define 
\eqb
\begin{split}
	\cT_\ell^{\eta} &:= \{ t\in\cJm\,:\,\ol\eta(t)\cap \ell\neq\emptyset \},\\
	\cT_\ell^{\op{C},d} &:= \{ t\in\cJm\,:\,\ \cT^{\op{C},d}_t \cap \cT_\ell^{\eta} \neq\emptyset \},\qquad\qquad
	\cT_\ell^{\op{D},d} := \{ t\in\cJm\,:\,\ \cT^{\op{D},d}_t \cap \cT_\ell^{\eta} \neq\emptyset \},\\
	\cT_\ell^{\op{cell}} &:= \{ t\in\cJLm\,:\,\ \cT^{\op{cell}}_t \cap \cT_\ell^{\eta} \neq\emptyset \}
	=(\cT^{\op{C},1}_\ell \cap \cJLm )\cup \{t\in\cJLm \,:\,\exists s\in \cT^{\op{C},1}_\ell\text{\,such\,that\,}\{s,t\}\in\ET \},\\
	\cT_\ell^{\op{large}}&:= \{ t\in\cJLm\,:\,\ \cT^{\op{large}}_t \cap \cT_\ell^{\eta} \neq\emptyset \}
	=(\cT^{\op{C},3}_\ell \cap \cJLm)\cup \{t\in\cJLm \,:\,\exists s\in \cT^{\op{C},3}_\ell\text{\,such\,that\,}\{s,t\}\in\ET \}.
\end{split}
\label{eq:85}
\eqe
See Figure \ref{fig:Tell} for an illustration of $\cT^{\eta}_\ell$.

The following lemma will be used to check that the tuples $(\cH^{\op{cell}},\sim_{\op{cell}},\frk c_{\op{cell}})$ and $(\cH^{\op{large}},\sim_{\op{large}},\frk c_{\op{large}})$ defined below are cell configurations in the sense that if two cells are connected by an edge then their intersection (viewed as a subset of the complex plane) is non-empty; see \eqref{eq:cell-conf-adjacency}. 
\begin{lem}
	If $\{u,s \}\in\EM$ then $u\in \cT^{\op{cell}}_s$ and $s\in {\cT}^{\op{cell}}_u$.
	\label{prop:em-adj}
\end{lem}
\begin{proof}
	By symmetry in $u$ and $s$ it is sufficient to prove the first assertion, namely that $u\in \cT^{\op{cell}}_s$.
	If $\{u,s \}\in\EM$ then there is a face $f\in \cF(M)$ such that both $\frak v_u$ and $\frak v_s$ are incident to $f$. We have $\{u,\frk t_f \},\{s,\frk t_f \}\in\ET$, where we identify $f$ with the associated vertex of the triangulation $T$ so that $\frk t_f$ is well-defined. Applying Assertion 3.\ of Proposition \ref{prop:time-sets} and using $\{u,\frk t_f \}\in\ET$ we have $u\in \cT^{\op{C},1}_{\frk t_f}$. Using this and $\{s,\frk t_f \}\in\ET$ along with \eqref{eq:wtT} we get $u\in \cT^{\op{cell}}_s$. 
\end{proof}

The following lemma will be used to verify that our cell configurations satisfy Assumption (iii) of Theorem \ref{thm-general-clt-uniform}. For a set  $\cT\subset\cJLm$ we say that an M-path is in $\cT$ if all vertices on the path are contained in $\cT$.
\begin{lem}
	Let $\ell$ be an arbitrary line segment and let $s,\wt s\in \cT^{\op{cell}}_\ell$. Then we can find an M-path from $s$ to $\wt s$ in $\cT^{\op{large}}_\ell$.
	\label{prop:lines}
\end{lem}
\begin{figure}
	\centering
	\includegraphics{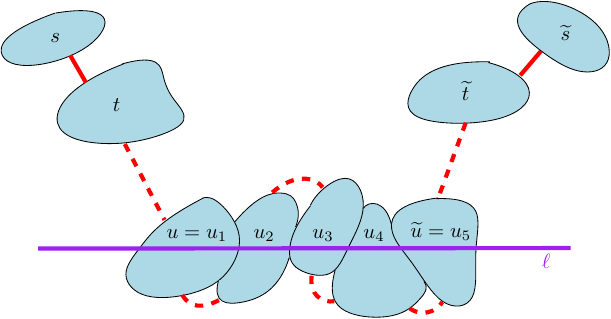}
	\caption{Illustration of the proof of Lemma \ref{prop:lines0}. The compact sets with light blue background and numbers $s,t,\dots$ are schematic drawings of the sets  $\ol{\eta}(s),\ol{\eta}(t),\dots$. The two red solid lines represent the edges $\{s,t \},\{\wt s,\wt t \}\in\ET$. The red dashed lines represent T-paths in $\cT^{\op{D},2}_\ell$.}
	\label{fig-line-conn}
\end{figure}
The following lemma is used in the proof of Lemma \ref{prop:lines}.
\begin{lem}
	In the context of Lemma \ref{prop:lines} we can find a T-path from $s$ to $\wt s$ such that all vertices on the path except possibly $s$ and $\wt s$ are in $\cT^{\op{D},2}_\ell$. 
	\label{prop:lines0}
\end{lem}
\begin{proof}
	The proof is illustrated in Figure \ref{fig-line-conn}. 
	By the definition of $\cT^{\op{cell}}_\ell$ 
	there are $t\in\cJm$ and $u\in\cJm$ such that $u\in\cT^{\eta}_\ell$, $t\in\cT^{\op{C},1}_u$, and either $\{s,t \}\in\ET$ or $t=s$.  Similarly, we can find a $\wt t\in\cJm$ and $\wt u\in\cJm$ such that the same holds with $\wt t,\wt u,\wt s$ instead of $t,u,s$.
	
	Next we argue that the set $\cT^{\eta}_\ell$ is G-connected. By \cite[Proposition 3.1]{gms-tutte}, we know that $\cT^{\eta}_\ell\cap \cJLm$ (resp.\ $\cT^{\eta}_\ell\cap \cJRm$) is G-connected 
	since by the cited proposition, for all any $t,t'\in\cT^{\eta}_\ell\cap \cJLm$ (resp.\ $t,t'\in\cT^{\eta}_\ell\cap \cJRm$) we can find a path of vertices in $\cT^{\eta}_\ell\cap\cJLm$ (resp.\ $\cT^{\eta}_\ell\cap\cJRm$) connecting $t$ and $t'$ such that \eqref{eq:line-under} (resp.\ the same condition with $R$) holds for any adjacent $s,u$ in the path. We get that $\cT^{\eta}_\ell$ is G-connected by choosing 
	$t\in\cT^{\eta}_\ell\cap \cJLm$ 
	and $t'\in\cT^{\eta}_\ell\cap \cJRm$ such that $I_t\cap I_{t'}\neq\emptyset$, so $t$ and $t'$ are G-adjacent.

	By the result of the previous paragraph we can find a G-path $u_1,\dots,u_k$ in $\cT^{\eta}_\ell$ such that $u_1=u$ and $u_k=\wt u$. In order to conclude the proof it is sufficient to prove the following three results since this allows us to construct a T-path from $t$ to $\wt t$ in $\cT^{\op{D},2}_\ell$ by concatenation: 
	(a) there is a T-path from $t$ to $u$ in $\cT^{\op{D},2}_\ell$,  
	(b) there is a T-path from $u_j$ to $u_{j+1}$ in $\cT^{\op{D},2}_\ell$ for $j=1,\dots,k-1$, and
	(c) there is a T-path from $\wt t$ to $\wt u$ in $\cT^{\op{D},2}_\ell$. We omit the proof of (c) since it is identical to the proof of (a).
	
	We have $t\in\cT^{\op{C},1}_u\subset \cT^{\op{D},2}_{u}$ by Proposition \ref{prop:time-sets}, Assertion 1. By using that  $\cT^{\op{D},2}_u$ is T-connected by Proposition \ref{prop:time-sets}, Assertion 5, we get that there is a T-path from $t$ to $u$ in $\cT^{\op{D},2}_u$. This implies (a) since $\cT^{\op{D},2}_u\subset \cT^{\op{D},2}_\ell$ by Proposition \ref{prop:time-sets}, Assertion 2. 
	In order to get (b), first observe that since $\{u_{j}, u_{j+1}\}\in\EG$ we have $u_{j+1}\in \cT^{\op{D},2}_{u_j}$ by Proposition \ref{prop:time-sets}, Assertion 4. By using that $\cT^{\op{D},2}_{u_j}\subset \cT^{\op{D},2}_\ell$ is T-connected we get (b).
\end{proof}

\begin{figure}
	\centering
	\includegraphics[scale=1]{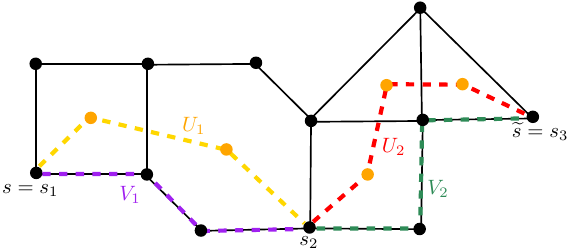}
	\caption{Illustration of the proof of Lemma \ref{prop:lines0}. The planar map drawn in black is a subset of the planar map $M$ encoded by $\cZ$ and the orange dots represent dual vertices. Recall that the primal (resp.\ dual) vertices are in bijective correspondence with $\cJLm$ (resp.\ $\cJRm$). The (four) dashed lines represent the paths $U_1,\dots,U_{k-1},V_1,\dots,V_{k-1}$ considered in the proof.}
	\label{fig-Tpath-to-Mpath}
\end{figure}

\begin{proof}[Proof of Lemma \ref{prop:lines}]
	The proof is illustrated in Figure \ref{fig-Tpath-to-Mpath}. 
	Consider the T-path of Lemma \ref{prop:lines0}, which consists of vertices from $\cJLm$ and $\cJRm$. We can write the path as $s_1,U_1,s_2,U_2,\dots,s_k$, where $s_j\in\cJLm$ and $U_j$ is a $T$-path from $s_j$ to $s_{j+1}$ consisting only of elements of $\cJRm$ except at its end-points. We will modify this path by replacing each $U_j$ by an M-path $V_j$ from $s_j$ to $s_{j+1}$ such that each vertex $t$ on the M-path $V_j$ is T-adjacent to some $t'$ on $U_j$ (i.e., $\{t,t' \}\in\ET$). This is sufficient to conclude the proof of the proposition since $\cT^{\op{D},2}_\ell\subset \cT^{\op{C},3}_\ell$ (Proposition \ref{prop:time-sets}, Assertion 1), since each path $U_j$ is in $\cT^{\op{D},2}_\ell$ (Lemma \ref{prop:lines0}) and by the definition of $\cT^{\op{large}}_\ell$ in \eqref{eq85}. The existence of appropriate paths $V_j$ follow since the subset of $\cJLm$ consisting of vertices that are T-connected to a vertex in $U_j$ is M-connected.

\end{proof}

For concreteness we assume in the remainder of this section that $h$ has the circle average embedding. However, the results we present would also have been true for an arbitrary other embedding chosen independently of $(h^{\op{IG}},\frk U)$. Recall the definition of $\ol\eta$ in \eqref{eq:ol-eta}. For $A\subset\C$ let $\partial A$ denote the boundary of $A$. 
\begin{lem}
	There exists a collection of connected closed sets with disjoint interior $\mcl A=\{A_{u,s}\,:\,s\in\cJLm,u\in\cT^{\op{cell}}_s \}$ such that for each $s\in\cJLm$ and $u\in\cT^{\op{cell}}_s$ we have   $\partial(\ol\eta(u))\subset A_{u,s}\subset\ol\eta(u)$ and
	\eqb
	2\cdot
	\#{\cT}_u^{\op{cell,inv}} 
	\cdot \op{area}(A_{u,s}) = \op{area}(\ol\eta(u)).
	\label{eq:split-cores}
	\eqe
	Furthermore, we may assume that for a measurable function $g$ we have $\mcl A=g(h,h^{\op{IG}},\frk U)$ and a.s.,
	\eqb
	C^{-1}(\mcl A-z)=
	g(h(C\cdot+z)+Q\log C,h^{\op{IG}}(C\cdot+z),\frk U-\tau_z ),\qquad\forall C>0,z\in\BB C.
	\label{eq:transl2}
	\eqe
	and that for measurable functions $g_n$, $n\in\N$, and any fixed $R>0$,
	\eqb
	\lim_{n\rta\infty}\P[\mcl A|_{B(0,R)} = g_n(h|_{B(0,n)}, h^{\op{IG}}|_{B(0,n)}, \frk U)|_{B(0,R)}]=1,
	\label{eq:local}
	\eqe
	where for any $\wt{\mcl A}$ given by a collection of subsets of $\C$, the restriction $\mcl A|_{B(0,R)}$ is defined by $\mcl A|_{B(0,R)}=\{ A\in\wt{\mcl A}\,:\,A\cap B(0,R)\neq\emptyset \}$.
	\label{prop:Uus}
\end{lem}
\begin{figure}
	\centering
	\includegraphics{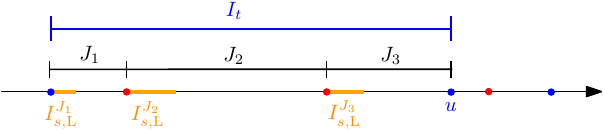}
	\caption{Illustration of the construction of $A_{u,s}$ given in the proof of Lemma \ref{prop:Uus}. The blue (resp.\ red) dots represent elements of $\cJLm$ (resp.\ $\cJRm$) so in particular we have $u\in\cJLm$. 
	}
	\label{fig:V}
\end{figure}
\begin{proof} 
	We present a possible explicit construction of the sets $A_{u,s}$. See Figure \ref{fig:V}. 
	Consider a connected component $J$ of $\R\setminus \cJm$ and let $u_{\op{L}}\in \cJLm$ (resp.\ $u_{\op{R}}\in \cJRm$) be such that $J\subset I_{u_{\op{L}}}$ (resp.\ $J\subset I_{u_{\op{R}}}$). Write $\eta(J)$ as the union of sets of the form 
	$\eta(I^J_{s,\op{L}})$ for $s\in \cT_{u_{\op{L}}}^{\op{cell,inv}}$ 
	and 
	$\eta(I^J_{s,\op{R}})$ for $s\in \cT_{u_{\op{R}}}^{\op{cell,inv}}$, where $I^J_{s,\op{L}},I^J_{s,\op{R}}\subset J$ are closed intervals. 
	Suppose the set of intervals 
	$\{I^J_{s,\op{L}} \,:\, s\in \cT_{u_{\op{R}}}^{\op{cell,inv}} \}
	\cup 
	\{I^J_{s,\op{R}} \,:\, s\in \cT_{u_{\op{R}}}^{\op{cell,inv}} \}$ have pairwise disjoint interior and union $\ol J$ such that 
	$2\cdot\# \cT_{u_{\op{L}}}^{\op{cell,inv}} \cdot \op{area}(\eta(I^J_{s,\op{L}})) = \op{area}(\eta(J))$ for each $s\in \cT_{u_{\op{L}}}^{\op{cell,inv}}$ and
	$2\cdot\# \cT_{u_{\op{R}}}^{\op{cell,inv}} \cdot \op{area}(\eta(I^J_{s,\op{R}})) =\op{area}(\eta(J))$ for each $s\in \cT_{u_{\op{R}}}^{\op{cell,inv}}$; such a collection of intervals can for example be uniquely determined by requiring that $I^J_{s_1,\op{R}}$ is to the left of $I^J_{s_2,\op{R}}$ on the real line if $s_1<s_2$.

	For $s\in\cJLm$ and $u\in \cT_{s}^{\op{cell}}\cap\cJLm$ let $J_1,\dots,J_k$ denote the connected components of $I_{u}\setminus \cJm$ and set 
	$$
	A_{u,s}=\bigcup_{j=1}^k \Big(\eta\Big(I^{J_j}_{s,\op{L}}\Big) 
	\cup 
	\bigcup_r \partial\Big(\eta\Big(I^{J_j}_{r,\op{L}}\Big)\Big)
	\cup
	\bigcup_r \partial\Big(\eta\Big(I^{J_j}_{r,\op{R}}\Big)\Big)
	\Big),
	$$ 
	where $r\in\cJLm$ ranges over all possible values such that $I^{J_j}_{r,\op{L}}$ (resp.\ $I^{J_j}_{r,\op{R}}$) is defined. The reason we are including the two unions over $r$ is to ensure that $A_{u,s}$ is connected and that $\partial\ol\eta(u)\subset A_{u,s}$.
	Similarly, if $s\in\cJLm$ and $u\in \cT_{s}^{\op{cell}}\cap\cJRm$ let $J_1,\dots,J_k$ denote the connected components of $I_{u}\setminus \cJm$ and define $A_{u,s}$ by the same formula except that we consider  $I^{J_j}_{s,\op{R}}$ instead of $I^{J_j}_{s,\op{L}}$. It is immediate from the construction that the desired properties as stated in the lemma are satisfied. 
\end{proof}
For $s\in\cJLm$ define
$$
H^{\op{cell}}_s := \bigcup_{u\in{\cT}_s^{\op{cell}}} A_{u,s}.
$$
Define $\sim_{\op{cell}}$ such that $H_t\sim_{\op{cell}} H_{t'}$ if and only if $\{t,t' \}\in\EM$. Set $\frk c(H_t,H_{t'})=1$ for all such $t$ and $t'$.

Define
$(\cH^{\op{large}},\sim_{\op{large}},\frk c_{\op{large}})$ in exactly the same way as $(\cH^{\op{cell}},\sim_{\op{cell}},\frk c_{\op{cell}})$, except that we use ${\cT}^{\op{large}}_{t}$ 
instead of ${\cT}^{\op{cell}}_{t}$. 

\begin{lem}
	The tuples $(\cH^{\op{cell}},\sim_{\op{cell}},\frk c_{\op{cell}})$ and $(\cH^{\op{large}},\sim_{\op{large}},\frk c_{\op{large}})$ are cell configurations and satisfy \eqref{eq:cell-conf}. Furthermore, for each $t\in\cJLm$ the set $\cT_t^{\op{cell}}$ is G-connected.
	\label{prop25}
\end{lem}
\begin{proof}
	For each $t\in\cJLm$ the set $\cT_t^{\op{cell}}$ is G-connected since the set $\cT^{\op{C},1}_u$ is G-connected for any $u\in\cJm$ (Proposition \ref{prop:time-sets}, Assertion 5) and since each set $\cT^{\op{C},1}_u$ in the union defining $\cT^{\op{cell}}_s$ contains $s$ (Proposition \ref{prop:time-sets}, Assertion 3). This implies that each set $H^{\op{cell}}_t$ is connected as a subset of $\C$. 
	To prove that $(\cH^{\op{cell}},\sim_{\op{cell}},\frk c_{\op{cell}})$ defines a cell configuration
	we also need to argue that if $H_t\sim_{\op{cell}} H_{t'}$ (equivalently, $\{t,t' \}\in\cE_{\op{M}}$) then $H_t\cap H_{t'}\neq\emptyset$. This follows from Lemma \ref{prop:em-adj}. It is immediate from the construction that \eqref{eq:cell-conf} is satisfied.  The proof for $(\cH^{\op{large}},\sim_{\op{large}},\frk c_{\op{large}})$ 
	is identical. 
\end{proof}

Our next goal is to prove Lemma \ref{lem:transl-inv}, which says that our cell configurations is translation invariant modulo scaling. We start by proving that $\mcl H^{\op{cell}}$ is defined from $(h,h^{\op{IG}},\frk U)$ in a translation and scale invariant way which is also local (Lemma \ref{prop26}). Then we describe a sense in which $(h,h^{\op{IG}},\frk U)$ is translation invariant modulo scaling (Lemma \ref{prop30}). Combining these two results we get Lemma \ref{lem:transl-inv}.

The next lemma is an immediate consequence of Lemma  \eqref{prop:Uus} and the proof is therefore omitted.
\begin{lem}
	There is a measurable function $g$ such that a.s.,
	\eqb
	C^{-1}(\mcl H^{\op{cell}}-z)=
	g(h(C\cdot+z)+Q\log C,h^{\op{IG}}(C\cdot+z),\frk U-\tau_z )\qquad\forall z\in\C,C>0.
	\label{eq:transl4}
	\eqe
	There are also measurable functions $g_n$ taking their value in the set of cell configurations such that for any fixed $R>1$,
	\eqb
	\lim_{n\rta\infty}\P[\mcl H^{\op{cell}}|_{B(0,R)} = g_n(h|_{B(0,n)}, h^{\op{IG}}|_{B(0,n)}, \frk U)|_{B(0,R)}]=1,
	\label{eq:almost-local}
	\eqe
	where the restriction of a cell configuration to $B(0,R)$ means the collection of cells intersecting $B(0,R)$, along with information about adjacency and conductances for these cells.
	The above two results also hold for $(\cH^{\op{large}},\sim_{\op{large}},\frk c_{\op{large}})$ instead of $(\cH^{\op{cell}},\sim_{\op{cell}},\frk c_{\op{cell}})$.
	\label{prop26}
\end{lem}

Next we will describe in Lemma \ref{prop30} a sense in which $(h,h^{\op{IG}},\frk U)$ is translation invariant modulo scaling. In order to describe this notion of translation invariance we need to introduce the concept of uniform dyadic systems on $\C$, which were earlier considered in  \cite{gms-random-walk}. 
A dyadic system on $\C$ is a two-dimensional counterpart of the dyadic systems defined in Section \ref{sec:kmt1}. If $S\subset\C$ is a square, a dyadic child of $S$ is one of the four squares $S'\subset S$ with half the side length of $S$ whose corners include one corner of $S$ and the center of $S$. A dyadic parent of $S$ is one of the four squares with $S$ as a dyadic child. A square $S'$ is a dyadic ancestor (resp.\ descendant) of $S$ if we can find squares $S_1,\dots,S_k$ such that $S=S_1$, $S'=S_k$ and $S_{j+1}$ is a dyadic parent (resp.\ child) of $S_j$ for $j=1,\dots,k-1$. A dyadic system on $\C$ is a collection $\frk D^{\C}$ of closed squares with the following properties:
\begin{enumerate}
	\item If $S\in\frk D^{\C}$ then all the four dyadic children of $S$ are also in $\frk D^{\C}$,
	\item if $S\in\frk D^{\C}$ then exactly one of the dyadic parents of $S$ is also in $\frk D^{\C}$, and
	\item any two squares in $\frk D^{\C}$ have a common dyadic ancestor. 
\end{enumerate}
A uniform dyadic system is a random dyadic system defined as follows. Let $s$ be sampled uniformly from $[0,1]$ and conditional on $s$ let $w$ be sampled uniformly from  $[0,2^s]\times[0,2^s]$. Set $S_0=[0,2^s]\times[0,2^s]-w$. For $k\in\N$ inductively let $S_k$ be sampled uniformly from the four dyadic parents of $S_{k-1}$. Let $\frk D^{\C}$ be the set of all dyadic descendants of $S_k$ for $k\in\N$.
\begin{lem}
	Suppose $(h,h^{\op{IG}},\frk U)$ is sampled from $\PE$ and let $\frk D^{\C}$ be an independent uniform dyadic system on $\C$. Let $S\in\frk D^{\C}$ be such that $0\in S$ and such that for a measurable function $g$ we have $S=g(h,h^{\op{IG}},\frk U,\frk D^{\C})$. Also suppose that for measurable functions $g_n$,
	\eqb
	\lim_{n\rta\infty}\PE[ S = g_n(h|_{B(0,n)},h^{\op{IG}}|_{B(0,n)},\frk U,\frk D^{\C}) ] = 1
	\label{eq51}
	\eqe 
	and that 
	\eqb
	C^{-1}(S-z) = g( h(C\cdot + z)+Q\log C,h^{\op{IG}}(C\cdot +z),\frk U-\tau_z,C^{-1}(\frk D^{\C}-z)),\qquad\forall z\in S, C>1.
	\label{eq52}
	\eqe

	Then if $z$ is a point sampled from Lebesgue measure in $S$,  for an appropriately chosen random constant $C>0$,
	\eqb
	(h(C\cdot+z)+Q\log C, h^{\op{IG}}(C\cdot+z), \frk U-\tau_z) \eqD (h,h^{\op{IG}},\frk U).
	\label{eq57}
	\eqe
	\label{prop30}
\end{lem}
\begin{proof}
	For concreteness we will assume that $h$ has the circle average embedding. 
	Let $R_n\rta\infty$ and let $z_n$ be sampled uniformly from Lebesgue measure on $B(0,R_n)$. By \eqref{eq50}, conformal invariance of the whole plane GFF, and the definition of $\frk D^{\C}$ as a uniform dyadic system on $\C$, there are random $C_n>0$ such that the following convergence holds in total variation distance for any fixed $R>0$

	$$
	\big(h(C_n\cdot+z_n)+Q\log C_n|_{B(0,R)}, h^{\op{IG}}(C_n\cdot+z_n)|_{B(0,R)}, \frk U-\tau_{z_n},  C_n^{-1}(\frk D^{\C}-z_n)\big) \rta (h|_{B(0,R)},h^{\op{IG}}|_{B(0,R)},\frk U,\frk D^{\C}).
	$$
	For $z\in\C$ and $C>0$ define $S(z,C)=g(h(C\cdot+z)+Q\log C, h^{\op{IG}}(C\cdot+z),\frk U-\tau_z,C^{-1}(\frk D^{\C}-z))$. By \eqref{eq51}, we can add the cell $S$ to the tuple above
	\eqb
	\begin{split}
		&(h(C_n\cdot+z_n)+Q\log C_n|_{B(0,R)}, h^{\op{IG}}(C_n\cdot+z_n)|_{B(0,R)}, \frk U-\tau_{z_n}, C_n^{-1}(\frk D^{\C}-z_n),S(z_n,C_n)) \\
		&\qquad\qquad\qquad\qquad\qquad\qquad\rta (h|_{B(0,R)},h^{\op{IG}}|_{B(0,R)},\frk U,\frk D^{\C},S).
	\end{split}
	\label{eq54}
	\eqe
	Let $w_n$ be sampled uniformly from Lebesgue measure on $S(z_n,C_n)$ and let $\wh w_n=z_n+C_nw_n$ be the associated point in the coordinates of $h$. By \eqref{eq54} and with $z$ as in \eqref{eq57}, for appropriate random constants $C_n,C$,

	\eqbn
	\begin{split}
		&(h(C_n\cdot+\wh w_n) + Q\log C_n |_{B(0,R)}, 
		h^{\op{IG}}(C_n\cdot+\wh w_n)|_{B(0,R)}, \frk U-\tau_{\wh w_n} 
		) \\
		&\qquad\qquad\qquad\qquad\rta
		((h(C\cdot+z)+Q\log C)|_{B(0,R)}, h^{\op{IG}}(C\cdot+z)|_{B(0,R)}, \frk U-\tau_z 
		).
	\end{split}
	\eqen
	Observe that $\P[z_n+C_nS(z_n,C_n)\subset B(0,R_n)]\rta 1$ as $n\rta\infty$. This implies that with probability $1-o_n(1)$, the total variation distance between $z_n$ and $\wh w_n$, conditioned on the fields and the dyadic system, goes to zero as $n\rta\infty$. This implies further that the total variation distance between the left sides of the last two displays goes to zero as $n\rta\infty$ (where, here and in the rest of this paragraph, we ignore the last two elements of the tuple in \eqref{eq54})
	and then that the total variation distance between the right sides of the last two displays 
	goes to zero. By using that the right sides of the last two displays do not depend on $n$ we get that these two right sides are equal in law, which concludes the proof since $R$ was arbitrary.

\end{proof}

\begin{lem}
	The cell configurations  $(\cH^{\op{cell}},\sim_{\op{cell}},\frk c_{\op{cell}})$ and $(\cH^{\op{large}},\sim_{\op{large}},\frk c_{\op{large}})$ satisfy Assumption \ref{item-resampling} (translation invariance modulo scaling). %
	\label{lem:transl-inv}
\end{lem}
\begin{proof}
	Let $\frk D^{\C}$ be an independent uniform dyadic system on $\C$. For $z\in\C$ and $m>0$ define
	$\wh S_m^z$ to be the largest square $S\in\frk D^{\C}$ such that $z\in S$ and 
	$$
	\sum_{H\in \cH^{\op{cell}}\,:\, H\cap S\neq\emptyset } \frac{\op{Area}(H\cap S)}{\op{Area}(H)} \leq m,
	$$
	i.e., $\wh S_m^z$ is the largest square containing $z$ which contains at most $m$ cells, where cells partly contained in the square are partially accounted for. 
	By Condition 5 of \cite[Lemma 2.3]{gms-random-walk}, $\cH^{\op{cell}}$ is translation invariant modulo scaling if and only if $(\cH^{\op{cell}}-z_0,\frk D^{\C}-z_0)$ is equal in law to $(\cH^{\op{cell}},\frk D^{\C})$ modulo scaling for $z_0$ sampled uniformly from Lebesgue measure on $\wh S_m^0$. 
	By \eqref{eq:transl4}, in order to conclude the proof of the lemma it is therefore sufficient to argue that \eqref{eq57} holds with $z$ sampled from Lebesgue measure on $\wh S_m^0$. In order to get this it is sufficient to argue that the random square $S:=\wh S^0_m$ satisfies the assumption of Lemma \ref{prop30}, which is an immediate consequence of Lemma \ref{prop26} and the definition of $\wh S^0_m$.

	The proof for $\cH^{\op{large}}$ is identical.
\end{proof}

\begin{lem}
	The cell configurations  $(\cH^{\op{cell}},\sim_{\op{cell}},\frk c_{\op{cell}})$ and $(\cH^{\op{large}},\sim_{\op{large}},\frk c_{\op{large}})$ satisfy Assumption \ref{item-ergodic}
	(ergodicity modulo scaling). 
	\label{prop:erg}
\end{lem}	
\begin{proof}
	We only do the proof for $\cH^{\op{cell}}$ since the proof for $\cH^{\op{large}}$ is identical. Let $F=F(\cH^{\op{cell}})$ be a bounded measurable function of $\cH^{\op{cell}}$ satisfying $F(C(\cH^{\op{cell}}-z_0))=F(\cH^{\op{cell}})$ for all $z_0\in\C$ and $C>0$. We need to show that $F$ is a.s.\ equal to a deterministic constant. Since $\cH^{\op{cell}}$ is a measurable function of $(h,h^{\op{IG}},\frk U)$ we have $F(\cH^{\op{cell}})=F_1(h,h^{\op{IG}},\frk U)$ for some bounded measurable function $F_1$. By translation invariance modulo scaling for $F$ and \eqref{eq:transl4}, a.s., 
	\eqb
	F_1\big( h(C\cdot+z_0)+Q\log C, h^{\op{IG}}(C\cdot+z_0), \frk U-\tau_{z_0} \big)=F_1(h,h^{\op{IG}},\frk U),\quad\forall C>0,z_0\in\C.
	\label{eq:transl}
	\eqe
	By arguing exactly as in the proof of \cite[Proposition 3.1]{gms-tutte} we get that $F_1$ is a.s.\ equal to a function only of $\frk U$, i.e., $F_1(h,h^{\op{IG}},\frk U)=F_2(\frk U)$ a.s.\ for some function $F_2$. 
	Furthermore, by the cited proof we may assume that for all $z\in\C$ and all $\frk U$ we have $F_2(\frk U-\tau_z) = F_2(\frk U)$. Defining $B=\{\tau_z\,:\, z\in\C \}$, for any fixed $t\in\R$ it holds a.s.\ that $t\in B$, the set $B$ is random and independent of $\frk U$, and 
	\eqb
	F_2(\frk U-t) = F_2(\frk U)\qquad\forall t\in B.
	\label{eq:transl5}
	\eqe
	Recall from Remark \ref{rmk:U-determines-D} that there are two sources of randomness determining $\frk U$, namely the randomness of $\frk D$ and the randomness of the standard normal random vectors indexed by $\frk D$. In this proof it will be convenient to separate these two sources of randomness and we will therefore work with standard normal random vectors indexed by $\Z\times\Z$ instead of $\frk D$. 
	Explicitly, for $i,j\in\Z$ we let $U_{i,j}=U_{J}$ if $J\in\frk D$, $|J|=|I(j)|$, and if the left end-point of $J$ minus the left end-point of $I(j)$ is $i|J|$; in particular we will have $U_{0,j}=U_{I(j)}$ for all $j\in\Z$. There is a measurable and bijective function $\beta$ such that $\frk U = \beta(\cA)$ for $\cA=(\theta, (V_j)_{j\in\N}, (U_{i,j})_{i,j\in\Z} )$, where $\theta$ and $(V_j)_{j\in\N}$ are as in the definition of a dyadic system. Define $F_3=F_2\circ\beta$. 
	For any $t\in\R$ define $\cA-t$ such that $\frk U-t=\beta(\cA-t)$ and observe that by \eqref{eq:transl5}, $F_3(\cA-t)=F_3(\cA)$ for almost all $\cA$ and Lebesgue almost all $t\in\R$. 
	Define $\alpha_n(\cA):=(\theta, (V_j)_{j\leq n}, (U_{i,j})_{|i|,|j|\leq n} )$ and then define a function $F^n_3$ such that a.s.,
	$$
	F^n_3(\alpha_n(\cA)) := \E[F_3(\cA)  \,|\,\alpha_n(\cA)].
	$$
	Then $F^n_3(\alpha_n(\cA))\rta F_3(\cA)$ a.s.\ by martingale convergence. Using this along with the fact that $\cA-t\eqD \cA$ for any $t$ chosen independently of $\cA$, we get that for such $t$ and some (deterministic) $\eps_n\rta 0$,
	$$
	\P[ |F_3(\cA)-F^n_3(\alpha_n(\cA)) |>\eps_n/2]<\eps_n/2,\qquad 
	\P[ |F_3(\cA-t)-F^n_3(\alpha_n(\cA-t)) |>\eps_n/2]<\eps_n/2.
	$$
	This and $F_3(\cA-t)=F_3(\cA)$ imply 
	$$
	\P[ |F^n_3(\alpha_n(\cA))-F^n_3(\alpha_n(\cA-t)) |>\eps_n]<\eps_n.
	$$
	If $t$ is chosen uniformly at random from $[(2n+3) 2^{n+1},(2n+4)2^{n+1}]$ then the tuple $\alpha_n(\cA-t)$ is independent of the tuple $\alpha_n(\cA)$. Therefore, for independent copies $\wt\theta\eqD\theta$, $\wt V_j\eqD V_j$ and $\wt U_{i,j}\eqD U_{i,j}$,
	$$
	\P[ |F^n_3(\theta,(V_j)_{j\leq n}, (U_{i,j})_{|i|,|j|\leq n})-F^n_3(\wt\theta, (\wt V_j)_{j\leq n}, (\wt U_{i,j})_{|i|,|j|\leq n}) |>\eps_n]<\eps_n.
	$$
	This gives further that for deterministic constants $a_n\in\R$,
	$$
	\P[ |F^n_3(\theta,(V_j)_{j\leq n}, (U_{i,j})_{|i|,|j|\leq n} )-a_n|>\eps_n]<\eps_n,
	$$
	which implies by the a.s.\ convergence $F^n_3(\alpha_n(\cA))\rta F_3(\cA)$ that $F_3(\cA)$ is a.s.\ equal to a constant, which concludes the proof.
\end{proof}
The following lemma is used to verify that Assumption (iii) in Theorem \ref{thm-general-clt-uniform} is satisfied for our cell configurations. Unless $\ell$ is short, the event in (b) has a small probability and we will see later that this case can be disregarded when verifying Assumption (iii) of the cited theorem. In other words the lemma is saying that, disregarding the low probability event in (b), we have $t\in\cT_\ell^{\op{cell}}$ if and only if  $H^{\op{cell}}_t\cap \ell\neq\emptyset$.
\begin{lem}
	It holds with probability one that for any line segment $\ell$ and $t\in\cJLm$ we have
	$t\in\cT_\ell^{\op{cell}}$ if and only if at least one of the following two properties hold: (a) $H^{\op{cell}}_t\cap \ell\neq\emptyset$ or (b) $\ell\subset\ol\eta(u)$ for some $u\in\cT^{\op{cell}}_t$.
	The same results hold with ``$\mathrm{cell}$'' replaced by ``$\mathrm{large}$''. 
	\label{prop:lines2}
\end{lem}
\begin{proof}
	If (b) holds then by the definition of $\cT^{\op{cell}}_t$ we have $t\in\cT^{\op{cell}}_\ell$. Suppose now that (a) holds.
	Since $A_{u,t}\subset\ol\eta(u)$ for each $t\in\cJLm$ and $u\in\cT^{\op{cell}}_t$ (Lemma \ref{prop:Uus}) we have $H^{\op{cell}}_t\subset \cup_{u\in\cT^{\op{cell}}_t}\ol\eta(u)$, so there is some $u\in \cT^{\op{cell}}_t$ such that $\ol\eta(u)\cap \ell\neq\emptyset$ (equivalently, $u\in\cT^\eta_\ell$). This implies that $t\in\cT^{\op{cell}}_\ell$.
	
	Conversely, suppose now that $t\in\cT^{\op{cell}}_\ell$. Then it follows from the definition of $\cT^{\op{cell}}_\ell$ that $\ol\eta(u)\cap \ell\neq\emptyset$ for some $u\in \cT^{\op{cell}}_t$. At least one of the following must occur: (a)' $\ell\cap\partial(\ol\eta(u))\neq\emptyset$ or (b)' $\ell\subset\ol\eta(u)$. We now conclude by observing that (b)' is identical to (b) while (a)' implies (a) since  $\partial(\ol\eta(u))\subset A_{u,t}\subset H_t^{\op{cell}}$.
\end{proof}

The next result gives a general condition under which a cell configuration has no macroscopic cells. It is a corollary of \cite[Lemma 2.9]{gms-random-walk}.
\begin{lem}
	Suppose $(\cH,\sim,\frk c)$ is a cell configuration which is translation invariant modulo scaling, ergodic modulo scaling and satisfies the moment condition in \eqref{eqn-moment}. Then a.s., for each $\eps\in(0,1)$   it holds for large enough $r>1$ that for 
	all $H\in\cH$ intersecting $[-r,r]^2$ we have $\op{diam}(H)\leq\eps r$.
	\label{prop:no-macro-cells}
\end{lem}
\begin{proof}
	Let $\frk D^{\C}$ be a uniform dyadic system on $\C$ and let $S_0,S_1,\dots$ be as defined in the definition of such a system. %
	By \cite[Lemma 2.9]{gms-random-walk}, for each $\wt\eps\in(0,1)$ it holds a.s.\ for large enough $k\in\N$ that for all $H\in\cH$ intersecting $S_k$ we have $\op{diam}(H)\leq\wt\eps 2^{s+k}$. For $j=1,2,3$ and $k\in\N$ let $S^j_k$ be the clockwise rotation of $S_k$ around 0 by $\pi j/2$. Then $(S^j_1,S^j_2,\dots)$ and $(S_1,S_2,\dots)$ are equal in law for each $j\in\{1,2,3 \}$, so the result of the previous paragraph still holds with $S_k$ replaced by $S^j_k$. Observe that $[-2^{k-1},2^{k-1}]^2 \subset S_k\cup S^1_k\cup S^2_k\cup S^3_k$. Combining these observations we see that for large enough $k\in\N$, for all $H\in\cH$ intersecting $[-2^{k-1},2^{k-1}]^2$ we have $\op{diam}(H)\leq \wt \eps 2^{s+k}$. We conclude the proof by choosing $\wt\eps=0.1\eps$.
\end{proof}

The proof of the following moment bound will be deferred to Section \ref{sec:moment}.
\begin{prop}
	The cell configurations $(\cH^{\op{cell}},\sim_{\op{cell}},\frk c_{\op{cell}})$ and $(\cH^{\op{large}},\sim_{\op{large}},\frk c_{\op{large}})$ defined in Section \ref{sec:cells3} satisfy the moment condition \ref{item-moment}.
	\label{prop:moment}
\end{prop}

\begin{lem}
	Assumptions (i)-(iii) of Theorem \ref{thm-general-clt-uniform} are satisfied.
	\label{prop:i-iii}
\end{lem}
\begin{proof}
	Assumption (i) is immediate since both $\cH^{\op{cell}}$ and $\cH^{\op{large}}$ are both isomorphic to $M$. 
	Assumption (ii) follows since $\cT^{\op{cell}}_s\subset \cT^{\op{large}}_s$ for each $s\in\cJLm$ and by the requirement $\partial(\ol\eta(u))\subset A_{u,s}\subset\ol\eta(u)$ in Lemma \ref{prop:Uus}. To prove Assertion (iii) let $r>1$,  $\ell\subset\C$ be a horizontal or vertical line segment of length at least $r$ and at distance at most $r$ from the origin, and suppose that $H^{\op{cell}}_t$ and $H^{\op{cell}}_{t'}$ both intersect $\ell$. Then by Lemma \ref{prop:lines2} we have $t,t'\in\cT^{\op{cell}}_\ell$, so by Lemma \ref{prop:lines} we can find an M-path 
	from $t$ to $t'$ in $\cT^{\op{large}}_\ell$. To conclude the proof it is sufficient to argue that for any $s\in \cT^{\op{large}}_\ell$ on this M-path it holds that $H^{\op{large}}_s\cap \ell\neq\emptyset$. By Lemma \ref{prop:lines2} we know that either $H^{\op{large}}_s\cap \ell\neq\emptyset$ or $\ell\subset\ol\eta(u)$ for some $u\in\cT^{\op{large}}_s$. By Lemma \ref{prop:no-macro-cells} and Proposition \ref{prop:moment} the second scenario does not occur if $r$ is sufficiently large, which concludes the proof.

\end{proof}

\subsection{Cell moment bound}
\label{sec:moment}
The goal of this section is to prove Proposition \ref{prop:moment}. 
In our proof of the proposition we will only consider  $(\cH^{\op{cell}},\sim_{\op{cell}},\frk c_{\op{cell}})$ 
since the proof for 
$(\cH^{\op{large}},\sim_{\op{large}},\frk c_{\op{large}})$ is the same. Recall from \eqref{eq:ol-eta} that $\ol\eta(t)=\eta(I_t)$ for $t\in\cJm$, where $I_t$ is  
the interval defined at the beginning of Section \ref{sec:cells1}, and set
$$
H^{\eta}_t := \ol\eta(t),\qquad 
R_t := \frac{\op{diam}(\ol\eta(t))^2}{\op{area}(\ol\eta(t))}.
$$

Notice that $H_t^{\op{cell}}$ is defined for $t\in\cJLm$, whereas $H_t^\eta$ is defined for every $t\in\cJm$. In the moment condition \ref{item-moment} we consider the origin-containing cell, which motivates the following definition. 
Let $\op{cell}(0)\in\cJLm$ be the a.s.\ unique element of $\cJLm$ such that $0\in H_{\op{cell}(0)}^{\op{cell}}$. Below we generalize this definition to $\op{cell}(s)$ for $s\in\R$.

Since ${\cT}^{\op{cell}}_{\op{cell}(0)}$ is G-connected by Lemma \ref{prop25}, by \eqref{eq:cell-conf}, and by the property of G-connectedness discussed right below \eqref{eq:ol-eta}, we have
$$
\diam(H^{\op{cell}}_{{\op{cell}(0)}}) \leq \sum_{s\in {\cT}^{\op{cell}}_{\op{cell}(0)}}\diam(\ol\eta(s)).
$$
By  \eqref{eq:cell-conf},
$$
\op{area}(H_{\op{cell}(0)}^{\op{cell}})
= 
\sum_{u\in {\cT}^{\op{cell}}_{\op{cell}(0)}} \frac{1}{2\cdot\#{\cT}_u^{\op{cell,inv}}}\op{area}( \ol\eta(u) ).
$$
Combining the last two indented equations we get
\eqbn
\begin{split}
	&\EE\Bigg[ \frac{\diam\big(H^{\op{cell}}_{{\op{cell}(0)}}\big)^2}{\op{area}\big(H^{\op{cell}}_{\op{cell}(0)}\big)} \deg(H^{\op{cell}}_{\op{cell}(0)}) \Bigg]
	\leq
	2\EE\Bigg[ \frac{\Big(\sum_{t\in {\cT}^{\op{cell}}_{\op{cell}(0)}}\diam(\ol\eta(t) ) \Big)^2}{\sum_{t\in {\cT}^{\op{cell}}_{\op{cell}(0)}}\op{area}(\ol\eta(t))/\# {\cT}^{\op{cell,inv}}_t} \deg(H^{\op{cell}}_{\op{cell}(0)}) \Bigg]\\
	&\qquad\qquad\qquad\qquad\leq 
	2\EE\Bigg[ \frac{  \sum_{t\in {\cT}^{\op{cell}}_{\op{cell}(0)}} \op{area}(\ol\eta(t) ) R_t }
	{\sum_{t\in {\cT}^{\op{cell}}_{\op{cell}(0)}}\op{area}(\ol\eta(t))} \deg(H^{\op{cell}}_{\op{cell}(0)}) \#{\cT}^{\op{cell}}_{\op{cell}(0)} \max_{t\in {\cT}^{\op{cell}}_{\op{cell}(0)} } \# {\cT}^{\op{cell,inv}}_t \Bigg]\\
	&\qquad\qquad\qquad\qquad\leq 
	2\EE\Bigg[  \max_{t\in {\cT}^{\op{cell}}_{\op{cell}(0)}}  R_t 
	\cdot
	\Big(\# {\cT}^{\op{cell}}_{\op{cell}(0)} \Big)^2
	\cdot
	\Big(\max_{t\in {\cT}^{\op{cell}}_{\op{cell}(0)}} \# {\cT}^{\op{cell,inv}}_t \Big) \Bigg],
\end{split}
\eqen
where in the last step we use that $\deg(H^{\op{cell}}_{\op{cell}(0)})\leq \# {\cT}^{\op{cell}}_{\op{cell}(0)}$ by Lemma \ref{prop:em-adj}.

The following two propositions, which will be proved in the next two subsections, imply Proposition \ref{prop:moment} in combination with the last indented equation.
\begin{prop}
	$\displaystyle \EE\Big[(%
	\# {\cT}^{\op{cell}}_{\op{cell}(0)})^p 
	\cdot 
	\Big(\max_{t\in {\cT}^{\op{cell}}_{\op{cell}(0)}} \# {\cT}^{\op{cell,inv}}_t\Big)^p\Big]
	<\infty$ for all $p>1$.
	\label{prop:moment-T}
\end{prop}

\begin{prop}
	$\displaystyle \EE\Big[\max_{t\in {\cT}_{\op{cell}(0)}^{\op{cell}}}(R_t)^p\Big]<\infty$ for all $p>1$.
	\label{prop:moment-R}
\end{prop}

\subsubsection{Proof of Proposition \ref{prop:moment-T}}
	For $s\in\cJm$ define
	\eqb
	\wt\cT^{\op{C},1}_s = \bigcup_{ u\in \cT^{\op{C},1}_s } \cT^{\op{C},1}_u. 
	\label{eq:wtcelldef}
	\eqe
	We will prove that $\max_{s\in\wt{\cT}^{\op{C},1}_{\op{cell}(0)}} \#\wt\cT^{\op{C},1}_s$ has finite moments of all orders. This is sufficient to conclude the proof of  Proposition \ref{prop:moment-T} by the following lemma.
	\begin{lem}
		For any $s\in\cJLm$ we have $\cT^{\op{cell}}_s\subseteq \wt\cT^{\op{C},1}_s$ and for any $s\in\cJm$ we have $\cT^{\op{cell,inv}}_s\subseteq \wt\cT^{\op{C},1}_s$.
		\label{prop24}
	\end{lem}
\begin{proof}
	We get $\cT^{\op{cell}}_s\subseteq \wt\cT^{\op{C},1}_s$ by using the definitions in equations \eqref{eq:wtT} and \eqref{eq:wtcelldef} along with Assertion 3.\ of Proposition \ref{prop:time-sets}. 
	We get $s\in \cT^{\op{cell}}_s$ for $\cJLm$ by using \eqref{eq:wtT} and Assertion 1.\ of Proposition \ref{prop:time-sets}. 
	To prove the second assertion we first observe that by the definitions of $\cT^{\op{cell}}_s$ and $\cT^{\op{cell,inv}}_s$ in \eqref{eq:wtT} and \eqref{eq74}, respectively, we have
	$$
	\cT^{\op{cell,inv}}_s = \bigcup_{t\in \cT^{\op{C},1}_s} \{u\in\cJLm \,:\,u=t\text{\,\,or\,\,}\{u,t\}\in\ET \}.
	$$
	This gives the second assertion in combination with Assertion 3.\ of Proposition \ref{prop:time-sets}.
\end{proof}
	Before explaining our proof strategy, let us define $\op{up}(s,J)\in J$ and $\op{cell}(s)\in\cJLm$ for $s\in\R$ and $J\subset \R$ a locally finite subset of $\R$ which is unbounded above, where we remark that the latter definition generalizes the definition of $\op{cell}(0)$ from the second paragraph below Proposition \ref{prop:moment} to general $s\in\R$. 
	First let $\op{up}(s,J)=\inf (J\cap[s,\infty))$. In other words, we choose the smallest element of $J$ which is greater than or equal to $s$. Typically we will use this definition with $J$ equal to $\cJ$, $\cJm$, $\cJLm$ or $\cJRm$. 
	For $s\in\R$ let $\op{cell}(s)\in\cJLm$ be such that $\eta(s)\in H^{\op{cell}}_{\op{cell}(s)}$. For fixed $s\in\R$ there is a.s.\ a unique choice of $\op{cell}(s)\in\cJLm$ satisfying this. In order for $\op{cell}(s)$ to be a.s.\ well-defined simultaneously for \emph{all} $s\in\R$, in case there are multiple $t\in\cJLm$ such that $\eta(s)\in H^{\op{cell}}_t$, we let $\op{cell}(s)=\op{up}(s,\cJLm)$. 

	Let $h>1$. For $s\in\cJm$ and $k\in\N\cup\{\infty \}$ define
	\eqb
	\wt\cT^{h,k}_s = \bigcup_{ u\in \cT^{h,k}_s } \cT^{h,k}_u. 
	\label{eq85}
	\eqe
	Then define 
	$$
	\op{cell}_{h,\infty}(t)
	:=
	\wt{\cT}^{h,\infty}_{\op{up}(t,\cJLm)}\cup \wt{\cT}^{h,\infty}_{\op{up}(t,\cJRm)},\qquad
	\op{cell}_{h,k}(t)
	:=
	\wt{\cT}^{h,k}_{\op{up}(t,\cJLm)}\cup \wt{\cT}^{h,k}_{\op{up}(t,\cJRm)}.
	$$ 
	For $A\subset\cJm$ set $\wt{\cT}^{h,\infty}_{A}=\cup_{t\in A}\wt{\cT}^{h,\infty}_{t}$.  
	In order to prove Proposition \ref{prop:moment-T} we will apply Corollary \ref{prop:eq0} with 
	$$
	G_h(\cZ,\frk D)=\1\Big(\max_{s\in \wt{\cT}^{h,\infty}_{ \op{cell}_{h,\infty}(0) } } \# \wt{\cT}^{h,\infty}_s>h^{18}\Big), \qquad
	G_{h,k}(\cZ,\frk D)=\1\Big(\max_{s\in \wt{\cT}^{h,k}_{\op{cell}_{h,k}(0)}} \# \wt{\cT}^{h,k}_s>h^{18}\Big),
	$$
	where we write $G_h$ and $G_{h,k}$ instead of $G$ and $G_k$, respectively, in order to emphasize the dependence on $h$, and we use the convention that the maximum of an empty set is 0. We get
	\eqbn
	\begin{split}
		A_h= \{\exists t\in I(0) \text{ such that}\max_{s\in \wt{\cT}^{h,\infty}_{\op{cell}_{h,\infty}(t)}} \# \wt{\cT}^{h,\infty}_s>h^{18} \},\\
		A_{h,k}= \{\exists t\in I(0) \text{ such that}\max_{s\in \wt{\cT}^{h,k}_{\op{cell}_{h,k}(t)}} \# \wt{\cT}^{h,k}_s>h^{18} \},
	\end{split}
	\eqen
	where we write $A_h$ and $A_{h,k}$ instead of $A$ and $A_k$, respectively.
	
	For $s\in\cJm$ and $d\in\N$ define
	\eqb
	\wt\cT^{\op{D},d}_s = \bigcup_{ u\in \cT^{\op{D},d}_s } \cT^{\op{D},d}_u. 
	\label{eq:wtDdef}
	\eqe
	The next lemma says that various time sets we will consider in the proof are contained in $I(m)$ with high probability for large $m$.
\begin{lem} There is a $c\in(0,1)$ such that for any $h>1$ and $m\in\N$,
	\eqb
	\PE[(E_{h,m})^c] 
	\leq c^{-1} h^{c^{-1}}2^{-c m},
	\qquad\text{where}\qquad 
	E_{h,m}:=\Big\{ \bigcup_{t\in I(0)
	} \bigcup_{
	s\in \wt{\cT}^{h,\infty}_{\op{cell}_{h,\infty}(t)} 
	} \wt{\cT}^{h,\infty}_s \subset I(m)\Big\}.
	\label{eq26}
	\eqe
	Moreover, the same holds with $\PQ$ instead of $\PE$. Finally, with $d\in\N$ and $c$ depending on $d$,
	$$
	\PE[(\wt E_{h,\op{D},d,m})^c] 
	\leq c^{-1} h^{c^{-1}} 2^{-c m},
	\qquad\text{where}\qquad 
	\wt E_{h,\op{D},d,m}:=\Big\{ \bigcup_{t\in I(0)}\bigcup_{s\in \wt{\cT}^{\op{D},d}_{\op{cell}_{h,\infty}(t)}} \wt{\cT}^{\op{D},d}_s \subset I(m)\Big\}.
	$$
		\label{T-tail}
\end{lem}
\begin{proof}
	In order to prove \eqref{eq26} we will apply Assertion (iii) of Proposition \ref{prop:Text-len-tail} repeatedly. 
	By the definition of $\wt{\cT}^{h,\infty}_s$, the time set appearing in the definition of $E_{h,m}$ is equal to
	\eqb
	\begin{split}
		\{ s\in\cJm\,:\, 
		\exists t_2,\dots,t_6\in \cJm, t_1\in I(0) \text{\,\,such\,\,that\,\,}&
		t_2\in \op{cell}_{h,\infty}(t_1), \\
		&s=t_6 \text{\,and\,} t_{j+1}\in \cT^{h,\infty}_{t_j}\text{ for }j=2,\dots,5 \}.
	\end{split}
	\label{eq25}
	\eqe
	Let $J^0$ be the shortest closed interval such that 
	$$
	\{\op{up}(s,\cJLm)\,:\, s\in I(0) \}\cup \{\op{up}(s,\cJRm)\,:\, s\in I(0) \} \subseteq J^0.
	$$ 
	Then define $J^1,J^2,\dots$ such that for $j=1,2,\dots$, $J^j$ is the shortest closed interval such that
	$$
	\bigcup_{s\in J^{j-1}\cap\cJm }\cT^{h,\infty}_s \subset J^j.
	$$
	Applying Assertion (iii) of Proposition \ref{prop:Text-len-tail} six times with $J=J^0,J^1,\dots,J^4,J^5$ we get that $\PE[J^6\not\subset I(m)]\leq c^{-1}h^{c^{-1}}2^{-cm}$. This concludes the proof since by \eqref{eq85} we know that if $t_2,\dots,t_6$ are as in \eqref{eq25} then $t_j\in J^j$ for $j=2,\dots,6$. 

	The second and third assertions of the lemma follow by the exact same argument using 
 	Assertion (iv) and Assertion (i), respectively, of Proposition \ref{prop:Text-len-tail} instead of Assertion (iii). The details are therefore omitted. 

\end{proof}

The next lemma says that with high probability under $\PQ$, no two points in $I(\lceil \wt h \rceil)\cap\cJ$ are very close. This is an immediate consequence of the fact that $\cJ$ is a unit rate Poisson point process. The second assertion is an immediate consequence of the definition of $\op{up}(s,J)$ for $s\in\R$ and $J\subset\R$ locally finite.
\begin{lem}
	If 
	$$
	\wh E_{\wt h}=\{\forall s\in\cJ\cap [-2^{\wt h}, 2^{\wt h}], t\in\cJ 
	\text{ such that } s\neq t \text{ we have }|s-t|>2^{-2\wt h} \}
	$$ 
	then $\PQ[(\wh E_{\wt h})^c]\leq 4\cdot 2^{-\wt h}$. Furthermore, the event $\wh E_{\wt h}$ is contained in the following event for $J=\cJ$ and $J=\cJm$
	$$
	\big\{ \{ \op{up}(k2^{-2\wt h},J)\,:\, k\in\Z \}\cap [-2^{\wt h}, 2^{\wt h}] = J\cap [-2^{\wt h}, 2^{\wt h}] \big\}.
	$$
	\label{prop:cover}
\end{lem}

\begin{lem}
	Given any $\wt C>1$ it holds that $\PQ[\# \cT^{h,\infty}_{\op{up}(0,\cJm)} >h^{9}]\leq 2^{-\wt Ch^{9/8}}$ for sufficiently large $h>1$. 
	\label{prop:Tsize}
\end{lem}
\begin{proof}
	The majority of the proof will be devoted to showing that for sufficiently large $h>1$,
	\eqb
	\PQ[\# \cT^{h,\infty}_{\op{up}(0,\cJm)}>h^{8}]\leq 2^{-\wt Ch}.
	\label{eq30}
	\eqe
	Once this has been established the lemma follows since with $r:=h^{9/8}>1$,
	$$
	\PQ[\# \cT^{h,\infty}_{\op{up}(0,\cJm)}>h^{9}]
	\leq 
	\PQ[\# \cT^{r^{8/9},\infty}_{\op{up}(0,\cJm)}>r^{8}]
	\leq 
	\PQ[\# \cT^{r,\infty}_{\op{up}(0,\cJm)}>r^{8}]
	= 2^{-\wt Cr} = 2^{-\wt Ch^{9/8}},
	$$
	where we get the second inequality by using that $\cT^{r^{8/9},\infty}_{\op{up}(0,\cJm)}\subseteq \cT^{r,\infty}_{\op{up}(0,\cJm)}$.
	
	Let us now prove \eqref{eq30}. By Lemma \ref{T-tail}, and since $\wt{\cT}^{h,\infty}_{\op{up}(0,\cJm)}$ is contained in the time set appearing in the definition of the event $E_{h,m}$ in \eqref{eq26}, for sufficiently large $C_0>1$,
	$$ 
	\PQ[ \wt\cT^{h,\infty}_{\op{up}(0,\cJm)}\not\subset I(C_0 h) ] \leq 2^{-2\wt Ch}.
	$$
	We claim that in order to prove \eqref{eq30} it is sufficient to prove that
	\eqb
	\PQ[\exists t\in I(C_0 h)\cap\cJ\,s.t.\,\#\cT^{h,\infty}_{t,\op{L}}>h^4/2] \leq 2^{-2\wt Ch}.
	\label{eq76}
	\eqe
	Suppose the events in the latter two indented equations do not occur, and suppose the event in the latter indented equation with ``R'' instead of ``L'' does not occur either. Then we have $\#\cT^{h,\infty,0}_{t}\leq h^4$ for all $t\in I(C_0h)\cap\cJ$ by \eqref{eq27} and further $\#\cT^{h,\infty}_{{\op{up}(0,\cJm)}}\leq (h^4)^2$ by \eqref{eq28}. This proves the claim. The remainder of the proof will be devoted to establishing \eqref{eq76}.

	Recalling notation from Lemma \ref{prop:cover}, a union bound gives 
	\eqb
	\begin{split}
		&\PQ[\exists t\in I(C_0 h)\cap\cJ\,s.t.\,\#\cT^{h,\infty}_{t,\op{L}}>h^4/2]%
		\leq 
		\PQ[(\wh E_{\lceil C_0h\rceil})^c] + \sum_{t\in \cI } \PQ[\#\cT^{h,\infty}_{{\op{up}(t,\cJ)},\op{L}}>h^4/2],
	\end{split}
\label{eq87}
	\eqe
	where $\cI=\{ t=k2^{-2\lceil C_0h\rceil}\,:\,k\in \Z, t\in[-2^{\lceil C_0h\rceil}-1,2^{\lceil C_0h\rceil}] \}$. The first probability on the right side is bounded above by $2^{-2\wt Ch}$ for sufficiently large $C_0>1$ by Lemma \ref{prop:cover}. Since $\cZ$ is translation invariant in law the sum on the right side is bounded above by $\#\cI \cdot\PQ[\#\cT^{h,\infty}_{{\op{up}(0,\cJ)},\op{L}}>h^4/2]$, where we have $\#\cI\leq 2^{3\lceil C_0h\rceil+2}$. 
	The lemma now follows by using that by the definition \eqref{eq29} of $\cT^{h,\infty}_{0,\op{L}}$, for a symmetric simple random walk $(W_j)_{j\in\N_0}$ on $\Z$ started from 0 and for some $c\in(0,1)$, 
	\eqbn
	\begin{split}
		&\PQ[\#\cT^{h,\infty}_{{\op{up}(0,\cJ)},\op{L}}>h^4/2]\\ 
		&\,\,\,\leq 
		\P[\#\{ j\in\{ 1,\dots,\tau\}\,:\, W_{j}=W_{j-1}-1 \text{\,and\,} \{W_{j-1},W_{j}\}\cap[-(h+1),h+1]\neq\emptyset \}>h^4/2]
		\leq 2^{-c h^2},
	\end{split}
	\eqen
	where $\tau=\inf\{j\in\N_0\,:\,W_j<-(h+1) \}$.
\end{proof}

\begin{lem}
	$\PQ[A_{h}]\leq 2^{-ch}$.
	\label{prop:Pah}
\end{lem}
\begin{proof}
	Recalling the definition of the event $E_{h,m}$ in \eqref{eq26} and notation from Lemma \ref{prop:cover}, we have 
	\eqbn
	\begin{split}
		\PQ[A_h]&=
		\PQ\Big[\exists t\in I(0) \text{ such that }\max_{s\in \wt{\cT}^{h,\infty}_{\op{cell}_{h,\infty}(t) }} \# \wt{\cT}^{h,\infty}_s>h^{18}\Big]\\
		&\leq\PQ[(E_{h,h})^c] + \PQ[\exists t\in I(h)\cap\cJm \,s.t.\, \#\cT^{h,\infty}_t>h^9]\\
		&\leq 2^{-ch} +\PQ[(\wh E_h)^c]
		+2^{3\lceil h\rceil+1}\PQ\big[\# \cT^{h,\infty}_{\op{up}(0,\cJm)}>h^{9}\big]\leq 2^{-c h},
	\end{split}
	\eqen
	where in the second step we used a union bound and the definition of $\wt{\cT}^{h,\infty}_s$, 
	in the third step we used 
	Lemma \ref{T-tail}, a union bound, and translation invariance in law of $\cZ$ in a similar fashion as in \eqref{eq87}, and in the last step we used Lemmas \ref{prop:cover} and \ref{prop:Tsize}.
\end{proof}

\begin{lem}
	There is a $c\in(0,1)$ such that for any $h>1$ and $k\in\N$ we have
	$\PE[\1_{A_h}\neq \1_{A_{h,k}}]\leq h^{c^{-1}} 2^{-ck}$.
	\label{prop:Ah-neq-Ahk}
\end{lem}
\begin{proof}
	We can write $I(k)\setminus I(k/2)$ as the union of two disjoint intervals $I_k^{\op{L}}$ and $I_k^{\op{R}}$ such that $I_k^{\op{L}}$ (resp.\ $I_k^{\op{R}}$) is to the left (resp.\ right) of $I(k/2)$. 
	By Lemma \ref{prop:localT}(iii) 
	the considered probability is bounded above by $$\PE[(E_{h,k/2})^c] + \PE[I_k^{\op{L}}=\emptyset]+ \PE[\cJLm\cap I_k^{\op{L}}=\emptyset;I_k^{\op{L}}\neq\emptyset] + \PE[\cJRm\cap I_k^{\op{L}}=\emptyset;I_k^{\op{L}}\neq\emptyset],$$ where $E_{h,k/2}$ is defined in \eqref{eq26}. The first term on the right side is bounded by $h^{c^{-1}}2^{-ck}$ by Lemma \ref{T-tail}. The second term on the right side is bounded above by $2^{-ck}$ by the definition of a uniform dyadic decomposition of $\R$. By Theorem \ref{prop:zero-vs-gamma-walk} the third term on the right side is bounded above by $2^{C_0k}\PQ[\cJLm\cap I_k^{\op{L}}=\emptyset;I_k^{\op{L}}\neq\emptyset]^{1/4}+2^{-k/C_0}$. We have $\PQ[\cJLm\cap I_k^{\op{L}}=\emptyset;I_k^{\op{L}}\neq\emptyset]\leq 2^{-c2^{k/2}}$ since $I_k^{\op{L}}$ has length at least $2^{k/2}$ on the event that $I_k^{\op{L}}\neq\emptyset$ and it follows that the third term on the right side of the indented equation is also bounded by $2^{-ck}$. The fourth term is bounded similarly.
\end{proof}

\begin{lem}
	There is a $c\in(0,1)$ such that $\PE[(\ol E_h)^c]\leq \exp(-c h)$ for all $h>1$, where 
	$$
	\ol E_h:=\{\wt{\cT}^{\op{C},1}_{s}\subset \wt{\cT}^{h,\infty}_{s}\,\forall s\in \wt{\cT}^{\op{C},1}_{\op{cell}_{h,\infty}(0) }\}
	\cap
	\{\op{cell}(0)\in \op{cell}_{h,\infty}(0) \}.
	$$
	\label{prop:cell-neq-Th}
\end{lem}
\begin{proof}
By first using Assertion 1.\ of Proposition  \ref{prop:time-sets} and then using Lemma \ref{prop:localT}(iv), for any $\eps\in(0,1)$, 
	\eqbn
	\PE[(\ol E_h)^c] 
	\leq 
	\PE\big[\exists s\in\wt\cT^{\op{D},2}_{\op{cell}_{h,\infty}(0)} \text{\,s.t.\,}\wt\cT^{\op{D},2}_s\not\subseteq \wt\cT^{h,\infty}_s\big].
	\eqen
	By the first part of \eqref{eq:cell-conf}  we have $\op{cell}(0)\in \cT^{\op{cell,inv}}_{\op{up}(0,\cJLm)}\cup \cT^{\op{cell,inv}}_{\op{up}(0,\cJRm)}$. 
	Furthermore, by Lemma \ref{prop24} and Assertion 1.\ of Proposition \ref{prop:time-sets}, $\cT^{\op{cell,inv}}_{u} \subseteq \wt{\cT}^{\op{C},1}_{u}  \subseteq 
	\wt{\cT}^{\op{D},2}_{u}$ for any $u\in\cJm$. Using this and the definition $\op{cell}_{h,\infty}(u)=
	\wt{\cT}^{h,\infty}_{\op{up}(u,\cJLm)}\cup 
	\wt{\cT}^{h,\infty}_{\op{up}(u,\cJRm)}$ for $u\in\R$,
	$$
	\PE[ \op{cell}(0)\not\in\op{cell}_{h,\infty}(0) ]
	\leq 
	\PE\Big[
	\wt{\cT}^{\op{D},2}_{\op{up}(0,\cJLm)}  \cup  
	\wt{\cT}^{\op{D},2}_{\op{up}(0,\cJRm)} \not\subset \wt{\cT}^{h,\infty}_{\op{up}(0,\cJLm)}  \cup  
	\wt{\cT}^{h,\infty}_{\op{up}(0,\cJRm)}
	\Big] 
	$$
	By using that $\{\op{up}(0,\cJLm),\op{up}(0,\cJRm)\}\subset \op{cell}_{h,\infty}(0)$, the right sides of the last two displays are bounded above by 
	\eqbn
	\PE\big[\exists s\in\wt\cT^{\op{D},2}_{\op{cell}_{h,\infty}(0)} \text{\,s.t.\,}\wt\cT^{\op{D},2}_s\not\subseteq \wt\cT^{h,\infty}_s\big]
	\leq
	\PE[(\wt E_{h,\op{D},2,\eps h})^c] + \PE[(\check E_{h,\eps h})^c],
	\eqen
	where $\wt E_{h,\op{D},d,\eps h}$ is defined in Lemma \ref{T-tail} and $\check E_{h,m}$ is the event on the left side of \eqref{eq18}. By Lemma \ref{T-tail} the first probability on the right side is bounded above by $c^{-1}2^{-c\eps h}$ for some $c$. By Lemmas \ref{prop:eucl-error} and \ref{prop22} for $\eps$ sufficiently small the second probability on the right side is bounded above by $2^{-ch}$ for some $c$ depending on $\eps$, which concludes the proof.
\end{proof}

\begin{proof}[Proof of Proposition \ref{prop:moment-T}]
	By first doing a union bound and then using Lemma \ref{prop:cell-neq-Th}  and Corollary \ref{cor:zero-vs-gamma}, for some $c\in(0,1)$ independent of $h$,
	$$
	\PE\Big[  \max_{t\in \wt{\cT}^{\op{C},1}_{\op{cell}(0)}} \# \wt{\cT}^{\op{C},1}_t>h^{18}\Big] 
	\leq \PE[ A_h ] + \PE[(\ol E_h)^c]
	\leq h^{c^{-1}}c^{-1}\PQ[ A_h ]^c + 2^{-c h},
	$$
	where we see that Corollary \ref{cor:zero-vs-gamma} can be applied by Lemma \ref{prop:Ah-neq-Ahk} and since $A_h$ and $A_{h,k}$ are of the desired form \eqref{eq:Ak}, $A_{h,k}\in\sigma(\cZ|_{I(k)},I(k))$ by Assertion (ii) of Lemma \ref{prop:localT}, and $A_{h,k}\subset A_{h,k+1}\subset A_h$ by Assertion (i) of Lemma \ref{prop:localT}. Decreasing $c$ if necessary, the right side of the last display is bounded by $2^{-c h +1}$ by Lemma \ref{prop:Pah}. This concludes the proof since the expectation considered in the proposition is bounded above by $\PE\big[  \max_{t\in \wt{\cT}^{\op{C},1}_{\op{cell}(0)}} \big(\# \wt{\cT}^{\op{C},1}_t\big)^{2p}\big]$ by Lemma \ref{prop24}.
\end{proof}

\subsubsection{Proof of Proposition \ref{prop:moment-R}}
We fix $\wh C>1$ throughout this subsection and let $C_0>1$ be the constant appearing in Theorem \ref{prop:zero-vs-gamma-walk}. It is sufficient to argue that for all sufficiently large $M>1$,
$$
\PE\Big[\max_{t\in {\cT}_{\op{cell}(0)}^{\op{cell}}}R_t>M\Big] \leq 10M^{-\wh C}.
$$
\begin{lem}
	By choosing $C_1>C_0\wh C$ and $C_2>1$ sufficiently large, for all sufficiently large $M>1$ it holds that $\PE[B_1\cup B_2]<M^{-\wh C}$, where
	\eqb
	\begin{split}
		B_1&=\{\exists t\in {\cT}_{\op{cell}(0)}^{\op{cell}} \text{\,such\,that\,}  I_t\not\subset I(C_1\log M) \},\\
		B_2& = \{ \exists t\in \cJm\,:\, I_t\subset I(C_1\log M) \text{\,and\,}|I_t|\not\in (M^{-C_2},C_2\log M) \}.
	\end{split}
	\eqe
	\label{prop6-0}
\end{lem}
\begin{proof}
	Define $\wt B_1= \{\exists t\in {\cT}_{\op{cell}(0)}^{\op{cell}} \text{\,such\,that\,}  t\not\in I_{M}(C_1\log M) \}
	=\{ {\cT}_{\op{cell}(0)}^{\op{cell}}\not\subset I_{M}(C_1\log M) \}$, where $I_{M}(C_1\log M):=[a_1+M,a_2)$ for $I(C_1\log M)=[a_1,a_2)$. Define the event $B_\partial$ by
	$$
	B_\partial = \{\cJLm\cap[a_1,a_1+M]=\emptyset \}\cup \{\cJRm\cap[a_1,a_1+M]=\emptyset \},
	$$
	where we suppose $C_1$ is sufficiently large so that $2^{0.5C_1\log M}>M^{1.1}$.  
	Observe that $B_1\subset \wt B_1\cup B_\partial$, so in order to conclude it is sufficient to show that $\PE[\wt B_1\cup B_2\cup B_\partial]\leq M^{-\wh C}$.

	By choosing $C_1$ sufficiently large we get that for sufficiently large $M$
	$$
	\PE[\wt B_1] \leq \PE[{\cT}_{\op{cell}(0)}^{\op{cell}}\not\subset I(0.5C_1\log M)] + \PE[I(0.5C_1\log M)\not\subset I_{M}(C_1\log M) ] \leq M^{-\wh C-1},
	$$
	where we get the last inequality by using Lemma \ref{T-tail} and that, for sufficiently large $M$, we have $\{ I(0.5C_1\log M)\not\subset I_{M}(C_1\log M)\}$ only if $I(0.5C_1\log M)$ and $I(C_1\log M)$ have the same left end-point. 
	
	By Theorem \ref{prop:zero-vs-gamma-walk} with $k=C_1\log M$, 
	\eqbn
	\begin{split}
		\PE[B_2] \leq M^{C_0C_1} 
		\PQ[B_2]^{1/4}+ M^{-C_1/C_0}\leq M^{-\wh C-1},
	\end{split}
	\eqen
	where we used in the last step that $\cJLm$ and $\cJRm$ have the law of rate $0.25$ Poisson point processes so that $\PQ[B_2]\leq M^{-5\wh C-4C_0C_1-4}$ for sufficiently large $C_2$. By another application of Theorem \ref{prop:zero-vs-gamma-walk},
	\eqbn
	\begin{split}
		\PE[B_\partial] \leq M^{C_0C_1} 
		\PQ[B_\partial]^{1/4}+ M^{-C_1/C_0}\leq M^{-\wh C-1},
	\end{split}
	\eqen
	where we used in the last step that $\PQ[B_\partial]\leq 2e^{-0.25M}$.

\end{proof}

The next lemma is an immediate consequence of Lemma \ref{prop:contained1}.
\begin{lem}
	By choosing $C_3>C_0\wh C$ sufficiently large we can guarantee that
	$\PE[\eta( I(C_1\log M) )\not\subset B(M^{C_3}/2)]\leq M^{-\wh C}$, where we use the $M^{C_3}$-embedding of the 0-cone.
	\label{prop6}
\end{lem}
The next lemma is an immediate consequence of Lemmas \ref{prop:contained3} and \ref{prop:cone-moment} along with scale invariance and Markov's inequality.
\begin{lem}
	For $C_4>1$ sufficiently large, if the field has the $M^{C_3}$-embedding,
	$$
	\PQ\big[ \ol{B(0,2M^{-C_3/2})} \not\subset \eta([-M^{C_4}/2,M^{C_4}/2])\big] 
	+\PQ[\mu_h(B(0,M^{C_3})) > M^{C_4}/2 ]
	\leq M^{-3C_0C_3-3\wh C}.
	$$ 
	\label{prop7}
\end{lem}
The next lemma is proven in Appendix \ref{app3}.
\begin{lem}
	For any $C>0$ it holds that for $M>1$ sufficiently large, $\PQ[ B_3] \leq M^{-C}$, where
	$$ 
	B_3 =\Big\{ \max\Big\{ \frac{\op{diam}(\eta(J))^2}{\op{area}(\eta(J))} \,:\, J\subset [-M^{C_4},M^{C_4}] \text{ interval}, \,M^{-C_2}\leq|J|\leq C_2\log M\Big\} >M\Big\}.
	$$
	\label{prop:prop44}
\end{lem}

The next lemma is a variant of the former lemma for $\PE$ instead of $\PQ$.
\begin{lem}
	For all sufficiently large $M>1$ we have $\PE[B_4]\leq 2M^{-\wh C}$, where
	$$ 
	B_4 =\Big\{ \max\Big\{ \frac{\op{diam}(\eta(J))^2}{\op{area}(\eta(J))} \,:\, J\subset I(C_1\log M) \text{ interval}, \,M^{-C_2}\leq|J|\leq C_2\log M\Big\} >M\Big\}.
	$$
	\label{prop9}
\end{lem}
\begin{proof}
	Define
	$$
	\tau_-=\inf\{t<0\,:\,\eta([t,0])\subset B(0.5M^{C_3}) \},\qquad
	\tau_+=\sup\{t>0\,:\,\eta([0,t])\subset B(0.5M^{C_3}) \}.
	$$
	Most of the proof will be devoted to showing that $\PE[B_5]\leq M^{-\wh C}$ for all sufficiently large $M$, where
	$$ 
	B_5 =\Big\{ \max\Big\{ \frac{\op{diam}(\eta(J))^2}{\op{area}(\eta(J))} \,:\, J\subset [\tau_-,\tau_+] \text{ interval}, \,M^{-C_2}\leq|J|\leq C_2\log M\Big\} >M\Big\}.
	$$
	Once this has been proven the lemma follows by combining this estimate with Lemma \ref{prop6}.
	
	We will apply Proposition \ref{prop:zero-vs-gamma-ball} with $G=\1_{B_5}$ and $2^{2k}=M^{C_3}$. In the notation of the cited proposition, 
	$$
	\PE[B_5] = \PE[A(\{0 \})] \leq \PQ\big[A\big(\ol{B(0,M^{-C_3/2})}\big)\big]^{1/2} M^{C_0C_3/2}.
	$$
	By first doing a union bound and then applying Lemmas \ref{prop7} and \ref{prop:prop44},
	\eqbn
	\begin{split}
		\PQ\big[A\big(\ol{B(0,M^{-C_3/2})}\big)\big] 
		&\leq 
		\PQ\big[ \ol{B(0,2M^{-C_3/2})} \not\subset \eta([-M^{C_4}/2,M^{C_4}/2])\big] \\
		&\qquad\qquad+ \PQ[ \mu_h( B(0,M^{C_3}))>M^{C_4}/2 ]
		+ \PQ[B_3] 
		\leq M^{-2C_0C_3-2\wh C}.
	\end{split}
	\eqen
	By combining the last two indented equations we get $\PE[B_5]\leq M^{-\wh C}$ as desired. 
\end{proof}

\begin{proof}[Proof of Proposition \ref{prop:moment-R}]
	By a union bound and Lemmas \ref{prop6-0}, \ref{prop6}, and \ref{prop9},
	\eqb
	\begin{split}
		&\PE\Big[\max_{t\in \cT_{\op{cell}(0)}^{\op{cell}}}R_t>M\Big]
		\leq \PE[B_1]+ \PE[ B_2 ]+\PE[B_4]\leq 3M^{-\wh C}. 
	\end{split}
	\eqe
\end{proof}

\section{Convergence of the Tutte embedding}
\label{sec:tutte-conv}
	
In this section we conclude the proof of Theorem \ref{thm-tutte-conv0}. Recall from Proposition \ref{prop:mullin-wired-disk2} that the planar map appearing in the theorem can be bijectively encoded by a lattice walk $\cZ^n$ on $\N_0^2$ from $(\lceil \sqrt{n} \rceil-1,0)$ to $(0,0)$. Couple $\cZ^n$ with a Brownian cone excursion in the first quadrant from $(1,0)$ to $(0,0)$ such that $n^{-1/2}\cZ^n$ with appropriate time normalization converges in the scaling limit to $Z$. Letting $(\D,h,\eta)$ denote SLE-decorated quantum disk associated with $Z$, we get a natural embedding $\wt\phi:\cV(M)\to\ol{\D}$ of $M$ by setting $\wt\phi(v)=\eta(\frk t_v)$, where $\frk t_v$ is the time of $\cZ^n$ (appropriately rescaled) associated with the vertex $v$ in the bijection from Proposition \ref{prop:mullin-wired-disk2}.

A key step of the proof is to argue that simple random walk on $M$ embedded with $\wt\phi$ converges to 2D Brownian motion modulo reparametrization of time. Once we have this result, along with some consistency and boundary conditions, it follows from a general result (see Lemmas \ref{lem:rw-conv-implies-tutte-conv1} and \ref{lem:rw-conv-implies-tutte-conv2}, which are close to ideas in \cite{gms-tutte}) that $\wt\phi$ and the Tutte embedding $\phi$ are uniformly close. This concludes the proof since it follows from mating-of-trees theory that Theorem \ref{thm-tutte-conv0} holds with $\wt\phi$ instead of $\phi$. We refer to Section \ref{sec:tutte-conv2} for details on how to deduce Theorem \ref{thm-tutte-conv0} from the invariance principle for random walk on $M$ embedded with $\wt\phi$.

We will obtain the required invariance principle in two steps. We get a first variant of the invariance principle (see Proposition \ref{prop:bm-conv-E}) by using Propositions \ref{prop:conditions} and \ref{prop:moment} along with Theorem \ref{thm-general-clt-uniform}. This invariance principle is for the SLE-decorated 0-quantum cone and the associated walk (as constructed in Theorem \ref{prop:zero-vs-gamma-walk}) and planar map. In Section \ref{sec:tutte-conv0} we transfer this invariance principle to the setting of the SLE-decorated $\gamma$-quantum cone (Proposition \ref{prop0}) and in Section \ref{sec:tutte-conv1} we transfer the invariance principle further to the desired setting of the SLE-decorated quantum disk (Lemma \ref{prop5}).

Section \ref{sec:tutte-conv0} is the hardest part of Section \ref{sec:tutte-conv}. By using Lemma \ref{lem:rooted}, which gives an absolute continuity result for the 0-quantum cone and the (recentered at another point) $\gamma$-quantum cone restricted to bounded sets, it is relatively easy to compare the continuum boundary length process $Z$ restricted to compact sets in the two settings. 
However, the associated walks $\cZ$ are hard to compare since they are determined from $Z$ in a non-local way. 
Another difficulty is that we know very little about the law of the walk $\cZ$ in the 0-quantum cone setting and to get Lemma \ref{prop5} we need to consider a walk with exactly i.i.d.\ increments. See Section \ref{sec:tutte-conv00} for a discussion of our proof strategy in Section \ref{sec:tutte-conv0}. On the other hand, the arguments in Section \ref{sec:tutte-conv1} are fairly standard and the key idea is that the 2D Brownian cone excursion and the associated walk, restricted to an interval away from the initial and terminal times, are absolutely continuous with respect to the corresponding processes in the $\gamma$-quantum cone case.

\subsection{Random walk on planar map encoded by i.i.d.\ walk converges to 2D Brownian motion}
\label{sec:tutte-conv0}

The goal of this section is to get an invariance principle for the simple random walk on a random planar map encoded by a walk with i.i.d.\ steps and coupled with a $\gamma$-quantum cone, namely Proposition \ref{prop0} below.
	Let $I\subset\R$ be a bounded and possibly random interval of positive length. For each $n\in\N$ let
	\begin{compactitem}
		\item $\cZ^n:I\to\R^2$ be a continuous-time random walk with steps $(\pm 1,0)$ and $(0,\pm 1)$  (Definition \ref{def:cont-time-rw}) and 
		\item $Z^n:I\to\R^2$ be a continuous process with values in $\R^2$ which, conditioned on $I$, is absolutely continuous with respect to the boundary length process  restricted to $I$ under one of the measures $\PQ$ or $\PE$ (Definition \ref{def:PQ-PE}) or is absolutely continuous with respect to an uncorrelated Brownian cone excursion in the first quadrant from $(1,0)$ to $(0,0)$.
	\end{compactitem}
	 Consider the SLE-decorated LQG surface encoded by $Z^n$ (see Section \ref{sec:mot}) and embed it as $(\D,h_n,\eta_n)$ in the unit disk such that $\eta_n(t_{\op{mid}})=0$ and $\eta_n$ starts at $1$, where $t_{\op{mid}}$ is the mid-point of the interval $I$. The walk $\cZ^n$ encodes a tree-decorated planar map $(M_n,\frk T_n)=\Phi_{\op{M}}(\cZ^n)$; see Section \ref{sec:bijection-disk}. Embed $M_n$ in $\ol{\BB D}$ such that a vertex $v$ is sent to $\eta_n(\frk t_v)$, where we refer to Section \ref{sec:bijection} for the notation $\frk t_v$. 

For $v\in \cV(M_n)$ let $W_n: \{0,1,\dots,T \}\to \cV(M_n)$ denote the simple random walk on $M_n$ started from $v$ and stopped at the random time $T$ at which the walk hits a vertex on the boundary of $M_n$. Let $\wt W_n: [0,T]\to\ol{\D}$ be the embedding of $W_n$, which is defined by $\wt W_n(t):=\eta_n(\frk t_{W_n(t)})$ for $t\in\{0,1,\dots,T \}$ and by linear interpolation otherwise. Let $\sigma^{v,n}$ denote the law of $\wt W_n$ conditioned on $(\cZ^n,Z^n)$. Let $\sigma^z$ denote the law of Brownian motion started from $z$ run until it first exits $\D$. We view $\sigma^{v,n}$ and $\sigma^z$ as probability measures on the space of curves in $\ol{\BB D}$ viewed modulo reparametrization of time (see \eqref{eq:curves-mod-reparam}), where $\sigma^z$ is a deterministic probability measure while $\sigma^{v,n}$ is a random probability measure determined by $(\cZ^n,Z^n)$. 

The terminology in the following definition will be used frequently throughout this section and Section \ref{sec:tutte-conv1}. Roughly speaking, it says that simple random walk on the planar map encoded by $\cZ^n$ and embedded in $\ol{\D}$ as in the previous paragraph, converges to a 2D Brownian motion in a quenched sense. The function $d_{\op{Prok}}$ denotes the Prokhorov metric.
\begin{definition}
	We say that {\emph{random walk converges to Brownian motion on $(\cZ^n,Z^n)_{n\in\N}$}} if 
	\eqb
	\sup_{v\in\cV(M_n)} d_{\op{Prok}}\big(\sigma^{v,n},\sigma^{\eta_n(\frk t_{v})}\big)
	\label{eq:rw-conv}
	\eqe
	converges to 0 in probability as $n\rta\infty$.
	\label{def:rw-conv}
\end{definition}

\begin{prop}
	For each $n\in\N$ there is a probability measure on pairs $(\cZ^n,Z^n)$ such that the following hold.
	\begin{compactitem}
		\item[(i)] $\cZ^n=(\cZ^n_t)_{t\in[0,1)}$ is a rate $2^n$ continuous-time random walk with values in $\R^2$ such that $\cZ^n_0=0$ and, conditioned on the jump times, the jumps are uniformly and independently distributed on $\{(1,0), (-1,0),(0,1),(0,-1)\}$.
		\item[(ii)] $Z^n=(Z^n_t)_{t\in[0,1)}$ has the law of a variance $\frac 12$ uncorrelated 2D Brownian motion restricted to $[0,1)$.
		\item[(iii)] Random walk converges to Brownian motion on $(\cZ^n,Z^n)_{n\in\N}$ (see Definition \ref{def:rw-conv}).
		\item[(iv)] $\|Z^n-2^{-n/2}\cZ^n\|_\infty$ converges in probability to 0 as $n\rta\infty$.
	\end{compactitem}
	\label{prop0}
\end{prop}
Using Propositions \ref{prop:conditions} and \ref{prop:moment} along with Theorem \ref{thm-general-clt-uniform} we immediately get a variant of Proposition \ref{prop0} where the processes $\cZ^n$ and $Z^n$ have a different law. Namely, the law is equal to the law of $(Z,\cZ)$ under the measure $\PE$, after applying a certain scaling and restriction of the processes; see Proposition \ref{prop:bm-conv-E} below. We will take this result as a starting point and then show that the invariance principle in (iii) still holds as we modify the law of $(Z^n,\cZ^n)$ step by step until we reach the desired law as described in (i) and (ii) above. The intermediary measures will be called $\PW$ and $\PI$, and the associated invariance principles are stated below in Lemmas \ref{prop:bm-conv-whQ} and \ref{prop10}, respectively, while important properties of $\PI$ are stated in Lemma \ref{prop8}. The most difficult partial results are Lemmas \ref{prop10} and \ref{prop8} and the proofs of these are deferred to Sections \ref{sec:proof-prop10} and \ref{sec:proof-prop8}, respectively. We refer to the beginning of those sections for proof outlines.

In Sections \ref{sec:kmt3} and \ref{sec:kmt-0-cone} we considered couplings of processes $Z$ and $\cZ$ such that the rate at which $\cZ$ jumps is of order 1. As you can see in Proposition \ref{prop0}, we will now consider variants of these couplings where the rate of jumps is of order $2^n$ instead. 
Let $(h,\eta,\frk U)$ be sampled from one of the measures $\PQ$ or $\PE$, 
let $Z$ be the associated boundary length process, and
let $\cZ$ be the process coupled with $Z$ as in Theorem \ref{prop-kmt} (for $\PQ$) or \ref{prop:zero-vs-gamma-walk} (for $\PE$) in the setting of \eqref{eq:p-v}. Recalling the definition \eqref{eq:X-def} of $X_J$ for intervals $J\subset\R$, define $\frk X$ by
\eqb
\frk X = \{X_{I(k)} \,:\,k\geq 2 \}.
\label{eq:X}
\eqe
Define the process $\wt Z: I(0)\to\R^2$ by 
\eqb
\wt Z=\br(Z,I(0)),
\label{eq:wtZ}
\eqe 
where we recall the definition of $\br$ from Section \ref{sec:kmt0}. For a dyadic system $\frk D$ and $n\in\N$ let $2^n\frk D$ be the dyadic system given by $\{ 2^n I\,:\,I\in\frk D \}$. For a set $\frk U=(U_I\,:\,I\in\frk D)$ indexed by a dyadic system $\frk D$ let $2^n\frk U$ be the set $\{ U^n_{I}\,:\,I\in 2^n\frk D \}$ indexed by the dyadic system $2^n\frk D$ such that $U^n_I=U_{2^{-n} I}$.
\begin{lem}
	Consider the setting above under one of the probability measures $\PQ$ or $\PE$. For each $n\in\N$ there is a process $\cZ^n:\R\to\R^2$ such that
	\eqb
	((2^{n/2}Z_{2^{-n}t})_{t\in \R},
	2^n \frk U,
	(\cZ^n_{2^{-n}t})_{t\in \R})
	\eqD
	(Z,\frk U,\cZ),
	\label{eq:def-cZn}
	\eqe
	and the process $\cZ^n$ is measurable with respect to $(Z,\frk U)$.
	
	Define
	\eqb\begin{split}
		&\wh{\cS}^n:=(\wh S^n_\ell)_{\ell\in[\dZ]}\quad\text{and}\quad
		\wh S^n_\ell
		:=
		\#\{ t\in I(1)\,:\,\Delta_t\cZ^{n}= v_\ell \},\\
		&\cS^n:=( S^n_\ell)_{\ell\in[\dZ]}\quad\text{and}\quad S^n_\ell:=\#\{ t\in I(0)\,:\,\Delta_t\cZ^{n}= v_\ell \}.
	\end{split}\eqe
	Then 
	\begin{compactitem} 
		\item[(i)]\label{item1}  
		$\wh{\cS}^n$ is measurable with respect to $(\frk X,(U_{I(j)})_{j\geq 2},(V_j)_{j\geq 2} )$, where $V_1,V_2,\dots\in\{-1,1 \}$ are the random variables used to define $\frk D$ in Section \ref{sec:kmt1}, 
		\item[(ii)]\label{item2} $\cS^n$ is measurable with respect to $(X_{I(1)},U_{I(1)},V_1,\wh{\cS}^n)$; more explicitly, $(X^1_{I(1)},\dots,X^{\dZ}_{I(1)})=\wt\cU(X_{I(1)},U_{I(1)})$ and, if the functions $\tau$ and $\wt\cU$ are as in Lemmas \ref{prop:tus} and \ref{prop:split-rand}, respectively, then $\cS^n_\ell=\tau(\wh{\cS}^n_\ell,2X^\ell_{I(1)})$ (resp.\ $\cS^n_\ell=\wh{\cS}^n_\ell-\tau(\wh{\cS}^n_\ell,2X^\ell_{I(1)})$) if $V_1=1$ (resp.\ $V_1=-1$), and 
		\item[(iii)]\label{item3} $\cZ^n|_{I(0)}$ is measurable with respect to $(\wt Z, (U_{I})_{I\in\frk D, I\subseteq I(0)}, \cS^n)$.
	\end{compactitem}
	\label{prop:def-cZn}
\end{lem}
\begin{proof}
	First observe that $((2^{n/2}Z_{2^{-n}t})_{t\in \R},
	2^n \frk U)
	\eqD
	(Z,\frk U)$. The walk $\cZ$ is determined by $(Z,\frk U)$ under both $\PQ$ and $\PE$, see Remark \ref{rmk:split-rand} and Theorem \ref{prop:zero-vs-gamma-walk}. This implies \eqref{eq:def-cZn} and that $\cZ^n$ is measurable with respect to $((2^{n/2}Z_{2^{-n}t})_{t\in \R},
	2^n \frk U)$, which is further measurable with respect to $(Z,\frk U)$. It remains to prove Assertions (i), (ii) and (iii). 
	The formula for $(X^1_{I(1)},\dots,X^{\dZ}_{I(1)})$ in Assertion (ii) follows from Assertion (iii) Lemma \ref{prop:split-rand}. 
	The rest of Assertion (ii) follows from Assertion (i) of Lemma \ref{prop:coupling-properties}. 
	Assertion (i) follows from Assertion (ii) of Lemma \ref{prop:coupling-properties} and Assertion (iii) of Lemma \ref{prop:split-rand}. 
	To get Assertion (iii), we first observe that $(\br(Y^\ell,I(0)))_{\ell\in[\dZ]}$ is measurable with respect to $(\wt Z, (U_{I})_{I\in\frk D, I\subseteq I(0)}, \cS^n)$. 
	Using this, Assertion (iii) of Lemma \ref{prop:coupling-properties}, and \eqref{eq97}, we get that $\cZ^n|_{I(0)}$ is also measurable with respect to $(\wt Z, (U_{I})_{I\in\frk D, I\subseteq I(0)}, \cS^n)$.
\end{proof}

\begin{prop}
	Let $(h,\eta,\frk U)$ be sampled from the measure $\PE$ in Definition \ref{def:PQ-PE}, let $Z$ be the associated boundary length process, and for each $n\in\N$ let $(\cZ^n)_{n\in\N}$ be as in Lemma \ref{prop:def-cZn}. Then random walk converges to Brownian motion on $(\cZ^n|_{I(0)},Z|_{I(0)})_{n\in\N}$. Furthermore, $\|(Z-2^{-n/2}\cZ^n)|_{I(0)}\|_\infty$ converges a.s.\ to 0 as $n\rta\infty$.
	\label{prop:bm-conv-E}
\end{prop}
\begin{proof}
	The last assertion of the proposition is immediate from \eqref{eq:def-cZn} and Theorem \ref{prop:zero-vs-gamma-walk}.
	By Propositions \ref{prop:conditions} and \ref{prop:moment} along with Theorem \ref{thm-general-clt-uniform}, there is a deterministic covariance matrix $\Sigma$ such that for any fixed compact set $A\subset\C$ it is a.s.\ the case that as $\ep \rta 0$, the maximum over all $z\in A$ of the Prokhorov distance between the conditional law of $\ep Y^{z/\ep}$ given $\cH^{\op{cell}}$ and the law of Brownian motion started from $z$ with covariance matrix $\Sigma$ tends to 0, where $\cH^{\op{cell}}$ is the cell configuration in Proposition \ref{prop:conditions} and $Y^z$ is the simple random walk on this cell configuration started from $z\in\C$ as in Theorem \ref{thm-general-clt-uniform}. We get that $\Sigma$ is a scalar multiple of the identity matrix by using rotational invariance in law of the cell configuration. The proposition now follows from \eqref{eq:def-cZn} along with a similar argument as in the proof of \cite[Lemma 3.5]{gms-tutte}. 
\end{proof}
For $k\in\N$ suppose now that the 0-quantum cone in $\PE$ has the $2^k$-embedding. Then define the probability measure $\PW$ on tuples $(h,\eta,\frk U)$ by 
\eqb
d\PW = \EE[\mu_h(B(0,2^{-k}))]^{-1}\mu_h(B(0,2^{-k}))\,d\PE,
\label{eq:defPW}
\eqe
where we use this notation since ``W'' is the first letter in ``weighted''. 
Our motivation for defining $\PW$ is Lemma \ref{lem:rooted}. 
The boundary length process $Z$ and the processes $\cZ^n$ as in Lemma \ref{prop:def-cZn} can still be defined under this measure by absolute continuity. Also note that the law of the Radon-Nikodym derivative (namely, $\mu_h(B(0,2^{-k}))$) does not depend on $n$. These observations give the following lemma.
\begin{lem}
	The conclusion of Proposition \ref{prop:bm-conv-E} still holds if we replace $\PE$ by $\PW$ for any fixed $k\in\N$. It also holds if we let $k$ depend on $n$ and send $k\rta\infty$ sufficiently slowly as $n\rta\infty$. In other words, if for each $n\in\N$ we sample $(h,\eta,\frk U)$ from $\P_{\op{W}}^{k(n)}$ for suitably chosen $k(n)$ diverging slowly to $\infty$, we let $Z^n$ denote the boundary length process, and we define the continuous-time random walk $\cZ^n$ as in Lemma \ref{prop:def-cZn}, then random walk converges to Brownian motion on $(\cZ^n|_{I(0)},Z^n|_{I(0)})_{n\in\N}$ and $\|(Z^n-2^{-n/2}\cZ^n)|_{I(0)}\|_\infty$ converges to 0 as $n\rta\infty$.
	\label{prop:bm-conv-whQ}
\end{lem}

The measures $\PQ$ and $\PW$ each induce a law on a tuple $(\wt Z,X_{I(1)},\frk U,\frk X)$. Let $\PI$ be the law on such tuples such that $(\wt Z,X_{I(1)},\frk U)$ has the same law as under the measure $\PQ$ while $\frk X$ is independent and has the same law as under the measure $\PW$. We use the letter ``I'' in this notation since it is the first letter of the word ``independent''.

By \eqref{eq:Z-brZ-X} and Lemma \ref{prop:split-rand} (see also Remark \ref{rmk:split-rand}), 
$\Delta_{I(0)} Z$ can be expressed as a function of $(X_{I(1)},\frk U,\frk X)$ and $\Delta_{I(1)}Y^1, \dots, \Delta_{I(1)}Y^{\dZ}$ can be expressed as a function of $(\frk U,\frk X)$. This allows us further to express  $Z|_{I(0)}$ as a function of $(\wt Z,X_{I(1)},\frk U,\frk X)$. %
Next we define $\cZ^n:I(0)\to\R^2$ 
as in (i)-(iii) of Lemma \ref{prop:def-cZn}. We conclude that we have natural processes $Z,\cZ^n:I(0)\to\R^2$ and random variables $\Delta_{I(1)}Y^1,\dots,\Delta_{I(1)}Y^{\dZ}$ under the measure $\PI$.
\begin{lem} 
The conclusion of Proposition \ref{prop:bm-conv-E} still holds if we replace $\PE$ by $\PI$ and send $k\rta\infty$ sufficiently slowly as $n\rta\infty$. 
More precisely, random walk converges to Brownian motion on $(\cZ^n|_{I(0)},Z|_{I(0)})_{n\in\N}$ and $\|(Z-2^{-n/2}\cZ^n)|_{I(0)}\|_\infty\rta 0$ under the measures $\P_{\op{I}}^{k(n)}$ if $k(n)\rta\infty$ sufficiently slowly.
Furthermore, the total variation distance between the law of $\op{sh}(Z|_{I(0)})$ under $\PI$ and $\PQ$ goes to 0 as $k\rta\infty$, and $Y^1|_{I(1)},\dots,Y^{\dZ}|_{I(1)}$ under $\PI$ converges in law to $Y^1|_{I(1)},\dots,Y^{\dZ}|_{I(1)}$ under $\PQ$ as $k\rta\infty$. 
\label{prop10}
\end{lem}
The proof is deferred to Section \ref{sec:proof-prop10}.
\begin{lem}
	The total variation distance between law of $\cZ^n$ under the measure $\PI$ and the law of $\cZ^n|_{I(0)}$ under the measure $\PQ$ goes to zero if $k\rta\infty$ sufficiently slowly as $n\rta\infty$.
	\label{prop8}
\end{lem}
The proof is deferred to Section \ref{sec:proof-prop8}.
\begin{proof}[Proof of Proposition \ref{prop0}]
	For $n\in\N$ let $k(n)$ be an increasing sequence with values in $\N$ such that $k(n)\rta\infty$ slowly as $n\rta\infty$ and the conclusion of Lemma \ref{prop10} is satisfied for this sequence. For any fixed $n\in\N$ consider the probability space associated with $\PI$ for $k=k(n)$. Enlarge the probability space by adding new processes $\wh Z: I(0)\to\R^2$ and $\wh{\cZ}^n:I(0)\to\R^2$ which satisfy the following: 
	\begin{compactitem}
		\item $\wh Z=Z$ except on an event with probability going to 0 as $n\rta\infty$, and the process $\wh Z$ has the law of a variance $\frac 12$ two-sided 2D Brownian motion restricted to $I(0)$,
		\item $\wh{\cZ}^n=\cZ^n$ except on an event with probability going to 0 as $n\rta\infty$, and $\wh{\cZ}^n$ is a rate $2^n$ continuous-time random walk such that, conditioned on the jump times, the jumps are uniformly and independently distributed on $\{(1,0),(-1,0),(0,1),(0,-1) \}$.
	\end{compactitem}
It is possible to define such processes by Lemma \ref{prop10} (for $\wh Z$) and Lemma \ref{prop8} (for $\wh{\cZ}^n$). Since $(\wh{\cZ}^n,\wh Z)$ is close in total variation distance to $(\cZ^n,Z)$, it follows from the first assertion of Lemma \ref{prop10} that random walk converges to Brownian motion on $(\wh{\cZ}^n,\wh Z)_{n\in\N}$ (see Definition \ref{def:rw-conv}). We conclude by using that $(\wh{\cZ}^n,\wh Z)$ has exactly the same law as the processes appearing in the statement of Proposition \ref{prop0} after pre-composition with a translation map (so the processes are defined on $[0,1)$) and adding a random constant (so the processes take the value 0 at 0), and these operations do not change the properties described in Assertions (iii) and (iv).

\end{proof}

\subsubsection{Discussion of proof in Section \ref{sec:tutte-conv0}}
\label{sec:tutte-conv00}
In this section we discuss our proof strategy for deriving Proposition \ref{prop0} from Proposition \ref{prop:bm-conv-E}.  
In order to prove Proposition \ref{prop0}, we need to understand random walk on a planar map encoded by a walk $\cZ^n$ with exactly i.i.d.\ increments, while we a priori have very little information about the law of the (non-i.i.d.) random walk $\cZ^n$ from Proposition \ref{prop:bm-conv-E}. 
In Section \ref{sec:kmt-0-cone} we derived several general comparison results for the 0-quantum cone and the $\gamma$-quantum cone, see e.g.\ Theorem \ref{prop:zero-vs-gamma-walk} and Corollary \ref{cor:zero-vs-gamma}. 
One could have hoped to use these results to transfer our invariance principle for simple random walk from the setting of the 0-quantum cone (Proposition \ref{prop:bm-conv-E}) to the setting of $\gamma$-quantum cone (Proposition \ref{prop0}), but it turns out that our earlier results are of little help. The cited theorem and corollary are helpful for bounding the probability of very unlikely events and the estimates are quite crude due to the large Radon-Nikodym derivative term on the right side. In order to transfer the invariance principle we need to understand the exact asymptotic law of random walk on certain embedded planar maps, which cannot be understood only by estimating the probability of very unlikely events. Furthermore, the results from Section \ref{sec:kmt-0-cone} which say that $\cZ$ is ``almost'' locally determined (see e.g.\ \eqref{eq:same2}) are not sufficient; we would have needed some variant of \eqref{eq:same2} where the probability on the left side goes to 0 as $m\rta\infty$ while $n$ is of constant order, and this result is not true. 

Let us describe the strategy we follow in order to prove Proposition \ref{prop0}. For any $J\in\frk D$ it holds that $\cZ|_J$ is measurable with respect to $\br(Z,J)$ (recall \eqref{eq:br}), $\frk U$, and the number of jumps $S(\cY^{\ell},J)$ that $\cZ|_J$ makes of type $v_\ell\in\{(\pm 1,0), (0,\pm 1) \}$ (recall Lemma \ref{prop43}); see Assertion (iii) of Lemma \ref{prop:def-cZn}. By using absolute continuity (after appropriate recentering) of $\br(Z,J)$ for the $\gamma$-quantum cone and the 0-quantum cone, we get a variant of Proposition \ref{prop0} under a probability measure $\PW$ where $\br(Z,J)$ has (almost) the same law as in Proposition \ref{prop0}; see Lemma \ref{prop:bm-conv-whQ}. A major difficulty in deriving Proposition \ref{prop0} from Lemma \ref{prop:bm-conv-whQ} is that the random variables $\br(Z,J)$ and $\{ S(\cY^{\ell},J)\,:\,\ell\in[\dZ]\}$ are \emph{not} independent in the context of Proposition \ref{prop0} or \ref{prop:bm-conv-whQ}, which makes it hard to understand the law of $\cZ|_J$ in these lemmas. This motivates us to consider a probability measure $\PI$ where these two random variables are exactly independent for some suitably chosen $J\in\frk D$, namely $J=I(1)$; see Lemma \ref{prop10}. We prove in Section \ref{sec:proof-prop10} that our random walk invariance principle also holds under $\PI$. The argument is quite involved, and we consider various conditional laws for the GFF, along with the geometric interpretation of the random variables $X_J$ in terms of LQG lengths, in order to argue that the law of $\{ S(\cY^{\ell},I(1))\,:\,\ell\in[\dZ]\}$ given $\br(Z,I(1))$ under $\PW$ is not too different from its unconditional law; see the beginning of Section \ref{sec:proof-prop10} for a more detailed outline. Then we argue in Section \ref{sec:proof-prop8} that the walk $\cZ^n:I(0)\to\Z^2$ under the measure $\PI$ is close in total variation distance to a walk with exactly i.i.d.\ increments. When proving this, we use the precise law of various random variables associated with $(Z,\cZ)$ under $\PQ$.

\subsubsection{Proof of Lemma \ref{prop10}}
\label{sec:proof-prop10}

Let us first give an outline of the proof of Lemma \ref{prop10}. 
Random walk converges to Brownian motion on $(\cZ^n|_{I(0)}, Z|_{I(0)})$ under $P_{\op{W}}^{k(n)}$ by Lemma \ref{prop:bm-conv-whQ} 
and 
$(\cZ^n|_{I(0)},Z|_{I(0)})$ (at least up to a shift) is measurable with respect to  $(\sh(Z|_{I(0)}), X_{I(1)}, \frk X, \frk U)$ by Lemma \ref{prop:tuple}.
Therefore, in order to establish convergence of random walk to Brownian motion on $(\cZ^n|_{I(0)},Z|_{I(0)})$ under $P_I^{k(n)}$, it is sufficient to prove that $(\sh(Z|_{I(0)}), X_{I(1)}, \frk X, \frk U)$ under $\PI$ is absolutely continuous on a high probability event with respect to the law of the associated tuple under $\PW$. 
It is not hard to establish such absolute continuity for $(\sh(Z|_{I(0)}), X_{I(1)}, \frk U)$ due to Lemma 4.5, but adding $\frk X$ to the tuple requires significant work. By the definition of $\PI$ and due to the absolute continuity result for $(\sh(Z|_{I(0)}), X_{I(1)}, \frk U)$, in order to conclude we need to prove that the conditional law of $\frk X$ given $(\sh(Z|_{I(0)}), X_{I(1)}, \frk U)$ on a high probability event is absolutely continuous with respect to its law given (only) $\Delta_{I(1)} Z$ under the measure $\PW$; the reason we condition on $\Delta_{I(1)} Z$ here is that under $\PI$, the random variables $(\sh(Z|_{I(0)}), X_{I(1)}, \frk U)$ and $\frk X$ are independent conditioned on the value of $\Delta_{I(1)} Z$.

In Lemma \ref{lem:far} we prove the desired absolute continuity result with $\frk X_{m+1,\infty}$ instead of $\frk X$, where $\frk X_{m+1,\infty}$ is defined in \eqref{eq95} and is measurable with respect to
$(h,h^{\op{IG}})$ 
restricted to $\C\setminus \eta(I(m))$. We prove this lemma by considering the conditional law of $(h,h^{\op{IG}})$ given its value restricted to 
$\eta(I(1))$ (Lemma 6.12) along with a general result on the law of the total mass of a Liouville measure (Lemma 6.13). 
In order to strengthen the absolute continuity result to concern $\frk X$ instead of $\frk X_{m+1,\infty}$, a natural proof strategy would be to prove that conditioned on $\frk U$ and $(h,h^{\op{IG}})|_{\eta(I(1))\cup (\C\setminus \eta(I(m)))}$  (which determine $(\sh(Z|_{I(0)}), X_{I(1)}, \frk U, \frk X_{m+1,\infty})$), the conditional law of $\frk X_{2,m}$ on a high probability event is absolutely continuous with respect to its law conditioned on $\Delta_{I(1)}Z$ only. However, this result is not true; the four boundary lengths 
$\partial^{\pm,\op{L}}_{I(m)},\partial^{\pm,\op{R}}_{I(m)}$ (see \eqref{eq96}), which are measurable with respect to $(h,h^{\op{IG}})$ restricted to $\C\setminus \eta(I(m))$, puts constraints on what kind of values 
$\frk X_{2,m}$ can achieve; by a recursive argument using the geometric interpretation of $\frk X$ in terms of LQG boundary lengths, one can check that the smaller the four considered boundary lengths are, the smaller is the set of possible values $\frk X_{2,m}$.  
This motivates us to introduce the event $\cE(m,C_1)$, which says that all four boundary lengths $\partial^{\pm,\op{L}}_{I(m)},\partial^{\pm,\op{R}}_{I(m)}$ are larger than $C_1$. In Proposition \ref{prop36} we prove that $\cE(m,C_1)$ happens with uniformly positive probability given a particular $\sigma$-algebra $\cG_0$ such that in particular $\frk X_{m+1,\infty}$ is measurable with respect to $\cG_0$. Once Proposition \ref{prop36} has been proven, we can conclude the proof rather quickly, but proving the proposition requires significant work.

One technical difficulty in the proof of Proposition \ref{prop36} is that we would like to consider the conditional law of $h$ and $h^{\op{IG}}$ given their restriction to sets which are not local sets \cite{ss-contour}. In some cases we work out the explicit conditional law (see e.g.\ Lemma \ref{prop31}(i)), while in other cases slightly enlarging the non-local set allows us to understand the explicit conditional law (which motivates the introduction of the parameter $s>0$ in \eqref{eq89}). Via an involved geometric argument, we work out explicitly the set of achievable values for $\partial^{\pm,\op{L}}_{I(m)},\partial^{\pm,\op{R}}_{I(m)}$, conditioned on $\frk X_{m+1,\infty}$ and the LQG boundary length of certain flow line segments.

If $z:[t_1,t_2)\to\R^2$ for $t_1<t_2$ let $\op{sh}(z)$ denote the function from $[0,t_2-t_1)$ to $\R^2$ given by $t\mapsto z(t+t_1)-z(t_1)$. In other words, $\op{sh}(z)$ is equal to $z$ composed with a translation, plus some constant which ensures that $\op{sh}(z)$ takes the value 0 at 0.
\begin{lemma}
	In the setting of Definition \ref{def:rw-conv},
	random walk converges to Brownian motion on $(\cZ^n,Z^n)_{n\in\N}$ if and only if random walk converges to Brownian motion on $(\op{sh}(\cZ^n),\op{sh}(Z^n))_{n\in\N}$.
	\label{prop:sh}
\end{lemma}
\begin{proof}
	The lemma is an immediate consequence of the fact that $Z^n$ encodes the same SLE-decorated LQG surface as $\op{sh}(Z^n)$ and that $\cZ^n$ encodes the same tree-decorated planar map as $\op{sh}(\cZ^n)$.
\end{proof}
Several places below we work with $\sh(Z|_{I(1)})$ and $\sh(Z|_{I(0)})$ instead of $Z|_{I(1)}$ and $Z|_{I(0)}$. The reason for this is Assertion (i) of the following lemma, which does not hold if we consider $Z|_{I(1)}$ due to different  behavior close to 0 under the two measures. This different behavior is due to the fact that $h$ has no log singularity at the origin under $\PW$ while it has a $\gamma$-log singularity at the origin under $\PQ$.

\begin{lem}
	\begin{compactitem}
		\item[(i)] The total variation distance between $\op{sh}(Z|_{I(1)})$ under the measures $\PW$ and $\PQ$ converges to 0 as $k\rta\infty$.
		\item[(ii)]  $(\Delta_{I(1)}Y^1,\dots,\Delta_{I(1)}Y^{\dZ})$ under $\PI$ converges in law to $(\Delta_{I(1)}Y^1,\dots,\Delta_{I(1)}Y^{\dZ})$ under $\PQ$ as $k\rta\infty$.
		\item[(iii)] The total variation distance between $Z|_{I(0)}$ under the measures $\PI$ and $\PQ$ converges to 0 as $k\rta\infty$.
	\end{compactitem}
	\label{prop20}
\end{lem}
\begin{proof}
	Sample $z_0$ from the LQG area measure restricted to  $B(0,2^{-k})$. By Lemmas \ref{lem:rooted} and \ref{lem:trans-inv}, $(h(\cdot+z_0)|_{B(0,2^k-2^{-k})},h^{\op{IG}}(\cdot+z_0)|_{B(0,2^k-2^{-k})},\frk U-\tau_{z_0} )$ under the measure $\PW$ is equal in law to $(h|_{B(0,2^k-2^{-k})},h^{\op{IG}}|_{B(0,2^k-2^{-k})},\frk U )$ under the measure $\PQ$ if we consider the $2^k$-embedding of the cone. Observe that $\op{sh}(Z|_{I(1)})$ is measurable with respect to these tuples on the following event $A$
	$$
	A=\{B(0,2^{-k})\subset\eta(I(0))\} \cap \{\eta(I(1))\subset B(0,2^{k-1})\}.
	$$
	By using that the probability of $A$ under both measures converges to 1 as $k\rta\infty$ we get Assertion (i) of the lemma.
	
	To get Assertion (ii) we use that by \eqref{eq:Z-brZ-X}, for any $m\in\N$ (to be chosen below) and for $\ell=1,\dots,\dZ$,
	$$
	\Delta_{I(1)}Y^\ell = \frk S^\ell_1+\frk S^\ell_2,
	\quad\op{where}\quad
	\frk S^\ell_1=\sum_{j=2}^m  2^{1-j/2}V_j X^\ell_{I(j)},\quad
	\frk S^\ell_2=\sum_{j=m+1}^\infty  2^{1-j/2}V_j X^\ell_{I(j)}.
	$$
	The random variables $\frk S^\ell_1$ and $\frk S^\ell_2$ are measurable with respect to $\frk X$ and $V_2,V_3,\dots$, which have the same joint law under $\PI$ and $\PW$, so we can work with $\PW$ instead of $\PI$ when establishing Assertion (ii). In order to conclude the proof it is sufficient to argue that for any $\eps>0$ we can choose $m$ sufficiently large such that

	\begin{compactitem}
		\item[(a)] the total variation distance between $(\frk S^1_1,\dots,\frk S^{\dZ}_1)$ under $\PW$ and $\PQ$ is smaller than $\eps$ for all $k$ sufficiently large, and
		\item[(b)] for $\ell=1,\dots,\dZ$ it holds that $\PW[|\frk S^\ell_2|>\eps]<\eps$ for all sufficiently large $k\in\N$ and $\PQ[|\frk S^\ell_2|>\eps]<\eps$. 
	\end{compactitem}
	We will first choose $m$ sufficiently large such that (b) holds for all  sufficiently large $k$ and then we argue that (a) holds with this choice of $m$ for all sufficiently large $k$. 
	
	It is immediate that for all sufficiently large $m$ we have $\PQ[|\frk S^\ell_2|>\eps]<\eps$. In order to establish (b) we only need to prove that we can find $m$ such that $\PW[|\frk S^\ell_2|>\eps]<\eps$ for all  sufficiently large $k\in\N$. We have
	$$
	\frk S^\ell_2 = \frk S^\ell_3 + \frk S^\ell_4,
	\quad\text{where}\quad
	\frk S^\ell_3=\sum_{j=m+1}^{k}  2^{1-j/2}V_j X^\ell_{I(j)},\quad
	\frk S^\ell_4=\sum_{j=k+1}^{\infty}  2^{1-j/2}V_j X^\ell_{I(j)}.
	$$
	By Lemmas \ref{lem:rooted} and \ref{prop:contained1}, and since, under the $2^k$-embedding, the field of the $\gamma$-quantum cone restricted to $B(0,2^k)$ stochastically domainates the field of the 0-quantum cone restricted to $B(0,2^k)$, it holds that 
	 $\PW[\eta(I(k))\subset B(0,2^{k-1})]\rta 1$ as $k\rta\infty$. Using this, similar considerations as in the proof of Assertion (i), and $X_{I(j)} \in\sigma(Z|_{I(k)})$ for $j=m+1,\dots,k$, the total variation distance between $\frk S^\ell_3$ under $\PW$ and $\PQ$ goes to 0 as $k\rta\infty$. 
	 It follows that by choosing $m$ sufficiently large, for all sufficiently large $k$,
	\eqb
	\PW[|\frk S^\ell_3|>\eps/2] < 
	\PQ[|\frk S^\ell_3|>\eps/2] + \eps/3
	<\eps/2,
	\label{eq:frkS3}
	\eqe 
	where we first pick $m$ sufficiently large such that the second inequality holds and we then observe that the first inequality holds for all sufficiently large $k$. 
	
	In order to bound $\frk S^\ell_4$, 
	by first using the definition of $\PW$ along with Lemma \ref{lem:rooted} and H\"older's inequality, then using Lemma \ref{prop:zero-vs-gamma0}, and finally using that $X_{I(j)}^\ell$ are independent and centered normal random variables of variance 0.25 under $\PQ$, for any $n\geq k$ and constants $C,C_0>1$, 
	\eqbn
	\begin{split}
		&\PW[|X^\ell_{I(j)}|> 0.1\eps 2^{0.1(j+k)} ] 
		\leq  
		2^{Ck} \PE[|X^\ell_{I(j)}|>0.1\eps 2^{0.1(j+k)}]^{0.1}\\
		&\quad\leq 
		2^{Ck} (2^{C_0n} \PQ[|X^\ell_{I(j)}|>0.1\eps 2^{0.1(j+k)}] + 2^{-n/C_0})^{0.1}
		\leq 2^{Ck} (2^{C_0n} e^{-C^{-1}\eps 2^{0.2(j+k)}} + 2^{-n/C_0})^{0.1}.
	\end{split}
	\eqen
	Letting $n=\lceil 2^{0.2(j+k)} \rceil$ and possibly increasing $C$, we get that the right side is smaller than $2^{-C^{-1}\eps 2^{0.2(j+k)}}$. By the definition of $\frk S^\ell_4$ and a union bound it follows that for sufficiently large $k$,

	\eqbn
	\begin{split}
		\PW[|\frk S^\ell_4|>\eps/2]
		&\leq 
		\PW\Big[ |\frk S^\ell_4| > 2 \sum_{j=k}^\infty 0.1\eps 2^{0.1k-0.4j} \Big] 
		< 
		\sum_{j=k}^{\infty}
		\PW[ |X_{I(j)}^\ell| > 0.1\eps 2^{0.1(j+k)}  ]\\   
		&< 
		\sum_{j=k}^{\infty} 2^{-C^{-1}\eps 2^{0.2(j+k)}}
		< \eps/2,
	\end{split}
	\eqen
	which implies (b) combined with \eqref{eq:frkS3}. Then (a) follows by similar considerations as for $\frk S^\ell_3$ above, in particular we apply Lemma \ref{lem:rooted}.

	To prove Assertion (iii), recall that $Z|_{I(0)}$ is equal to $\wt Z$ plus a linear function chosen such that $Z|_{I(0)}$ takes the value 0 at 0 and such that $Z|_{I(0)}$ has the correct increment $\Delta_{I(0)}Z$ during $I(0)$. 
	It is sufficient to prove the following in order to conclude: 
	(a) $\wt Z$ is equal in law under the two measures, 
	(b) $\wt Z$ and $\Delta_{I(0)}Z$ are independent under both measures, and
	(c) $\Delta_{I(0)}Z$ are close in total variation distance under the two measures. 
	The four random variables 
	\eqb
	\wt Z,\qquad 
	V_1,\qquad 
	X_{I(1)},\qquad
	(\Delta_{I(1)}Y^1,\dots,\Delta_{I(1)}Y^{\dZ})
	\label{eq58}
	\eqe 
	are independent under both $\PI$ and $\PQ$, and $(\wt Z,V_1,X_{I(1)})$ is equal in law under the two measures by the definition of $\PI$, which immediately implies (a).  Furthermore, (c) follows from Lemma \ref{prop:sum-tv} below, Assertion (ii), and
	$$
	\Delta_{I(0)}Z = A+B,
	\quad\text{where}\quad 
	A=2^{-1/2} V_1 X_{I(1)},\quad 
	B=\frac 12\Delta_{I(1)}Z=\frac 12 \sum_{\ell=1}^{\dZ} v_\ell\Delta_{I(1)}Y^\ell.
	$$
	Finally, we deduce (b) from the last display and independence of the four random variables in \eqref{eq58}.
\end{proof}

\begin{lem}
	Let $d\in\N$ and let 
	$(A_k,B_k)_{k\in\N}$ be a sequence of random variables such that 
	$A_k$ and $B_k$ take values in $\R^d$ for each $k$, 
	$A_k$ and $B_k$ are independent for each $k\in\N$, 
	$A_k\eqD \wt A$ for some normal vector $\wt A\sim\cN(\mu_1,\Gamma_1)$ and each $k\in\N$ with $\Gamma_1$ positive definite, and 
	$B_k$ converges in law to a normal vector $\wt B\sim\cN(\mu_2,\Gamma_2)$ as $k\rta\infty$. 
	Then $A_k+B_k$ converges in total variation distance to the normal vector $\wt A+\wt B\sim\cN(\mu_1+\mu_2,\Gamma_1+\Gamma_2)$, where $\wt A$ and $\wt B$ are independent.
	\label{prop:sum-tv}
\end{lem}
Recall the definition of $\frk X$ from \eqref{eq:X} and for $m_1,m_2\in\Z$ set 
\eqb
\frk X_{m_1,m_2} = (X_{I(m_1)},X_{I(m_1+1)},\dots,X_{I(m_2)}),\qquad
\frk X_{m_1,\infty} = (X_{I(m_1)},X_{I(m_1+1)},\dots).
\label{eq95}
\eqe
Observe that $\frk X$ can be obtained by concatenating $\frk X_{2,m-1}$, $X_{I(m)}$ and $\frk X_{m+1,\infty}$.

\begin{lem}
	Let $k,m\in\N$.  
	Working with the probability measure $\PW$, let $\cG$ be the $\sigma$-algebra generated by $(h|_{\eta(I(1))}, h^{\op{IG}}|_{\eta(I(1))})$, and let $E_{k,m},\wt E_{k,m}$ be events of $\PW$-probability $1-o_{k}(1)-o_m(1)$ such that $E_{k,m}$ (resp.\ $\wt E_{k,m}$) is measurable with respect to $\Delta_{I(1)}Z$ (resp.\ $\cG$). Consider the following two ways of sampling $\frk X_{m+1,\infty}$:
	\begin{compactitem}
		\item[(a)] Sample $\frk X_{m+1,\infty}$ from $\PW$ conditioned on $\Delta_{I(1)}Z$ and on the event $E_{k,m}$.
		\item[(b)] Sample $\frk X_{m+1,\infty}$ from $\PW$ conditioned on $\cG$ and on the event $\wt E_{k,m}$.
	\end{compactitem}
	It is possible to choose $E_{k,m},\wt E_{k,m}$ such that the law of $\frk X_{m+1,\infty}$ in (a) conditioned on an event of probability $1-o_{k}(1)-o_m(1)$ is absolutely continuous with respect to the law in (b). 
	\label{lem:far}
\end{lem}

\begin{lem}
	In the setting of Lemma \ref{lem:far} let $\check E_{k,m}$ be an event of $\PW$-probability $1-o_{k}(1)-o_m(1)$ such that $\check E_{k,m}$ is measurable with respect to $\cG$. Consider the following two ways of sampling $(h|_{\BB C \setminus B(0,\frk R(2^{0.5m}))}, h^{\op{IG}}|_{\BB C \setminus B(0,\frk R(2^{0.5m}))})$:
	\begin{compactitem}
		\item[(a)] Sample $(h|_{\BB C \setminus B(0,\frk R(2^{0.5m}))}, h^{\op{IG}}|_{\BB C \setminus B(0,\frk R(2^{0.5m}))})$ from $\PW$ conditioned on  $\Delta_{I(1)} Z$.
		\item[(b)] Sample $(h|_{\BB C \setminus B(0,\frk R(2^{0.5m}))}, h^{\op{IG}}|_{\BB C \setminus B(0,\frk R(2^{0.5m})) })$ from $\PW$ conditioned on $\cG$ and on the event $\check E_{k,m}$.
	\end{compactitem}
	It is possible to choose $\check E_{k,m}$ such that the law in (a) conditioned on an event with probability $1-o_{k,m}(1)$ is absolutely continuous with respect to the law in (b). 
	\label{lem:far0}
\end{lem}

\begin{proof}
	First we will give a partial description of the conditional law of $(h,h^{\op{IG}})$ given $\cG$ and assuming the $\cG$-measurable probability $1-o_k(1)$ event $\wt E=\{ B(0,2^{-k})\subset \eta(I(1)) \}$ occurs. Conditioned on $\cG$, the law of $h^{\op{IG}}$ in the complement of $\eta(I(1))$ is the harmonic extension of the field, plus an independent Dirichlet GFF. The conditional law of $h$ given $\cG$ is less explicit, but one can sample from the law as follows, first observing that we may work under the measure $\PE$ instead of the measure $\PW$ since the Radon-Nikodym derivative of $\PW$ with respect to $\PE$ is $\cG$-measurable on the event $\wt E$: 
	\begin{compactitem}
		\item[(i)] Sample $h$ restricted to $B(0,r)$ conditioned on $\cG$, where $r:=\op{diam}(\eta(I(1)))$. In particular, this determines $A_{-\log r}$, where $A_t:=h_{e^{-t}}(0)$ is defined as in Section  \ref{sec:lqg}. 
		\item[(ii)] Since we work with the $2^k$-embedding (not the circle average embedding) here, we know that $(A_{t})_{t\geq -k\log 2}$ has the law of a standard Brownian motion started from 0. We sample $(A_{t})_{t\in[-k\log 2,-\log r ]}$ conditioned on its value at $t=-\log r$ from the previous bullet point.
		\item[(iii)] Then sample $(A_t)_{t\leq -k\log 2}$ from its marginal law under $\PE$, i.e., independently of $\cG$ and of the random variables sampled in (i) and (ii).
		\item[(iv)] Finally sample $h-h_{|\cdot|}(0)$, which is independent of $A$ and has the same law as the analogous process for a whole-plane GFF. This means in particular that $h-h_{|\cdot|}(0)$ is given by the harmonic extension of $h-h_{|\cdot|}(0)$ from $B(0,r)$ to $\C\setminus B(0,r)$, plus an independent Gaussian field in $\C\setminus B(0,r)$ with 0 boundary conditions.
	\end{compactitem}

	Observe that $h|_{\C\setminus B(0,2r)}$ given $h|_{B(0,r)}$ is mutually absolutely continuous with respect to its unconditional law. 
	Also observe that the event $\frk R(2^{0.5m})>2r$ is measurable with respect to $h|_{\C\setminus B(0,2r)}$ and that this event has probability $1-o_m(1)$. Combining these observations gives the lemma.
\end{proof}

Let $D\subset\C$ be a domain, let $\alpha\in(0,2)$, and let $\frk m$ be a finite and non-trivial Radon measure on $D$ of dimension strictly larger than $\alpha^2/2$, where we refer to e.g.\ \cite[equation (1.3)]{berestycki-gmt-elementary} for the definition of the dimension of $\frk m$. By \cite{berestycki-gmt-elementary}, it is possible to define a measure $\wt\mu$ given by $e^{\alpha h}\,d\frk m$ for $h$ some variant of a Gaussian free field in $D$. Typically $\frk m$ is either Lebesgue area measure on $D$ or the natural measure (defined via Minkowski content) on an SLE-type curve. If $\frk m$ is Lebesgue area measure in $D$ and $\alpha=\gamma$, then $\mu$ is simply the Liouville area measure $\mu_h$ from Section \ref{sec:lqg}. If $\frk m$ is the natural measure on an SLE$_\kappa$-type curve and $\alpha=\gamma/2$, then $\wt\mu$ is the LQG length measure on the curve.  
\begin{lem}
	In the context above, suppose $h$ is the sum of a zero boundary GFF and a deterministic continuous function $\frk h:D\to\R$ such that $\wt\mu(D)<\infty$ a.s. Then the random variable $\wt\mu(D)$ is an absolutely continuous random variable such that its probability density function is strictly positive on $(0,\infty)$. 
	\label{prop32}
\end{lem}
Below we denote by $\rho_{\frk h}:(0,\infty)\to[0,\infty)$ the probability density function of the Liouville area measure $\mu(D)$ if $h$ is the sum of $\frk h$ and a zero-boundary GFF. In other words, $\rho_{\frk h}$ is the probability density function discussed in the lemma when $\frk m$ is Lebesgue area measure and $\alpha=\gamma$.
\begin{proof}
	For $\eps>0$ let $D_\eps=\{z\in D\,:\,\op{dist}(z,\C\setminus D)>\eps \}\cap B(0,\eps^{-1})$ be the set of points in $D$ that have distance at least $\eps$ from $\C\setminus D$ which are contained in $B(0,\eps^{-1})$. By choosing $\eps$ sufficiently small we may suppose $\frk m(D_\eps)>0$. Let $\phi_0:D\to\R$ be a smooth and compactly supported function which is constant in $D_\eps$ and has Dirichlet energy 1. Let $\phi_0,\phi_1,...$ be a smooth orthonormal basis for the Hilbert space closure for the Dirichlet inner product of the set of smooth compactly supported functions in $D$ and let $\beta_0,\beta_1,\dots$ be i.i.d.\ standard normal random variables such that $h=\frk h+\beta_0 \phi_0+\beta_1 \phi_1+\dots$. We have
	$\wt\mu_{h}(D)=f(\beta_0)$, where $f:\R\to\R$ is the random function defined by 
	 $f(\beta):=\int e^{\alpha\beta \phi_0}\,d\check\mu$ and $\check\mu$ is the measure $e^{\alpha \check h}d\frk m$ with $\check h:=\frk h+\beta_1\phi_1+\beta_2\phi_2+\dots$. We have $f'(\beta)=\int e^{\alpha\beta \phi_0}\alpha\phi_0(z)\,d\check\mu(z)$. By non-triviality of $\wt\mu|_{D_\eps}$ (see e.g.\ \cite{ghss-negative-moments}), $f^{-1}:(f(-\infty),\infty)\to\R$ is continuously differentiable with a strictly positive derivative. Letting $\rho$ denote the standard normal probability density function, the probability density of $\wt\mu(D)$ conditioned on $\check h$ is equal to $\rho(f^{-1}(x))(f^{-1})'(x)$ for $x\in(f(-\infty),\infty)$ and 0 otherwise. We conclude by averaging over the randomness of $\check h$, sending $\eps\rta 0$, and using that $\rho(x)\rta 0$ as $x\rta\pm\infty$.
	
\end{proof}

\begin{figure}
	\begin{center}
		\includegraphics{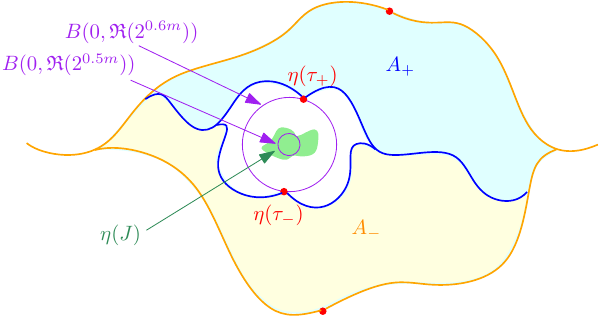}
		\caption{Illustration of the statement and proof of Lemma \ref{prop:far2}.}
	\end{center}
\label{fig-AA}
\end{figure}
See Figure \ref{fig-AA} for an illustration of the statement and proof of the next lemma.
\begin{lem}
	Define 
	\eqbn
	\begin{split}
	\tau_+&=\inf\{t\geq 0\,:\,\eta((t,\infty))\cap B(0,\frk R(2^{0.6m}) )=\emptyset \},\\
	\tau_- &=\sup\{t\leq 0\,:\,\eta((-\infty,t))\cap B(0,\frk R(2^{0.6m}) )=\emptyset \}.	
	\end{split}
	\eqen 
	Let $\wt{\cG}$ be the $\sigma$-algebra generated by $\frk U$, $h^{\op{IG}}$, $h|_{\C\setminus B(0,\frk R(2^{0.6m}) )}$ and $h|_{\eta(I(1))}$. 
	Then the sets 
	\eqbn
	\begin{split}
	A_+&:=\eta\Big(\Big[
	\tau_+, 
	\tau_+ + x_m^+ - x_1^+ -  
	\mu\big(\eta([0,\tau_+]) \setminus B(0,\frk R(2^{0.6m}) ) \big) 
	\Big]\Big) \subset\C,\\
	A_-&:=\eta\Big(\Big[
	\tau_- + x_m^- - x_1^- +  
	\mu\Big(\eta([\tau_-,0]) \setminus B(0,\frk R(2^{0.6m}) ) \Big), 
	\tau_-
	\Big]\Big)\subset\C	
	\end{split}
	\eqen
	are measurable with respect to $\wt{\cG}$, where we set $\eta([t_1,t_2])=\emptyset$ if $t_2<t_1$. There is an event $E$ measurable with respect to $\wt{\cG}$ and of $\PW$-probability $1-o_{k}(1)-o_m(1)$ such that the following two measures are a.s.\ mutually absolutely continuous under $\PW$:
	\begin{compactitem}
		\item[(a)] the law of $(\eta(x_m^-),\eta(x_m^+))$ conditioned on $\wt\cG$ and the event $E$, restricted to the event of probability $1-o_{k}(1)-o_m(1)$ that $(\eta(x_m^-),\eta(x_m^+))\in A_-\times A_+$, 
		\item[(b)] $\mu_h|_{A_- } \times \mu_h|_{A_+ }$ conditioned on $\wt{\cG}$ and the event $E$.
	\end{compactitem}

	The same conclusion holds if we instead let $\wt{\cG}$ be the $\sigma$-algebra generated by $\frk U$, $h^{\op{IG}}$, $h|_{\C\setminus B(0,\frk R(2^{0.6m}) )}$, i.e., $h|_{\eta(I(1))}$ is not measurable with respect to $\wt{\cG}$.
	\label{prop:far2}
\end{lem}

\begin{proof}
	We work on the event $B(0,2^{-k})\subset\eta(I(1))$ of probability $1-o_k(1)$ throughout the proof. This allows us to work with the probability measure $\PE$ instead of $\PW$ since the Radon-Nikodym derivative relating the two measures is determined by $h|_{B(0,2^{-k})}$.
	
	Observe that the set $\eta((-\infty,\tau_-])$ and $\eta(\tau_-)$ are measurable with respect to $\wt{\cG}$ since it $\eta(\tau_-)$ is the point at which $\eta$ hits $B(0,\frk R(2^{0.6m}))$. Similarly, $\eta([\tau_+,\infty))$ and $\eta(\tau_+)$ are also $\wt{\cG}$-measurable. This means that the two LQG areas appearing in the definition of $A_\pm$ are $\wt{\cG}$-measurable. By using that $x_m^\pm$ and $x_1^\pm$ are also $\wt{\cG}$-measurable, we conclude that $A_\pm$ are $\wt{\cG}$-measurable.
	
	Let $J$ be the shortest interval such that $B(0,\frk R(2^{0.5m}))\subset\eta(J)$. Let $E$ denote the probability $1-o_{m}(1)$ event that 
	\eqbn
	\begin{split}
		E = \{\eta(I(1))\subset B(0,\frk R(2^{0.5m}) ) \} 
		\cap 
		\{  \eta(J)\subset B(0,\frk R(2^{0.6m})) \} \cap \{A_+\neq\emptyset \}\cap \{A_-\neq\emptyset \}.
	\end{split}
	\eqen
	The set $\eta([0,\tau_+])$ is measurable with respect to $\wt{\cG}$. Furthermore, on the event $\{x_m^+>\tau_+\}$ of probability $1-o_{m}(1)$, the point $\eta(x_m^+)$ is measurable with respect to $\wt{\cG}$ and the LQG area $\mu(\eta([0,\tau_+]))$ since $\eta(x_m^+)$ is uniquely characterized by the property that $\mu(\eta([0,x_m^+]))=x_m^+$. Similarly, on the event $\{x_m^-<\tau_-\}$ of probability $1-o_{m}(1)$, the point $\eta(x_m^-)$ is measurable with respect to $\wt{\cG}$ and the LQG area $\mu(\eta([\tau_-,0]))$. By Lemma \ref{prop32}, the law of $\big(\mu(\eta([\tau_-,x_1^-])),\mu(\eta([x_1^+,\tau_+])) \big)$ conditioned on $\wt{\cG}$ is a.s.\ mutually absolutely continuous with respect to Lebesgue measure on $(0,\infty)^2$. This implies a.s.\ mutual absolute continuity of the measures in (a) and (b).

	The proof of the last assertion of the lemma is similar.
\end{proof}

\begin{proof}[Proof of Lemma \ref{lem:far}]
	Recall that $\frk X_{m+1,\infty}$ is defined by considering the LQG boundary length of a various flow line segments for $h^{\op{IG}}$ which are contained in the set $\eta((-\infty,x_m^-])\cup \eta([x_m^+,\infty))$. In particular, if we are given $h$ and $h^{\op{IG}}$ restricted to the set $\eta((-\infty,x_m^-])\cup \eta([x_m^+,\infty))$ we can determine $\frk X_{m+1,\infty}$. With probability $1-o_{m}(1)$, the set $\eta((-\infty,x_m^-])\cup \eta([x_m^+,\infty))$ has a positive distance from $B(0,\frk R(2^{0.6m}) )$. By using these observations, we see that there is a measurable function $G$ such that with probability $1-o_{m}(1)$ we have $\frk X_{m+1,\infty}=G(
	h^{\op{IG}}|_{\C\setminus B(0,\frk R(2^{0.5m}) )}, 
	h|_{\C\setminus B(0,\frk R(2^{0.6m}) )},
	\eta(x_m^-), \eta(x_m^+)
	)$.
	We can sample $\frk X_{m+1,\infty}$ in the following steps: (i) Sample $h^{\op{IG}}|_{\C\setminus B(0,\frk R(2^{0.5m}) )}, 
	h|_{\C\setminus B(0,\frk R(2^{0.6m}))}$, (ii) sample $\eta(x_m^-), \eta(x_m^+)$, (iii) define $\frk X_{m+1,\infty}$ via the function $G$. We get the lemma by first using Lemma \ref{lem:far0} to compare the law in (i) in the settings (a) and (b), and then using Lemma \ref{prop:far2} to compare the law in (ii) (conditional on (i)) in the settings (a) and (b). 
\end{proof}

Our next goal is to establish Proposition \ref{prop36} below. Before stating it, we need to introduce some notation. For a possibly random interval $I\subset\R$ with left (resp.\ right) end-point given by $t_-\in\R$ (resp.\ $t_+\in\R$), define the following four LQG boundary lengths: 
\eqb
\partial_I^{\pm,\op{L}}=\nu_h\Big(\eta_{t_\pm}^{\op{L}}\cap\eta(I)\Big),\qquad
\partial_I^{\pm,\op{R}}=\nu_h\Big(\eta_{t_\pm}^{\op{R}}\cap\eta(I)\Big).
\label{eq96}
\eqe
Notice that these segments are pairwise disjoint if we exclude the end-points and that the union of the four boundary segments equals the boundary of $\eta(I)$. For $C_1>0$ and $m\in\Z$ let $\cE(m,C_1)$ denote the event that the LQG-length of all the four segments for the interval $I(m)$ is larger than $C_1$, i.e.,
$$
\cE(m,C_1) = \{ 
\partial_{I(m)}^{\pm,\op{L}}>C_1, 
\partial_{I(m)}^{\pm,\op{R}}>C_1
\}.
$$
We let $o_{k\gg m}(1)$ denote a quantity which converges to 0 as $k,m\rta\infty$ assuming that $k\gg m$, i.e., there exists some diverging function $g$ such that $m\leq g(k)$.
\begin{prop}
	For $m\in\{3,4,\dots \}$ let $\cG_0$ be the $\sigma$-algebra generated by $(h,h^{\op{IG}})|_{\eta(I(1))}$ and $\frk X_{m+1,\infty}$. Then under the probability measure $\PW$ and for any fixed $C_1>1$, we have $\PW[\cE(m,C_1)\,|\,\cG_0]>0$ with probability $1-o_{k\gg m}(1)$.
	\label{prop36}
\end{prop}

In the next lemma we describe the conditional law of the projection of $h,h^{\op{IG}}$ onto $H(\op{int}(\eta(I)))$ given the projection of the same distributions onto $H^\perp(\op{int}(\eta(I)))$, where $I$ is an interval and $\op{int}(A)$ denotes the interior of $A$ for any $A\subset\C$. In Assertion (a) of the lemma we sample a GFF conditioned on inducing a certain LQG area. Let us describe how to make precise sense of this conditioning. Let $U$ be a bounded domain, let $a>0$, let $\frk h$ be a continuous function in $U$, let $h$ be a zero-boundary GFF in $U$ independent of $\frk h$, and let $\rho_{\frk h}:(0,\infty)\to(0,\infty)$ be the density of $\mu_{\frk h}(U)$ as defined below Lemma \ref{prop32}. The law of $h$ conditioned on the probability zero event that $\mu_{h+\frk h}(U)=a$ can be sampled as follows. Let $\phi_0,\dots,\phi_n$ be an orthonormal basis for $H(U)$ such that $h=\sum_{k=0}^\infty\beta_k\phi_k$ for i.i.d.\ standard normal random variables $\beta_k$. Let $f(\beta):=\int e^{\gamma\beta\phi_0}\,d\check\mu$, where $\check\mu$ is the Liouville area measure induced by $\check h=\frk h+\beta_1\phi_1+\beta_2\phi_2+\dots$. Conditioning on $\check h$, let $\rho_{\check h}$ denote the probability density function of $\mu_h(U)$ given $\check h$.  
By Bayes' rule, we may sample $h$ conditioned on the probability zero event that $\mu_{h+\frk h}(U)=a$ by first sampling $\check h$ from its marginal law reweighted by $\rho_{\check h}(a)/\rho_{\frk h}(a)$, and then choosing $\beta_0\in\R$ such that $\mu_{h+\frk h}(U)=a$.

\begin{lem}
	Let $(h,h^{\op{IG}})$ be sampled from $\PW$ and let $I$ be a bounded interval not containing 0 which is independent of $(h,h^{\op{IG}})$. Let $\cG_0$ be the $\sigma$-algebra generated by $I$, the projection of $h$ and $h^{\op{IG}}$ onto $H^\perp(\op{int}(\eta(I)) )$. Conditioned on $\cG_0$ and  the event that $\eta(I)\subset B(0,2^k)\setminus B(0,2^{-k})$, the projections of $h$ and  $h^{\op{IG}}$ onto $H(\op{int}(\eta(I)) )$ independently have the following conditional laws:
	\begin{compactitem}
		\item[(a)] The projection of $h$ onto $H(\op{int}(\eta(I)))$ has the law of a zero boundary GFF in $\eta(I)$ conditioned on the event that $\mu_h(\eta(I))=|I|$.
		\item[(b)] The projection of $h^{\op{IG}}$ onto $H(\op{int}(\eta(I)))$ has the law of a zero boundary GFF in $\eta(I)$.
	\end{compactitem}
\label{prop31}
\end{lem}
\begin{proof}
	Assertion (b) follows from the Markov property of imaginary geometry flow lines (see e.g.\ \cite[Theorem 1.1]{ig4}) and independence of $h$ and $h^{\op{IG}}$. To prove Assertion (a) and get conditional independence of the random variables in (a) and (b), we define $\cG_1$ to be the $\sigma$-algebra generated by $\cG_0$ and $h^{\op{IG}}$ and we want to prove Assertion (a) conditioned on $\cG_1$ and the event that $\eta(I)\subset B(0,2^k)\setminus B(0,2^{-k})$. Conditioning on the latter event allows us to assume that the $h$ is a whole-plane GFF with circle average 0 on $B(0,2^k)$. 
	Let $t_-$ (resp.\ $t_+$) be equal to the left (resp.\ right) end-point of $I$, and let $\delta>0$. It is sufficient to argue that the projection of $h$ onto $H(\op{int}(\eta(I)))$ restricted to $\op{int}(\eta(I))\setminus \ol{B(\eta(t_+),\delta)}$ has the desired law.

	We can write $\R=I_-\cup I\cup I_+$ for $I_-, I, I_+$ pairwise disjoint such that $I_-$ (resp.\ $I_+$) are to the left (resp.\ right) of $I$. We assume $0\in I_-$; the other case $0\in I_+$ can be treated similarly. The set  $K:={\ol{\eta(I_-)}\cup B(0,2^{k})^c} $ is a local set for $h$, so given the projection of $h$ onto $H^\perp( \C\setminus K )$, the conditional law of the projection of $h$ onto $H( \C\setminus K )$ is a zero-boundary GFF. We cannot argue in quite the same way if we replace $K$ by $\ol{\eta(I_-)}\cup\ol{\eta(I_+)}$ since $\mu_{h}(\eta(I))=|I|$, which implies that $\ol{\eta(I_-)}\cup\ol{\eta(I_+)}$ is not a local set. To get around this issue, let $\eps>0$ and let $S$ be an $\eps\times\eps$ box intersecting $B(0,2^k)\setminus\ol{\eta(I_-)}$ which is equal to a complementary component of the grid $\eps\Z^2$. Consider the union of $\ol{\eta(I_-)}\cup B(0,2^k)^c$ and the flow lines started from a point inside $S$, and let $U$ be the complementary component of this set with $\eta(t_-)$ on its boundary. Then $\C\setminus U$ is a local set for $h$, so conditioning on the $\sigma$-algebra $\wt{\cG}_\eps$ generated by $h^{\op{IG}}$ and on the projection of $h$ onto $H^\perp(U)$, the conditional law of the projection of $h$ onto $H(U)$ is a zero boundary GFF in $U$. 
	If we condition on the positive probability event that $\eta(t_+)\in S$, the projection of $h$ onto $H(U)$ will have the law of a zero boundary GFF conditioned on the event $\mu_h(U)\in A$, where $A$ is the $\wt{\cG}_\eps$-measurable set of values of $\mu_h(U)$ which give $\eta(t_+)\in S$. Letting $E$ denote any event which is measurable with respect to the projection of $h$ onto $H(\op{int}(\eta(I)))$ restricted to $\op{int}(\eta(I))\setminus \ol{B(\eta(t_+),\delta)}$, it will hold that $\P[E\,|\,\wt{\cG}_\eps]$ converges to the associated probability for a zero-boundary GFF in $\eta(I)$ conditioned on the event that $\mu_h(\eta(I))=|I|$. On the other hand, $\P[E\,|\,\wt{\cG}_\eps]\rta \P[E\,|\,\cG_1]$ by backward martingale convergence since  $\wt{\cG}_\eps\downarrow\cG_1$ and Assertion (a) follows.

\end{proof}

\begin{figure}
	\centering
	\includegraphics{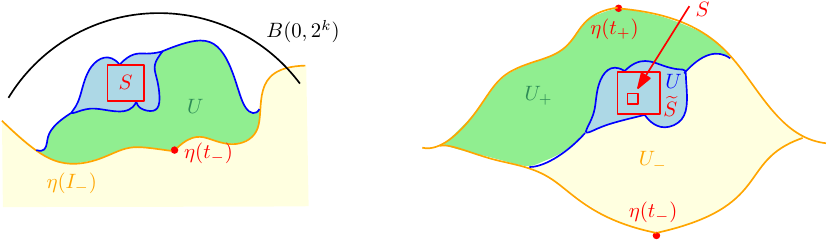}
	\caption{Illustration of the proofs of Lemmas \ref{prop31} (left) and \ref{prop33} (right). In the left (resp.\ right) figure the blue domain and the blue curves represent the east- and west-going flow lines started from a point in the box $S$ (resp.\ $\wt S$) and the points separated from $\infty$ by such flow lines.}
\end{figure}

\begin{lem}
	Consider the setting of Lemma \ref{prop31} and let $t_-$ (resp.\ $t_+$) denote the left (resp.\ right) end-point of the interval $I$. Let $b$ be a possibly random number in the interior of $I$ which is independent of $(h,h^{\op{IG}})$. 
	Let $\cG_0$ be the $\sigma$-algebra generated by $I$, $b$, and the projection of $h$ and $h^{\op{IG}}$ onto $H^\perp(\op{int}(\eta(I)) )$. 
	Suppose the $\cG_0$-measurable event $\eta(I)\subset B(0,2^k)\setminus B(0,2^{-k})$ occurs. Then the conditional law of the random variable $(\Delta_{[t_-,b]} L,\Delta_{[t_-,b]} R)$ given $\cG_0$ is a.s.\ absolutely continuous with respect to Lebesgue measure on $\R^2$ and its probability density function is a.s.\ strictly positive a.e.\ on the set $(-\partial^{-,\op{L}}_{I},\infty)\times
	 (-\partial^{-,\op{R}}_{I},\infty)$ and zero otherwise.
	\label{prop33}
\end{lem}
\begin{proof}
	We work under the same conditioning as in the statement of the lemma throughout the proof although we might not explicitly state it, and in part of the proof we will also condition on additional random variables and events. 
	A first attempt to prove the lemma would be to consider the east-going and the west-going flow line started from $\eta(b)$ and then apply the second assertion of Lemma \ref{prop32}. This strategy does not quite work since the conditional law of $h$ given $\cG_1$ and the point $\eta(b)$ is not that of a GFF in $\eta(I)$; in fact, this conditional law will be singular with respect to a GFF since it will have $\gamma$-log singularity at the point $\eta(b)$. 
	
	We modify the strategy as follows. For $\eps>0$ let $S$ be an $\eps\times\eps$ box in the grid $\eps\Z^2$ and for $s>2\eps$ let $\wt S$ be the $s\times s$ box in the grid $s\Z^2$ containing $S$. Let $K\subset\C$ be the union of all east- and west-going flow lines started from a point in $\wt S$, in addition to points separated from infinity by such flow lines. Then $K$ can a.s.\ be written as the union of two infinite curves $\eta_{\op{L}},\eta_{\op{R}}$ and $\ol{U}$ for some simply connected domain $U\subset\C$. Let $U_-$ (resp.\ $U_+$) be the complementary component of $\eta(I)\setminus K$ which has $\eta(t_-)$ (resp.\ $\eta(t_+)$) on its boundary. We can sample $h,h^{\op{IG}}$ projected onto $H(D)$ with $D=\op{int}(\eta(I))$ in the following steps, where we get the law of the various fields and sets from Lemma \ref{prop31}:
	\begin{compactitem}
		\item[(a)] Sample $U$ and the projection of $h$ and $h^{\op{IG}}$ onto $H^\perp(D\setminus \ol U)$ according to the law in Lemma \ref{prop31}.
		\item[(b)] Sample $\eta_{\op{L}}$ and $\eta_{\op{R}}$ and the projection of $h,h^{\op{IG}}$ onto $H^\perp(D\setminus K)$.
		\item[(c)] Let the projection of $h$ onto $H(U_-)$ and $H(U_+)$ be given by independent zero-boundary GFFs conditioned on the event that $\mu_h(D)=|I|$.
	\end{compactitem}

	Now let $S$ and $\wt S$ denote the random boxes as above such that $\eta(b)\in S\subset\wt S$. Let $\cG_{s,\eps}$ denote the $\sigma$-algebra generated by $S,\wt S$ and the randomness in (a), where if $\eps=0$ we condition on the point $\eta(b)$. We will argue that conditioned on $\cG_{s,0}$, with probability $1-o_s(1)$ and for appropriate $\delta_s\rta 0$, the conditional law of the random variable $(\Delta_{[t_-,b]} L,\Delta_{[t_-,b]} R)$ is a.s.\ absolutely continuous with respect to Lebesgue measure on $\R^2$ and its probability density function is a.s.\ strictly positive a.e.\ on the set $(-\partial^{-,\op{L}}_{I}+\delta_s, \infty) \times
		(-\partial^{-,\op{R}}_{I} + \delta_s, \infty)$ and zero a.e.\ on the set $\R^2 \setminus (-\partial^{-,\op{L}}_{I},\infty) \times
		(-\partial^{-,\op{R}}_{I},\infty)$. This is sufficient to complete the proof of the lemma upon sending $s\rta 0$.
		
	Conditioned on $\cG_{s,\eps}$ for $\eps>0$, the random variables in (b) and (c) are absolutely continuous with respect to their laws conditioned only on the randomness in (a), and the Radon-Nikodym derivative can be expressed in terms of the areas $\mu_h(U_-)$ and $\mu_h(U_+)$ since we need to guarantee that $\eta(b)\in S$. When $\eps\rta 0$ these conditional laws converge such that in the limit we condition on $\mu_h(U_-)$ and $\mu_h(U_+)$ to take explicit $\cG_{s,0}$-measurable values. 
	In the remaining argument, we work under this limiting law where we are conditioning on $\eta(b)$ and the random variables in (a). By Lemma \ref{prop32}, the LQG lengths of
	$\eta_{\op{L}}\cap D,\eta_{\op{R}}\cap D$ are absolutely continuous with respect to Lebesgue measure on $\R^2$ with a strictly positive density on $(0,\infty)^2$. Notice that Lemma \ref{prop32} is only stated for a single GMC measure, while we consider the joint law of two GMC measures here, but our conclusion still holds since we can first apply the lemma with $\eta_{\op{L}}\cap D$ and then (conditioned on the LQG length of $\eta_{\op{L}}\cap D$) with $\eta_{\op{R}}\cap D$. 
	
	This implies our desired result (as stated at the end of the previous paragraph) since $\Delta_{[t_-,b]} L$ is equal to $\partial^{+,\op{L}}_{I([t_-,b])}$ 
	(which is greater than or equal to the LQG length of $\eta_{\op{L}}\cap D$ plus a $\cG_{s,0}$-measurable  random variable converging in law to 0 as $s\rta 0$) 
	minus the LQG length  
	$\partial^{-,\op{L}}_{I([t_-,b])}$
	(which is between 0 and $\partial^{-,\op{L}}_{I}$), and the same holds with R instead of L. 
	
\end{proof}

For each $\ell\in\R$ and $q\in\N$ define 
$$
\frk L(\ell,q) = \{ (d_1,\dots,d_q)\in\R^q\,:\, d_j>-(\ell+d_1+\dots+d_{j-1}) \,\text{for}\, j=1,\dots,q \}.
$$ 
For $\ell,r\in\R$ and $q\in\N$ define
\eqb
\frk L(\ell,r,q) = \{ (d^{\op{L}}_1,\dots,d^{\op{L}}_q, 
 d^{\op{R}}_1,\dots,d^{\op{R}}_q)
 \in\R^{2q}\,:\,
(d^{\op{L}}_1,\dots,d^{\op{L}}_q)\in \frk L(\ell,q),
(d^{\op{R}}_1,\dots,d^{\op{R}}_q)\in \frk L(r,q)
 \}.
 \label{eq:frkL}
\eqe
Observe that $\frk L(\ell,q)$ is monotone in $\ell$ in the sense that $\frk L(\ell,q)\subsetneq\frk L(\wt\ell,q)$ if $\ell<\wt\ell$. Similarly, $\frk L(\ell,r,q)\subseteq \frk L(\wt\ell,\wt r,q)$ if $\ell\leq \wt\ell$ and $r\leq \wt r$. 
\begin{lem}
	Consider the setting of Lemma \ref{prop33}, except that instead of $b\in I$ we have $b_1<\dots<b_q$ in the interior of $I$ for some $q\in\N$, such that  $(b_1,\dots,b_q)$ is independent of $(h,h^{\op{IG}})$.
	Condition on $I$, $(b_1,\dots,b_q)$, the projection of $h$ and $h^{\op{IG}}$ onto $H^\perp(\op{int}(\eta(I)) )$, and the event that $\eta(I)\subset B(0,2^k)\setminus B(0,2^{-k})$. Then the following hold: 
	\begin{compactitem}
		\item[(a)] The conditional law of the random variable 
		$$
		(
		\Delta_{[t_-,b_1]} L, 
		\Delta_{[b_1,b_2]} L, \dots,
		\Delta_{[b_{q-1},b_q]} L, 
		\Delta_{[t_-,b_1]} R,
		\Delta_{[b_1,b_2]} R, \dots,
		\Delta_{[b_{q-1},b_q]} R)
		$$ 
		is absolutely continuous with respect to Lebesgue measure on $\R^{2q}$. Furthermore, its density is a.s.\ strictly positive a.e.\ on the set $\frk L(\partial_I^{-,\op{L}},\partial_I^{-,\op{R}},q)\subset\R^{2q}$ and zero a.e.\ on its complement.
		\item[(b)] The conditional law of the random variable 
		$$
		(
		-\Delta_{[b_q, t_+]} L,
		-\Delta_{[b_{q-1}, b_q]} L, \dots,
		-\Delta_{[b_1, b_2]} L, 
		-\Delta_{[b_q, t_+]} R,
		-\Delta_{[b_{q-1}, b_q]} R, \dots,
		-\Delta_{[b_1, b_2]} R)
		$$ 
		is absolutely continuous with respect to Lebesgue measure on $\R^{2q}$. Furthermore, its density is a.s.\ strictly positive a.e.\ on the set $\frk L(\partial_I^{+,\op{L}},\partial_I^{+,\op{R}},q)\subset\R^{2q}$ and zero a.e.\ on its complement.
	\end{compactitem}
	\label{prop34} 
\end{lem}

\begin{lem}
	For $q\in\N$ consider real-valued random variables $U^{\op{L}}_0,\dots,U^{\op{L}}_q,U^{\op{R}}_0,\dots,U^{\op{R}}_q$ on some probability space. Let $\cG_0\subset\cG_1\subset\cG_2\subset\dots$ be increasing $\sigma$-algebras such that $U^{\op{L}}_j,U^{\op{R}}_j$ are measurable with respect to $\cG_j$. Suppose that the conditional law of $(U^{\op{L}}_j, U^{\op{R}}_j)$ given $\cG_{j-1}$ is a.s.\ absolutely continuous with respect to Lebesgue measure on $\R^2$ such that its density is a.s.\ strictly positive a.e.\ on the set $(-(U^{\op{L}}_0+U^{\op{L}}_1+\dots+U^{\op{L}}_{j-1}),\infty)\times(-(U^{\op{R}}_0+U^{\op{R}}_1+\dots+U^{\op{R}}_{j-1}),\infty)\subset\R^2$ and is a.s.\ zero elsewhere.
	
	Then the conditional law of $(U^{\op{L}}_1,\dots,U^{\op{L}}_q,U^{\op{R}}_1,\dots,U^{\op{R}}_q)$ given $\cG_0$ is a.s.\ absolutely continuous with respect to $\R^{2q}$. Furthermore, the density is a.s.\ positive a.e.\ on $\frk L(U^{\op{L}}_0,U^{\op{R}}_0,q)$ and is zero a.e.\ on $\R^{2q}\setminus \frk L(U^{\op{L}}_0,U^{\op{R}}_0,q)$.
	\label{lem:abs-cont-R}
\end{lem}
\begin{proof}
	We may assume without loss of generality that $\cG_0$ is the $\sigma$-algebra generated by $\emptyset$. 
	First observe that $(U_j^{\op{L}},U_j^{\op{R}})$ is almost surely absolutely continuous with respect to Lebesgue measure on $\R^2$ conditioned on $\cG_0$ and $(U^{\op{L}}_1, \dots, U^{\op{L}}_{j-1}, U^{\op{R}}_1, \dots, U^{\op{R}}_{j-1})$. Denote the conditional density by $f_{(U_j^{\op{L}},U_j^{\op{R}})|(U^{\op{L}}_1, \dots, U^{\op{L}}_{j-1}, U^{\op{R}}_1, \dots, U^{\op{R}}_{j-1})} (\cdot | \cdot)$ and observe that the density of $(U^{\op{L}}_1,\dots,U^{\op{L}}_{q},U^{\op{R}}_1,\dots,U^{\op{R}}_{q})$ is given by 
	$$
	\prod_{j=1}^q f_{(U_j^{\op{L}},U_j^{\op{R}})|(U^{\op{L}}_1, \dots, U^{\op{L}}_{j-1}, U^{\op{R}}_1, \dots, U^{\op{R}}_{j-1})} (u_j^{\op{L}},u_j^{\op{R}} | u^{\op{L}}_1, \dots, u^{\op{L}}_{j-1}, u^{\op{R}}_1, \dots, u^{\op{R}}_{j-1}),
	$$
	which is strictly positive on $\frk L (U^{\op{L}}_0, U^{\op{R}}_0, q)$
\end{proof}
\begin{proof}[Proof of Lemma \ref{prop34}]
	We only prove Assertion (a) since the proof of Assertion (b) is identical. We apply Lemma \ref{lem:abs-cont-R} with $(U_0^{\op{L}},U_0^{\op{R}})=(\partial^{-,\op{L}}_I,\partial^{-,\op{R}}_I)$,   $U_j^{\op{L}}=\Delta_{[b_{j-1},b_{j}]}L$, $U_j^{\op{R}}=\Delta_{[b_{j-1},b_{j}]}R$, $\cG_0$ is the $\sigma$-algebra generated by $I$, $(b_1,\dots,b_q)$, the projection of $h$ and $h^{\op{IG}}$ onto $H^\perp(\op{int}(\eta(I)))$, and the event $\eta(I)\subset B(0,2^k)\setminus B(0,2^{-k})$, where $b_0=t_-$. For $j=1,\dots,q$ let $\cG_j$ be the $\sigma$-algebra generated by $\cG_0$ and the projection of $h$ and $h^{\op{IG}}$ onto $H^{\perp}(\op{int}(\eta([b_{j},t_+])))$. The assumptions of Lemma \ref{lem:abs-cont-R} are satisfied by Lemma  
	\ref{prop33}. 
\end{proof}

\begin{figure}
	\centering
	\includegraphics{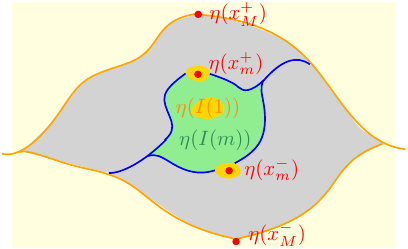}
	\caption{Illustration of the sets $K_1,K_2,K_3\subset\C$. 
	The set $K_1$ is shown in light yellow and orange, 
	the set $K_2$ is the union of $K_1$ and the four blue curves, and 
	the set $K_3$ is the union of $K_2$ and the green set.}
\end{figure}
For any $j\in\Z$ write $I(j)=[x_{j}^-,x_{j}^+)$. Let $m\in\{3,4,\dots \}$, $M\in\{m+5, m+6, \dots \}$ and $s=0.1$. Define set sets $K_1,K_2,K_3\subset\C$ by
\eqb
\begin{split}
	K_1 &= \eta(\ol{I(1)})
	\cup \eta( (-\infty,x_{M}^-] ) 
	\cup \eta( [x_{M}^+,\infty))
	\cup \eta( [x_m^+ -s,x_m^+ +s] )
	\cup \eta( [x_m^- -s,x_m^- +s] ),\\
	K_2&=K_1\cup \eta^{\op{L}}_{\eta(x_m^+)}
	\cup \eta^{\op{L}}_{\eta(x_m^-)}
	\cup \eta^{\op{R}}_{\eta(x_m^+)}
	\cup \eta^{\op{R}}_{\eta(x_m^-)},\\
	K_3 &= K_2\cup\eta(\ol{I(m)}). 
\end{split}
\label{eq89}
\eqe
Let $\cG_1$ be the $\sigma$-algebra generated by $\frk D$, 
the projection $\frk h_1$ of $h$ onto $H^\perp(\C\setminus K_1)$,  
the projection of $h^{\op{IG}}$ onto $H^\perp(\C\setminus K_3)$.  
Let $\cG_2$ be the $\sigma$-algebra generated by $\cG_1$ and the projection $\frk h_2$ of $h$ onto $H^\perp(\C\setminus K_2)$. Set $m_*=\inf\{ j\in\Z\,:\,
x_j^+=x_{M}^+ 
\text{\,or\,}
x_j^-=x_{M}^- \}$ and 
\eqb
\frk X'=(X_{I(m+1)}, X_{I(m+2)}, \dots, X_{I(m_*-1)}, X_{I(m_*+1)}, \dots, X_{I(M-2)}, X_{I(M-1)}).
\label{eq:Xm}
\eqe
Let $o_{k\gg M\gg m}(1)$ denote a quantity converging to 0 as $k,M,m\rta\infty$, where we require that $k$ is sufficiently large as compared to $M$ and that $M$ is sufficiently large as compared to $m$. 
Let $\wt\cE$ be the following event, which is measurable with respect to $\cG_1$ and has probability $1-o_{k\gg M\gg m}(1)$:
\eqb
\begin{split}
	\wt\cE = &\{ \eta(I(M))\subset B(0,2^k) \}
	\cap \{B(0,2^{-k})\subset \eta(I(1)) \}
	\cap \{\op{the\,sets\,}
	\eta(\ol{I(1)}),
	\eta( (-\infty,x_{M}^-] )\cup \eta( [x_{M}^+,\infty)),\\	&\qquad\qquad\qquad\qquad\qquad\qquad
	\eta( [x_m^+ -s,x_m^+ +s] ),
	 \eta( [x_m^- -s,x_m^- +s]  \op{\,are\,pairwise\,disjoint}\}.
\end{split}
\label{eq88}
\eqe
\begin{lem}
	Working under the probability measure $\PW$, the  conditional law of $\frk X'$ given $\cG_2$ and the event $\wt{\cE}$ is absolutely continuous with respect to Lebesgue measure on $\R^{2(M-m-2)}$. Letting $\rho_2$ denote the conditional density of $\frk X'$ given $\cG_2$ and assuming $x_{m_*}^-=x_M^-$, the set $A:=\{x\in\R^{2(M-m-2)}\,:\,\rho_2(x)>0 \}$ depends only on $\frk D$ and the following random variables, each of which describe either an LQG boundary length or an LQG boundary length difference: 
	\eqb
	\begin{split}
		&\partial_{[x_m^+ + s,x_{M}^+]}^{-,\op{L}},\quad
		\partial_{[x_m^+ + s,x_{M}^+]}^{-,\op{R}},\quad
		\partial_{[x_{M}^-,x_{m}^- - s]}^{+,\op{L}},\quad
		\partial_{[x_{M}^-,x_{m}^- - s]}^{+,\op{R}},\\
		&\Delta_{[x_m^--s,x_m^++s]}Z,\quad
		\Delta_{[x_M^-,x_m^--s]}Z,\quad  \Delta_{[x_m^--s,x_m^-]}Z,\quad 
		\Delta_{[x_m^+,x_m^++s]}Z.
	\end{split}
	\label{eq:4-delta}
	\eqe
	If $x_{m_*}^+=x_M^+$ instead of $x_{m_*}^-=x_M^-$ the same statement holds with $\Delta_{[x_m^++s,x_M^+]}Z$ instead of $\Delta_{[x_M^-,x_m^--s]}Z$.
	
	Given $\cG_2$ and $\frk X'$, and if $x_{m_*}^-=x_M^-$, we can write $X_{I(m_*)}$ and $X_{I(M)}$ as explicit linear functions of 
	$\Delta_{I(m)}Z$, $\frk X', \Delta_{I(M)}Z, \Delta_{[x_M^-,x_m^-]}Z$, where the coefficients of the linear functions are determined by $\frk D$ and the coefficient of $\Delta_{[x_M^-,x_m^-]}Z$ is non-zero. If instead $x_{m_*}^+=x_M^+$, the same holds with $\Delta_{[x_m^+,x_M^+]}Z$ instead of $\Delta_{[x_M^-,x_m^-]}Z$.
	\label{prop35}
\end{lem}
\begin{proof}
	Let 
	$$
	U_+=\op{int}( \eta([x_m^+ + s, x_{M}^+]) ),\qquad
	U_-=\op{int}( \eta([x_{M}^-, x_m^- - s]) ).
	$$
	Given $\cG_2$, the conditional law of $h|_{U_+\cup U_-}$ is given by $\frk h|_{U_+\cup U_-}$ plus an independent zero-boundary GFF in each of the domains $U_+$ and $U_-$, conditioned on the events $\mu_h(U_+)=x_{M}^+-x_m^+-s$ and $\mu_h(U_-)=x_{m}^- - x_{M}^- -s$. Choose $q_+,q_-\in\N_0$ and $b^+_1<\dots<b^+_{q_+}$ and $b^-_1<\dots<b^-_{q_-}$ such that
	$$
	\{b^+_1,\dots,b^+_{q_+} \} = \{ x^+_{m+1},\dots,x^+_{M-1} \}\setminus\{  x^+_m,x^+_M\},\quad
	\{b^-_1,\dots,b^-_{q_-} \} = \{ x^-_{m+1},\dots,x^-_{M-1} \}\setminus\{x^-_m,  x^-_M\}.
	$$
	Applying Lemma \ref{prop34} twice with the intervals $[x_m^+ +s, x_{M}^+]$ and $[x_{M}^-, x_m^- - s]$, we get absolute continuity of the following tuples with respect to Lebesgue area measure on $\R^{2(q_+ + q_-)}=\R^{2(M-m-2)}$ 
	\eqb
	\big(
	\Delta_{[x_m^++s,b^+_1]} Z, \dots,
	\Delta_{[x_m^++s,b^+_{q_+}]} Z, 
	\Delta_{[b^-_1,x_m^--s]} Z, \dots,
	\Delta_{[b^-_{q_-},x_m^--s]} Z\big).
	\label{eq:4-delta-LR}
	\eqe
	Furthermore, the set at which the probability density function is positive, is a function of 
	$\partial_{[x_m^+ + s,x_{M}^+]}^{-,\op{L}},
	\partial_{[x_m^+ + s,x_{M}^+]}^{-,\op{R}},
	\partial_{[x_{M}^-,x_{m}^- - s]}^{+,\op{L}},
	\partial_{[x_{M}^-,x_{m}^- - s]}^{+,\op{R}}$. 
	We suppose in the remainder of the proof that $x_{m_*}=x_M^-$, but the case $x_{m_*}=x_M^+$ can be handled similarly by considering $\Delta_{[x_m^++s,x_M^+]}Z$ instead of $\Delta_{[x_M^-,x_m^--s]}Z$ below.
	Observe that $\frk X'$ can be expressed as linear function of \eqref{eq:4-delta-LR} and the $\cG_2$-measurable random variables $\Delta_{[x_m^--s,x_m^++s]}Z$,
	$\Delta_{[x_M^-,x_m^--s]}Z$,  $\Delta_{[x_m^--s,x_m^-]}Z$, and $\Delta_{[x_m^+,x_m^++s]}Z$ such that the coefficients of the linear function are measurable with respect to $\frk D$; one can see this by first determining $Z_{b_j^\pm}-Z_{x_m^+}$ for all $j$ from \eqref{eq:4-delta-LR} and  $\Delta_{[x_m^--s,x_m^++s]}Z$, and then further determining $\frk X'$ by also using 	$\Delta_{[x_M^-,x_m^--s]}Z$,  $\Delta_{[x_m^--s,x_m^-]}Z$, and $\Delta_{[x_m^+,x_m^++s]}Z$. 
	Moreover, the function from \eqref{eq:4-delta-LR} to $\frk X'$ is invertible conditioned on $\cG_2$ and one can check from the construction in the previous sentence that the claimed monotonicity holds. 
	It follows that the conditional law of $\frk X'$ given $\cG_2$ and the event $\wt\cE$ is absolutely continuous with respect to Lebesgue measure and that the set $A$ in the statement of the lemma depends only on $\frk D$ and \eqref{eq:4-delta}.
	 
	It remains to prove the assertion about $X_{I(m_*)}$ and $X_{I(M)}$. We will assume that $x_{m_*}^-=x_M^-$, but the case $x_{m_*}^+=x_M^+$ can be done similarly. First observe that for any $j\in\Z$, by the definition of $X_{I(j+1)}$,
	\eqb
	\Delta_{I(j+1)}Z \text{ is a linear function of } \Delta_{I(j)}Z\text{ and } X_{I(j+1)},
	\label{eq:j}
	\eqe 
	where, here and in the remainder of this paragraph, all linear functions we consider have coefficients determined by $\frk D$. By using \eqref{eq:j} iteratively for $j=m,m+1,\dots,m_*-2$, we see that each of $\Delta_{I(m+1)} Z,\dots,\Delta_{I(m_*-1)}Z$ can be expressed as a linear function of 
	\eqb
	\Delta_{I(m)}Z, X_{I(m+1)},\dots,X_{I(m_*-1)}.
	\label{eq59}
	\eqe
	The increment $\Delta_{[x_m^-,x_{m_*-1}^-]}Z$ is a linear function of the increments $\Delta_{I(m)} Z,\Delta_{I(m+1)} Z,\dots,\Delta_{I(m_*-1)}Z$. 
	Combining the last two results we see that  $\Delta_{[x_m^-,x_{m_*-1}^-]}Z$ can be expressed as a linear function of \eqref{eq59}. 
	Since $\Delta_{[x_{m_*-1}^-,x_{m}^-]}Z=\Delta_{I(m_*-1)}Z - \Delta_{[x_{m}^-,x_{m_*-1}^+]}Z$, it follows that $\Delta_{[x_{m_*-1}^-,x_{m}^-]}Z$ is a linear function of \eqref{eq59}. 
	Now using 
	$\Delta_{[x_{M}^-,x_{m_*-1}^-]}Z=\Delta_{[x_M^-,x_m^-]}Z - \Delta_{[x_{m_*-1}^-,x_{m}^-]}Z$, we get that 
	$\Delta_{[x_{M}^-,x_{m_*-1}^-]}Z$ is a linear function of
	\eqb
	\Delta_{I(m)}Z,X_{I(m+1)},\dots,X_{I(m_*-1)},\Delta_{[x_M^-,x_m^-]}Z,
	\label{eq60}
	\eqe 
	where, both here and in similar statements below, the coefficient of $\Delta_{[x_M^-,x_m^-]}Z$ is non-zero. 
	Since $x_{m_*-1}^+=x_{m_*}^+$ by the definition of $m_*$, we have $\Delta_{I(m_*)}Z=\Delta_{[x_{m_*}^-,x_{m_*}^+]}Z=\Delta_{[x_{M}^-,x_{m_*-1}^+]}Z=\Delta_{[x_{M}^-,x_{m_*-1}^-]}Z+\Delta_{I(m_*-1)}Z$, which implies that   $\Delta_{I(m_*)}Z$ is a linear function \eqref{eq60}.  
	Finally we get that $X_{I(m_*)}$ is a linear function of \eqref{eq60} by using that $X_{I(m_*)}$ is a linear function of $\Delta_{I(m_*)}Z$ and $\Delta_{I(m_*-1)}Z$.
	  
	Using \eqref{eq:j} iteratively for $j=m_*,m_*+1,\dots,M-2$, we express  $\Delta_{I(m_*+1)}Z,\dots,\Delta_{I(M-1)}Z$ as linear functions of \eqref{eq60} and $X_{I(m_*+1)},\dots,X_{I(M-1)}$. We get that $X_{I(M)}$ is a linear function of \eqref{eq60}, $X_{I(m_*+1)},\dots,X_{I(M-1)}$, and $\Delta_{I(M)}Z$ by using that $X_{I(M)}$ is a linear function of $\Delta_{I(M)}Z$ and $\Delta_{I(M-1)}Z$. 
\end{proof}

\begin{proof}[Proof of Proposition \ref{prop36}]
	Let $\cG_3$ be the $\sigma$-algebra generated by $\cG_1$ and the four LQG boundary lengths $\partial_{[x_m^+ + s,x_{M}^+]}^{-,\op{L}},
	\partial_{[x_m^+ + s,x_{M}^+]}^{-,\op{R}},
	\partial_{[x_{M}^-,x_{m}^- - s]}^{+,\op{L}},
	\partial_{[x_{M}^-,x_{m}^- - s]}^{+,\op{R}}$. 
	Notice that the random variables in \eqref{eq:4-delta} are measurable with respect to $\cG_3$. 
	Let $\cG_4$ be the $\sigma$-algebra generated by $\cG_3$ and $\frk X'$. Notice that $\frk X_{m+1,\infty}$ is measurable with respect to $\cG_4$ by the last two assertions of Lemma \ref{prop35} and since $\Delta_{I(m)}Z,\Delta_{I(M)}Z,\Delta_{[x_M^-,x_m^-]}Z$ are $\cG_4$-measurable. This implies 
	 $\cG_0\subset\cG_4$. Therefore it is sufficient to prove the lemma with $\cG_0$ replaced by $\cG_4$, i.e., we want to prove that for a suitable choice of $M$ we have $\PW[\cE(m,C_1)\,|\,\cG_4]>0$ with probability  $1-o_{k\gg m}(1)$. We condition on the probability $1-o_{k\gg M\gg m}(1)$ event $\wt{\cE}$ from \eqref{eq88} throughout the proof.
	
	Let $U_\pm\subset\C$ be the domains defined in the proof of Lemma \ref{prop35} and also define the following $\cG_1$-measurable domains
	$$
	V_+ = \op{int}(\eta([x_1^+,x^+_m-s])),\qquad
	V_- = \op{int}(\eta([x_m^-+s,x^-_1])).
	$$
	The law of $h$ given $\cG_3$ is $\frk h_1$ plus a zero boundary GFF in each of the domains $U_\pm,V_\pm$, conditioned such that $\partial_{[x_m^+ + s,x_{M}^+]}^{-,\op{L}},
	\partial_{[x_m^+ + s,x_{M}^+]}^{-,\op{R}},
	\partial_{[x_{M}^-,x_{m}^- - s]}^{+,\op{L}},
	\partial_{[x_{M}^-,x_{m}^- - s]}^{+,\op{R}}$ take the required values and
	$$
		\mu_{h} (U_+)= x_{M}^+ -x_m^+    -s,\,\,
		\mu_{h} (U_-) = x_{m}^-  -x_{M}^- -s,\,\,
		\mu_{h} (V_+)=x_{m}^+  -x_1^+    -s,\,\,
		\mu_{h} (V_-)= x_{1}^-  -x_m^-    -s.
	$$
	One can sample from this law in a similar fashion as described above Lemma \ref{prop31}. When sampling $h$ given $\cG_3$ we are given the following four LQG lengths
	\eqb
	\partial_{[x_m^+ ,x_{M}^+]}^{-,\op{L}},\quad
	\partial_{[x_m^+ ,x_{M}^+]}^{-,\op{R}},\quad
	\partial_{[x_{M}^-,x_{m}^- ]}^{+,\op{L}},\quad
	\partial_{[x_{M}^-,x_{m}^- ]}^{+,\op{R}}.
	\label{eq:4-delta-2}
	\eqe
	but we are not given the associated four LQG boundary measures, except that we know the restriction of the boundary measures to $\eta([x_m^+-s,x_m^++s])\cup \eta([x_m^--s,x_m^-+s])$. Notice that the 
	SLE-type curve $\eta_{\eta(x_m^-)}^{\op{L}}\cap \eta(I(m))$ 
	with LQG length $\partial_{I(m)}^{-,\op{L}}$ 
	is contained in the 
	SLE-type curve $\eta_{\eta(x_m^+)}^{\op{L}}\cap \eta(I([x_m^+ ,x_{M}^+]))$
	with LQG length $\partial_{[x_m^+ ,x_{M}^+]}^{-,\op{L}}$ and when we sample $h$ conditioned on $\cG_3$ the LQG length of $\eta_{\eta(x_m^-)}^{\op{L}}\cap \eta(I(m))$  may take a value arbitrarily close to $\partial_{[x_m^+ ,x_{M}^+]}^{-,\op{L}}$. A similar result holds for the other curves defining the four LQG lengths appearing in $\cE(m,C_1)$. Therefore it holds with positive conditional probability given $\cG_3$ that $\cE(m,C_1)$ occurs as long as each of the four random variables in \eqref{eq:4-delta-2} is strictly larger than $C_1$. Furthermore, it follows from Lemma \ref{prop35} and the fact that $\cG_3\subset\cG_2$ that the

	conditional law of $\frk X'$ given $\cG_3$ is absolutely continuous with respect to Lebesgue measure. Let $\rho_3$ denote the conditional density.
	
	Recall that $\cG_4$ is generated by $\cG_3$ and $\frk X'$ and notice that $\cG_2$ is generated by $\cG_3$ and $\frk h_2$. Using this and Bayes' rule, the conditional law of $\frk h_2$ given $\cG_4$ is equal to its conditional law given $\cG_3$, weighted by $\rho_2(\frk X')/\rho_3(\frk X')$, where $\rho_2$ is as in Lemma \ref{prop35}. Notice that $\rho_2(\frk X')$ is almost surely positive since Lemma \ref{prop35} says that the set $A$ in the lemma statement is measurable with respect to $\cG_3$. Therefore the conditional law of $\frk h_2$ given $\cG_4$ is mutually absolutely continuous with respect to its conditional law given $\cG_3$. By using that $\cE(m,C_1)$ is measurable with respect to $\frk h_2$, it follows that $\cE(m,C_1)$ occurs with positive conditional probability given $\cG_4$ as long as each of the four random variables in \eqref{eq:4-delta-2} is strictly larger than $C_1$. When we send $M\rta\infty$ such that $m\ll M\ll k$, the probability that the random variables in \eqref{eq:4-delta-2} are larger than $M$ goes to 1.
\end{proof}

\begin{lem}
	The law of $\frk X_{2,m-1}$ under $\PW$ conditioned on $Z|_{I(1)}$ and $\frk X_{m+1,\infty}$ is  mutually absolutely continuous with respect to Lebesgue measure on $\R^{2(m-2)}$ with probability $1-o_{k\gg m}(1)$.
	\label{lem:mid2} 
\end{lem}

\begin{proof}
	Define the probability $1-o_{k\gg m}(1)$ event $E$ by
	$$
	E=\{ B(0,2^{-k})\subset \eta(I(1))\}\cap\{\eta(I(m))\subset B(0,2^k)\}.
	$$ 

	It is sufficient to show that on the event $E$, for arbitrary $K>1$ the law of $\frk X_{2,m-1}$ under $\PW$ conditioned on $Z|_{I(1)}$, $\frk X_{m+1,\infty}$ and $\frk D$ is absolutely continuous and has a strictly positive density on $(-K,K)^{2(m-2)}$. Define 
	$m_{**}=\inf\{ j\in\Z\,:\,
	x_j^+=x_{m}^+ 
	\text{\,\,or\,\,}
	x_j^-=x_{m}^- \}$ .  Suppose $m_{**}>1$; the case $m_{**}\leq 1$ can be treated in a similar way. Then define  
	$$
	\frk X'' = (X_{I(2)}, X_{I(3)}, \dots, X_{I(m_{**}-1)}, X_{I(m_{**}+1)}, \dots, X_{I(m-1)}).
	$$ 
	Arguing exactly as in the proof of Lemma \ref{prop35}, we get that $X_{I(m_{**})}$ is an explicit linear function of $\Delta_{I(1)}Z$, $\frk X''$ and $\Delta_{[x_m^-,x_1^-]}Z$ (resp.\ $\Delta_{[x_1^+,x_m^+]}Z$) with coefficients determined by $\frk D$ if $x_{m_{**}}^-=x_{m}^-$ (resp.\ $x_{m_{**}}^+=x_{m}^+$), 
	such that the coefficient of $\Delta_{[x_m^-,x_1^-]}Z$ (resp.\ $\Delta_{[x_1^+,x_m^+]}Z$) in the linear function is non-zero. Therefore, in order to conclude the proof it is sufficient to show that for any $\wt K>1$, on the event $E$, the tuple given by the concatenation of $\frk X''$ and $\Delta_{[x_m^-,x_1^-]}Z$ (resp.\ $\Delta_{[x_1^+,x_m^+]}Z$) is absolutely continuous and has a strictly positive density on $(-\wt K,\wt K)^{2(m-2)}$ conditioned on $Z|_{I(1)}$, $\frk X_{m+1,\infty}$ and $\frk D$. In the remainder of the proof we assume $x_{m_{**}}^-=x_{m}^-$ so that we consider $\Delta_{[x_m^-,x_1^-]}Z$ instead of $\Delta_{[x_1^+,x_m^+]}Z$.
	  
	Define
	$$
	K_1=\eta(I(1))
	\cup \eta((-\infty,x_m^-])
	\cup \eta([x_m^+,\infty)),\qquad
	K_2 = K_1\cup\eta_{\eta(x_1^-)}^{\op{L}}
	\cup\eta_{\eta(x_1^-)}^{\op{R}}
	\cup\eta_{\eta(x_1^+)}^{\op{L}}
	\cup\eta_{\eta(x_1^+)}^{\op{R}}.
	$$
	Then sample $\frk X'',\Delta_{[x_m^-,x_1^-]}Z$ in the following steps, conditioned on $Z|_{I(1)}$, $\frk X_{m+1,\infty}$ and $\frk D$:
	\begin{compactitem}
		\item[(i)] Sample the projection of $h$ and $h^{\op{IG}}$ onto $H^\perp(\op{int}(K_1))$ conditioned on $Z|_{I(1)}$, $\frk X_{m+1,\infty}$ and $\frk D$. We condition on the event $E$ in the remainder of the proof; it is convenient to do this conditioning since it allows us to work with the measure $\PE$ instead of $\PW$, and since $h$ restricted to $B(0,2^k)$ has the law of a whole-plane GFF restricted to $B(0,2^k)$. By Proposition \ref{prop36}, for any $C_1>2\wt K$ the event $\cE(m,C_1)$ occurs with positive probability.

		\item[(ii)] Sample the projection of $h$ and $h^{\op{IG}}$ onto $H^\perp(\op{int}(K_2))$ conditioned on the random variables from (i). By Lemma \ref{prop32}, the law of $\Delta_{[x_m^-,x_1^-]}Z$ is absolutely continuous with respect to Lebesgue measure and has a strictly positive density on $[-\wt K,\wt K]^2$. Furthermore, this also holds if we condition on the positive probability event $\check{\cE}(m,C_1)$ that the point at which 
		$\eta^{\op{L}}_{\eta(x_1^-)}$ 
		(resp.\ $\eta^{\op{R}}_{\eta(x_1^-)}$) hits the boundary of $\eta(I(m))$ has distance of at least $C_1/2$ (as measured along the boundary of $\eta(I(m))$) from the points $\eta(x_m^-)$ and $\eta(x_m^+)$. 
		\item[(iii)] Sample $\frk X''$ conditioned on the random variables from (i) and (ii).
	\end{compactitem}

	In order to conclude the proof it is sufficient to argue that by choosing $C_1$ sufficiently large, when we sample $\frk X''$ in (iii) then the following hold
	\begin{itemize}
		\item[(a)] the law of $\frk X''$ is absolutely continuous with respect to Lebesgue measure on $\R^{2(m-3)}$, and
		\item[(b)] conditioned on the positive probability event $\cE(m,C_1)\cap \check{\cE}(m,C_1)\cap \{ \Delta_{[x_m^-,x_1^-]}Z\in[-\wt K,\wt K]^2\}$, the random variable $\frk X''$ is absolutely continuous with a strictly positive density on $[-\wt K,\wt K]^{2(m-3)}$. 
	\end{itemize}
	We will proceed similarly as in the proof of Lemma \ref{prop35}. Choose $q_+,q_-\in\N_0$ and $b^+_1<\dots<b^+_{q_+}$ and $b^-_1>\dots>b^-_{q_-}$ such that
	$$
	\{b^+_1,\dots,b^+_{q_+} \} = \{ x^+_{2},\dots,x^+_{m-1} \}\setminus\{  x^+_1,x^+_m\},\qquad
	\{b^-_1,\dots,b^-_{q_-} \} = \{ x^-_{2},\dots,x^-_{m-1} \}\setminus\{  x^-_1,x^-_m\}.
	$$
	Applying Lemma \ref{prop34} twice with the intervals $[x_1^+, x_{m}^+]$ and $[x_{m}^-, x_1^-]$, we get absolute continuity of the following tuples with respect to Lebesgue area measure on $\R^{2q_+ + 2q_-}=\R^{2(m-3)}$ 
	\eqb
	(
	\Delta_{[x_1^+,     b^+_1]}    Z,
	\Delta_{[b_1^+,     b^+_2]}    Z, \dots,
	\Delta_{[b_{q_+-1}^+,     b^+_{q_+}]}Z, 
	\Delta_{[b^-_1,     x_1^-]}          Z,
	\Delta_{[b^-_2,     b_1^-]}          Z, \dots,
	\Delta_{[b^-_{q_-}, b_{q_--1}^-]}          Z).
	\label{eq:4-delta-LRp}
	\eqe
	Furthermore, the set at which the probability density function is positive, is a function of 
	\eqb
	\partial_{[x_1^+,   x_{m}^+]}^{-,\op{L}},\quad
	\partial_{[x_1^+,   x_{m}^+]}^{-,\op{R}},\quad
	\partial_{[x_{m}^-, x_{1}^-]}^{+,\op{L}},\quad
	\partial_{[x_{m}^-, x_{1}^-]}^{+,\op{R}}.
	\label{eq:4-delta-p}
	\eqe
	The tuple $\frk X''$ can be expressed as a linear function of \eqref{eq:4-delta-LRp}, $\Delta_{[x_1^-,x_1^+]}Z$ and $\Delta_{[x_m^-,x_1^-]}Z$ such that the coefficients of the linear function are measurable with respect to $\frk D$. Conditioned on $\Delta_{[x_1^-,x_1^+]}Z$ and $\Delta_{[x_m^-,x_1^-]}Z$, the map from \eqref{eq:4-delta-LRp} to $\frk X''$ is invertible. 
	It follows that the conditional law of $\frk X''$ given  $Z|_{I(1)}$, $\frk X_{m+1,\infty}$, $\frk D$ and the randomness sampled in (i) and (ii), is absolutely continuous with respect to Lebesgue measure. 
	Furthermore, by using the monotonicity of $\frk L$ as stated below \eqref{eq:frkL} we can find $\check K>1$ such that if the four random variables  \eqref{eq:4-delta-p} are all larger than $\check K$ then the conditional density of $\frk X''$ is strictly positive on $[-\wt K,\wt K]^{2(m-3)}$. 
	We can guarantee that the four random variables \eqref{eq:4-delta-p} are larger than $\check K$ on the event in (b) for $C_1$ sufficiently large, since, for example, 
	$\partial_{[x_m^-, x_{1}^-]}^{+,\op{L}} 
	= \Delta_{[x_m^-,  x_1^-]}L 
	+ 
	\partial_{[x_m^-,  x_{1}^-]}^{-,\op{L}} 
	\geq -\wt K+C_1/2$.
\end{proof}

\begin{lem}
	Under the measure $\PI$, the conditional law of $\frk X_{2,m-1}$ given $\Delta_{I(1)}Z$ and $\frk X_{m+1,\infty}$ is absolutely continuous with respect to Lebesgue measure on $\R^{2(m-2)}$ on an event of probability $1-o_{k\gg m}(1)$.
	\label{lem:mid1}
\end{lem}
\begin{proof}
	By the definition of $\PI$ and \eqref{eq:Z-brZ-X2}, the random variables we consider have the same joint law under $\PI$ and $\PW$, so we may work under the measure $\PW$ instead. The lemma is therefore an immediate consequence of Lemma \ref{lem:mid2}.

\end{proof}

\begin{lem}
	The tuple $(\sh(Z|_{I(0)}), X_{I(1)}, \frk X, \frk U)$ under $\PI$ conditioned on an event of probability $1-o_k(1)$ is absolutely continuous with respect to $(\sh(Z|_{I(0)}), X_{I(1)}, \frk X, \frk U)$ under $\PW$.
	\label{prop:PW-PI-ac}
\end{lem}
\begin{proof}
		We can sample $(\sh(Z|_{I(0)}), X_{I(1)}, \frk X,\frk U)$ from the measure $\PW$ in the following way for any $m\in\{3,4,\dots \}$:
	\begin{compactitem}
		\item[(i)] Sample $(h|_{\eta(I(1))},\eta|_{I(1)})$ and $\frk U$ from $\PW$, which in particular determines %
		$(\br(Z|_{I(0)}), X_{I(1)},\Delta_{I(1)}Z)$. 
		\item[(ii)] Sample $\frk X_{m+1,\infty}$ from $\PW$ conditioned on $(h|_{\eta(I(1))},\eta|_{I(1)})$.
		\item[(iii)] Sample $\frk X_{2,m-1}$ from $\PW$ conditioned on $\frk X_{m+1,\infty}$ and $(h|_{\eta(I(1))},\eta|_{I(1)})$.
		\item[(iv)] Define $X_{I(m)}$ as a function of $\Delta_{I(1)}Z,\frk X_{m+1,\infty},\frk X_{2,m-1}$ via \eqref{eq:Z-brZ-X2}.
	\end{compactitem}
	Recall that under the measure $\PI$, $(\br(Z|_{I(0)}), X_{I(1)})$ has the same law as under $\PQ$ while $\frk X$ is independent and has the same law as under $\PW$. By \eqref{eq:Z-brZ-X2}, $\frk X$ is measurable with respect to $\Delta_{I(1)}Z$ and $(\frk X_{2,m-1},\frk X_{m+1,\infty})$. 
	Furthermore, $(\br(Z|_{I(0)}), X_{I(1)},\Delta_{I(1)}Z)$ has total variation distance $o_k(1)$ to its law under $\PQ$, since under $\PQ$ we have independence of the three random variables in the tuple and the law of $\Delta_{I(1)}Z$ under $\PW$ has total variation distance $o_k(1)$ its law under $\PQ$. By using these observations, sampling from the measure $\PI$ can be done in the following way, up to a $o_k(1)$ error in total variation distance
	\begin{compactitem}
		\item[(i)'] Sample $(\sh(\br(Z|_{I(0)})), X_{I(1)},\Delta_{I(1)}Z,\frk U)$ from the measure $\PQ$.
		\item[(ii)'] Sample $\frk X_{m+1,\infty}$ independently from $\PW$, except that we condition on $\frk U$ and on the value of $\Delta_{I(1)}Z$ from (i)'. 
		\item[(iii)'] Sample ${\frk X}_{2,m-1}$ from $\PW$ conditioned on $\frk U$, $\frk X_{m+1,\infty}$ and $\Delta_{I(1)}Z$. 
		\item[(iv)'] Define $X_{I(m)}$ as a function of $\Delta_{I(1)}Z,\frk X_{m+1,\infty},\frk X_{2,m-1},\frk U$ via \eqref{eq:Z-brZ-X2}.
	\end{compactitem}

	By Assertion (i) of Lemma \ref{prop20}, the tuples $(\sh\br((Z|_{I(0)})), X_{I(1)},\Delta_{I(1)}Z,\frk U)$ in (i) and (i)' have total variation distance $o_k(1)$. 
	By Lemma \ref{lem:far}, the tuple $(\sh(\br(Z|_{I(0)})), X_{I(1)}, \Delta_{I(1)}Z, \frk U, \frk X_{m+1,\infty})$ from (i)' and (ii)' conditioned on an event of probability $1-o_{k,m}(1)$ is absolutely continuous with respect to the same tuple from (i) and (ii). By Lemmas \ref{lem:mid1} and \ref{lem:mid2}, the same statement is true for the tuple $(\sh(\br(Z|_{I(0)})), X_{I(1)}, \Delta_{I(1)}Z, \frk U, \frk X_{m+1,\infty},\frk X_{2,m-1} )$, which concludes the proof since $\sh(\br(Z|_{I(0)}))$ and $\Delta_{I(1)}Z$ determine $\sh(Z|_{I(0)})$.

\end{proof}

\begin{lem}
	The pair $\big(\sh(Z|_{I(0)}),\sh(\cZ^n|_{I(0)})\big)$ is measurable with respect to $(\sh(Z|_{I(0)}), X_{I(1)}, \frk X,\frk U)$
	\label{prop:tuple} 
\end{lem}
\begin{proof}
	Lemma \ref{prop:def-cZn} implies that $\big(\sh(Z|_{I(0)}), \sh(\wt Z), \frk U, X_{I(1)}, \frk X\big)$ determines  $\big(\sh(Z|_{I(0)}),\sh(\cZ^n|_{I(0)})\big)$. This completes the proof since $\sh(Z|_{I(0)})$ determines $\sh(\wt Z)$. 
\end{proof}

\begin{proof}[Proof of Lemma \ref{prop10}]
	Let $\eps>0$ and let $A_n$ be the event that the supremum in \eqref{eq:rw-conv} is larger than $\eps$ if $(\cZ^n,Z)$ is replaced by $\big(\sh(\cZ^n|_{I(0)}),\sh(Z|_{I(0)})\big)$. By Lemmas \ref{prop:PW-PI-ac} and \ref{prop:tuple} there is an event $E_k$ of $\PI$-probability $1-o_k(1)$ such that 
	$\big(\sh(Z|_{I(0)}),\sh(\cZ^n|_{I(0)})\big)$ under $\PI$ conditioned on $E_k$ is absolutely continuous with respect to the same tuple under $\PW$. Let $\xi_k$ denote the associated Radon-Nikodym derivative. We have

	\eqb
	\PI[A_n] \leq \PI[E^c_k] + \E_{\op{W}}^k[\1_{A_n} \xi_k].
	\label{eq75}
	\eqe
	For any fixed $k$, it follows from Lemmas \ref{prop:bm-conv-whQ} that $\1_{A_n}\rta 0$ in $\PW$-probability. By sending $k=k(n)\rta\infty$ sufficiently slowly it follows from \eqref{eq75} that $\PI[A_n]\rta 0$ as $n\rta\infty$. 
	Since $\eps$ is arbitrary, it follows from Lemma \ref{prop:sh} that random walk converges to Brownian motion on $(\cZ^n|_{I(0)},Z|_{I(0)})_{n\in\N}$ under the measures $\PI$ if $k=k(n)\rta\infty$ sufficiently slowly. We get by a similar argument that  $\|(Z-2^{-n/2}\cZ^n)|_{I(0)}\|_\infty$ converges in law to 0 under $\PI$ as $n\rta\infty$ and $k=k(n)\rta\infty$ sufficiently slowly. The last assertion of the lemma follow from Assertion (ii) of Lemma \ref{prop20}.
\end{proof}

\subsubsection{Proof of Lemma \ref{prop8}}
\label{sec:proof-prop8}
By Assertion (iii) of Lemma \ref{prop:def-cZn}, and since $(\wt Z,(U_{I})_{I\in\frk D,I\subseteq I(0)} )$ has the same law under $\PI$ and $\PQ$, in order to complete the proof it is sufficient to argue the following:
	\begin{enumerate}
		\item[(i)] $(\wt Z,(U_{I})_{I\in\frk D,I\subseteq I(0)} )$ and $\cS^n$ are independent under both the measures $\PI$ and $\PQ$, where $\cS^n$ is as defined in Lemma \ref{prop:def-cZn}.
		\item[(ii)] The total variation distance between the laws of $\cS^n$ under the two measures $\PI$ and $\PQ$ goes to zero if $k\rta\infty$ sufficiently slowly as $n\rta\infty$.
	\end{enumerate}
In our proof of both (i) and (ii) we will work with the following three random variables:
$$
\theta_1 = (\frk X, (U_{I(j)})_{j\geq 2} , (V_j)_{j\geq 2} ), \quad
\theta_2 = (X_{I(1)},U_{I(1)},V_1),\quad
\theta_3 = (\wt Z, (U_{I})_{I\in\frk D,I\subseteq I(0)}),
$$
where $V_1,V_2\in \{-1,1 \},\dots$ are the random variables appearing in the definition of the dyadic system $\frk D$ in Section \ref{sec:kmt1}. 
The three random variables $\theta_1,\theta_2,\theta_3$ are independent for both $\PQ$ and $\PI$; for $\PQ$ this follows from Lemmas \ref{prop:br-indep} and \ref{prop:split-rand}, while for $\PI$ it follows by definition and from the result for $\PQ$. The law of $\theta_2$ (resp.\ $\theta_3$) is the same under the two measures, while the marginal law of $\theta_1$ (more specifically, the marginal law of $\frk X$) is different. 
It follows that $\theta_3$ and $\cS^n$ are independent since the latter random variable is a measurable function of $\theta_1$ and $\theta_2$ by Assertions (i) and (ii) of Lemma \ref{prop:def-cZn}. This gives (i) above.

It remains to prove (ii), which follows as an immediate consequence from the following claims, which will be proven below, along with Lemma \ref{prop:tv-poisson}:
\begin{compactitem}
	\item[(a)] Under $\PQ$ the random variables $\cS^n=(S^n_\ell)_{\ell\in[\dZ]}$ have the law of independent Poisson random variables with parameters $p_12^{n},\dots,p_{\dZ}2^{n}$.
	\item[(b)] With $\PI$-probability converging to 1 as $k=k(n)\rta\infty$ sufficiently slowly as $n\rta\infty$, we have $|\wh S^n_\ell-p_\ell 2^{n+1}-2^{n/2}\Delta_{I(1)}Y^\ell|<n^2$ for $\ell=1\dots,\dZ$.
	\item[(c)] The random variables $\Delta_{I(1)}Y^1,\dots,\Delta_{I(1)}Y^{\dZ}$ under $\PI$ converge jointly in law to the same random variables under $\PQ$, namely independent random variables with law $\cN(0,2p_{\ell})$ for $\ell=1,\dots,\dZ$.
	\item[(d)] Conditioned on $\wh{\cS}=(\wh S^n_\ell)_{\ell\in[\dZ]}$ it holds that $S_1,\dots,S_{\dZ}$ under $\PI$ are independent binomial random variables such that $S_\ell$ has law $B(\wh S_\ell,0.5)$ for $\ell=1,\dots,\dZ$.
\end{compactitem}

Claim (a) is immediate from the fact that the set of jumps of $\cY^\ell$ has the law of a rate $p_\ell$ Poisson point process, independently for each $\ell\in\{1,\dots,\dZ \}$. 
Claim (c) is immediate from the last assertion of Lemma \ref{prop10} and the fact that $Y^1,\dots,Y^{\dZ}$ have the law of independent variance $p_1,\dots,p_{\dZ}$ Brownian motions under $\PQ$. 
In order to prove claim (d), first observe that $\wh{\cS}^n$ and $\theta_2=(X_{I(1)},U_{I(1)},V_1)$ are independent since the former random variable is a function of $\theta_1$ by Lemma \ref{prop:def-cZn}(i). Furthermore, $X_{I(1)}$ is a two-dimensional Gaussian vector with covariance matrix given by $0.5$ times the identity matrix and $U_{I(1)}$ is a standard $\dZ$-dimensional normal vector. By using this and Lemma \ref{prop:def-cZn}(ii), conditioned on $\wh{\cS}^n$ it holds that $\cS^n$ is a collection of independent binomial $B(\wh S^n_\ell,0.5)$ random variables, where the randomness of the $B(\wh S^n_\ell,0.5)$ random variables depend on $\theta_2$. This proves claim (d).

In order to prove claim (b), first observe that $\wh S^n_\ell$ and $\Delta_{I(1)}Y^\ell$ are measurable with respect to $\theta_1$ by Assertion (i) of Lemma \ref{prop:def-cZn} and \eqref{eq:Z-brZ-X}, respectively. By using this along with the fact that $\theta_1$ has the same law under $\PI$ and $\PW$ by the definition of $\PI$ we get
\eqb
\begin{split}
	\PI[|\wh S^n_\ell-(2^{n+1}p_\ell + 2^{n/2} \Delta_{I(1)} Y^\ell )|>n^2] 
	&= \PW[|\wh S^n_\ell-(2^{n+1} p_\ell + 2^{n/2} \Delta_{I(1)} Y^\ell )|>n^2]\\
	&\leq C_k\PE[|\wh S^n_\ell-(2^{n+1}p_\ell +2^{n/2} \Delta_{I(1)} Y^\ell )|>n^2]^{0.1},
\end{split}
\eqe
where we use the definition \eqref{eq:defPW} of $\PW$ in the second step along with H\"older's inequality and Lemma \ref{lem:rooted}. We get claim (b) by using that the probability on the right side of the last display is bounded above by $2^{-n}$ for sufficiently large $n$ by \eqref{eq56} and \eqref{eq:def-cZn} since  
$$
|\wh S_\ell^n - (2^{n+1}p_\ell + 2^{n/2}\Delta_{I(n)}Y^{\ell} )|
= |\Delta_{I(n+1)}\cY^\ell - \Delta_{I(n)}Y^{\ell}| 
\leq 2 \sup_{t\in I(n+1)} |\cY^\ell_t-Y^\ell_t|.
$$

We conclude the proof of (ii) by using (a)-(d) above along with the following lemma with $\wh R_\ell=\wh S^n_\ell$ and $R_\ell=S^n_\ell$ under $\PI$ and $\wt p_\ell=p_\ell$. Claims (b) and (c) above imply that $|\wh R_\ell-P_\ell|\leq c2^{n/2}$ with probability $1-o_n(1)$.
\begin{lem}
	Let $n\in\N$, let $\wt p_1,\dots,\wt p_{\dZ}\in(0,1)$, 
	let $P_1,\dots,P_{\dZ}$ be a collection of independent  normal random variables such that $P_\ell\sim\cN(\wt p_\ell 2^{n+1},\wt p_\ell 2^{n+1})$, and let $\wh R_1,\dots,\wh R_{\dZ}$ be arbitrary random variables with value in $\N_0$ such that for any $c\in(0,1)$, $|\wh R_\ell-P_\ell|\leq c2^{n/2}$ for $\ell=1,\dots,\dZ$ with probability $1-o_n(1)$. Conditioned on these random variables, let $R_1,\dots,R_{\dZ}$ be independent binomial random variables such that $R_\ell$ is binomial with parameters $\wh R_\ell$ and 1/2. Then the total variation distance between $(R_1,\dots,R_{\dZ})$ and a vector of $\dZ$ independent Poisson random variables of parameters $\wt p_12^n,\dots, \wt p_{\dZ}2^n$, respectively, goes to 0 as $n\rta\infty$. 
	\label{prop:tv-poisson}
\end{lem}
\begin{proof}
	By Lemma \ref{prop:RN-poisson} we may assume instead that the $P_\ell$'s are independent Poisson random variables such that $P_\ell\sim\op{Pois}(\wt p_\ell 2^{n+1})$. Furthermore, if we had $\wh R_\ell=P_\ell$ instead of the approximate equality in the statement of the lemma, then the vector $(R_1,\dots,R_{\dZ})$ would have exactly the law of independent Poisson random variables of parameters $\wt p_12^n,\dots, \wt p_{\dZ}2^n$, respectively. Therefore, in order to include it is sufficient to argue that the total variation distance between $\op{Bin}(m,0.5)$ and $\op{Bin}(m+j,0.5)$ goes to zero as $m\rta\infty$ as long as $|j|m^{-1/2}\rta 0$. The latter result is immediate by the explicit formula for binomial probabilities. 

\end{proof}

\subsection{Random walk on tree-weighted disk converges to 2D Brownian motion}
\label{sec:tutte-conv1}

Our goal in this section is to prove the following lemma, which is the invariance principle that will be used in Section \ref{sec:tutte-conv2} to conclude the proof of Theorem \ref{thm-tutte-conv0}. The process $Z$ in (i) has exactly the law of the Brownian cone excursion described at the end of Section \ref{sec:sle}.
\begin{lem}
	There is a coupling of $( (\cZ^n)_{n\in\N},Z, (\frk v^n_{\op{U}})_{n\in\N} , \frk t_{\op{U}},A)$ such that 
	\begin{compactitem}
		\item[(i)] $Z:[0,A]\to\R_+^2$ has the law of a variance $\frac 12$ uncorrelated Brownian cone excursion from $(1,0)$ to $(0,0)$ of random duration $A>0$,
		\item[(ii)] $\cZ^n:[0,A]\to\R^2$ is a continuous-time random walk such that
		$(M_n,\frk T_n):=\phi_{\op{M}}(\cZ^n)$  
		has the same law as the tree-decorated planar map in Theorem \ref{thm-tutte-conv0},
		\item[(iii)] random walk converges to Brownian motion on $(\cZ^n,Z)_{n\in\N}$ (recall Definition \ref{def:rw-conv}),  	
		\item[(iv)] $\|Z-n^{-1/2}\cZ^n\|_\infty\rta 0$ in probability as $n\rta\infty$, 
		\item[(v)] %
		conditioned on $((\cZ^n)_{n\in\N},Z)$, on the event $\cV(M_n)\setminus\partial\cV(M_n)\neq\emptyset$ of probability $1-o_n(1)$, the vertex $\frk v^n_{\op U}$ is uniform on $\cV(M_n)\setminus\partial\cV(M_n)$, while on the complement of this event let $\frk v^n_{\op U}$ be the root vertex; furthermore, $\frk t_{\op{U}}$ has the law of a uniform random variable on $[0,A]$ conditioned on $((\cZ^n)_{n\in\N},Z)$, 
		\item[(vi)] $\frk t_{\frk v^n_{\op{U}}}$ converges in probability to $\frk t_{\op{U}}$, and
		\item[(vii)] the number of vertices of $\cV(M_n)\setminus\partial\cV(M_n)$ divided by $An$ converges in probability to a deterministic constant. 
	\end{compactitem}
	\label{prop5}
\end{lem}
One difference between this section and Section \ref{sec:tutte-conv0} is that in this section the processes for different $n$ are defined on the same probability space and, furthermore, the coupling is such that the continuum process $Z$ (denoted by $Z$ here and by $Z^n$ in Section \ref{sec:tutte-conv0}) is the same for all $n$. The main input to the proof of Lemma \ref{prop5} will be Proposition \ref{prop0}. We first prove a variant of this proposition, namely Lemma \ref{prop1}, where the rate of the continuous-time walk can be an arbitrary natural number, rather than a power of two, and the duration of the processes is some arbitrary fixed $a>0$. In Lemma \ref{prop2} we reweight the probability measure so that $Z|_{[\eps,a-\eps]}$ has the law of a Brownian cone excursion of duration $a$ restricted to $[\eps,a-\eps]$. We argue in the proof that under this reweighted law, the walk $\cZ^n|_{[\eps,a-\eps]}$ is close in total variation distance to a continuous-time walk on $[0,a]$ with i.i.d.\ increments conditioned to form a cone excursion and restricted to $[\eps,a-\eps]$. Then in Lemma \ref{prop12} we prove a variant of the latter lemma for cone excursions of random duration, which requires us to study the law of the number of steps of $\cZ^n$ in this setting. Finally we can conclude the proof of Lemma \ref{prop5} by sending $\eps\rta 0$ and argue that the invariance principle holds on the full planar map, not just the part of the planar map encoded by $\cZ^n|_{[\eps,a-\eps]}$ for $\eps>0$.
\begin{lem}
	Let $a>0$.
	There is a coupling of $Z,\cZ^1,\cZ^2,\cZ^3,\dots$ such that 
	\begin{compactitem}
		\item[(i)] $(Z_t)_{t\in[0,a]}$ has the law of a variance $\frac 12$ uncorrelated 2D Brownian motion restricted to $[0,a]$, 
		\item[(ii)] for each $n\in\N$ the process $(\cZ^n_t)_{t\in[0,a]}$ is a continuous-time random walk with rate $n$ starting from $\cZ^n_0=0$ and, conditioned on the set of jump times, the steps of $\cZ^n$ are uniformly and independently sampled from $\{(1,0),(-1,0),(0,1),(0,-1) \}$, 
		\item[(iii)] random walk converges to Brownian motion on $(\cZ^n,Z)$, and 	
		\item[(iv)] $\|Z-n^{-1/2}\cZ^n\|_\infty\rta 0$ in probability as $n\rta\infty$.
	\end{compactitem} 
	\label{prop1}
\end{lem}
\begin{proof} 
	Consider the setting of Proposition \ref{prop0} and denote the processes in that lemma by $(\check \cZ^1,\check Z^1),(\check \cZ^2,\check Z^2),\dots$. 
	For the finitely many $n$ with $an\leq1$, sample $\cZ^n$ with the marginal law in (ii), independently of $Z$.
	For the remaining $n$ choose $m=m(n)\in\N$ such that $2^{m-1}<an\leq 2^m$.  Set $\cZ^n_t=\check{\cZ}^m_{tn2^{-m}}$ and $Z^n_t=n^{-1/2}2^{m/2}\check{Z}^m_{tn2^{-m}}$ for $t\in[0,a]$. Finally couple the processes so they are defined on the same probability space with $Z^n=Z^1$ for all $n$.

\end{proof}

\begin{lem}
	Let $a,\eps>0$ be such that $3\eps<a$. There is a coupling of $(\cZ^n,Z)$ such that 
	\begin{compactitem}
		\item[(i)] $Z:[0,a]\to\R_+^2$ has the law of a variance $\frac 12$ uncorrelated Brownian cone excursion of duration $a$ in the first quadrant from $(1,0)$ to $(0,0)$,
		\item[(ii)] $\cZ^n$ has the same law as in Lemma \ref{prop1}, except that the walk starts from $\cZ^n_0=(\lceil \sqrt n\rceil-1,0)$ and we condition on the events that $\cZ^n_a=(0,0)$ and $\cZ^n_t\in\R_+^2$ for all $t\in[0,a]$,
		\item[(iii)] random walk on $(\cZ^n|_{[\eps,a-\eps]},Z|_{[\eps,a-\eps]})$ converges to Brownian motion, 
		\item[(iv)] $\|(Z-n^{-1/2}\cZ^n)|_{[\eps,a-\eps]}\|_\infty\rta 0$ in probability as $n\rta\infty$,
		\item[(v)] the number of vertices of $\cV(M_n)\setminus\partial\cV(M_n)$ divided by $an$ converges in probability to a deterministic constant independent of $a$. 
	\end{compactitem}	
	\label{prop2}
\end{lem}
\begin{proof}
	For any $\delta>0$ the number of downward jumps made by the first coordinate of the walk in (ii) will be in $((a-\delta)n/4,(a+\delta)n/4)$ with probability converging to 1 as $n\rta\infty$. Using this and the  correspondence between $\cV(M_n)$ and downward jumps of the first coordinate of $\cZ^n$ (see the text above Lemma \ref{prop-bij-steps}), we see that (v) follows from (ii). Therefore, in order to conclude the proof it is sufficient to find processes $\cZ^{n,\op{mid}},Z^{\op{mid}} : [\eps,a-\eps]\to\R^2$ such that

	\begin{compactitem}
		\item[(i)'] $Z^{\op{mid}}$ is equal in law to $Z|_{[\eps,a-\eps]}$, where $Z^n$ is as in (i),
		\item[(ii)'] the total variation distance between $\cZ^{n,\op{mid}}$ and $\cZ^n|_{[\eps,a-\eps]}$ goes to zero as $n\rta\infty$, where $\cZ^n$ is as in (ii),
		\item[(iii)'] random walk on $(\cZ^{n,\op{mid}},Z^{\op{mid}})$ converges to Brownian motion, and
		\item[(iv)'] $\|Z^{\op{mid}}-n^{-1/2}\cZ^{n,\op{mid}}\|_\infty\rta 0$ in probability as $n\rta\infty$.
	\end{compactitem}
	We will obtain processes satisfying (i)'-(iv)' by reweighting the law in Lemma \ref{prop1} by an appropriate Radon-Nikodym derivative. Let $(\cZ^{n,\op{u}},Z^{\op{u}})$ for $n\in\N$ have the law of the processes in Lemma \ref{prop1}, where we use the letter ``u'' in the notation since it is the first letter of ``unconditioned''. Letting $L^{\op{u}}$ (resp.\ $R^{\op{u}}$) denote the first (resp.\  second) coordinate of $Z^{\op{u}}$, define  $m=(m_1,m_2)\in\R^2$ by
	$$
	m_1 = \inf_{t\in[\eps,a-\eps]} L^{\op{u}}_t-L^{\op{u}}_{\eps},\qquad
	m_2 = \inf_{t\in[\eps,a-\eps]} R^{\op{u}}_t-R^{\op{u}}_{\eps}.
	$$

	Let $Y^0$ be a variance $\frac 12$ Brownian motion, 
	let $Y:\R_+\to\R$ be a variance $\frac 12$ Brownian motion conditioned to stay above $-1$ during $[0,\eps]$, and 
	let $Y^+$ be a variance $\frac 12$ Brownian motion conditioned to stay positive during $[0,\eps]$. Then let $\rho^0$ (resp.\ $\rho,\rho^+$) be the probability density of $Y^0_\eps$ (resp.\ $Y_\eps,Y^+_\eps$). If we reweight the law of $Z^{\op{u}}$ by 
	$$
	\xi=\xi(Z^{\op{u}}_\eps,Z^{\op{u}}_{a-\eps},m) := 
	c\rho(L^{\op{u}}_\eps-1) \rho^+(R^{\op{u}}_\eps) \rho^+(L^{\op{u}}_{a-\eps}) \rho^+(R^{\op{u}}_{a-\eps})
	\rho^0(L^{\op{u}}_{\eps})^{-1}
	\rho^0(R^{\op{u}}_{\eps})^{-1}
	\1_{m_1 >-L^{\op{u}}_\eps;m_2>-R^{\op{u}}_\eps}
	$$
	for a normalizing constant $c>0$, 
	then $Z^{\op{u}}|_{[\eps,a-\eps]}$ has exactly the law of the process $Z^{\op{mid}}$ in (i)' above. 
	
	Similarly for the discrete process, let $\cL^{n,\op{u}}$ (resp.\ $\cR^{n,\op{u}}$) denote the first (resp.\ second)  coordinate of $\cZ^{n,\op{u}}$ and define  $m^n=(m^n_1,m^n_2)\in\R^2$ by
	$$
	m^n_1 = \inf_{t\in[\eps,a-\eps]} \cL^{n,\op{u}}_t-\cL^{n,\op{u}}_{\eps},\qquad
	m^n_2 = \inf_{t\in[\eps,a-\eps]} \cR^{n,\op{u}}_t-\cR^{n,\op{u}}_{\eps}.
	$$
	Let $Y^{n,0}$ be a continuous time random walk with rate $0.5n$ and uniform steps $\pm 1$, 
	let $Y^n$ have the same law as $Y^{n,0}$ except that the process is conditioned to be $\geq-\lceil \sqrt{n}\rceil+1$ during $[0,\eps]$, and 
	let $Y^{n,+}$ have the same law as $Y^n$ except that we also condition the process to be non-negative on $[0,\eps]$. 
	Let $\rho^0_n$ (resp.\ $\rho_n,\rho_n^+$) be the probability mass function of $Y^{n,0}_\eps$ (resp.\ $Y^n_\eps,Y^{n,+}_\eps$). If we reweight the law of $\cZ^{n,\op{u}}$ by a constant multiple of the following for $c_n$ a normalizing constant
	\eqbn
	\begin{split}
	\xi_n=&\,\xi_n(\cZ^{n,\op{u}}_\eps,\cZ^{n,\op{u}}_{a-\eps},m^n) 
	:=\,  
	c_n \rho_n(\cL^{n,\op{u}}_\eps-(\lceil \sqrt{n}\rceil-1)) \rho^+_n(\cR^{n,\op{u}}_\eps) \rho_n^+(\cL^{n,\op{u}}_{a-\eps}) \rho_n^+(\cR^{n,\op{u}}_{a-\eps})\\
	&\qquad\qquad\qquad\qquad\qquad\qquad\qquad\qquad\qquad
	\cdot\rho^0_n(\cL^{n,\op{u}}_\eps)^{-1} \rho^0_n(\cR^{n,\op{u}}_\eps)^{-1} 
	\1_{m^n_1 \geq - \cL^{n,\op{u}}_\eps; m^n_2\geq -\cR^{n,\op{u}}_\eps},	
	\end{split}
	\eqen
	then $\cZ^{n,\op{u}}|_{[\eps,a-\eps]}$ has exactly the law of the process $\cZ^{n,\op{mid}}$ in (ii)' above.

	The function $\xi_n$ (viewed as a function of $(n^{-1/2}\cZ^{n,\op{u}}_\eps,n^{-1/2}\cZ^{n,\op{u}}_{a-\eps},n^{-1/2}m^n)$) converges to the function $\xi$ (viewed as a function of $(Z^{\op{u}}_\eps,Z^{\op{u}}_{a-\eps},m)$), which is continuous except for the indicator in its definition
	Therefore $\xi_n$ converges in probability to $\xi$ as $n\rta\infty$. Furthermore, $\sup_n\E[\xi_n^p]<\infty$ and $\E[\xi^p]<\infty$ for all $p>1$. This implies that if we reweight the law of $((\cZ^{n,\op{u}})_{n\in\N},Z^{\op{u}})$ by $\xi$ and we set $\cZ^{n,\op{mid}}=\cZ^{n,\op{u}}|_{[\eps,a-\eps]}$ and $Z^{\op{mid}}=Z^{\op{u}}|_{[\eps,a-\eps]}$ then we get processes satisfying (i)'-(iv)' above.

\end{proof}

The next lemma is a variant of Lemma \ref{prop2} where the Brownian path $Z$ has a random duration $A>0$. Notice that $Z$ and $\cZ$ have the same law as in Lemma \ref{prop5}.
\begin{lem}
	Let $\eps\in (0,0.1)$. There is a coupling of $(Z,\cZ^1,\cZ^2,\cZ^3,\dots)$ such that 
		\begin{compactitem}
			\item[(i)] $Z:[0,A]\to\R_+^2$ has the law of a variance $\frac 12$ uncorrelated Brownian cone excursion from $(1,0)$ to $(0,0)$ of random duration $A>0$,
			\item[(ii)] for each $n\in\N$, $\cZ^n:[0,A]\to\R^2$ is a continuous-time random walk and 
			$(M_n,\frk T_n):=\phi_{\op{M}}(\cZ^n)$  
			has the same law as the tree-decorated planar map in Theorem \ref{thm-tutte-conv0},
			\item[(iii)] random walk converges to Brownian motion on $\big(\cZ^n|_{[\eps,A-\eps]},Z|_{[\eps,A-\eps]}\big)_{n\in\N}$ on the event that $A>3\eps$,  	
			\item[(iv)] $\1_{A>3\eps}\|(Z-n^{-1/2}\cZ^n)|_{[\eps,A-\eps]}\|_\infty\rta 0$ in probability as $n\rta\infty$, 
			\item[(v)] the number of vertices of $\cV(M_n)\setminus\partial\cV(M_n)$ divided by $An$ converges in probability to a deterministic constant. 
	\end{compactitem}
	\label{prop12}
\end{lem}

\begin{proof}%
	As mentioned above, the difference between the current lemma and Lemma \ref{prop2} is that the latter lemma considers processes with a fixed duration $a$ while this lemma considers processes with a random duration $A$. To describe this difference more precisely, first observe that the law of $A$ has a continuous density $\rho:\R_+\to[0,\infty)$ with respect to Lebesgue measure on $\R_+$ given by $\rho(t)=ct^{-3}\exp( -1/t)$ for a suitably chosen constant $c>0$; see e.g.\ \cite[Theorem 1.2]{ag-disk}. Sample a random variable $A$ of this law and, conditioned on $A$, let $\cZ^n$ and $Z$ be as in Lemma \ref{prop2} with $a=A$. Then $Z$ is equal in law to the process $Z$ in Lemma \ref{prop5}.

	For each $n\in\N$ let $\check{\cZ}^n$ be a process with the same law as the process $\cZ^n$ in Lemma \ref{prop5}, assuming the set of jumps is uniformly distributed on $[0,A]$ conditioned on $A$ and the ordered set of jumps, such that $\check{\cZ^n}$ is coupled together with $(Z,\cZ^1,\cZ^2,\dots)$ in an arbitrary way. Let $K_n$ (resp.\ $\check K_n$) denote the number of jumps made by the process ${\cZ}^n$ (resp.\ $\check{\cZ}^n$). The random variables $K_n$ and $\check K_n$ have the same support $\scr S=\{\lceil \sqrt{n}\rceil-1,\lceil \sqrt{n}\rceil+1,\lceil \sqrt{n}\rceil+3,\dots   \}$ and for any $k\in\scr S$, 
	the law of ${\cZ}^n$ conditioned on $K_n=k$ is equal to 
	the law of $\check{\cZ}^n$ conditioned on $\check K_n=k$.

	We will prove that the total variation distance between $\check K_n$ and $K_n$ is $o_n(1)$, which implies that the total variation distance between $\cZ^n$ and $\check\cZ^n$ (equivalently, between $\cZ^n$ and the process described in Assertion (ii) of the lemma) is $o_n(1)$. Assertions (iii), (iv) and (v) are satisfied by Lemma \ref{prop2} and since the law of $(Z,\cZ^n)$ here is a mixture of the associated laws in Lemma \ref{prop2} for different $a$. We conclude that the lemma follows once we have proven that the total variation distance between $\check K_n$ and $K_n$ is $o_n(1)$.
	
	Let $b>0$ be such that $ bn\in\scr S$. 
	By \cite[equation (12)]{guy-walk-in-cone}, the number of walks of length $ bn$ in $\N_0^2$ from $(\lceil\sqrt{n}\rceil-1,0)$ to $(0,0)$ is 
	\eqb
	\frac{\lceil\sqrt{n}\rceil}{( bn+1)( bn+2)}
	\binom{ bn+2}{0.5( bn-\lceil\sqrt{n}\rceil+1)}
	\binom{ bn+2}{0.5( bn+\lceil\sqrt{n}\rceil+1)} = c_0 
	4^{ bn}n^{-5/2}\rho(b)(1+o_n(1))
	\label{eq:guy}
	\eqe
	for some $c_0>0$, where the asymptotics are obtained via Stirling's formula and the $o_n(1)$ is uniform for all $b$ in compact subsets of $(0,\infty)$.

	It follows that
	$$
	\P[\check K_n= bn ] = \frac 2n \rho(b) (1+o_n(1)),
	$$
	where again the $o_n(1)$ is uniform for all $b$ in compact subsets of $(0,\infty)$.
	Now let $\wt{\P}$ be a probability measure such that the law of $\cZ^n$ under $\wt{\P}$ is that of an (unrestricted) continuous-time nearest-neighbor random walk $\cZ$ on $\Z^2$ with rate $n$ started from $\cZ^n_0=(\lceil \sqrt{n} \rceil-1,0)$, and let $E$ denote the following event with $Q_1=\{ (x,y)\in\R^2\,:\,x\geq 0,y\geq 0 \}$ being the closed first quadrant
	$$
	E = \{ \cZ^n_a = (0,0),\cZ^n_t \in Q_1\,\forall t\in[0,a] \}.
	$$
	Then by Bayes' rule, \eqref{eq:guy}, and using that the law of $K_n$ under $\wt{\P}$ is that of an intensity $an$ Poisson random variable, if $r\in\R$ is such that $an+r\sqrt{an}\in \scr S$,
	\eqbn
	\begin{split}
	&\P[K_n = an+r\sqrt{an}\,|\,A=a] = \wt{\P}[K_n=an+r\sqrt{an}\,|\,E] \\
	&\qquad=
	\frac{\wt{\P}[E\,|\,K_n=an+r\sqrt{an}]   }
	{\sum_{m\,:\,an+m\in\scr S}  \wt{\P}[E\,|\,K_n=an+m] \wt{\P}[K_n=an+m]} \wt{\P}[K_n=an+r\sqrt{an}]\\
	&\qquad =
	\frac{c_0 n^{-5/2} \rho( a+ra^{1/2}n^{-1/2} )  }
	{\sum_{m\,:\,an+m\in\scr S}  c_0 n^{-5/2} \rho( a+m/n) \wt{\P}[K_n=an+m]}\cdot  (an)^{-1/2}\sigma(r)(1+o_n(1)) \\
	&\qquad= 2(an)^{-1/2}\sigma(r)(1+o_n(1)),	
	\end{split}
	\eqen 
	where $\sigma$ denotes the density of a standard normal random variable and $o_n(1)$ is uniform for all $r$ in compact subsets of $\R$. Furthermore, $\P[|K_n-an|/\sqrt{an}>C\,|\,A=a]=o_C(1)$, where the $o_C(1)$ is uniform over all $n\in\N$ and $a$ restricted to a compact subset of $(0,\infty)$. 
	It follows that
	\eqbn
	\begin{split}
		\P[K_n = bn]
		&=
		\int_{0}^{\infty} \rho(a) \P[K_n = bn\,|\,A=a]\,da
		=
		\int_{0}^{\infty} 2\rho(a) (an)^{-1/2} \sigma(a^{-1/2} (b-a)\sqrt{n} )(1+o_n(1)) \,da\\
		&= \frac 2n \rho(b) (1+o_n(1)).
	\end{split}
	\eqen

\end{proof}

\begin{proof}[Proof of Lemma \ref{prop5}]
	Let $\eps_n\rta 0$ be a sequence of numbers in $(0,0.1)$ converging slowly to 0. Let $Z$ have the law described in Assertion (i) of the lemma and then, for each $n\in\N$, let $\cZ^n$ have the same conditional law given $Z$ as in Lemma \ref{prop12} with $\eps=\eps_n$. Then Assertions (ii) and (iv) are satisfied, where for (iv) we need to use tightness $\cZ^n$ near the end-points. 

	In order to prove Assertions (v) and (vi), notice that 
	the existence of appropriate $(\frk v^n_{\op{U}})_{n\in\N}$ and $\frk t_{\op{U}}$ would have 
	been immediate if we did not condition $\cZ^n$ to form a cone excursion since the set of $(-1,0)$ jumps of $\cZ^n$ have the law of a Poisson point process. Then notice that the assertions still hold under our conditioning since the event we are conditioning on (namely, $\cZ^n$ forming a cone excursion) has polynomial probability in $n$. 
	Let $\delta\in(0,0.1)$. In order to conclude the proof it is sufficient to argue that 
	\eqb
	\limsup_{n\rta\infty}\sup_{v\in\cV(M_n)}
	d_{\op{Prok}}(\sigma^{\eta_n(\frk t_v)},\sigma^{v,n})=o_\delta(1).
	\label{eq66}
	\eqe

	Set $D_n=\eta([0,\eps_n])\cup\eta([A-\eps_n,A])$ and observe that $D_n\subset B(1,\delta^2)$ for sufficiently large $n$ a.s. Therefore, Brownian motion started on $\ol{\D}\setminus B(1,\delta)$ is hits $D_n$ before existing $\ol{\D}$ with probability $o_\delta(1)$ and, by Assertion (iii) of Lemma \ref{prop12}, the same holds for the random walk $\wt W_n$ if it is started on $\ol{\D}\setminus B(1,\delta)$. It follows from Assertion (iii) of Lemma \ref{prop12} that
	\eqb
	\limsup_{n\rta\infty}
	\sup_{v\in\cV(M_n)\,:\,\eta_n(\frk t_v)\in\ol{\D}\setminus B(1,\delta)}
	d_{\op{Prok}}(\sigma^{v,n},\sigma^{\eta_n(\frk t_v)}) = o_\delta(1).
	\label{eq55}
	\eqe

	Letting $\wt\sigma$ denote the law of the deterministic stochastic process $(\wt\gamma_t)_{t\in[0,1]}$ given by $\wt\gamma_t=1$ for all $t\in[0,1]$, we have
	\eqb
	\sup_{z\in\ol{\D}\cap B(1,\delta) }d_{\op{Prok}}(\wt\sigma,\sigma^z)=o_\delta(1),
	\label{eq67}
	\eqe
	since a Brownian motion started in $\ol{\D}\cap B(1,\delta)$ exits $\D$ before exiting $B(1,\sqrt{\delta})$ with probability $1-o_\delta(1)$.

	Let $\tau=\inf\{t\in [0,T]\cap\Z\,:\, \wt W_n(t)\not\in B(1,\delta) \}$, where we recall that $T$ is the lifetime of $\wt W_n$. If $\tau=\infty$ then we have $d_{\op{CMP}}(\wt W,\wt\gamma)=o_\delta(1)$. If $\tau<\infty$ and $n$ is large, then we will most likely also have  $d_{\op{CMP}}(\wt W_n, \wt\gamma)=o_\delta(1)$ since, by Assertion (iii) of Lemma \ref{prop12}, $\wt W_n$ is likely to hit a boundary vertex of the map before reaching distance more than $o_\delta(1)$ away from the point 1. This implies
	$$
	\limsup_{n\rta\infty}\sup_{v\in\cV(M_n)\,:\,\eta_n(\frk t_v)\in\ol{\D}\cap B(1,\delta)}
	d_{\op{Prok}}(\wt\sigma,\sigma^{v,n})=o_\delta(1),
	$$

	which gives \eqref{eq66} in combination with \eqref{eq55} and \eqref{eq67}.
\end{proof}

\subsection{Convergence of Tutte embedded planar map to LQG}
\label{sec:tutte-conv2}

The following lemma gives a general condition which guarantees that some embedding $\wt\phi$ is close to the Tutte embedding. Namely, it says that if the random walk on the embedded map is close in law to 2D Brownian motion and some natural boundary and consistency conditions (i)-(iii) are satisfied, then the two embeddings are close.
\begin{lem}
	Let $M$ be a planar map with disk topology and distinguished vertices $v_0\in \cV(M)\setminus\partial \cV(M)$ and $v_1\in\partial \cV(M)$. Let $\phi:\cV(M)\to\ol{\D}$ be the Tutte embedding of $M$ centered at $v_0$ such that $\phi(v_1)=1$, and let $\wt\phi:\cV(M)\to\ol{\D}$ be an embedding of $M$ such that the following holds.
	\begin{compactitem}
		\item[(i)] Embedding fixed as for Tutte: $\wt\phi(v_0)=0$ and $\wt\phi(v_1)=1$,
		\item[(ii)] Boundary embedded to boundary: $\wt\phi(v)\in\partial \D$ for all $v\in \partial \cV(M)$,
		\item[(iii)] Cyclic ordering of boundary preserved: The map $v\mapsto \wt\phi(v)$ restricted to $\partial \cV(M)$ preserves the cyclic ordering of vertices, i.e., if $w_1,w_2,w_3,\dots$ are ordered counterclockwise on $\partial \cV(M)$ then $\wt\phi(w_1),\wt\phi(w_2),\wt\phi(w_3),\dots$ are ordered counterclockwise on $\partial\ol{\D}$.
	\end{compactitem}
	Then for some universal constant $C>1$,
	$$
	\|\wt\phi-\phi\|_\infty<C\sup_{v\in\cV(M)}d_{\op{Prok}}(\wt\sigma^{v},\sigma^{\wt\phi(v)}),
	$$
	where $\wt\sigma^{v}$ is defined exactly as $\sigma^{v,n}$ at the beginning of Section \ref{sec:tutte-conv0} except that $\eta_n(\frk t_\cdot)$ is replaced by $\wt\phi(\cdot)$ everywhere. 
	\label{lem:rw-conv-implies-tutte-conv1}
\end{lem}
\begin{proof}
	First let us introduce some notation. For $v\in \cV(M)$ let $W^v:\N_0\to \cV(M)$ denote simple random walk on $M$ satisfying $W^v(0)=v$ and let $\tau:=\inf\{t\in\N_0\,:\,W^v(t)\in\partial \cV(M) \}$ denote the hitting time of the boundary. For $w_1,w_2\in\partial \cV(M)$ let $\op{arc}_M(w_1,w_2)$ denote the boundary vertices between $w_1$ and $w_2$ in counterclockwise direction from $w_1$ to $w_2$. For $z\in\D$ let $B^z:[0,\infty)\to\C$ denote a 2D Brownian motion satisfying $B^z(0)=z$ and let $\sigma:=\inf\{t\geq 0\,:\, B^z(t)\not\in\D \}$ denote the exit time from the unit disk. For $z_1,z_2\in\partial\D$ let $\op{arc}_{\D}(z_1,z_2)$ denote the counterclockwise boundary arc of $\D$ from $z_1$ to $z_2$.
	
	Define $\delta:=\sup_{v\in\cV(M)}d_{\op{Prok}}(\wt\sigma^{v},\sigma^{\wt\phi(v)})$, 
	where for arbitrary finite measures $\mu_1,\mu_2$ on a metric space $(\frk M,d)$ equipped with its Borel $\sigma$-algebra, 
	$$
	d_{\op{Prok}}(\mu_1,\mu_2):=\inf\{\eps>0\,:\,
	\mu_1(A)\leq \mu_2(A_\eps)+\eps
	\text{ and }
	\mu_2(A)\leq \mu_1(A_\eps)+\eps
	\text{ for all Borel sets }
	 A\subset \frk M
	 \}
	$$ 
	and $A_\eps:=\{ x\,:\,\exists y\in A, d(x,y)<\eps\}$ is the $\eps$-neighborhood of $A$ for any $A\subset \frk M$ and $\eps>0$.
	Throughout the proof we let $C>1$ denote some universal constant which may change from line to line. We start by arguing that
	\eqb
	\|(\wt\phi-\phi)|_{\partial \cV(M)}\|_\infty<C\delta.
	\label{eq:bdy-error}
	\eqe
	By the definition of $d_{\op{Prok}}$ and $d_{\op{CMP}}$ (see \eqref{eq:curves-mod-reparam}), Assumption (i), and since the law of $B^0(\sigma)$ is the uniform measure on $\partial\D$,
	\eqb
	|\P[B^0(\sigma) \in \op{arc}_{\D}(1,\wt\phi(v))]
	- \P[\wt\phi\circ W^{v_0}(\tau) \in \op{arc}_{\D}(1,\wt\phi(v)) ] | < C\delta.
	\label{eq:hitting}
	\eqe
	
	By first using Assumptions (ii) and (iii), and then using the definition of the Tutte embedding, the second probability on the left side differs by at most $C\delta$ from
	$$
	\P[W^{v_0}(\tau) \in \op{arc}_{M}(v_1,v) ] = \P[B^0(\sigma) \in \op{arc}_{\D}(1,\phi(v))],
	$$
	so
	$$
	|\P[B^0(\sigma) \in \op{arc}_{\D}(1,\wt\phi(v))]-\P[B^0(\sigma) \in \op{arc}_{\D}(1,\phi(v))]|<C\delta.
	$$
	Since the law of $B^0(\sigma)$ is the uniform measure on $\partial\D$ it follows that $|\wt\phi(v)-\phi(v)|<C\delta$. 

	Next we will treat vertices in the bulk and let $v\in \cV(M)\setminus\partial \cV(M)$. The following is proven exactly as \eqref{eq:hitting} 
	$$
	|\E[ B^{\wt\phi(v)}(\sigma) ] - \E[\wt\phi\circ W^v(\tau)] | \leq C\delta. %
	$$
	The first expectation equals $\wt\phi(v)$ by the optional stopping theorem while the second expectation satisfies $|\E[\wt\phi\circ W^v(\tau)]-\E[\phi\circ W^v(\tau)]|\leq C\delta$ by \eqref{eq:bdy-error}. By the definition of the Tutte embedding, $\phi(v)=\E[\phi\circ W^v(\tau)]$. Combining these results we get $|\wt\phi(v)-\phi(v)|\leq C\delta$.
\end{proof}

We will be using the following variant of Lemma \ref{lem:rw-conv-implies-tutte-conv1}. %
\begin{lem}
	Consider the setting of Lemma \ref{lem:rw-conv-implies-tutte-conv1}, except that Assertions (i)-(iii) are required to hold for some embedding $\check\phi:\cV(M)\to\ol{\D}$ instead of $\wt\phi$. 
	Then there is a universal constant $C>1$ such that 
	$$
	\|\wt\phi-\phi\|_\infty<
	C\|\check\phi-\wt\phi\|_\infty
	+
	C\sup_{v\in\cV(M)}d_{\op{Prok}}(\wt\sigma^{v},\sigma^{\wt\phi(v)}).
	$$
	\label{lem:rw-conv-implies-tutte-conv2}
\end{lem}
\begin{proof}
	
	Arguing just as for \eqref{eq:hitting}, for any $v\in\partial\cV(M)$,
	\eqb
	|\P[B^0(\sigma) \in \op{arc}_{\D}(1,\check\phi(v))]
	- \P[\check\phi\circ W^{v_0}(\tau) \in \op{arc}_{\D}(1,\check\phi(v)) ] | < C\|\check\phi-\wt\phi\|_\infty
	+
	C\sup_{v\in\cV(M)}d_{\op{Prok}}(\wt\sigma^{v},\sigma^{\wt\phi(v)}).
	\eqe
	Arguing as in the proof of Lemma \ref{lem:rw-conv-implies-tutte-conv1}, this implies that $|\check\phi(v)-\phi(v)|$ is bounded by the right side of the last display, so, after increasing $C$, $|\wt\phi(v)-\phi(v)|$ is also bounded by the right side of the last display. The bound for interior vertices follows from the bound of boundary vertices by the same argument as in Lemma \ref{lem:rw-conv-implies-tutte-conv1}.
\end{proof}

In the next lemma we will apply the previous lemma to the planar map appearing in Theorem \ref{thm-tutte-conv0} and an embedding $\wt\phi_n$ coming from mating of trees. The hardest part of the proof is to argue that the embedding $\wt\phi_n$ almost satisfies Assertions (ii) and (iii) from Lemma \ref{lem:rw-conv-implies-tutte-conv1}. We will prove this by using Assertion (iv) of Lemma \ref{prop5} along with properties of Brownian motion.

\begin{figure}
	\includegraphics[scale=1]{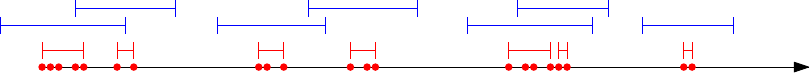}
	\caption{Illustration of the proof of Lemma \ref{prop:MOT-emb-close-to-Tutte}. The red dots represent the fractal set $\cA$ given by the running infima of $L$. The red (resp.\ blue) intervals are the intervals $[t_j,s_j]$ (resp.\ $\wt{\cA}_j$).}
	\label{A-tilde}
\end{figure}

Let $( (\cZ^n)_{n\in\N},Z, (\frk v^n_{\op{U}})_{n\in\N} , \frk t_{\op{U}},A)$ %
be as in Lemma \ref{prop5}, let $(M_n,\frk T_n):=\PhiM(\cZ^n)$ be the associated tree-decorated planar maps, and let $(\BB D,h,\eta)$ be the associated SLE-decorated LQG surface embedded such that $\eta(0)=1$ and $\eta(\frk t_{\op{U}})=0$. Let $\phi_n:\cV(M_n)\to\ol{\D}$ denote the Tutte embedding of $M_n$ centered at $\frk v^n_{\op{U}}$. Define another embedding $\wt\phi_n:\cV(M_n)\to\ol{\D}$ by $\wt\phi_n(v):=\eta(\frk t_v)$ for $v\in \cV(M_n)$.
\begin{lem}
	In the setting described above, $\|\wt\phi_n-\phi_n\|_\infty\rta 0$ in probability as $n\rta\infty$.
	\label{prop:MOT-emb-close-to-Tutte}
\end{lem}
\begin{proof}
	Define $\check\phi_n:\cV(M_n)\to\ol{\D}$ by $\check\phi_n(\frk v^n_{\op{U}})=0$,  $\check\phi_n(v)=\wt\phi_n(v)$ for $v\in \cV(M_n)\setminus(\partial \cV(M_n)\cup\{\frk v^n_{\op{U}}\})$,  
	set $\check\phi_n(\frk u)=\phi_n(\frk u)=1$ for $\frk u$ equal to the starting vertex of the root edge $\frk e$, and
	define $\check\phi_n$ as follows for other boundary vertices. 
	For $v\in\partial \cV(M_n)\setminus\{\frk u \}$ let $\tau_v=\inf\{ t\geq 0\,:\,L_t=n^{-1/2}\cL^n_{\frk t_v} \}$ and set $\check\phi_n(v)=\eta(\tau_v)$. Then $\check\phi_n$ satisfies Assertions (i)-(iii) of Lemma \ref{lem:rw-conv-implies-tutte-conv1}.  
	By Assertion (iii) of Lemma \ref{prop5} and Lemma \ref{lem:rw-conv-implies-tutte-conv2}, in order to conclude the proof of the lemma it is sufficient to argue that $\|\check\phi_n-\wt\phi_n\|_\infty\rta 0$ as $n\rta\infty$. We have $\check\phi_n(\frk v^n_{\op{U}})\rta 0$ by Assertion (vi) of Proposition \ref{prop5}. Notice that $\tau_v$ and $\frk t_v$ are running infima for $L$ and $\cL^n$, respectively. Let $\eps>0$ and pick a (random) $\delta_1>0$ such that $|\eta(t)-\eta(s)|<\eps$ for $|t-s|<3\delta_1$. Let $T>0$ denote the duration of $L$ and let $\cA=\{t\geq 0\,:\, L_t<L_s\,\,\forall s< t \}\cap[\delta_1,T-\delta_1]$ denote the set of running infima for $L$ away from the end-points of $L$. See Figure \ref{A-tilde} for an illustration. Define $t_j$ for $j\in\N_0$ inductively by $t_0=0$ and $t_j=\inf\{ t\in\cA \,:\, t\geq t_{j-1}+\delta_1 \}$ for $j\in\N$, where the infimum of an empty set is $\infty$. Let $m=\max\{ j\in\N \,:\,t_j<\infty \}$. Set $s_j=\sup (\cA\cap[t_j,t_j+\delta_1])$ and $\wt\cA_j=[t_j-\delta_1,s_j+\delta_1]$. We have $\eta(s_j)=\eta(t_{j+1})$ when $t_{j+1}<\infty$ since $\min_{t\in[s_j,t_{j+1}]}L_t=L_{s_j}=L_{t_{j+1}}$ for such $j$.

	There exists a (random) $\delta_2\in(0,\delta_1/2)$ such that $L_t>L_{s_j}+\delta_2=L_{t_{j+1}}+\delta_2$ for all $j=0,\dots,m-1$ and $t\in[s_j+\delta_1,t_{j+1}-\delta_1]$ since $L$ makes an excursion above $L_{s_j}=L_{t_{j+1}}$ during $[s_j,t_{j+1}]$. Furthermore, by possibly decreasing $\delta_2$ it also holds that $L_{s_j}<L_{t_j}-\delta_2$ for $j=0,\dots,m$. If $\|L-n^{-1/2}\cL^n\|_\infty<\delta_2/3$, $v\in\partial\cV(M_n)\setminus\{\frk u \}$ 
	and $\tau_v\in\wt{\cA}_j$, then we have $\frk t_v\in \wt{\cA}_{j-1}\cup \wt{\cA}_j\cup\wt{\cA}_{j+1}$ by the definition of $\tau_v$ and $\delta_2$. 
	In this setting, if we assume $\frk t_v\in \wt{\cA}_{j-1}$,
	$$
	|\wt\phi_n(v)-\check\phi(v)|
	= 
	|\eta(\frk t_v)-\eta(\tau_v)|
	\leq 
	|\eta(\frk t_v)-\eta(s_{j-1})|
	+ |\eta(s_{j-1})-\eta(t_j) |
	+ |\eta(t_{j})-\eta(\tau_v)|<\eps+0+\eps.
	$$ 
	By a similar argument we also have $|\wt\phi_n(v)-\check\phi_n(v)|<2\eps$ if $\frk t_v\in \wt{\cA}_{j}$ or $\frk t_v\in \wt{\cA}_{j+1}$. 
	This concludes the proof for $v$ such that $\frk t_v\in[\delta_1,T-\delta_1]$. %
	On the other hand, if $\frk t_v\in[0,\delta_1]\cup[T-\delta_1,T]$ we have $|\wt\phi_n(v)-\check\phi_n(v)|\leq |\wt\phi_n(v)-1|+|1-\check\phi_n(v)|=o_\eps(1)$, where the bound for $\check\phi_n$ follows by using that $\check\phi_n$ satisfies Assumption (iii) of Lemma \ref{lem:rw-conv-implies-tutte-conv1} and our bound on $|\wt\phi_n(w)-\check\phi_n(w)|$ for $w\in \partial\cV(M_n)$ with $\frk t_w\in [\delta_1,T-\delta_1]$. 
\end{proof}

\begin{proof}[Proof of Theorem \ref{thm-tutte-conv0}]
	Define the measure $\mu^{\op{MOT}}_n$ exactly as the measure $\mu_n$ in \eqref{eq:mu-ell}, except that we use the embedding $\wt\phi_n$ instead of the Tutte embedding $\phi_n$. Lemma  \ref{prop:MOT-emb-close-to-Tutte} implies that  $d_{\op{Prok}}(\mu^{\op{MOT}}_n,\mu_n)\rta 0$ in probability as $n\rta\infty$, while Assertions (v) and (vii) of Lemma \ref{prop5} implies that $d_{\op{Prok}}(c\mu^{\op{MOT}}_n,\mu)\rta 0$ in probability as $n\rta\infty$ for a deterministic constant $c>0$.

	In order to conclude the proof of the first assertion of the theorem, we need to argue that $d_{\op{Prok}}(\nu_n,\nu)\rta 0$. Define $\check\nu_n$ just as $\nu_n$, but with the embedding $\check\phi_n$ from the proof of Lemma \ref{prop:MOT-emb-close-to-Tutte}. For any $v\in\partial \cV(M_n)$, by the interpretation of $\cL$ as a discrete boundary length process measuring the height in the spanning tree, which corresponds to the length of an appropriate boundary segment when we are on a boundary vertex %
	$$
	\check \nu_n(\op{arc}_{M_n}(1,\check\phi_n(v)))
	= n^{-1/2}\cL_{\frk t_v}+o_n(1)
	\rta L_{\tau_v}
	= \nu(\op{arc}_{\partial\D}(1,\eta(\tau_v) ) ),
	$$ 

	so $d_{\op{Prok}}(\check\nu_n,\nu)\rta 0$ in probability as $n\rta\infty$. It follows from the proof of Lemma \ref{prop:MOT-emb-close-to-Tutte} that $\|\check\phi_n-\wt\phi_n\|_\infty\rta 0$ and we argued above that  
	$\|\phi_n-\wt\phi_n\|_\infty\rta 0$. Therefore $\|\check\phi_n-\phi_n\|_\infty\rta 0$, which implies 
	$d_{\op{Prok}}(\check\nu_n,\nu_n)\rta 0$. We conclude from the triangle inequality that $d_{\op{Prok}}(\nu_n,\nu)\rta 0$ as desired.
	
	The second assertion of the theorem follows by using Assertion (iii) of Lemma \ref{prop5} along with Lemma  \ref{prop:MOT-emb-close-to-Tutte}.
\end{proof}

\appendix

\section{LQG and SLE estimates}
\label{app2}
This appendix collects a number of estimates and other results on LQG and SLE. 

\subsection{Absolute continuity for 0-quantum cone and $\gamma$-quantum cone}
The following result will be used in the proof of Lemma \ref{lem:rooted}. We omit the proof since it is identical to the one of \cite[Corollary 3.17]{bp-lqg-book}, which proves the same result for $U=D$.
	\begin{lem}
		Let $D\subseteq\C$ be a bounded domain and let $U\subseteq D$ be a bounded domain. Let $\wt h$ be a centered Gaussian field in $D$ with covariance function $G:D\times D\to\R$ of the form $G(x,y)=-\log|x-y|+g(x,y)$, where $g:\ol D\times\ol D\to\R$ is continuous. Let $z_0$ be sampled from $\mu_{\wt h}|_{U}$ renormalized to be a probability measure. Conditioned on $z_0$, the law of $\wt h$ weighted by $\mu_{\wt h}(U)\E[\mu_{\wt h}(U)]^{-1}$ is equal to the law of $\wt h+\gamma G(\cdot,z_0)$.
	\label{lem:rooted0}
	\end{lem}
	
	\begin{proof}[Proof of Lemma \ref{lem:rooted}]
		Throughout the proof let $\wh h$ be a whole-plane GFF with additive constant chosen such that the circle of radius $2^R$ has average 0. Then $\wh h$ is a centered Gaussian field in $\C$ with covariance 
		$$
		(z_1,z_2)\mapsto\log \frac{|2^{-R}z_1|_+|2^{-R}z_2|_+}{|2^{-R} (z_1-z_2)|},\qquad |\cdot|_+:=\max\{|\cdot|,1 \},
		$$
		see e.g.\ \cite[Equation (2.5)]{cw-liouville}. 
		Recall from Section \ref{sec:mot} that $h|_{B(0,2^R)}\eqD \wh h|_{B(0,2^R)}$ and that the field of the $2^R$-embedding of a $\gamma$-quantum cone restricted to $B(0,2^R)$ is equal in law to $(\wh h-\gamma\log|2^{-R}\cdot |)|_{B(0,2^R)}$.
		
		By Lemma \ref{lem:rooted0}, conditioned on $z_0$ the field $h|_{B(2^R)}$ reweighted by $\mu_h(U)\E[\mu_h(U)]^{-1}$ is equal in law to $(\wh h-\gamma \log|2^{-R}(\cdot-z_0)|)|_{B(2^R)}$. 
		This implies that conditioned on $z_0$ the field $h(\cdot+z_0)|_{B(2^R-2^{-k})}$ reweighted by $\mu_h(U)\E[\mu_h(U)]^{-1}$ is equal in law to $(\wh h(\cdot+z_0)-\gamma \log|2^{-R}(\cdot)|)|_{B(2^R-2^{-k})}$. 
		Using that $\wh h(\cdot+z_0)|_{B(2^R-2^{-k})}\eqD \wh h|_{B(2^R-2^{-k})}$ we get 
		the first assertion of the lemma.
		
		We get the upper bound in \eqref{eq:annulus-moment} by applying \eqref{eq:ball-upper} with exponent $p+9$ 
		\eqbn
		\begin{split}
			&\EE[\mu_h(U)^{-p}] 
			= 
			\int_0^\infty \PE[\mu_h(U)^{-p}>y]\,dy
			= 
			\int_0^\infty p2^{k (2+\gamma^2/2) p} x^{p-1}\PE[\mu_h(U)<2^{-k(2+\gamma^2/2)}/x]\,dx\\
			&\quad\leq 
			C'_p \exp(C'_p (R+k))
			+ 
			C'_p\exp(C'_p k)\int_{\exp(\sqrt{R+k})}^\infty x^{p-1} \exp(-c(\log x)^2/(R+k))+x^{-10}\,dx\\ 
			&\quad\leq 
			C_p \exp(C_p (R+k)),
		\end{split}
		\eqen
		where $C_p,C'_p>1$ are constants depending on $p$. When proving the last assertion we may assume that $U=B(0,2^R)$ since making the radius smaller can only decrease the area. The result now follows by scaling and by using that GMC measures have finite moments up to power $4/\gamma^2$.
\end{proof}

\subsection{Estimates for LQG areas and SLE}
Lemmas \ref{prop-mu_h-inf}, \ref{prop-mu_h-sup}, and \ref{prop:cone-moment} give upper and lower bounds for the LQG area of a ball. In the first two lemmas we consider a whole-plane GFF while in the third lemma we consider a $\gamma$-quantum cone. In Lemmas \ref{prop:contained1}, \ref{prop:contained2} and \ref{prop:contained3} we bound the probability that a segment of a space-filling SLE parametrized by LQG area contains or is contained in given balls.
\begin{lem} \label{prop-mu_h-inf}
	For any $p>1$ there is a $C$ such that the following holds.
	Let $h$ be a whole-plane GFF, let $\mu_h$ be its $\gamma$-quantum area measure, and let $(h_r(z))_{z\in\C,r>0}$ be its circle average process.  For each $z\in \BB C$, $r>0$, and $\delta \in (0,1)$,  
	\eqb \label{eqn-mu_h-inf}
	\BB P\!\left( \mu_h(B(z,r)) < \delta r^{2+\frac{\gamma^2}{2} }  e^{\gamma h_r(z)} \right)  \leq C\delta^p.
	\eqe 
	Furthermore, if the additive constant of the field has been chosen such that $h_r(z)$ is a centered Gaussian of variance $\sigma^2$ then
	\eqb 
	\BB P\!\left( \mu_h(B(z,r)) < \delta r^{2+\frac{\gamma^2}{2} }  \right)   \leq
	C\delta^p +\exp\Big(-\frac{(\log \delta^{-1})^2}{C\sigma^2 }\Big), 
	\label{eq:ball-upper}
	\eqe  
\end{lem}
\begin{proof}
	The first assertion of the lemma is \cite[Lemma 3.12]{ghm-kpz}. The cited lemma assumes $r\in(0,1)$, but this assumption is not needed. The second assertion of the lemma follows from the first assertion with exponent $2p$ and parameter $\delta^{1/2}$, along with the fact that $\P[e^{\gamma h_r(z)}<\delta^{1/2}]$ is bounded above by the second term on the right side of \eqref{eq:ball-upper}, where we might need to increase the constant $C$.
\end{proof}

\begin{lem}\label{prop-mu_h-sup}
	For any $p\in(0,4/\gamma^2)$ there is a $C>1$ such that with notation as in the previous lemma, 
	\eqb \label{eqn-mu_h-sup}
	\BB P\!\left( \mu_h(B(z,r)) > \delta^{-1} r^{2+\frac{\gamma^2}{2} }  e^{\gamma h_{r}(z)} \right)   \leq C\delta^{p}. 
	\eqe  
	Furthermore, if the additive constant of the field has been chosen such that $h_r(z)$ is a centered Gaussian of variance $\sigma^2$ then for any $\delta_1,\delta_2\in(0,1)$,
	\eqb
	\BB P\!\left( \mu_h(B_r(z)) > (\delta_1\delta_2)^{-1} r^{2+\frac{\gamma^2}{2} } \right)   \leq
	C\delta^{p}_1 +
	\exp\Big(-\frac{(\log \delta^{-1}_2)^2}{C\sigma^2 }\Big). 
	\label{eq:ball-lower}
	\eqe  
\end{lem}
\begin{proof}
	The first assertion follows from the same argument as \cite[Lemma 3.12]{ghm-kpz}, except that we use that $\E[  \mu_h(B(0,1))^p ]<\infty$ for any $p\in(0,4/\gamma^2)$ if $h$ has the circle average embedding. The second assertion follows from the first assertion by a similar argument as for \eqref{eq:ball-upper}.
\end{proof}

The next lemma considers the $\gamma$-quantum cone instead of the whole-plane GFF, which requires us to consider annuli centered at the origin due to the drift of the field.
\begin{lem}
	Suppose $h$ is the field of circle average embedding of the $\gamma$-quantum cone. Then for any $p\in(0,\op{min}\{1, 2Q/\gamma  - 2\})$ we have $\E[\mu(B(0,1))^p]<\infty$. Furthermore, there is a $C>1$ depending only on $p$ such that for any $s\in(0,1]$, $x>10$ and $y\geq 1$, 
	$$
	\P[\mu(B(0,s )) >  s^{Q\gamma  - \gamma^2} x y^{\sqrt{\log s^{-1}}} ] 
	\leq Cx^{-p} + \exp(-(\log y)^2/(2\gamma^2)).
	$$ 
	\label{prop:cone-moment}
\end{lem}
\begin{proof}
	We have $\mu(B(0,s))\eqD s^{Q\gamma-\gamma^2}e^{\gamma \sqrt{\log s^{-1}}\xi}\mu(B(0,1))$ for $\xi$ an independent  standard normal random variable, where we get $s^{\gamma Q}$ from \eqref{eq:coord-change}, we get $s^{-\gamma^2}$ from the drift of the process $A$ in the definition of the $\gamma$-quantum cone in Section \ref{sec:lqg}, and we get $\xi$ from the randomness of $A$. 	Letting $A_j$ be the annulus with radii $2^{-j}$ and $2^{-j-1}$, if $p<\op{min}\{1, 2Q/\gamma  - 2\}$ and by a similar scaling argument for the annulus as for the ball,
	$$
	\E[\mu(B(0,1))^p] = \E\Big[ \Big(\sum_{j=0}^\infty \mu(A_j)\Big)^p \Big] \leq \sum_{j=0}^\infty \E[\mu(A_j)^p] \leq C_p \E[\mu(A_0)^p] < \infty.
	$$
	Furthermore, $\E[\mu(B(0,s))^p]= s^{Q\gamma p - \gamma^2p - \gamma^2p^2/2} \E[\mu(B(0,1))^p]$, so 
	\eqbn
	\begin{split}
		\P[\mu(B(0,s )) >  s^{Q\gamma  - \gamma^2} x y^{\sqrt{\log s^{-1}}} ] &= \P\big[\mu(B(0,1 )) e^{-\gamma\sqrt{\log s^{-1}}\xi} > x y^{\sqrt{\log s^{-1}}} \big]\\
		&\leq \P[\mu(B(0,1 )) > x ] 
		+ \P[e^{-\gamma\xi} > y]
		\leq Cx^{-p} + \exp(-(\log y)^2/(2\gamma^2)).
	\end{split} 
	\eqen
\end{proof}

\begin{lem}
	For any $C_1>1$ there is a $C_2>1$ such that for all $r\in\R$ and $s>0$, if we assume the field $h$ has the $2^r$-embedding,
	\eqb
	\PE[ \eta([-2^{\gamma Q r - s},2^{\gamma Q r - s}]) \not\subset B(0,2^r) ] \leq C_2 2^{-C_1 s}.
	\eqe
	\label{prop:contained1}
\end{lem}
\begin{proof}
	By \eqref{eq:coord-change} and \cite[Proposition 4.13(i)]{wedges} the considered probability is equal to $\PE[ \eta([-1,1]) \not\subset B(0,2^{s/(\gamma Q)}) ]$ if we suppose the field $h$ has the $2^{s/(\gamma Q)}$-embedding. The desired bound now follows from \cite[Lemma 4.3]{gms-tutte}.
\end{proof}

\begin{lem}
	For any $\eps\in(0,0.1)$ there is a $c>0$ such that the following hold for any $r>0$, $s\in (-r/\eps,(1-\eps)r]$, and $m\in ( s\gamma Q+\eps r ,\infty)_\Z$ if we assume the field $h$ has the $2^r$-embedding
	\eqb
	\PE[ B(0,2^s) \not\subset \eta( I(m)) ] \leq c^{-1} 2^{-cr}.
	\eqe	
	\label{prop:contained2}
\end{lem}
\begin{proof}
	Pick $u\in (s+0.1\eps r,s+0.2\eps r)$ and define
	\eqb
	\tau^+ = \inf\{t\geq 0\,:\,\eta([0,t])\not\subset B(0,2^{u}) \},\qquad
	\tau^- = \sup\{t\leq 0\,:\,\eta([t,0])\not\subset B(0,2^{u}) \}.
	\label{eq:tau-pm}
	\eqe
	Then
	\eqbn
	\begin{split}
		\PE[B(0,2^{s})\not\subset \eta[I(m)]] 
		&\leq 
		\PE[B(0,2^{s})\not\subset \eta[\tau^-,\tau^+]]
		+\PE[\eta[\tau^-,\tau^+] \not \subset\eta(I(m))]\\
		&\leq 2^{-cr} 
		+ \PE[\mu(B(0,2^{u}))>2^{m-0.1\eps r}] 
		+ \PE[ [-2^{m-0.1\eps r},2^{m-0.1\eps r}] \not \subset I(m) ],
	\end{split}
	\eqen
	where we use \cite[Proposition 6.2]{hs-euclidean} to bound  $\PE[B(0,2^{s})\not\subset \eta[\tau^-,\tau^+]]$. By the definition of a uniform dyadic system, we get that the third term on the right side is also bounded above by $2^{-cr}$. Finally, the second term on the right side of the last display is also bounded by $2^{-cr}$ since
	$$
	\PE[\mu(B(0,2^{u}))>2^{m-0.1\eps r}] 
	\leq \PE[\mu(B(0,2^{s+0.2\eps r}))>2^{s\gamma Q +0.9\eps r}] \leq 2^{-cr}, 
	$$
	where we use \eqref{eq:ball-lower} and $\gamma Q=2+\gamma^2/2<4$ in the last step.
\end{proof}

We also need a variant of the above lemma where the probability decays like an arbitrarily large power of $2^r$ when we work with the $2^r$-embedding of the field.
\begin{lem}
	For any $\wt C_0>1$ one can find $\wt C_1>1$ such that for any $r>0$, if the field $h$ has the $2^r$-embedding,
	\eqbn
	\PQ[ B(0,2^{-r/2}) \not\subset \eta( [-2^{\wt C_1 r},2^{\wt C_1 r}]) ] \leq 2^{-\wt C_0 r}.
	\eqen	
	\label{prop:contained3}
\end{lem}
In the proof we will use the following estimate from \cite{gms-tutte} for the radii $\frk R(b)$ in~\eqref{eqn-mass-hit-time}.
\begin{lem}[\cite{gms-tutte}] \label{lem-cone-hit-tail}
	Consider the circle-average embedding of an $\alpha$-quantum cone for $\alpha<Q$. Then there is a constant $a=a(\alpha,\gamma)>0$ such that for each $b_2>b_1>0$ and each $C > 1$,
	\begin{align} \label{eqn-cone-hit-tail}
		&\BB P\left[ C^{-1} (b_2/b_1)^{ \tfrac{1 }{\gamma (Q-\alpha) }} \leq \frk R(b_2)/\frk R(b_1) \leq C (b_2/b_1)^{ \tfrac{1 }{\gamma (Q-\alpha)}} \right] \notag\\
		&\qquad\geq 1  -   3 \exp\left( - \frac{  a  (\log C)^2 }{   \log (b_2/b_1) +  \log C   } \right)  .
	\end{align}
\end{lem}

\begin{proof}[Proof of Lemma \ref{prop:contained3}]
	We assume throughout the proof that $r$ has been chosen sufficiently large.  
	Setting $u=\wt C_2r$ for some $\wt C_2>1$ and defining $\tau^\pm$ as in \eqref{eq:tau-pm}, by similar considerations as in the proof of Lemma \ref{prop:contained2},
	\eqbn
	\begin{split}
		\PQ[B(0,2^{-r/2})\not\subset \eta( [-2^{\wt C_1 r},2^{\wt C_1 r}]) ] 
		&\leq 
		\PQ[B(0,2^{-r/2})\not\subset \eta[\tau^-,\tau^+] ] 
		+ \PQ[\mu(B(0,2^{\wt C_2 r}))>2^{\wt C_1 r}].
	\end{split}
	\eqen
	By \cite[Proposition 6.2]{hs-euclidean}, if we choose $\wt C_2$ sufficiently large, the first probability on the right side is smaller than $2^{-2\wt C_0 r}$. Observe that by the definition of the $2^r$-embedding as given in Section \ref{sec:lqg} the second probability on the right side is bounded above by the following if we use the circle average embedding instead of the $2^r$-embedding,
	\eqbn
	\begin{split}
	\PQ[\mu(B(0,\frk R(2^{r\gamma Q}) 2^{\wt C_2 r}))>2^{\wt C_1 r}]
	\leq 
	\PQ[\mu(B(0,\frk R(2^{\wt C_3 r }) ) )>2^{\wt C_1 r}]
	+ \PQ[\frk R(2^{\wt C_3 r }) < \frk R(2^{r\gamma Q}) 2^{\wt C_2 r} ],
	\end{split}
	\eqen
	where we get the inequality by a union bound for arbitrary $\wt C_3$. 
	By choosing $\wt C_3$ sufficiently large and applying Lemma \ref{lem-cone-hit-tail}, the second probability is smaller than $2^{-2\wt C_0 r}$. 
	By scaling, the first probability is equal to $\PQ[\mu(\D) >2^{\wt C_1 r-\wt C_3 r}]$, which by Lemma \ref{prop:cone-moment} with $p\in(0,\min\{1,2Q/\gamma-2\})$ is smaller than a constant multiple of $2^{-p\wt C_1 r+p\wt C_3 r}$. Choosing $\wt C_1$ sufficiently large this is smaller than $2^{-2\wt C_0 r}$. We conclude that the probability in the statement of the lemma is bounded above by $3\cdot 2^{-2\wt C_0 r}$.
\end{proof}

\subsection{Maximal diameter of SLE cell}
	The following result, which is essentially proven in \cite{gms-harmonic}, will be used in the proof of Lemma \ref{prop:max-cell-diam} below.
	\begin{lem}\label{prop:max-cell-diam0}
		Suppose $h$ is the circle-average embedding (equivalently, $1$-embedding) of a 0-quantum cone. For each $q \in ( 0 , 2(2+\gamma)^{-2} )$, each $\rho \in (0,1)$, and each $\ep \in (0,1)$,
		\eqb \label{eqn-max-cell-diam0}
		\BB P\left[ \op{diam}\left(\eta(I) \right) \leq \ep^{q} ,\: \forall I \text{ s.t.\ } |I|=\ep/2,\eta(I)\subset B(0,\rho)  \right] \geq 1 -  \ep^{\alpha(q,\gamma) + \wt g(\eps)}  ,
		\eqe   
		where $\wt g$ is a function depending only on $q,\rho,\gamma$ such that $\wt g(x)\rta 0$ as $x\rta 0$, and
		\eqb
		\alpha(q,\gamma) := \frac{q}{2\gamma^2  } \left( \frac{1}{q} - 2 - \frac{\gamma^2}{2} \right)^2 - 2q>0.
		\label{eq:alpha}
		\eqe
	\end{lem}
	\begin{proof}
		From the proof of \cite[Lemma 2.7]{gms-harmonic} we know that if $\theta$ is a uniform number between 0 and 1 and 
		\eqbn \label{eqn-max-cell-diam00}
		E(\theta)=\{\op{diam}\left(\eta(I) \right) \leq \ep^{q} ,\: \forall I \text{ s.t.\ } I=[\ep (\theta+n),\ep(\theta+n+1)],n\in\Z,\eta(I)\cap B(0,\rho)\neq\emptyset\}
		\eqen
		then $\BB P\left[ E(\theta)  \right] \geq 1 -  \ep^{\alpha(q,\gamma) + \wt g(\eps) }$. Since $(\theta+n)_{n\in\Z}\eqD(\theta+1/2+n)_{n\in\Z}$ we have further
		$$
		\BB P\left[ E(\theta)\cap E(\theta+1/2)  \right] \geq 1 -  \ep^{\alpha(q,\gamma) + o_\ep(1)}.
		$$
		The lemma now follows since the event in the latter indented equation for any fixed $\theta$ implies the event in \eqref{eqn-max-cell-diam0}.
	\end{proof}
	
Combining Lemma \ref{prop:max-cell-diam0} and a scaling argument gives the following.
	\begin{lem}
		For $r>1$ suppose $h$ is the $2^r$-embedding of a 0-quantum cone. For each $q \in ( 0 , 2(2+\gamma)^{-2} )$, each $\rho \in (0,1)$, and each $m<\gamma Q r-1$, 
		\eqb
		\begin{split}
			&\BB P\left[ \max\{ \op{diam}(\eta([t,t+2^{m}]))\,:\, \eta([t,t+2^{m}])\subset B(0,\rho 2^r),t\in\R \} \leq 2^{r-\gamma Q r q+qm+q}   
			\right] \\
			&\qquad\geq 1 -  2^{(-\gamma Q r+m+1)\alpha(q,\gamma) + g(\gamma Q r-m)}  ,
		\end{split}
		\label{eq:maxdiam}
		\eqe  
		where $g$ is a function depending only on $q,\rho,\gamma$ such that $g(x)/x\rta 0$ as $x\rta\infty$, and $\alpha(q,\gamma)$ is given by \eqref{eq:alpha}. Furthermore,
		\eqb
		\begin{split}
			&\E[\max\{ \op{diam}(\eta([t,t+2^{m}]))^2\,:\, \eta([t,t+2^{m}])\subset B(0,\rho 2^r),t\in\R  \}]\\ 
			&\qquad\qquad\leq 2^{2(r-\gamma Q r q+qm+q)} + 2^{2(r+1) + 
				(-\gamma Q r + m + 1)\alpha(q,\gamma)+ g(\gamma Q r-m) }.
		\end{split}
		\label{eq72}
		\eqe
		\label{prop:max-cell-diam}
	\end{lem}
	\begin{proof}
		Let $h^0$ and $h^r$ be whole-plane GFFs normalized to have mean 0 on $\partial B(0,1)$ and $\partial B(0,2^r)$, respectively. We will work with $h^r$ instead of $h$ throughout the proof, which is justified since $h|_{B(0,2^r)}\eqD h^r|_{B(0,2^r)}$. By conformal invariance of the GFF we may couple $h^0$ and $h^r$ such that $h^0=h^r\circ\phi$, where $\phi:\C\to\C$ is the conformal map given by $\phi(z)=2^rz$. Let $\eta_0$ and $\eta_r$ be whole-plane space-filling SLE$_\kappa$ parametrized by $\mu_{h^0}$ and $\mu_{h^r}$, respectively, and suppose the curves are coupled together such that $\eta_0$ modulo parametrization is equal to $2^{-r}\eta_r$. By \eqref{eq:coord-change}, for any $A\subset\C$ we have
		$$
		\mu_{h^r}(A)=2^{\gamma Q r}\mu_{h^0}(\phi^{-1}(A)).
		$$
		Using this we get that for any $\eps\in(0,1)$,
		\eqbn
		\begin{split}
			&\{ \max\{ \op{diam}(\eta_0([t,t+\eps/2]))\,:\, \eta_0([t,t+\eps/2])\subset B(0,\rho),t\in\R  \}\leq \eps^q\} \\
			&\qquad= 
			\{\max\{ \op{diam}(\eta_r([t,t+2^{\gamma Q r-1}\eps]))\,:\, \eta_r([t,t+2^{\gamma Q r-1}\eps])\subset B(0,\rho 2^r),t\in\R  \} \leq 2^r\eps^q \}.
		\end{split}
		\eqen
		Setting $\eps=2^{-\gamma Q r+m+1}$ gives the first assertion of the lemma. The second assertion \eqref{eq72} of the lemma follows from the first assertion since if the event in \eqref{eq:maxdiam} fails to hold then the considered diameter is still smaller than $2^{r+1}$. The first (resp.\ second) term on the right side of \eqref{eq72} comes from the case when the event in \eqref{eq:maxdiam} holds (resp.\ does not hold).
\end{proof}

\subsection{Roundness of cells}
\label{app3}
In this appendix we will prove Lemma \ref{prop:prop44}. Most of the proof will be devoted to showing the following lemma, which is identical to Lemma \ref{prop:prop44} except that we require $0\in J$ in the event $\wt B^0_3$ and that we write $C_5$ instead of $C$. 
\begin{lem}
	For any given $C_2,C_4,C_5>0$ it holds that for $M>1$ sufficiently large, $\PQ[ \wt B^0_3] \leq M^{-C_5}$, where
	$$ 
	\wt B^0_3 =
	\Big\{ \max
	\Big\{ \frac{\op{diam}(\eta(J))^2}{\op{area}(\eta(J))} \,:\, J\subset [-M^{C_4},M^{C_4}] \text{ interval}, \,M^{-C_2}\leq|J|\leq C_2\log M, 0\in J 
	\Big\} >M
	\Big\}.
	$$
	\label{prop:prop44-0}
\end{lem}
Before giving the proof of Lemma \ref{prop:prop44-0} let us explain why it implies Lemma \ref{prop:prop44}.
\begin{proof}[Proof of Lemma \ref{prop:prop44} assuming Lemma \ref{prop:prop44-0}]
	Defining $\cI=\{0.5kM^{-C_2}\,:\,k\in\Z \}\cap [-M^{C_4},M^{C_4}]$ and letting $\wt B^{t_0}_3$ be defined exactly as in the statement of Lemma \ref{prop:prop44-0} except that we require $t_0\in J$ instead of $0\in J$,
	$$
	\PQ[B_3] 
	\leq \PQ\bigg[ \bigcup_{t_0 \in\cI} \wt B^{t_0}_3 \bigg]
	\leq \sum_{t_0\in\cI} \PQ[\wt B^{t_0}_3] = \#\cI\cdot \PQ[\wt B^0_3],
	$$

	where we use translation invariance in law in the last step. We conclude by using $\#\cI\leq 4M^{C_4+C_2}$ and Lemma \ref{prop:prop44-0} with $C_5>C_2+C_4+C$.
\end{proof}
Recall $\frk R(b)$ from \eqref{eqn-mass-hit-time} and  observe that for any bounded interval $I\subset\R$ and $a>0$,
\eqb
\PQ[\eta(I)\subset B(0,\frk R(b))] = \PQ[\eta(aI)\subset B(0,\frk R(ab))].
\label{eq:Rb-scaling}
\eqe

\begin{lem} \label{lem-radius-segment}  
	For any $C>1$ and $\ep\in(0,1)$, for sufficiently large $M$ and if we work with the circle-average embedding of the $\gamma$-quantum cone, 
	\eqb \label{eqn-radius-segment-upper}
	\PQ\left[ \eta([ 0 ,1]) \not\subset B_{\frk R(1/M ) }(0) \right] \geq 1 - M^{-\min\{1,2Q/\gamma-2\} +\eps}
	\eqe
	and 
	\eqb\label{eqn-radius-segment-lower}
	\PQ\left[  \eta([-1, 1])  \subset B_{  \frk R(M) / 8  }(0) \right] \geq 1 - M^{-C}. 
	\eqe  
\end{lem}
\begin{proof}
	This was proven for a 0-quantum cone in \cite[Lemma 4.3]{gms-tutte} (see also Lemma \ref{prop:contained1}). The second assertion for a $\gamma$-quantum cone is immediate from the 0-quantum cone case since the field of a $\gamma$-quantum cone stochastically dominates the field of a 0-quantum cone if both fields have the $2^r$-embedding for some $r$ and we consider the fields restricted to $B(0,2^r)$. Following the proof of \cite[Lemma 4.3]{gms-tutte}, in order to prove the first assertion we need to argue that $\EQ[\mu_h(\D)^\beta]<\infty$ for any $\beta<\op{min}\{1, 2Q/\gamma  - 2\}$, where $h$ is the field of a $\gamma$-quantum cone with the 1-embedding. This follows from Lemma \ref{prop:cone-moment}.
\end{proof}

\begin{proof}[Proof of Lemma \ref{prop:prop44-0}]
	We proceed similarly as for the proof of \cite[Proposition 4.4]{gms-tutte}. We will apply~\cite[Lemma 3.6]{ghm-kpz} (which says that $\eta$ is unlikely to cross an annulus without absorbing a Euclidean ball of radius comparable to the distance across the annulus) at several scales and use Lemma~\ref{lem-radius-segment} to argue that with high probability, for $J$ as in the lemma, $\eta|_J$ has to cross one of these scales. 
	
	Let $C_6 >2C_5/\min\{1,2Q/\gamma-2\}+C_2$. The random variables $\frk R(M)$ and $\frk R(M^{-C_6})$ for $M > 1$ depend only on $ h$, so are independent from the curve $ \eta$, viewed modulo monotone re-parameterization.
	For $r>0$, write $\tau(r)=\inf\{t\geq0:|\eta(t)|\geq r\}$ and $\tau^-(r)=\sup\{t\leq0:|\eta(t)|\geq r\}$.  By~\cite[Lemma~3.6]{ghm-kpz} and a union bound over $k$, we have 
	\alb
	&\BB P\left[ \op{area}\left( \eta([ 0 ,\tau(  e^k ) ]  ) \right) \geq 4e^2M^{-1} e^{2k} ,\: \forall k \in \left[  \log \frk R(M^{-C_6}) -1  ,   \log \frk R(M)     +1 \right]_{\BB Z} \,|\, h \right]  \\
	&\qquad \geq 1 - \left( \log \left( \frk R(M) / \frk R(M^{-C_6}) \right)+2\right) e^{-c_0 M^{1/2}}            .
	\ale
	for $c_0 > 0$ a constant depending only on $\gamma$. By combining this with Lemma~\ref{lem-cone-hit-tail} (applied with $\alpha = \gamma$ and, e.g., $e^{M^{1/4}}$ in place of $M$), we get that except on an event of probability $o_M^\infty(M)$ (at a rate depending only on $C_6$ and $\gamma$), 
	\eqb \label{eqn-exp-scale-area}
	\op{area}\left( \eta([ 0 ,\tau(  e^k ) ]) \right) \geq M^{-1} e^{2k} ,\quad \forall k \in \left[  \log \frk R(M^{-C_6}) -1   ,   \log \frk R(M)   +1  \right]_{\BB Z}.
	\eqe 
	By the invariance of the law of $\eta$ under time reversal, it also holds except on an event of probability $o_M^\infty(M)$ that
	\eqb \label{eqn-exp-scale-area'}
	\op{area}\left( \eta([\tau^-(e^k), 0  ]) \right) \geq M^{-1} e^{2k} ,\quad \forall k \in \left[  \log \frk R(M^{-C_6}) -1   ,   \log \frk R(M)  +1  \right]_{\BB Z}. 
	\eqe
	By Lemma~\ref{lem-radius-segment}, \eqref{eq:Rb-scaling}, and the invariance of the law of $\eta$ under time reversal, it holds with probability $1-M^{-\min\{1,2Q/\gamma-2\}(C_6-C_2) + o_M(1)}\geq 1-M^{-2C_5}$ for sufficiently large $M$ that 
	\eqb \label{eqn-use-free-gff-mass}
	\begin{split}
		\eta([0, M^{-C_2}/2]) ,\eta([-M^{-C_2}/2 , 0]) &\not\subset B(0,\frk R(M^{-C_6}))  \\   
		\eta([ -C_2\log M ,C_2\log M]) &\subset B(0,\frk R(M)) .
	\end{split}
	\eqe  
		
	Suppose now that \eqref{eqn-exp-scale-area}, \eqref{eqn-exp-scale-area'}, and~\eqref{eqn-use-free-gff-mass} all hold (which happens with probability at least $1-M^{-2C_5}$ of  sufficiently large $M$). We seek to bound the maximum appearing in the lemma statement. 
	To this end, suppose $J$ is an interval as in the lemma statement. Then
	\eqbn
	\begin{split}
		\eta(J) \subset \eta([- C_2\log M , C_2\log M] ),  
		\quad \op{and}\quad
		\eta([-M^{-C_2}/2, 0])\subset \ol{\eta(J)}
		\op{or}
		\eta([0, M^{-C_2}/2]) \subset \ol{\eta(J)}.
	\end{split}
	\eqen
	By~\eqref{eqn-use-free-gff-mass}, if we let $k$ be the smallest integer for which $\eta(J) \subset B_{e^{k+1}}(0)$, then 
	\eqbn
	k\in \left[  \log \frk R(M^{-C_6}) -1   ,    \log \frk R(M)    +1  \right]_{\BB Z} .
	\eqen
	Furthermore, by the definition of $k$ either $ \eta([  0 , \tau(  e^k ) ]) $ or $ \eta([ \tau^-(e^k) , 0  ]) $ is contained in $\ol{\eta(J)}$. By~\eqref{eqn-exp-scale-area} and~\eqref{eqn-exp-scale-area'},
	\eqb \label{eqn-area-diam-implication}
	\op{area}\left(  \eta(J) \right) 
	\geq 4e^{2}M^{-1}e^{2k}
	\geq  M^{-1} \op{diam}\left( \eta(J)  \right)^2 ,
	\eqe 
	which implies that the event in the lemma does not occur.
\end{proof}

\bibliography{cibib} 

\def\cprime{$'$}
\begin{thebibliography}{KMSW19}

\bibitem[AFS20]{afs-metric-ball}
M.~Ang, H.~Falconet, and X.~Sun.
\newblock Volume of metric balls in {L}iouville quantum gravity.
\newblock {\em Electron. J. Probab.}, 25:Paper No. 160, 50, 2020. \MR{4193901}

\bibitem[AG19]{ag-disk}
M.~{Ang} and E.~{Gwynne}.
\newblock {L}iouville quantum gravity surfaces with boundary as matings of
  trees.
\newblock {\em ArXiv e-prints}, Mar 2019, \arxiv{1903.09120}.

\bibitem[APPS22]{apps-central-charge}
M.~Ang, M.~Park, J.~Pfeffer, and S.~Sheffield.
\newblock Brownian loops and the central charge of a {L}iouville random
  surface.
\newblock {\em Ann. Probab.}, 50(4):1322--1358, 2022. \MR{4420421}

\bibitem[Ber07]{bernardi-maps}
O.~Bernardi.
\newblock Bijective counting of tree-rooted maps and shuffles of parenthesis
  systems.
\newblock {\em Electron. J. Combin.}, 14(1):Research Paper 9, 36 pp.
  (electronic), 2007, \arxiv{math/0601684}. \MR{2285813 (2007m:05125)}

\bibitem[Ber17]{berestycki-gmt-elementary}
N.~Berestycki.
\newblock An elementary approach to {G}aussian multiplicative chaos.
\newblock {\em Electron. Commun. Probab.}, 22:Paper No. 27, 12, 2017,
  \arxiv{1506.09113}. \MR{3652040}

\bibitem[Ber25]{berestycki-spectral}
N.~Berestycki.
\newblock On the spectral geometry of liouville quantum gravity, 2025,
  2512.02538.

\bibitem[BG22]{berestycki-gwynne-lbm}
N.~Berestycki and E.~Gwynne.
\newblock Random walks on mated-{CRT} planar maps and {L}iouville {B}rownian
  motion.
\newblock {\em Comm. Math. Phys.}, 395(2):773--857, 2022. \MR{4487526}

\bibitem[BGS25]{bgs-smith}
F.~Bertacco, E.~Gwynne, and S.~Sheffield.
\newblock Scaling limits of planar maps under the {S}mith embedding.
\newblock {\em Ann. Probab.}, 53(3):1138--1196, 2025. \MR{4908099}

\bibitem[BHS23]{bhs-site-perc}
O.~Bernardi, N.~Holden, and X.~Sun.
\newblock Percolation on triangulations: a bijective path to {L}iouville
  quantum gravity.
\newblock {\em Mem. Amer. Math. Soc.}, 289(1440):v+176, 2023. \MR{4651497}

\bibitem[BK26]{bk-spectral}
N.~Berestycki and J.~Klein.
\newblock Spectral expansion of lqg heat trace and kpz scaling, 2026,
  2606.26958.

\bibitem[BP24]{bp-lqg-book}
N.~Berestycki and E.~Powell.
\newblock Gaussian free field and liouville quantum gravity.
\newblock {\em arXiv preprint arXiv:2404.16642}, 2024.

\bibitem[BRG24]{bg-idla}
A.~Bou-Rabee and E.~Gwynne.
\newblock Internal {DLA} on mated-{CRT} maps.
\newblock {\em Ann. Probab.}, 52(6):2173--2237, 2024. \MR{4815971}

\bibitem[BS01]{benjamini-schramm-topology}
I.~Benjamini and O.~Schramm.
\newblock Recurrence of distributional limits of finite planar graphs.
\newblock {\em Electron. J. Probab.}, 6:no. 23, 13 pp. (electronic), 2001,
  \arxiv{0011019}. \MR{1873300 (2002m:82025)}

\bibitem[BW24]{bw-weyl}
N.~Berestycki and M.~D. Wong.
\newblock Weyl's law in liouville quantum gravity, 2024, 2307.05407.

\bibitem[Che17]{chen-fk}
L.~Chen.
\newblock Basic properties of the infinite critical-{FK} random map.
\newblock {\em Ann. Inst. Henri Poincar\'e D}, 4(3):245--271, 2017,
  \arxiv{1502.01013}. \MR{3713017}

\bibitem[Che18]{chelkak-icm}
D.~Chelkak.
\newblock Planar ising model at criticality: state-of-the-art and perspectives.
\newblock In {\em Proceedings of the International Congress of Mathematicians:
  Rio de Janeiro 2018}, pages 2801--2828. World Scientific, 2018.

\bibitem[Che24]{chelkak-s-emb-periodic}
D.~Chelkak.
\newblock Ising model and s-embeddings of planar graphs, 2024. \MR{4885806}

\bibitem[CLR23]{clr-t-emb}
D.~Chelkak, B.~Laslier, and M.~Russkikh.
\newblock Dimer model and holomorphic functions on t-embeddings of planar
  graphs.
\newblock {\em Proc. Lond. Math. Soc. (3)}, 126(5):1656--1739, 2023.
  \MR{4588721}

\bibitem[Cur15]{curien-glimpse}
N.~Curien.
\newblock A glimpse of the conformal structure of random planar maps.
\newblock {\em Comm. Math. Phys.}, 333(3):1417--1463, 2015, \arxiv{1308.1807}.
  \MR{3302638}

\bibitem[CW25]{cw-liouville}
S.~Chatterjee and E.~Witten.
\newblock Liouville theory: an introduction to rigorous approaches.
\newblock {\em Journal of High Energy Physics}, 2025(2):1--41, 2025.

\bibitem[Dav88]{david-conformal-gauge}
F.~David.
\newblock Conformal field theories coupled to {2-D} gravity in the conformal
  gauge.
\newblock {\em {M}od. {P}hys. {L}ett. {A}}, (3), 1988.

\bibitem[DK89]{dk-qg}
J.~Distler and H.~Kawai.
\newblock Conformal field theory and {2D} quantum gravity.
\newblock {\em {N}ucl.{P}hys. {B}}, (321), 1989.

\bibitem[DKRV16]{dkrv-lqg-sphere}
F.~David, A.~Kupiainen, R.~Rhodes, and V.~Vargas.
\newblock Liouville quantum gravity on the {R}iemann sphere.
\newblock {\em Comm. Math. Phys.}, 342(3):869--907, 2016, \arxiv{1410.7318}.
  \MR{3465434}

\bibitem[DMS21]{wedges}
B.~Duplantier, J.~Miller, and S.~Sheffield.
\newblock Liouville quantum gravity as a mating of trees.
\newblock {\em Ast\'erisque}, (427):viii+257, 2021. \MR{4340069}

\bibitem[DS11]{shef-kpz}
B.~Duplantier and S.~Sheffield.
\newblock Liouville quantum gravity and {KPZ}.
\newblock {\em Invent. Math.}, 185(2):333--393, 2011, \arxiv{1206.0212}.
  \MR{2819163 (2012f:81251)}

\bibitem[Dur10]{durrett}
R.~Durrett.
\newblock {\em Probability: theory and examples}.
\newblock Cambridge Series in Statistical and Probabilistic Mathematics.
  Cambridge University Press, Cambridge, fourth edition, 2010. \MR{2722836
  (2011e:60001)}

\bibitem[GGJN20]{gjn-dirichlet}
O.~Gurel-Gurevich, D.~C. Jerison, and A.~Nachmias.
\newblock The {D}irichlet problem for orthodiagonal maps.
\newblock {\em Adv. Math.}, 374:107379, 53, 2020. \MR{4157579}

\bibitem[GHM15]{ghm-kpz}
E.~{Gwynne}, N.~{Holden}, and J.~{Miller}.
\newblock {An almost sure KPZ relation for SLE and Brownian motion}.
\newblock {\em {A}nnals of {P}robability}, to appear, 2015, \arxiv{1512.01223}.

\bibitem[GHMS17]{kappa8-cov}
E.~Gwynne, N.~Holden, J.~Miller, and X.~Sun.
\newblock Brownian motion correlation in the peanosphere for {$\kappa>8$}.
\newblock {\em Ann. Inst. Henri Poincar\'e Probab. Stat.}, 53(4):1866--1889,
  2017, \arxiv{1510.04687}. \MR{3729638}

\bibitem[GHS20]{ghs-map-dist}
E.~Gwynne, N.~Holden, and X.~Sun.
\newblock A mating-of-trees approach for graph distances in random planar maps.
\newblock {\em Probab. Theory Related Fields}, 177(3-4):1043--1102, 2020.
  \MR{4126936}

\bibitem[GHS23]{ghs-mating-survey}
E.~Gwynne, N.~Holden, and X.~Sun.
\newblock Mating of trees for random planar maps and {L}iouville quantum
  gravity: a survey.
\newblock {\em Panor. Synth\`eses}, 59:41--120, 2023. \MR{4619311}

\bibitem[GHSS18]{ghss-negative-moments}
C.~{Garban}, N.~{Holden}, A.~{Sep{\'u}lveda}, and X.~{Sun}.
\newblock {Negative moments for Gaussian multiplicative chaos on fractal sets}.
\newblock {\em ArXiv e-prints}, May 2018, \arxiv{1805.00864}.

\bibitem[GJN19]{gjn-macroscopic-circles}
O.~{Gurel-Gurevich}, D.~C. {Jerison}, and A.~{Nachmias}.
\newblock {A combinatorial criterion for macroscopic circles in planar
  triangulations}.
\newblock {\em ArXiv e-prints}, Jun 2019, \arxiv{1906.01612}.

\bibitem[GMS19]{gms-harmonic}
E.~{Gwynne}, J.~{Miller}, and S.~{Sheffield}.
\newblock {Harmonic functions on mated-{CRT} maps}.
\newblock {\em Electron. J. Probab.}, 24:no. 58, 55, 2019, \arxiv{1807.07511}.

\bibitem[GMS20]{gms-poisson-voronoi}
E.~Gwynne, J.~Miller, and S.~Sheffield.
\newblock The {T}utte embedding of the {P}oisson-{V}oronoi tessellation of the
  {B}rownian disk converges to {$\sqrt{8/3}$}-{L}iouville quantum gravity.
\newblock {\em Comm. Math. Phys.}, 374(2):735--784, 2020. \MR{4072229}

\bibitem[GMS21]{gms-tutte}
E.~Gwynne, J.~Miller, and S.~Sheffield.
\newblock The {T}utte embedding of the mated-{CRT} map converges to {L}iouville
  quantum gravity.
\newblock {\em Ann. Probab.}, 49(4):1677--1717, 2021. \MR{4260465}

\bibitem[GMS22]{gms-random-walk}
E.~Gwynne, J.~Miller, and S.~Sheffield.
\newblock An invariance principle for ergodic scale-free random environments.
\newblock {\em Acta Math.}, 228(2):303--384, 2022. \MR{4448682}

\bibitem[GP19]{gp-dla}
E.~{Gwynne} and J.~{Pfeffer}.
\newblock {External diffusion limited aggregation on a spanning-tree-weighted
  random planar map}.
\newblock {\em ArXiv e-prints}, January 2019, \arxiv{1901.06860}.

\bibitem[GR13]{gill-rohde-type}
J.~T. Gill and S.~Rohde.
\newblock On the {R}iemann surface type of random planar maps.
\newblock {\em Rev. Mat. Iberoam.}, 29(3):1071--1090, 2013, \arxiv{1101.1320}.
  \MR{3090146}

\bibitem[Guy00]{guy-walk-in-cone}
R.~K. Guy.
\newblock Catwalks, sandsteps and {P}ascal pyramids.
\newblock {\em J. Integer Seq.}, 3(1):Article 00.1.6, 1 HTML document, 2000.
  \MR{1750743}

\bibitem[HG25]{cp-gwynne-curvature}
A.~C. Hip and E.~Gwynne.
\newblock Gaussian curvature on random planar maps and liouville quantum
  gravity, 2025+, ArXiv e-print 2406.08674.

\bibitem[HS18]{hs-euclidean}
N.~Holden and X.~Sun.
\newblock S{LE} as a mating of trees in {E}uclidean geometry.
\newblock {\em Comm. Math. Phys.}, 364(1):171--201, 2018, \arxiv{1610.05272}.
  \MR{3861296}

\bibitem[HS23]{hs-cardy-embedding}
N.~Holden and X.~Sun.
\newblock Convergence of uniform triangulations under the {C}ardy embedding.
\newblock {\em Acta Math.}, 230(1):93--203, 2023. \MR{4567714}

\bibitem[HY26]{hy-cp-ru}
N.~Holden and P.~Yu.
\newblock Circle packing and riemann uniformization of random planar maps in an
  ergodic scale-free environment, 2026, 2603.06528.

\bibitem[Kah85]{kahane}
J.-P. Kahane.
\newblock Sur le chaos multiplicatif.
\newblock {\em Ann. Sci. Math. Qu\'ebec}, 9(2):105--150, 1985. \MR{829798
  (88h:60099a)}

\bibitem[KLRR22]{klrr-circle-patterns}
R.~Kenyon, W.~Y. Lam, S.~Ramassamy, and M.~Russkikh.
\newblock Dimers and circle patterns.
\newblock {\em Ann. Sci. \'Ec. Norm. Sup\'er. (4)}, 55(3):865--903, 2022.
  \MR{4553657}

\bibitem[KMSW19]{kmsw-bipolar}
R.~Kenyon, J.~Miller, S.~Sheffield, and D.~B. Wilson.
\newblock Bipolar orientations on planar maps and {${\rm SLE}_{12}$}.
\newblock {\em Ann. Probab.}, 47(3):1240--1269, 2019, \arxiv{1511.04068}.
  \MR{3945746}

\bibitem[KMT76]{kmt}
J.~Koml\'os, P.~Major, and G.~Tusn\'ady.
\newblock An approximation of partial sums of independent {RV}'s, and the
  sample {DF}. {II}.
\newblock {\em Z. Wahrscheinlichkeitstheorie und Verw. Gebiete}, 34(1):33--58,
  1976. \MR{0402883}

\bibitem[KPZ88]{kpz-scaling}
V.~Knizhnik, A.~Polyakov, and A.~Zamolodchikov.
\newblock {Fractal structure of 2D-quantum gravity}.
\newblock {\em {Modern Phys. Lett A}}, 3(8):819--826, 1988.

\bibitem[{Le }13]{legall-uniqueness}
J.-F. {Le Gall}.
\newblock Uniqueness and universality of the {B}rownian map.
\newblock {\em Ann. Probab.}, 41(4):2880--2960, 2013, \arxiv{1105.4842}.
  \MR{3112934}

\bibitem[{Le }14]{legall-sphere-survey}
J.-F. {Le Gall}.
\newblock {Random geometry on the sphere}.
\newblock {\em Proceedings of the {ICM}}, 2014, \arxiv{1403.7943}.

\bibitem[LSW17]{lsw-schnyder-wood}
Y.~{Li}, X.~{Sun}, and S.~S. {Watson}.
\newblock {Schnyder woods, SLE(16), and Liouville quantum gravity}.
\newblock {\em ArXiv e-prints}, May 2017, \arxiv{1705.03573}.

\bibitem[Mie13]{miermont-brownian-map}
G.~Miermont.
\newblock The {B}rownian map is the scaling limit of uniform random plane
  quadrangulations.
\newblock {\em Acta Math.}, 210(2):319--401, 2013, \arxiv{1104.1606}.
  \MR{3070569}

\bibitem[MS16]{ig1}
J.~Miller and S.~Sheffield.
\newblock Imaginary geometry {I}: interacting {SLE}s.
\newblock {\em Probab. Theory Related Fields}, 164(3-4):553--705, 2016,
  \arxiv{1201.1496}. \MR{3477777}

\bibitem[MS17]{ig4}
J.~Miller and S.~Sheffield.
\newblock Imaginary geometry {IV}: interior rays, whole-plane reversibility,
  and space-filling trees.
\newblock {\em Probab. Theory Related Fields}, 169(3-4):729--869, 2017,
  \arxiv{1302.4738}. \MR{3719057}

\bibitem[MS21]{lqg-tbm3}
J.~Miller and S.~Sheffield.
\newblock Liouville quantum gravity and the {B}rownian map {III}: the conformal
  structure is determined.
\newblock {\em Probab. Theory Related Fields}, 179(3-4):1183--1211, 2021.
  \MR{4242633}

\bibitem[Mul67]{mullin-maps}
R.~C. Mullin.
\newblock On the enumeration of tree-rooted maps.
\newblock {\em Canad. J. Math.}, 19:174--183, 1967. \MR{0205882 (34 \#5708)}

\bibitem[Pol02]{pollard}
D.~Pollard.
\newblock {\em A user's guide to measure theoretic probability}, volume~8 of
  {\em Cambridge Series in Statistical and Probabilistic Mathematics}.
\newblock Cambridge University Press, Cambridge, 2002. \MR{1873379}

\bibitem[RV14]{rhodes-vargas-review}
R.~Rhodes and V.~Vargas.
\newblock Gaussian multiplicative chaos and applications: {A} review.
\newblock {\em Probab. Surv.}, 11:315--392, 2014, \arxiv{1305.6221}.
  \MR{3274356}

\bibitem[Sch00]{schramm0}
O.~Schramm.
\newblock Scaling limits of loop-erased random walks and uniform spanning
  trees.
\newblock {\em Israel J. Math.}, 118:221--288, 2000, \arxiv{math/9904022}.
  \MR{1776084 (2001m:60227)}

\bibitem[She07]{shef-gff}
S.~Sheffield.
\newblock Gaussian free fields for mathematicians.
\newblock {\em Probab. Theory Related Fields}, 139(3-4):521--541, 2007,
  \arxiv{math/0312099}. \MR{2322706 (2008d:60120)}

\bibitem[She16a]{shef-zipper}
S.~Sheffield.
\newblock Conformal weldings of random surfaces: {SLE} and the quantum gravity
  zipper.
\newblock {\em Ann. Probab.}, 44(5):3474--3545, 2016, \arxiv{1012.4797}.
  \MR{3551203}

\bibitem[She16b]{shef-burger}
S.~Sheffield.
\newblock Quantum gravity and inventory accumulation.
\newblock {\em Ann. Probab.}, 44(6):3804--3848, 2016, \arxiv{1108.2241}.
  \MR{3572324}

\bibitem[SS13]{ss-contour}
O.~Schramm and S.~Sheffield.
\newblock A contour line of the continuum {G}aussian free field.
\newblock {\em Probab. Theory Related Fields}, 157(1-2):47--80, 2013,
  \arxiv{math/0605337}. \MR{3101840}

\bibitem[WP21]{wp-gff}
W.~Werner and E.~Powell.
\newblock {\em Lecture notes on the {G}aussian free field}, volume~28 of {\em
  Cours Sp\'ecialis\'es [Specialized Courses]}.
\newblock Soci\'et\'e{} Math\'ematique de France, Paris, 2021. \MR{4466634}

\bibitem[Zai98]{zaitsev-kmt}
A.~Y. Zaitsev.
\newblock Multidimensional version of the results of {K}oml\'os, {M}ajor and
  {T}usn\'ady for vectors with finite exponential moments.
\newblock {\em ESAIM Probab. Statist.}, 2:41--108, 1998. \MR{1616527}

\end{thebibliography}
\bibliographystyle{hmralphaabbrv}

\end{document}